%% file: reshaping.tex
\documentclass[10pt,twoside,dvips]{article}
\usepackage[a5paper,layoutwidth=148.5mm,layoutheight=7.77817in,outer=10mm,inner=38pt,tmargin=10mm,bmargin=10mm,footskip=15pt]{geometry}
\usepackage[utf8]{inputenc}
\usepackage[T1]{fontenc}

\usepackage{lmodern}
\usepackage{fix-cm}
\usepackage{microtype}
\newcommand\textb[1]{{\fontseries{b}\selectfont #1}}

\usepackage{natbib}
\bibpunct{(}{)}{;}{a}{}{,}
\newcommand\urlprefix{}
\usepackage[french,german,english]{babel}

\newsavebox{\labelbox}
\newlength{\leftm}
\newenvironment{liste}[1]{\sbox{\labelbox}{#1\enspace}\setlength\leftm{\wd\labelbox}\setlength\parindent{0pt}\setlength\leftskip{\leftm}}{\par}
\newcommand\ligne[1][]{\par\leavevmode\llap{\makebox[\leftm]{#1\enspace}}}

\usepackage[normalem]{ulem}
\makeatletter
\def\uuuline{\leavevmode \bgroup
 \UL@setULdepth
 \advance\ULdepth1.5\p@
 \markoverwith{\lower\ULdepth\hbox
   {\kern-.03em\vbox{\hrule width.2em\kern1\p@\hrule width.2em\kern1\p@\hrule}\kern-.03em}}%
 \ULon}
\makeatother
\usepackage{soulutf8}
\sodef\so{}{.1em}{.5em plus.1em minus.1em}{.5em plus.1em minus.1em}

\usepackage[np]{numprint}

\makeatletter
\newcommand{\leqnomode}{\tagsleft@true\let\veqno\@@leqno}
\newcommand{\reqnomode}{\tagsleft@false\let\veqno\@@eqno}
\makeatother

\usepackage{amsmath}
\usepackage{amssymb}
\renewcommand\implies{\;\Rightarrow\;}
\renewcommand\iff{\;\Leftrightarrow\;}
\newcommand\afini{\{a_1,\allowbreak\dots,a_n\}}
\newcommand\bfini{\{b_1,\allowbreak\dots,b_m\}}
\newcommand\cfini{\{c_1,\allowbreak\dots,c_p\}}
\newcommand\dfini{\{d_1,\allowbreak\dots,d_q\}}
\usepackage{stmaryrd}
\newcommand\divides{\mathbin|}

\usepackage[scr]{rsfso}
\usepackage{cancel}
\usepackage[inline,shortlabels]{enumitem}
\usepackage{framed}
\newenvironment{Cadre}{
  \framed\selectlanguage{german}}{\endframed}
\newenvironment{EngCadre}{
  \framed\selectlanguage{english}}{\endframed}

\usepackage{tocloft}

\AtBeginDocument{}
\renewcommand\bibname{References.}
\AtBeginDocument{\renewcommand\bibname{References.}}
\usepackage[nottoc]{tocbibind}
\usepackage[raggedright]{titlesec}
\titlelabel{\thetitle.\quad}

\PassOptionsToPackage{hyphens}{url}
\usepackage{hyperref}
\hypersetup{
  pdfpagelabels=true,
  bookmarksdepth=2,
  bookmarksopen=false,
  breaklinks=true,
  plainpages=false,
  colorlinks=true,
  linkcolor=blue,
  citecolor=red,
  urlcolor=red,
  hypertexnames=false,
  pdfencoding=pdfdoc
}

\usepackage{breakurl}
\usepackage{booktabs}
\usepackage{doi}

\ifx\texorpdfstring\undefined\newcommand\texorpdfstring[2]{#1}\fi
\ifx\pdfbookmark[3]{}\fi

\usepackage[capitalise]{cleveref}

\crefname{property}{Property}{Properties}
\crefname{inequality}{Inequality}{Inequalities}
\crefname{equality}{Equality}{Equalities}
\creflabelformat{inequality}{(#2#1#3)}
\creflabelformat{equality}{(#2#1#3)}
\crefname{equation}{}{}
\renewcommand\cpageref[1]{page~\pageref{#1}}

\let\su\mathfrak

\newcommand*\overneg[1]{%
  \setbox0\hbox{$\mathaccent"0362{#1}^H$}%
  \setbox2\hbox{$\mathaccent"0362{\kern0pt#1}^H$}%
  \ifdim\ht0=\ht2 \overline{#1}\else \bar#1\fi
  }
\newcommand*\dng[1]{%
  \setbox0\hbox{$\mathaccent"0362{#1}^H$}%
  \setbox2\hbox{$\mathaccent"0362{\kern0pt#1}^H$}%
  \ifdim\ht0=\ht2 \overline{\overline{#1}}\else \bar{\bar#1}\fi
  }
\newcommand\bigland{\mathop{\textstyle\bigwedge}\limits}
\newcommand\biglor{\mathop{\textstyle\bigvee}\limits}

\newcommand\barre[2][]{#1}
\newcommand\insere[1]{#1}
\newcommand\ajout[1]{}
\newcommand\Bs{-} 

\newcommand\bruch{}

\newcommand\bs{} 
\newcommand\nbsp{~}
\newcommand\bsp{\space}
\newcommand\comment[1]{\footnote{\selectlanguage{english}#1}}

\title{Lorenzen's reshaping of Krull's Fundamentalsatz\\for integral domains (1938--1953)}
\author{Stefan Neuwirth}
\date{}

\begin{document}

\maketitle

\begin{abstract}
  Krull's Fundamentalsatz, the generalisation of the main theorem
of elementary number theory to integral domains, is the starting point of
Lorenzen’s career in mathematics. This article traces a conceptual
  history of Lorenzen's successive reformulations of the
  Fundamentalsatz on the basis of excerpts of his articles. An edition
  of the extant correspondence of Lorenzen with Hasse, Krull, and Aubert
  provides a better understanding of the context of these
  investigations.
\end{abstract}

\tableofcontents

\section{Introduction.}

In the period considered, Paul Lorenzen is a professional
mathematician whose research in abstract algebra leads to new insights
on the way algebraic objects are given to us. His genuine interest in
mathematical logic certainly triggers these insights; conversely, they also influence his conception of a consistency proof for elementary number theory as the construction of an embedding of a preordered set into a $\sigma$-complete pseudocomplemented semilattice \citep[see][]{coquandneuwirth19}.

The broader framework in which this happens is that of lattice
theory getting to the core of mathematical research, after the work of
the precursors Dedekind, Schröder, and Skolem. The introduction of
lattice theory into ideal theory is credited to Wolfgang
\citet{krull26}, whose assistant Lorenzen is from 1939
until his habilitation in 1946.

Herbert \citet[page~146]{mehrtens79}, in considering Steinitz’ work, describes the development of
modern algebra in the following terms:
\begin{enumerate*}
\item the algebraic structure of concrete objects is unveiled;
\item it is given an abstract axiomatic formulation;
\item the abstract concepts are subjected to a detailed investigation of their structure.
\end{enumerate*}
It is at this third stage that lattice theory becomes instrumental.

\begin{EngCadre}
  \noindent After the first step, where the algebraic fundamental structures had emerged in the investigation of numbers, equations, and functions, and the second phase, where the abstract-axiomatic formulation had been achieved, now the abstract concepts became the base of a careful structural investigation.

  Seen in this way, the third phase in the development of abstract algebra, the ``modern age'', starts with Steinitz. Dedekind's work belongs to the first, concrete phase; it is part of this developmental stage by virtue of containing the germ of the further development. This is also attested by the fact that, in one detail, he himself carried out all three phases -- in the investigation of the dual group [the lattice]. For abstract algebra, the third phase is marked by E. Noether, E. Artin, H. Hasse, W. Krull, and others, whose most productive time were the 20s and 30s.
  %
  %
\end{EngCadre}

This article proposes to follow Lorenzen's work in algebra by focusing
on a specific theorem, Krull's Fundamentalsatz for integral
domains. The theorem is introduced on the basis of excerpts of
articles by \citet{krull30,krull32,krull36} with a short comment. Thus
the starting point of Lorenzen's research is described. Then excerpts
of each of his four articles on the subject
\citep{Lor1939,Lor1950,Lor1952,Lor1953} are given and commented on, in
chronological order, so as to follow his enterprise of deepening the
understanding of the Fundamentalsatz. In his letter to Krull dated 6~June 1944 (page~\pageref{19440606}), Lorenzen describes his goal as follows.

\begin{EngCadre}
  And\label{19440606t} yet the simplification and clarification of the proof methods is
  precisely the proper main goal of my work. I have not tried to
  generalise the theorems of multiplicative ideal theory at any price,
  even at the price of complicating -- on the contrary, I rather care only to discern the
  basic ideas of the proof methods in their ultimate simplicity. If e.g.\ I replace the concept of valuation by
  ``homomorphism of a semilattice into a linearly preordered set'',
  then I see therein a conceptual simplification instead of a
  complication. For the introduction of the concept of valuation
  (e.g.\ the prima facie arbitrary triangle inequality) is justified
  only by the subsequent success, whereas the concept of
  homomorphism bears its justification in itself. I would say that the
  homomorphism into a linear preorder is the ``pure concept'' that
  underlies the concept of valuation. And the underlying, pure concept
  seems to me to be undeniably the simpler one.

  If I am mistaken in this point, then I urge you to tell
  it to me, because it has to date always been my endeavour in my whole mathematical
  work to bring
  to light these underlying, pure concepts themselves in their simple
  and transparent clarity.

  In this tendency toward conceptual clarification, the paper [\citealt{Lor1950}] also differs from my dissertation [\citealt{Lor1939}]. That this clarification, this understanding of the inner significance, as you call it once, is
  the most urgent duty has really become evident to me only in the last
  years.
\end{EngCadre}

This article allows for several different readings, according to the objectives of the readers and to their German language skills. Krull's formulation and analysis of his Fundamentalsatz and Lorenzen's reshaping of it are presented as a process that takes place in seven articles, of which the arguments relevant for our history are presented in a comprehensive way, sometimes accompanied by translations of selected paragraphs; then excerpts of the original article, chosen so as to allow for a direct and coherent reading and comparison with my presentation, are reproduced. The correspondence of Lorenzen with Hasse, Krull, and Aubert, together with a few documents, provides an insight into their respective conceptions as well as into the human and administrative component of this history.

\section{\texorpdfstring{\citealt{krull30}}{Krull 1930}: a first attempt at introducing valuations for an integral domain.}
\label{sec:1930}

\selectlanguage{english}

\citet{krull30} describes the goal of a  Fundamentalsatz for integral domains: it should  generalise the ``main theorem of elementary number theory'' to the case of an integral domain~$\mathfrak O$ with field of fractions~$\mathfrak K$. The theorem to be generalised answers the question of when a nonzero element~$\alpha$ of~$\mathfrak K=\mathbb Q$ is divisible by another nonzero element~$\beta$, i.e.\ of when the quotient~$\alpha\cdot\beta^{-1}$ lies in~$\mathfrak O=\mathbb Z$. It does so in terms of the difference between the exponents with which each prime number~$p$ appears in the numerator and the denominator of~$\alpha$ (the \emph{value of\/~$\alpha$ at}~$p$) and of~$\beta$: the value of~$\alpha$ must be greater than or equal to the value of~$\beta$ at each~$p$.

Krull proposes an abstract version of these differences of exponents by
defining valuations as functions~$\mathfrak K\to\mathbb R\cup\{+\infty\}$ satisfying three properties~1, 2~and~3, and investigates when valuations decide divisibility.

\subsection{Excerpt.}

Our excerpt is from pages~531--532 of ``Idealtheorie in unendlichen algebraischen Zahlkörpern~II'', \emph{Mathematische Zeitschrift}~31.

\begin{Cadre}
  \noindent\hfil  §~1.\nopagebreak\smallskip

  \noindent\hfil\textbf{Bewertung; Begriff der Zahlentheorie.}\nopagebreak\medskip

Das Grundproblem und der Hauptsatz der elementaren Zahlentheorie kann
folgendermaßen formuliert werden:\smallskip

\so{Frage:} Wann ist das Element~$\alpha$ aus dem
Körper~$\mathfrak K_r$ der rationalen Zahlen hinsichtlich der
Ordnung~$\mathfrak O_r$ der ganzen rationalen Zahlen durch das
gleichfalls zu~$\mathfrak K_r$ gehörige Element~$\beta$ teilbar, d.\
h.\ wann liegt der Quotient~$\alpha\cdot\beta^{-1}$
in~$\mathfrak O_r$?\smallskip

\so{Antwort:} Man ordne jeder Primzahl~$p$ aus~$\mathfrak O_r$ eine
„Stelle“~$P$ zu und verstehe für beliebiges $\alpha\ne0$
aus~$\mathfrak K_r$ unter dem Wert von~$\alpha$ in~$P$ die Differenz
der Exponenten, mit denen die zu~$P$ gehörige Primzahl~$p$ bei einer
Quotientendarstellung von~$\alpha$ durch ganze Zahlen in Zähler und
Nenner aufgeht.  Dann gilt der Satz: $\alpha$~ist dann und nur dann
durch~$\beta$ teilbar, wenn der Wert von~$\alpha$ in bezug auf eine
beliebige Stelle niemals kleiner ist als der Wert von~$\beta$ in bezug
auf die gleiche Stelle.\smallskip

Es bedeute von jetzt ab~$\mathfrak K$ einen beliebigen Körper, unter
einer \emph{Ordnung}~$\mathfrak O$ möge ein Unterring
von~$\mathfrak K$ verstanden werden, der das Einheitselement enthält,
von~$\mathfrak K$ selbst verschieden ist und die Eigenschaft besitzt,
daß jedes Element aus~$\mathfrak K$ als Quotient von Elementen
aus~$\mathfrak O$ geschrieben werden kann, daß also~$\mathfrak K$ den
Quotientenkörper von~$\mathfrak O$ darstellt.  Wir nennen das
Körperelement~$\alpha$ durch das Körperelement~$\beta$
\emph{hinsichtlich\/~$\mathfrak O$ teilbar}, wenn
$\alpha\cdot\beta^{-1}$ zu~$\mathfrak O$ gehört und fragen:\smallskip

\emph{Welche Eigenschaften muß\/ $\mathfrak O$~besitzen, damit über die
  Teilbarkeitsverhältnisse der Elemente von\/~$\mathfrak K$
  hinsichtlich\/~$\mathfrak O$ nach dem Vorbild der elementaren
  Zahlentheorie durch Einführung von „Stellen“ und zugehörigen
  „Bewertungen“ entschieden werden kann?}\smallskip

Um die aufgeworfene Frage zu präzisieren, haben wir vor allem genau
festzulegen, was unter einer Stelle mit zugehöriger Bewertung
verstanden werden soll. Von einer \emph{Bewertung von}~$\mathfrak K$
sprechen wir, wenn jedem Körperelement $\alpha\ne0$ eindeutig eine
reelle Zahl~$r$, der Null das Symbol~$+\infty$ als \emph{Wert}
zugeordnet ist, und wenn dabei die folgenden Bedingungen erfüllt sind:
1.~\emph{Der Wert von\/~$\alpha\cdot\beta$ ist gleich der Summe der
  Werte von\/ $\alpha$~und\/~$\beta$}.  2.~\emph{Der Wert
  von\/~$\alpha+\beta$ ist nie kleiner als der kleinere der Werte
  von\/ $\alpha$~und\/~$\beta$}.  3.~\emph{Nicht alle Elemente haben
  den Wert\/~$0$}.
\end{Cadre}

\section{\texorpdfstring{\citealt{krull32}}{Krull 1932}: the Fundamentalsatz for integral domains.}
\label{sec:1932}

\citet{krull32} discovers that the concept of valuation in
\citealt{krull30} is too narrow for deciding divisibility. He starts by
redefining a valuation as a surjective
function~$w\colon\mathfrak K^*\to\varGamma$, with $\varGamma$ an
abelian linearly preordered group, not necessarily archimedean,
such that
\begin{align}
  w(a\cdot b)&=w(a)+w(b)\text,\label[equality]{w(ab)}\\
  w(a+b)&\geq\min(w(a),w(b))\text.\label[inequality]{w(a+b)}
\end{align}
The
former concept corresponds to subgroups~$\varGamma$ of the additive
group~$\mathbb R$ and is from now on called ``special valuation''. A valuation
defines the \emph{valuation ring}~$\mathfrak B$ formed by~$0$ together with the elements where it is nonnegative.

E.g.\ if $\mathfrak K=\mathbb Q$ and $p$~is a prime number, then the
value of $a\in\mathbb Q^*$ at~$p$ may be defined by writing~$a=p^na'$
with $n\in\mathbb Z$ and $p$~dividing neither the numerator nor the
denominator of~$a'$, and by letting $\varGamma=\mathbb Z$ and
$w(a)=n$. Then $\mathfrak B$~is the ring of those fractions that can be written with a denominator not divisible by~$p$.

Krull starts by proving two characterisations of a valuation ring~$\mathfrak B$.
\begin{enumerate}[wide,nosep]
\item It is a subring of the field of fractions~$\mathfrak K$
  that preorders it linearly by divisibility: if $a\in\mathfrak K$, then
  $a\in\mathfrak B$ or $a^{-1}\in\mathfrak B$. It is here that he
  needs the freedom of constructing an arbitrary abelian linearly
  preordered group~$\varGamma$.
\item It is a local ring (i.e.\ its nonunits form an ideal) such that every
  ring~$\mathfrak B'$ with
  $\mathfrak B\subset\mathfrak B'\subseteq\mathfrak K$ contains the
  inverse of a nonunit of~$\mathfrak B$. Note that this property gives rise to nonconstructive instantiations of such an inverse of nonunit in arguments by contradiction.
\end{enumerate}

Then he focuses on the connection with integral closedness.

\begin{EngCadre}
  As usual, the element~$p$ is to be called \emph{integrally dependent on the ring}~$\mathfrak R$ if it satisfies an equation $p^n+a_1p^{n-1}+\dots+a_n=0$ with coefficients~$a_i$ in~$\mathfrak R$; $\mathfrak R$ will be termed \emph{integrally closed} if every element in the field of fractions~$\mathfrak K$ that is integrally dependent on~$\mathfrak R$ already belongs to~$\mathfrak R$.
\end{EngCadre}

He observes on the way that
valuation rings are integrally closed. This proves one direction of
the Fundamentalsatz for integral domains, viz.
\begin{EngCadre}
  \noindent\textbf{Theorem 7. }\emph{A proper subring\/~$\mathfrak R$
  of\/~$\mathfrak K$ may be represented as intersection of (finitely or infinitely many)
valuation rings if and only if it is integrally closed.}
\end{EngCadre}
The other direction follows from the fact that if
$a\in\mathfrak K\setminus\mathfrak R$ with $\mathfrak R$~integrally
closed, then $a^{-1}$~is a nonunit of~$\mathfrak R[a^{-1}]$: if
$a$~does not belong to~$\mathfrak R$, then it is not integrally
dependent on~$\mathfrak R$ either, and this may be expressed as the
absence of a relation of the form
$a^{-1}(a_1+a_2a^{-1}+\dots+a_n(a^{-1})^{n-1})=1$.  Then Zorn's lemma
provides a maximal
subring~$\mathfrak B$ of~$\mathfrak K$
containing~$\mathfrak R[a^{-1}]$ but not $a$; $\mathfrak B$ turns out
to be a valuation ring by the second characterisation of a valuation ring.

\subsection{Excerpt.}

Our excerpt is from pages~163--165 and 168--169 of ``Allgemeine Bewertungstheorie'', \emph{Journal für die reine und angewandte Mathematik}~167.

\begin{Cadre}
  \noindent\hfil\textbf{§~1. Definition und Grundeigenschaften der allgemeinen
    Bewertungen.}\nopagebreak\medskip

$\mathfrak K$~bedeutet stets einen Körper, $\varGamma$~eine Abelsche
Gruppe mit der Addition als Verknüpfungsrelation und der~$0$ als
Einheitselement. [\dots] 
$\varGamma$ bzw.~$\mathfrak K$ soll \emph{linear geordnet}
heißen, wenn zwischen den Elementen von~$\varGamma$ bzw.~$\mathfrak K$
eine Ordnungsbeziehung definiert ist mit den üblichen
charakteristischen Eigenschaften:

\dotfill\smallskip

Von einer \emph{Bewertung\/~$B$ von\/~$\mathfrak K$ mit der
  Wertgruppe\/~$\varGamma$} soll gesprochen werden, wenn jedem
Element~$a\ne0$ aus~$\mathfrak K$ eindeutig ein Element~$\alpha=w(a)$
aus~$\varGamma$ als \emph{Wert} zugeordnet ist, und wenn dabei folgende
Bedingungen erfüllt sind: 1.~$w(a\cdot b)=w(a)+w(b)$.
--~2.~$w(a+b)\geq\min(w(a),w(b))$. --~3.~Zu jedem~$\alpha$ aus~$\varGamma$
gibt es ein~$a$ aus~$\mathfrak K$, so daß $\alpha=w(a)$.

\dotfill\smallskip

Die Gesamtheit der Elemente~$a$, die in einer festen Bewertung~$B$
von~$\mathfrak K$ nichtnegative Werte haben, bildet zusammen mit dem
Nullelement einen echten Unterring~$\mathfrak B$ von~$\mathfrak K$,
den „\emph{zu\/~$B$ gehörigen Bewertungsring}“.

\dotfill\smallskip

\textbf{Satz 1. }\emph{Ein echter Unterring\/~$\mathfrak B$
  von\/~$\mathfrak K$ ist dann und nur dann Bewertungsring, wenn\/
  $\mathfrak B$~den Körper\/~$\mathfrak K$ zum Quotientenkörper hat
  und wenn in\/~$\mathfrak B$ von zwei Elementen\/~$a_1$~und\/~$a_2$
  stets (mindestens) eines durch das andre teilbar ist.}

\dotfill\medskip

\noindent\hfil\textbf{§~3. Kennzeichnung der ganz abgeschlossenen
  Ringe.}\nopagebreak\medskip

Das Element~$p$ soll wie üblich \emph{vom Ringe\/~$\mathfrak R$ ganz
  abhängig} heißen, wenn es einer Gleichung
$p^n+a_1p^{n-1}+\cdots+a_n=0$ mit Koeffizienten~$a_i$
aus~$\mathfrak R$ genügt; $\mathfrak R$ wird \emph{ganz abgeschlossen}
genannt, wenn jedes von~$\mathfrak R$ ganz abhängige Element aus dem
Quotientenkörper~$\mathfrak K$ bereits zu~$\mathfrak R$ selbst gehört.
Bezeichnen wir für beliebiges~$a$ mit $\mathfrak R[a]$ stets den Ring
aller Polynome in~$a$ mit Koeffizienten aus~$\mathfrak R$, so können
wir die von~$\mathfrak R$ ganz abhängigen Elemente auch folgendermaßen
charakterisieren:\smallskip

\emph{Das Element\/~$p$ hängt von\/~$\mathfrak R$ dann und nur dann ganz
  ab, wenn\/ $p^{-1}$ in~$\mathfrak R[p^{-1}]$ Einheit ist.}\smallskip

In der Tat, $p^n-a_1p^{n-1}-\cdots-a_n=0$ ist vollkommen gleichwertig
mit $1=p^{-1}\cdot(a_1+a_2p^{-1}+\cdots+a_n(p^{-1})^{n-1})$ und das
Bestehen einer Gleichung der letzteren Form mit Koeffizienten~$a_i$
aus~$\mathfrak R$ ist notwendig und hinreichend dafür, daß $p^{-1}$
in~$\mathfrak R[p^{-1}]$ Einheit ist. -- Aus unserm Kriterium für ganz
abhängige Elemente folgt sofort:\smallskip

\emph{Ist\/~$a$ Nichteinheit in\/~$\mathfrak R$, so kann\/ $a^{-1}$
  niemals von\/~$\mathfrak R$ ganz abhängen.}\smallskip

Wäre nämlich $a^{-1}$ von~$\mathfrak R$ ganz abhängig, so müßte
$a=(a^{-1})^{-1}$ in~$\mathfrak R[a]=\mathfrak R$ gegen Voraussetzung
Einheit sein.\smallskip

\textbf{Satz 5. }\emph{Ein echter Unterring\/~$\mathfrak B$
  von\/~$\mathfrak K$ ist dann und nur dann Bewertungsring, wenn
  in\/~$\mathfrak B$ die Gesamtheit der Nichteinheiten ein Ideal
  bildet, und wenn in\/~$\mathfrak K$ jeder echte Oberring
  von\/~$\mathfrak B$ mindestens ein Reziprokes einer Nichteinheit
  von\/~$\mathfrak B$ enthält.}\smallskip

Daß jeder Bewertungsring die in Satz 5 angegebenen Eigenschaften
besitzt, ist klar nach §~1. -- Ist umgekehrt~$\mathfrak B$ irgendein
Ring mit diesen Eigenschaften, so muß zunächst $\mathfrak B$~ganz
abgeschlossen sein.  Denn der Ring~$\mathfrak B^*$ aller
von~$\mathfrak B$ ganz abhängigen Elemente aus~$\mathfrak K$ muß
mit~$\mathfrak B$ zusammenfallen, weil er sicher kein einziges
Reziprokes einer Nichteinheit von~$\mathfrak B$ enthält.  Es seien
ferner $a_1$~und~$a_2$ zwei beliebige Elemente aus~$\mathfrak B$, und
es sei etwa $a_2$~nicht in~$\mathfrak B$ durch~$a_1$ teilbar, also
$a_2\cdot a_1^{-1}=p$ nicht Element von~$\mathfrak B$. Bilden wir
dann~$\mathfrak B[p]$, so muß in~$\mathfrak B[p]$ das
Reziproke~$a^{-1}$ einer gewissen Nichteinheit~$a$ aus~$\mathfrak B$
auftreten, es muß also eine Gleichung $a^{-1}=a_0+a_1p+\cdots+a_np^n$
mit Koeffizienten~$a_i$ aus~$\mathfrak B$ gelten.  Diese Gleichung kann
aber in die Form
$(a_0\cdot a-1)\cdot(p^{-1})^n+a_1\cdot
a\cdot(p^{-1})^{n-1}+\cdots+a_n\cdot a=0$ gebracht werden, und hier
muß wegen der besonderen Eigenschaften von~$\mathfrak B$ der
Koeffizient~$a_0\cdot a-1$ als Differenz einer Nichteinheit und einer
Einheit selbst Einheit sein. Es ist also $p^{-1}$ von~$\mathfrak B$ ganz
abhängig und somit nach dem bereits Bewiesenen Element
von~$\mathfrak B$. --~Wir haben jetzt gezeigt: Sind $a_1$~und~$a_2$
zwei beliebige Elemente aus~$\mathfrak B$, so ist in~$\mathfrak B$
entweder $a_1$~durch~$a_2$ oder $a_2$~durch~$a_1$ teilbar.  Nach
Satz~1 muß~$\mathfrak B$ daher Bewertungsring sein. -- Mit Hilfe von
Satz~5 beweisen wir weiter:\smallskip

\textbf{Satz 6. } \emph{Zu jedem echten Unterring\/~$\mathfrak R$
  von\/~$\mathfrak K$ gibt es (mindestens) einen
  Bewertungsring\/~$\mathfrak B$, der\/ $\mathfrak R$~enthält.}\smallskip

Es sei~$a$ eine beliebige Nichteinheit aus~$\mathfrak R$. Dann
existiert, wie aus trivialen Wohlordnungsschlüssen zu ersehen,
in~$\mathfrak K$ mindestens ein Ring~$\mathfrak B$, der~$\mathfrak R$,
aber nicht $a^{-1}$~enthält, und der außerdem die Eigenschaft hat, daß
$a^{-1}$~in jedem echten Oberring von~$\mathfrak B$
vorkommt. $\mathfrak B$~wird nach Satz~5 Bewertungsring sein, wenn die
Gesamtheit der Nichteinheiten in~$\mathfrak B$ ein Ideal bildet.  Es
sei nun $\mathfrak p$ irgendein Primidealteiler von~$a$
in~$\mathfrak B$, \ $\mathfrak B_{\mathfrak p}$ bedeute den Ring aller
der Elemente, die sich als Quotienten von Elementen aus~$\mathfrak B$
mit durch~$\mathfrak p$ unteilbarem Nenner schreiben lassen.  Dann
enthält~$\mathfrak B_{\mathfrak p}$ das Element~$a^{-1}$ nicht, und es
muß daher $\mathfrak B_{\mathfrak p}=\mathfrak B$ sein. Das ist
aber nur möglich, wenn $\mathfrak p$ gerade aus \emph{allen}
Nichteinheiten von~$\mathfrak B$ besteht.\smallskip

\textbf{Satz 7. }\emph{Ein echter Unterring\/~$\mathfrak R$
  von\/~$\mathfrak K$ läßt sich dann und nur dann als Durchschnitt von
  (endlich oder unendlich vielen) Bewertungsringen darstellen, wenn er
  ganz abgeschlossen ist.}\smallskip

Daß alle Bewertungsringe und damit auch alle Durchschnitte von
Bewertungsringen ganz abgeschlossen sind, wurde beim Beweise von
Satz~5 mitbewiesen.  Es sei jetzt umgekehrt $\mathfrak R$~irgendein
ganz abgeschlossener Ring, $a$~ein nicht in~$\mathfrak R$ vorkommendes
Element aus~$\mathfrak K$.  Dann haben wir nur zu zeigen, daß
mindestens ein Bewertungsring~$\mathfrak B$ existiert, der
zwar~$\mathfrak R$, aber nicht $a$~enthält. Wir bilden
$\mathfrak R[a^{-1}]$; in diesem Ringe kann~$a$ nicht vorkommen, denn
andernfalls wäre~$a^{-1}$ in~$\mathfrak R[a^{-1}]$ Einheit, und es
müßte daher~$a$ von~$\mathfrak R$ ganz abhängen und damit gegen
Voraussetzung in~$\mathfrak R$ enthalten sein.  Da also~$a^{-1}$
in~$\mathfrak R[a^{-1}]$ Nichteinheit ist, gibt es nach dem Beweise
von Satz~6 einen Bewertungsring~$\mathfrak B$, der Obermenge
von~$\mathfrak R[a^{-1}]$ also erst recht auch von~$\mathfrak R$ ist,
in dem aber $a$~nicht vorkommt.  Damit ist schon alles bewiesen. -- Um
die Bedeutung von Satz~7 noch klarer hervortreten zu lassen, führen
wir den Begriff der \emph{Hauptordnung} ein.  Der Ring~$\mathfrak R$
soll Hauptordnung heißen, wenn in~$\mathfrak R$ über die
Teilbarkeitsverhältnisse der Elemente durch Einführung von Bewertungen
entschieden werden kann, d.\ h.\ wenn sich (endlich oder unendlich
viele) Bewertungen~$\mathfrak B_\tau$ des
Quotientenkörpers~$\mathfrak K$ so definieren lassen, daß der folgende
Satz gilt: In~$\mathfrak R$ ist das Element~$a$ dann und nur dann
durch das Element~$b$ teilbar, wenn $a$~in keiner der
Bewertungen~$\mathfrak B_\tau$ einen kleineren Wert besitzt als~$b$.

Man sieht ohne Schwierigkeit, daß der Ring~$\mathfrak R$ dann und nur
dann Hauptordnung ist, wenn er als Durchschnitt von (endlich oder
unendlich vielen) Bewertungsringen~$\mathfrak B_\tau$ dargestellt
werden kann.  Aus Satz 7 ergibt sich daher:\smallskip

\textbf{Satz 7*. }\emph{Ein echter Unterring\/~$\mathfrak R$
  von\/~$\mathfrak K$ ist dann und nur dann Hauptordnung, wenn er ganz
  abgeschlossen ist.}\smallskip

Wie man in plausibler und naheliegender Weise zu dem Begriff der
Hauptordnung gelangt, habe ich bereits früher ausführlich
auseinandergesetzt. Doch konnte ich dort die Hauptordnungen nicht so
einfach und befriedigend charakterisieren, wie es hier durch Satz~7*
geschehen ist.  Der Grund für diesen Mangel meiner früheren Arbeit lag
einfach darin, daß ich damals nur Bewertungen mit archimedischer
Wertgruppe in den Kreis der Betrachtung zog.  Satz~7* zeigt also, daß
die Einführung der allgemeinen Bewertungen nicht nur naheliegend,
sondern auch für die naturgemäße Behandlung mancher arithmetischer
Probleme schlechtweg notwendig ist.
\end{Cadre}

\section{\texorpdfstring{\citealt{krull36}}{Krull 1936}: the computational content of the Fundamentalsatz.}
\label{sec:krull36}

\label{method}In the introduction to his \citeyear{krull36}, \citeauthor{krull36} explains that the Fundamentalsatz is the basis for a proof method: if a theorem is to be proved for an integrally closed integral domain, then it mostly suffices to prove it in a valuation ring, in which the property of linearity grants very simple computational arguments.

He also explains the defect of the Fundamentalsatz: its proof does not
construct the valuation rings, and the example of polynomial rings
shows that this goal is out of reach.

\begin{EngCadre}
\noindent\textb{A.\enspace The valuation rings are characterised by the fact that in them, of two arbitrary elements, one is always divisible by the other. --\enspace B.\enspace Every integrally closed integral domain~$\mathfrak I$ may be represented as intersection of valuation rings.}

A and B form the basis of an important proof method: if a theorem on~$\mathfrak I$ is to be derived which one knows to hold for the intersection of arbitrarily many rings provided that it holds for each individual component ring, then, according to~B, the theorem concerned needs to be proved only for valuation rings, and in this latter case, on account of~A, usually a quite elementary computation leads to the goal.

On the other hand, B~serves as a connecting link for novel questions, whose fruitfulness is to be substantiated in this contribution. B~is a pure existence theorem. For the ``proof method'', this does not constitute a drawback; nevertheless it seems theoretically desirable to go beyond the bare standpoint of existence. This goal would be reached if one were in a position to represent in some constructive manner, for an arbitrary given~$\mathfrak I$, the set of all those valuation rings out of the field of fractions~$\mathfrak K$ that contain the integral domain~$\mathfrak I$. Yet even the consideration of quite elementary special cases, e.g.\ of polynomial rings, indicates that one is there confronted with an, at least so far, unsolvable task. -- On the other hand, it becomes apparent that, surprisingly, already the attempt to turn the given task into the ideal-theoretic leads to very curious connexions with certain penetrating investigations by H. Prüfer that refer chief{}ly to Kronecker.
\end{EngCadre}

\subsection{Excerpt.}

Here is the source of our translation, from page~546 of ``Beiträge zur Arithmetik kommutativer Integritätsbereiche
  I: Multiplikationsringe, ausgezeichnete Idealsysteme und
  Kroneckersche Funktionalringe'', \emph{Mathematische Zeitschrift}~41.

\begin{Cadre}
\noindent\so{A.~Die Bewertungsringe sind dadurch gekennzeichnet, daß in
ihnen von zwei beliebigen Elementen stets eines durch das
andere teilbar ist.~-- B.~Jeder ganz abgeschlossene Integritätsbereich~$\mathfrak I$  läßt sich als Durchschnitt von Bewertungsringen
darstellen.}

A. und B. bilden die Grundlage für eine wichtige Beweismethode:
Soll ein Satz über~$\mathfrak I$ abgeleitet werden, von dem man weiß, daß er für
den Durchschnitt beliebig vieler Ringe richtig ist, sobald er nur für jeden
einzelnen Komponentenring gilt, so braucht der betreffende Satz nach B.
nur für Bewertungsringe bewiesen zu werden, und in diesem letzteren
Fall führt wegen A. meistens eine ganz einfache Rechnung zum Ziel.

Andererseits dient B. als Anknüpfungspunkt für neuartige Fragestellungen, deren Fruchtbarkeit in diesem Beitrag nachgewiesen werden
soll. B. ist ein reines Existenztheorem. Für die „Beweismethode“ bedeutet das keinen Nachteil; trotzdem scheint es theoretisch wünschenswert, über den bloßen Existenzstandpunkt hinauszukommen. Dieses Ziel
wäre erreicht, wenn man imstande wäre, bei beliebig gegebenem~$\mathfrak I$ die
Menge aller der Bewertungsringe aus dem Quotientenkörper~$\mathfrak K$, die den
Integritätsbereich~$\mathfrak I$ umfassen, irgendwie konstruktiv darzustellen. Doch
lehrt schon die Betrachtung ganz einfacher Spezialfälle, z.\ B. der Polynomringe, daß man da vor einer wenigstens vorläufig unlösbaren Aufgabe
steht.~-- Dagegen zeigt es sich überraschenderweise, daß bereits der Versuch, die gestellte Aufgabe ins Idealtheoretische zu wenden, zu sehr
merkwürdigen Zusammenhängen mit gewissen tiefeindringenden, vorwiegend
an Kronecker anknüpfenden Untersuchungen von H. Prüfer führt.
\end{Cadre}

\section{\texorpdfstring{\citealt{Lor1939}}{Lorenzen 1939}: the Fundamentalsatz for preordered cancellative monoids.}
\label{sec:1939}

\citealt{Lor1939} is an ``abstract foundation of multiplicative ideal
theory'' on the theory of preordered cancellative monoids. Lorenzen’s
letter to Krull dated 22 March 1938 (page~\pageref{19380322}) provides
a dictionary that translates the concepts of ring theory into a
framework that takes into account only the multiplicative structure of
an integral domain~$\mathfrak I$: its nonzero elements form a
cancellative monoid~$\mathfrak g$, and the nonzero elements of its
field of fractions form what is today called the Grothendieck
group~$\mathfrak G$ of~$\mathfrak g$. Divisibility is the equivariant
preorder~$\preccurlyeq$ on~$\mathfrak G$ whose positive cone
is~$\mathfrak g$: $a\preccurlyeq b\iff\mathfrak g\ni b\cdot a^{-1}$.

Lorenzen’s analysis of Krull's valuation theory proceeds in three steps: elaborating a theory of ideals that fits preordered monoids
; making a valuation-theoretic analysis of lattice-preordered groups
; establishing a transfer of valuations between a preordered
monoid and the positive cone of the Grothendieck lattice-preordered group
of its system of ideals
.

\subsection{Systems of ideals.}
\label{sec:systems-ideals}

\citet{Pru1932} proposes to define a system of $r$-ideals as an operation
$\afini\mapsto\afini_r$ on finite subsets of the field of fractions of
an integral domain. Lorenzen adapts this definition into an operation
on finite or infinite bounded-from-below subsets~$\mathfrak a$ of the group~$\mathfrak G$ (i.e.\
subsets such that there is a~$c\in\mathfrak G$ with
$\mathfrak g\supseteq c\mathfrak a$; Lorenzen's notation for this is
$\mathfrak a^{-1}\ne0$), taking values in
subsets of~$\mathfrak G$, and satisfying four conditions:
\begin{enumerate*}
\item $\mathfrak a_r\supseteq\mathfrak a$;
\item $\mathfrak b_r\supseteq\mathfrak a\implies\mathfrak b_r\supseteq\mathfrak a_r$;\label{si:2}
\item $\{a\}_r=a\cdot\mathfrak g$;\label{si:3}
\item $a\cdot\mathfrak a_r=(a\cdot\mathfrak a)_r$.
\end{enumerate*}
The variable letter~$r$ stands for variable systems of ideals.

E.g.\ the system of $t$-ideals\label{t-ideal}, used in \cref{sec:ideals-latt-order}, is defined by
\[\mathfrak a_t=\bigl\{\,b\in\mathfrak G\mid\exists a_1,\dots,a_n\in\mathfrak a\enskip\forall a\in\mathfrak G\quad a\preccurlyeq a_1,\dots,a\preccurlyeq a_n\implies a\preccurlyeq b\,\bigr\}\]
for finite or infinite bounded-from-below sets~$\mathfrak a$ as the set of those elements~$b$ that are greater than or equal to every lower bound~$a$ of some finite subset $\{a_1,\dots,a_n\}$ of~$\mathfrak a$.

A system of ideals has the structure of a meet-semilattice-preordered
monoid (a \emph{meet-monoid} for short) for the
multiplication
$\mathfrak a_r\cdot\mathfrak b_r=(\mathfrak a\cdot\mathfrak b)_r$ and
the meet operation
$\mathfrak a_r\land\mathfrak b_r=(\mathfrak a\cup\mathfrak
b)_r$.  \Cref{si:3} above means that the group~$\mathfrak G$ may be embedded into this meet-monoid by identifying~$a$ with~$a\cdot\mathfrak g$.

The system of Dedekind ideals corresponds to the operation that associates to~$\mathfrak a$
the set~$\mathfrak a_{d}$ of nonzero elements of the
module generated by~$\mathfrak a$ in the integral domain~$\mathfrak I$
hidden behind~$\mathfrak g$.  The definition of this operation therefore
takes into account the additive structure of~$\mathfrak I$, but
once the definition is done, it will be forgotten.  Let us show how
this works for the three following concepts:
\begin{enumerate*}
\item valuation ring;
\item integral closedness;
\item integral dependence.
\end{enumerate*}

\begin{enumerate}[wide,nosep]
\item Lorenzen separately analyses \cref{w(ab)} and \cref{w(a+b)} of
  the definition of a valuation given on \cpageref{w(ab)}.
  \begin{itemize}[wide,nosep]
  \item He says that a supermonoid~$\mathfrak B$ of~$\mathfrak g$
    $(\mathfrak G\supseteq\mathfrak B\supseteq\mathfrak g)$ is
    \emph{linear} if for $a\in\mathfrak G$ either $a\in\mathfrak B$ or
    $a^{-1}\in\mathfrak B$, i.e.\ if the preorder on~$\mathfrak G$
    having~$\mathfrak B$ as its positive cone is linear.
  \item \cref{w(a+b)} expresses that a valuation ring is closed under
    addition. He observes that this is equivalent to
    $\mathfrak B\supseteq\afini\implies\mathfrak
    B\supseteq\afini_{d}$.
  \end{itemize}
\item \citet{Pru1932} shows that an integral domain is integrally
  closed if
  \[\afini_{d}\supseteq a\cdot\afini_{d}\quad\text{implies}\quad\mathfrak g\supseteq a\cdot\mathfrak g\text,\] i.e.\
  $1\preccurlyeq a$. This means that $\afini_{d}$ can be
  cancelled in the containment
  $\afini_{d}\supseteq a\cdot\afini_{d}$.
\item He further shows that an element~$a$ is integrally dependent
  on~$\afini_{d}$ if and only if there is
  $\cfini$ such that
  $\afini_{d}\cdot\cfini_{d}\supseteq
  a\cdot\cfini_{d}$.
\end{enumerate}

On the basis of these observations, Lorenzen proposes the following definitions in the multiplicative setting of a system of $r$-ideals for the cancellative monoid~$\mathfrak g$.
\begin{enumerate}[wide,nosep]
\item The multiplicative counterpart of a valuation
  ring is a linear \emph{$r$-\hspace{0pt}supermonoid of}~$\mathfrak g$,
  i.e.\ a linear supermonoid~$\mathfrak B$ of~$\mathfrak g$ such that
  $\mathfrak B\supseteq\afini\implies\mathfrak
  B\supseteq\afini_r$.
\item The monoid~$\mathfrak g$ is $r$-\emph{closed} if
  $\afini_r\supseteq a\cdot\afini_r\implies1\preccurlyeq a$.
\item The $r_{a}$-operation is given by setting $\afini_{r_{a}}$ equal to
  \[\bigl\{\,a\in\mathfrak G\mid\exists\cfini\enskip\afini_r\cdot\cfini_r\supseteq a\cdot\cfini_r\,\bigr\}\text,\]
  i.e.\ to the set of elements that are \emph{$r$-dependent} on~$\afini$. Note that the $r_{a}$-operation
  satisfies 
  \cref{si:3} of the definition of a system of
  ideals 
  if and only if $\mathfrak g$~is $r$-closed.
\end{enumerate}

Lorenzen observes that the arguments of
\citet[page~18]{Pru1932} show that the cancellativity of the
meet-monoid of $r$-ideals may be forced by passing to the meet-monoid of
$r_{a}$-ideals. One may then construct the Grothendieck
lattice-preordered group~$\mathfrak h^*$ of this system of
$r_{a}$-ideals (see \citealt[§~3A]{coquandlombardineuwirth19} for this construction). Its positive cone~$\mathfrak h$ is formed by the
fractions~$\frac{\afini_{r_{a}}}{\bfini_{r_{a}}}$ such that
$\bfini_{r_{a}}\supseteq\afini_{r_{a}}$. The
group~$\mathfrak G$ may be embedded into~$\mathfrak h^*$ by
identifying~$a$ with~$\frac{a\cdot\mathfrak g}{\mathfrak g}$.

\subsection{Ideals in a lattice-preordered group.}
\label{sec:ideals-latt-order}

For developing a ``valuation theory'' of~$\mathfrak g$ when $\mathfrak G$~is a
lattice-preordered group (an \emph{$\ell$-group} for short), Lorenzen
uses the $t$-ideals, whose definition, thanks to the existence of finite meets, simplifies to
\[
  \mathfrak a_t=\bigl\{\,b\in\mathfrak G\mid a_1\land\dotsb\land
  a_n\preccurlyeq b\text{ for some }a_1,\dots,a_n\in\mathfrak
  a\,\bigr\}\text.
\]
\citet[page~470]{Cli1940} emphasises their similarity with the
semilattice-theoretic ideals.  A finite $t$-ideal $\afini_t$ is the
set $(a_1\wedge\dots\wedge a_n)\mathfrak g$ and may thus be thought of as $a_1\wedge\dots\wedge a_n$, and the set-theoretic clothing is
needed only for an infinite $t$-ideal, which is a colimit of finite
ones, in particular for a \emph{prime} $t$-ideal, i.e.\ a $t$-ideal~$\mathfrak p$
such that
$ab\in\mathfrak p\implies a\in\mathfrak p\text{ or }b\in\mathfrak p$.

Lorenzen shows that the linear $t$-\hspace{0pt}supermonoids of~$\mathfrak g$
are the localisations~$\mathfrak g_\mathfrak p$, i.e.\ the monoids
of fractions~$\frac ab$ with $a\in\mathfrak g$ and
$b\in\mathfrak g\setminus\mathfrak p$, where $\mathfrak p$~is a prime
$t$-ideal. In fact, the nonunits of~$\mathfrak g_\mathfrak p$ have the
form~$\frac ps$ with $p\in\mathfrak p$ and
$s\in\mathfrak g\setminus\mathfrak p$; if $\frac{p'}{s'}$~is another nonunit, then the meet $\frac ps\land\frac{p'}{s'}=\frac{s'p\land sp'}{ss'}$ is again a nonunit because $ss'\in\mathfrak g\setminus\mathfrak p$ by primeness, $s'p,sp'\in\mathfrak p$ by \cref{si:2,si:3} of the definition of a system of ideals, and $s'p\land sp'\in\mathfrak p$ by $t$-\hspace{0pt}idealness. But this shows that
$\mathfrak g_\mathfrak p$~is linear, as every
$a\in\mathfrak G$ satisfies $a=\frac{a_+}{a_-}$ with
$a_-,a_+\in\mathfrak g$ and $a_-\land a_+=1$, so that either~$a_-$
or~$a_+$ is a unit of~$\mathfrak g_\mathfrak p$ and either~$a$
or~$a^{-1}=\frac{a_-}{a_+}$ is in~$\mathfrak g_\mathfrak p$.

The Fundamentalsatz follows easily in this case, in the following two shapes.
\begin{EngCadre}
  \noindent\so{Theorem }11. \ \emph{Let\/~$N$ be a set of prime\/ $t$-ideals
  of\/~$\mathfrak g$ such that every nonunit
  of\/~$\mathfrak g$ lies in at least one\/~$\mathfrak p\in N$, then\/~$\mathfrak g$ is the intersection of the linear
  supermonoids\/~$\mathfrak g_{\mathfrak p}$\/ $(\mathfrak p\in N)$}
\[
  \mathfrak g=\bigcap_N\mathfrak g_{\mathfrak p}\text.
\]\vspace{-\belowdisplayshortskip}%
\end{EngCadre}
For if $a\in\mathfrak G\setminus\mathfrak g$, then $a_-$~is a nonunit
of~$\mathfrak g$, so that it lies in some~$\mathfrak p\in N$ and is a nonunit of~$\mathfrak g_\mathfrak p$; therefore $a=\frac{a_+}{a_-}\notin\mathfrak g_\mathfrak p$.

\begin{EngCadre}
  \noindent\so{Theorem }12. \ \emph{If\/ $N$ contains all maximal\/ prime $t$-ideals, then, for all\/ $t$-ideals\/~$\mathfrak a$
    of\/~$\mathfrak g$, the intersection representation\/
    $\mathfrak a=\bigcap\limits_N\mathfrak a\mathfrak g_{\mathfrak p}$
    holds.}
\end{EngCadre}
For\label{commenton12} if $a\notin\mathfrak a$, then $\mathfrak g$~contains properly the $t$-ideal~$a^{-1}\mathfrak a\cap\mathfrak g$. Zorn's lemma provides a maximal prime $t$-ideal~$\mathfrak p$ that contains it, so that its elements are nonunits of~$\mathfrak g_\mathfrak p$. Therefore $a^{-1}\mathfrak a\mathfrak g_\mathfrak p\not\ni 1$ and $\mathfrak a\mathfrak g_{\mathfrak p}\not\ni a$.

Krull's well-ordering argument is confined to the proof of this last theorem.

\subsection{Transfer to the system of \texorpdfstring{$t$}t-ideals.}
\label{sec:transfer-t-system}

To consider a linear $r$-\hspace{0pt}supermonoid~$\mathfrak B$ of~$\mathfrak g$
($\mathfrak G\supseteq\mathfrak B\supseteq\mathfrak g$) is to consider
a linear $t$-\hspace{0pt}supermonoid~$\mathfrak L$ of the positive cone~$\mathfrak h$ of the Grothendieck lattice-preordered group~$\mathfrak h^*$ ($\mathfrak h^*\supseteq\mathfrak L\supseteq\mathfrak h$).
\begin{enumerate}[wide,nosep]
\item If $\mathfrak L$~is given,\label{1939p545t} then
  $\mathfrak B=\mathfrak L\cap\mathfrak G$ (using the identification
  of~$a$ with~$\frac{a\cdot\mathfrak g}{\mathfrak g}$) is a linear
  supermonoid of~$\mathfrak g$. Furthermore, if
  $\mathfrak B\supseteq\afini$, then
  $\mathfrak L\supseteq\{\frac{a_1\cdot\mathfrak g}{\mathfrak
    g},\dots,\frac{a_n\cdot\mathfrak g}{\mathfrak g}\}$, so that by
  hypothesis
  $\mathfrak L\supseteq(\frac{a_1\cdot\mathfrak g}{\mathfrak
    g}\wedge\dots\wedge\frac{a_n\cdot\mathfrak g}{\mathfrak
    g})\mathfrak h=\frac{\afini_{r_{a}}}{\mathfrak g}\mathfrak h$ and
  $\mathfrak B\supseteq\afini_{r_{a}}\supseteq\afini_{r}$.
  At this point of his argument, Lorenzen compares it to
  \citealt[Satz~15]{krull36}.
\item\sloppy Conversely, if $\mathfrak B$~is given, then $\mathfrak B$~is even an $r_{a}$-supermonoid: if
  $\mathfrak B\supseteq\afini$ and
  $\afini_{r_{a}}\ni a$, consider~$\cfini$ such that
  $\afini_r\cdot\cfini_r\supseteq a\cdot\cfini_r$; one may suppose that
  $c_1$~is minimal w.r.t.\ the linear
  monoid~$\mathfrak B$ among $c_1,\dots,c_p$, i.e.\
  $\mathfrak B\supseteq\{1,\dots,c_1^{-1}c_p\}$, so that
  $\mathfrak B\supseteq\afini\cdot\{1,\dots,c_1^{-1}c_p\}$; by hypothesis, $\mathfrak B\supseteq\afini_r\cdot\{1,\dots,c_1^{-1}c_p\}_r\ni a$ 
  and we conclude that $\mathfrak
  B\supseteq\afini_{r_{a}}$. For every finite subset~$\afini$
  of~$\mathfrak G$, one may suppose that $a_1$~is minimal
  w.r.t.~$\mathfrak B$ among $a_1,\dots,a_n$, i.e.\
  $\mathfrak B\supseteq\{1,\dots,a_1^{-1}a_n\}$, and then
  $\mathfrak B\supseteq\{1,\dots,a_1^{-1}a_n\}_{r_{a}}$ and
  $a_1\mathfrak B=\afini_{r_{a}}\mathfrak B$. Consider
  $\mathfrak
  L=\bigl\{\frac{\afini_{r_{a}}}{\bfini_{r_{a}}}\mid\bfini_{r_{a}}\mathfrak
  B\supseteq\afini_{r_{a}}\mathfrak B\bigr\}$: as for all $a,b\in\mathfrak G$ either $b\mathfrak B\supseteq a\mathfrak B$ or vice versa, $\mathfrak L$~is a linear supermonoid of~$\mathfrak h$; it is a linear $t$-supermonoid of~$\mathfrak h$ because the meet of two elements of~$\mathfrak L$ is in~$\mathfrak L$: if $\bfini_{r_{a}}\mathfrak
  B\supseteq\afini_{r_{a}}\mathfrak B$ and $\dfini_{r_{a}}\mathfrak
  B\supseteq\cfini_{r_{a}}\mathfrak B$, then $\{b_1d_1,\dots,b_md_q\}_{r_{a}}\mathfrak B\supseteq\{a_1d_1,\dots,a_nd_q,b_1c_1,\dots,b_mc_p\}_{r_{a}}\mathfrak B$. Finally, note that $\mathfrak L\cap\mathfrak G=\mathfrak B$.
\end{enumerate}

This transfer proves at once
\begin{EngCadre}
\noindent\so{Theorem }14. \ \emph{Every\/ $r$-closed monoid is
  an intersection of linear\/ $r$-\hspace{0pt}supermonoids.}
\end{EngCadre}

\subsection{Excerpt.}

Our excerpt is from pages~535--537, 540, and 544--546 of ``Abstrakte Begründung der multiplikativen Idealtheorie'', \emph{Mathematische Zeitschrift}~45.

Note that in this article and in his 1938 letters to Krull, Lorenzen consistently makes a choice opposite to the now standard one: he considers $\mathfrak g$ as the cone of nonpositive elements of~$\mathfrak G$, so that his preorder~$\preccurlyeq$ satisfies $a\preccurlyeq b\iff a\cdot b^{-1}\in\mathfrak g$. Consistently, he writes $a\lor b$ and $\mathfrak a_r+\mathfrak b_r$ where we write $a\land b$ and $\mathfrak a_r\land\mathfrak b_r$, respectively. He probably does so to match the preorder relation~$\subseteq$ between subsets of~$\mathfrak G$.

\begin{Cadre}
  \noindent\hfil§~1.\medskip

  \noindent\hfil\textbf{Idealsysteme.}\medskip

Wir wollen eine Menge~$\mathfrak g$, für deren Elemente $a,b,c,\dots$ eine kommutative
und assoziative Multiplikation definiert ist, eine \emph{Halbgruppe} nennen, wenn

1.~aus $ac =  bc$ stets $a =  b$ folgt,

2.~$\mathfrak g$ ein Einselement~$1$ enthält.

Zu jeder Halbgruppe~$\mathfrak g$ gibt es eindeutig die
Gruppe~$\mathfrak G$, die aus den Quotienten~$\frac ab$
besteht. $\mathfrak G$~heißt die Quotientengruppe von~$\mathfrak
g$. So ist z.\ B. die Menge der von Null verschiedenen Elemente eines
Integritätsbereiches eine Halbgruppe, und die von Null verschiedenen
Elemente des Quotientenkörpers bilden die Quotientengruppe.

Viele Definitionen der Teilbarkeitstheorie in Integritätsbereichen
lassen sich unmittelbar auf Halbgruppen übertragen. Von den Elementen
von~$\mathfrak G$ nennen wir die Elemente von~$\mathfrak g$ die ganzen
Elemente, $b$~heißt ein Teiler von~$a$, und $a$~heißt Vielfaches
von~$b$, wenn $\frac ab$~ganz ist. Das Hauptideal~$(a)$ besteht aus
allen Vielfachen von~$a$. Ist eine Verwechslung ausgeschlossen, dann
wird das Hauptideal $(a)$ einfach mit~$a$ bezeichnet. Ein ganzes
Element~$a$ heißt Einheit, wenn $a^{-1}$ ganz ist.

Die Dedekindsche Definition der Ideale in Integritätsbereichen läßt
sich natürlich nicht übertragen, da in Halbgruppen ja keine Addition
erklärt zu sein braucht. Dagegen ist der Begriff der Idealsysteme, wie
er von Prüfer und Krull aufgestellt worden ist, fast unmittelbar auch
in Halbgruppen anzuwenden. Zur Formulierung der Definition wollen wir,
wenn $\mathfrak a, \mathfrak b$ beliebige Untermengen
von~$\mathfrak G$ sind, unter $\mathfrak a\mathfrak b$ die Menge der
Elemente~$ab$ mit $a\in\mathfrak a$, $b\in\mathfrak b$ und unter
$\mathfrak a : \mathfrak b$ die Menge der Elemente~$c$ mit
$c\mathfrak b\subseteq\mathfrak a$ verstehen.  Statt $(1): \mathfrak a$
schreiben wir $\mathfrak a^{-1}$.\smallskip

\so{De{f}inition 1}. \ Es seien zu jeder endlichen oder zu jeder
beliebigen Untermenge~$\mathfrak a$ von~$\mathfrak G$, für die
$\mathfrak a^{-1}\ne0$ ist, eine Untermenge~$\mathfrak a_r$
von~$\mathfrak G$ so definiert, daß die folgenden Bedingungen erfüllt
sind:\smallskip

1. $\mathfrak a\subseteq\mathfrak a_r$,

2. aus $\mathfrak a \subseteq\mathfrak b_r$ folgt
$\mathfrak a_r\subseteq\mathfrak b_r$,

3. für ein Element~$a$ ist $a_r =  (a)$,

4. $a\mathfrak a_r = (a\mathfrak a)_r$.\smallskip

\noindent $\mathfrak a_r$ heißt dann das \emph{aus\/ $\mathfrak a$
  erzeugte\/ $r$-Ideal}. Ist~$\mathfrak a$ endlich, so
heißt~$\mathfrak a_r$ ein endliches $r$-Ideal. Die Menge aller
$r$-Ideale nennen wir das \emph{$r$-\hspace{0pt}Idealsystem}. Ist sogar für jede
Untermenge~$\mathfrak a$ mit $\mathfrak a^{-1}\ne0$ ein $r$-Ideal
definiert, so nennen wir die Menge aller $r$-Ideale das \emph{totale\/
$r$-\hspace{0pt}Idealsystem}.\smallskip

Die Summe zweier Ideale wird durch
$\mathfrak a_r+\mathfrak b_r=(\mathfrak a\cup\mathfrak b)_r$, das
Produkt durch
$\mathfrak a_r\mathfrak b_r = (\mathfrak a\mathfrak b)_r$
definiert. Diese Definitionen sind unabhängig von der
Erzeugung der Ideale. Wegen~4.\ gilt das Distributivgesetz:
\[
  (\mathfrak a_r+\mathfrak b_r)\mathfrak c_r = \mathfrak a_r\mathfrak c_r+\mathfrak b_r\mathfrak c_r\text.
\]

\dotfill\smallskip

Für beliebige Halbgruppen konstruieren wir jetzt zwei wichtige
Idealsysteme. Wir setzen:
\[
  \mathfrak a_s=\bigcup_{a\in\mathfrak a}(a),\qquad\mathfrak
  a_v=\bigcap_{\mathfrak a\subseteq(a)}(a)\text.
\]

\dotfill\smallskip

\so{De{f}inition 2}. \ Eine Halbgruppe heißt \emph{$r$-\hspace{0pt}abgeschlossen},
wenn $\mathfrak c_r : \mathfrak c_r = 1$ für jedes endliche
$r$-Ideal~$\mathfrak c_r$ gilt.

\dotfill\smallskip

Die Bedeutung der $r$-\hspace{0pt}Abgeschlossenheit liegt darin, daß in jeder
$r$-\hspace{0pt}abgeschlossenen Halbgruppe die folgende Konstruktion eines neuen
Idealsystems möglich ist [\dots]:\smallskip

\so{De{f}inition 3}. \ Ist $\mathfrak a$ eine endliche Untermenge
von~$\mathfrak G$, und ist~$\mathfrak g$ $r$-\hspace{0pt}abgeschlossen, so soll
$\mathfrak a_{r_a}$ aus den Elementen~$a$ bestehen, für die es ein
endliches $r$-Ideal~$\mathfrak c_r$ gibt mit
$a\mathfrak c_r\subseteq\mathfrak a_r\mathfrak c_r$.\smallskip

Die Mengen $\mathfrak a_{r_a}$ bilden dann ein Idealsystem, das System
der $r_a$-Ideale.

Die $r_a$-Ideale haben nun die Eigenschaft, eine \emph{Halbgruppe} zu bilden,
d.\ h.\ für beliebige endliche $r_a$-Ideale
$\mathfrak a, \mathfrak b, \mathfrak c$ folgt aus
$\mathfrak a\mathfrak c = \mathfrak b\mathfrak c$ stets
$\mathfrak a = \mathfrak b$. [\dots]

\dotfill\medskip

\noindent\hfil §~2.\nopagebreak\medskip

\noindent\hfil \textbf{Totale Idealsysteme, Endlichartigkeit.}\nopagebreak\medskip

[\dots] Wir bezeichnen dazu mit~$\mathfrak e$ stets Mengen aus
endlich vielen Elementen von~$\mathfrak G$. Ist dann $\mathfrak a$~eine beliebige Untermenge
von~$\mathfrak G$ mit $\mathfrak a^{-1}\ne0$, so setzen wir
\[
  \mathfrak a_{r_s}=\bigcup_{\mathfrak e\subseteq\mathfrak a}\mathfrak e_r\text,\quad\mathfrak a_{r_v}=\bigcap_{\mathfrak a\subseteq\mathfrak e_r}\mathfrak e_r\text.
\]

\dotfill\smallskip

Wir nennen das $v_s$-System kürzer das $t$-System.

\dotfill\medskip

\noindent\hfil §~3.\nopagebreak\medskip

\noindent\hfil \textbf{Vollständige Halbgruppen.}\nopagebreak\medskip

Im §~4 werden wir zeigen, wie der Idealbegriff dazu dienen kann, viele
Fragen über Halbgruppen auf die Untersuchung der folgenden sehr
speziellen Klasse von Halbgruppen zurückzuführen:\smallskip

\so{De{f}inition 5}\label{1939def5}. \ Eine Halbgruppe heißt \emph{vollständig}, wenn je
zwei Elemente~$a,b$ einen größten gemeinsamen Teiler~$a\lor b$
besitzen.

\dotfill\smallskip

Außerdem müssen wir noch den Begriff der Quotientenringe eines
Integritätsbereiches auf Halbgruppen übertragen. Sei
dazu~$\mathfrak g$ eine beliebige Halbgruppe, $S$~eine Unterhalbgruppe
von~$\mathfrak g$, dann soll $\mathfrak g_S$ die Halbgruppe aller
Quotienten~$\frac ab$ ($a\in\mathfrak g$, $b\in S$) sein. \ $\mathfrak g_S$ heißt eine \emph{Quotientenhalbgruppe}
von~$\mathfrak g$ (vgl.\ Grell [1]).

Die Nichteinheiten von~$\mathfrak g_S$, die gleichzeitig in~$\mathfrak g$
liegen, bilden ein $s$-\hspace{0pt}Primideal~$\mathfrak p$ von~$\mathfrak
g$. Bezeichnen wir dann mit $\mathfrak g_\mathfrak p$ die
Quotientenhalbgruppe der Elemente~$\frac ab$ ($a\in\mathfrak g$,
$b\in\mathfrak g$, $b\notin\mathfrak p$) so wird
$\mathfrak g_{\mathfrak p}=\mathfrak g_S$. Die Quotientenhalbgruppen
von~$\mathfrak g$ entsprechen also umkehrbar eindeutig den
\emph{$s$-\hspace{0pt}Primidealen} von~$\mathfrak g$.

Ist nun $\mathfrak g$~vollständig, so auch jede
Quotientenhalbgruppe~$\mathfrak g_S$, denn sind $s,s'$ Elemente
von~$S$, so liegt mit $\frac as,\frac{a'}{s'}$ auch
\[\frac as\lor\frac{a'}{s'}=\frac{s'a\lor sa'}{ss'}\]
in~$\mathfrak g_S$ und liefert gleichzeitig den größten gemeinsamen
Teiler von $\frac as,\frac{a'}{s'}$ bezüglich der Teilbarkeit in
$\mathfrak g_S$.

Umgekehrt\label{1939p545} ist aber auch jede Oberhalbgruppe~$\overline{\mathfrak g}$
von~$\mathfrak g$ (natürlich
$\mathfrak g\subseteq\overline{\mathfrak g}\subseteq\mathfrak G$), die
mit zwei Elementen $a, b$ auch $a\lor b$ enthält, eine
Quotientenhalbgruppe von~$\mathfrak g$, denn für jedes Element~$a$ von
$\mathfrak G$ gilt $a=\frac{a(1\lor a)^{-1}}{(1\lor a)^{-1}}$, wobei
$a(1\lor a)^{-1}$ und $(1\lor a)^{-1}$ Elemente von~$\mathfrak g$
sind. Ist $a\in\overline{\mathfrak g}$, so folgt
$1\lor a\in\overline{\mathfrak g}$ und
damit~$\overline{\mathfrak g}=\mathfrak g_S$, wenn $S$~aus den
Einheiten von~$\overline{\mathfrak g}$ besteht, die in~$\mathfrak g$
liegen.

Die wichtigsten Quotientenhalbgruppen sind die
Quotientenhalbgruppen~$\mathfrak g_{\mathfrak p}$, bei denen $\mathfrak p$~ein
$t$-\hspace{0pt}Primideal von~$\mathfrak g$ ist, denn es gilt:\smallskip

\so{Satz }10\label{1939satz10}. \ \emph{Jede
  Quotientenhalbgruppe\/~$\mathfrak g_{\mathfrak p}$ nach einem\/
  $t$-\hspace{0pt}Primideal\/~$\mathfrak p$ ist linear.}\smallskip

Jede Nichteinheit von $\mathfrak g_{\mathfrak p}$
 läßt sich nämlich als Quotient~$\frac ps$ ($p\in\mathfrak p$, $s\notin\mathfrak p$)
darstellen, und sind $\frac ps,\frac{p'}{s'}$
zwei Nichteinheiten, so ist auch
\[
  \frac ps\lor\frac{p'}{s'}=\frac{s'p\lor sp'}{ss'}
\]
Nichteinheit, da $s'p\lor sp'\in\mathfrak p$.
Für jedes $a \in\mathfrak G$ folgt also aus
\[{\frac1{1\lor a}}\lor{\frac1{1\lor a^{-1}}}=1\text,\]
daß $\frac1{1\lor a}$ oder $\frac1{1\lor a^{-1}}$ Einheit
von~$\mathfrak g_{\mathfrak p}$ ist, d.\ h.\
$a\in\mathfrak g_{\mathfrak p}$ oder
$a^{-1}\in\mathfrak g_{\mathfrak p}$.\smallskip

\so{Satz }11. \ \emph{Es sei\/~$N$ eine Menge von\/ $t$-\hspace{0pt}Primidealen
  von\/~$\mathfrak g$ derart, daß jede Nichteinheit
  von\/~$\mathfrak g$ in mindestens einem\/~$\mathfrak p\in N$ liegt, dann
  ist\/~$\mathfrak g$ Durchschnitt der linearen
  Oberhalbgruppen\/~$\mathfrak g_{\mathfrak p}$\/ $(\mathfrak p\in N)$}
\[
  \mathfrak g=\bigcap_N\mathfrak g_{\mathfrak p}\text.
\]
Denn ist $a\notin\mathfrak g$, so ist $(1\lor a)^{-1}$ Nichteinheit
von~$\mathfrak g$.  Also gibt es ein $\mathfrak p\in N$ mit
$(1\lor a)^{-1}\in\mathfrak p$, und hieraus folgt
$a\notin\mathfrak g_{\mathfrak p}$.\smallskip

\so{Satz }12. \ \emph{Enthält\/ $N$ alle maximalen\/ $t$-\hspace{0pt}Primideale, so gilt für alle\/ $t$-Ideale\/~$\mathfrak a$
von\/~$\mathfrak g$ die Durchschnittsdarstellung\/ $\mathfrak a=\bigcap\limits_N\mathfrak a\mathfrak g_{\mathfrak p}$.}\smallskip

Denn ist $a\notin\mathfrak a$, so ist $a^{-1}\mathfrak a\cap\mathfrak g\subset\mathfrak g$. Es gibt
also ein $\mathfrak p\in N$ mit $a^{-1}\mathfrak a\cap\mathfrak g\subseteq\mathfrak p$, woraus wieder $a\notin\mathfrak a\mathfrak g_{\mathfrak p}$ folgt.

Diese Sätze gestatten zunächst für vollständige Halbgruppen Fragen
über die linearen Oberhalbgruppen von~$\mathfrak g$ in rein
idealtheoretische Fragen zu verwandeln. Diese Reduktion auch für
allgemeinere Halbgruppen durchzuführen, ist das Ziel von §~4.\medskip

\noindent\hfil§~4.\medskip

\noindent\hfil \textbf{Idealbrüche.}\medskip

Die Konstruktion der vollständigen Halbgruppe, auf die die
Strukturuntersuchung allgemeinerer Halbgruppen zurückgeführt werden
soll, ist z.\ B. für total abgeschlossene Halbgruppen~$\mathfrak g$
leicht geschehen, denn dann bilden ja die sämtlichen ganzen $v$-Ideale
von~$\mathfrak g$ ihrerseits eine vollständige Halbgruppe.  Ebenso,
wenn $\mathfrak g$~eine Multiplikationshalbgruppe bezüglich irgendeines
Idealsystems ist, d.\ h.\ also wenn jedes endliche Ideal umkehrbar
ist, so bilden die ganzen Ideale der
Form~$\mathfrak a\mathfrak b^{-1}$, wobei $\mathfrak a, \mathfrak b$
endliche Ideale sind, eine vollständige Halbgruppe. Gehört nämlich
noch~$\mathfrak c\mathfrak d^{-1}$ dazu, so auch der größte gemeinsame
Teiler von~$\mathfrak a\mathfrak b^{-1}$ und
$\mathfrak c\mathfrak d^{-1}$; denn dieser ist:
\[
  \mathfrak a\mathfrak b^{-1}+\mathfrak c\mathfrak d^{-1}=(\mathfrak
  a\mathfrak d+\mathfrak b\mathfrak c)(\mathfrak b\mathfrak d)^{-1}\text.
\]
Aber wir brauchen von unserer Halbgruppe, in der etwa das System der
$r$-Ideale gegeben sei, nur vorauszusetzen, daß sie $r$-\hspace{0pt}abgeschlossen
ist, um auf eine vollständige Halbgruppe zu kommen. Denn
ist~$\mathfrak g$ $r$-\hspace{0pt}abgeschlossen, so bilden die endlichen
$r_a$-Ideale eine Halbgruppe; wir können also deren Quotientengruppe
bilden. Diese besteht aus den „\emph{Idealbrüchen}“
$\frac{\mathfrak a}{\mathfrak b}$, wobei $\mathfrak a, \mathfrak b$
endliche $r_a$-Ideale sind. Wir definieren dann für zwei Idealbrüche
$\frac{\mathfrak a}{\mathfrak b}\subseteq\frac{\mathfrak c}{\mathfrak
  d}$, wenn $\mathfrak a\mathfrak d\subseteq\mathfrak b\mathfrak
c$. Es wird dadurch die Relation~$\subseteq$ fortgesetzt in die
Quotientengruppe.

Betrachten wir nun alle Idealbrüche
$\frac{\mathfrak a}{\mathfrak b}\subseteq1$. Diese bilden eine
Halbgruppe~$\mathfrak h$, und zwar eine vollständige Halbgruppe. Der
größte gemeinsame Teiler von $\frac{\mathfrak a}{\mathfrak b}$ und
$\frac{\mathfrak c}{\mathfrak d}$ ist natürlich durch
$\frac{\mathfrak a\mathfrak d+\mathfrak b\mathfrak c}{\mathfrak
  b\mathfrak d}$ gegeben. Ist
$\frac{\mathfrak a}{\mathfrak b}\subseteq1$, so nennen wir
$\frac{\mathfrak a}{\mathfrak b}$ einen „\emph{ganzen
  Idealbruch}“. Die vollständige Halbgruppe der ganzen Idealbrüche
ersetzt uns den Kroneckerschen Funktionalring, wie er in Prüfer~[1]
konstruiert worden ist. Das wichtige Ergebnis, das Krull~[5] über den
Funktionalring beweist, läßt sich auf~$\mathfrak h$ übertragen. Es
handelt sich dabei um die linearen Oberhalbgruppen von~$\mathfrak g$
einerseits, $\mathfrak h$~andererseits.  Ist $\mathfrak L$~eine
lineare Oberhalbgruppe von~$\mathfrak h$, so bilden die Elemente~$a$
von~$\mathfrak G$, deren Hauptideal~$(a)$ in~$\mathfrak L$ liegt,
ersichtlich eine ebenfalls lineare Oberhalbgruppe~$\mathfrak B$
von~$\mathfrak g$. Enthält~$\mathfrak L$ mit endlich vielen Elementen
auch deren größten gemeinsamen Teiler, so hat~$\mathfrak B$ noch
folgende Eigenschaft:\smallskip

$\mathfrak B$ enthält mit endlich vielen Elementen $a_1,\dots,a_n$
auch das $r$-Ideal $(a_1,\dots,a_n)_r$.\smallskip

Nennen wir Halbgruppen mit dieser Eigenschaft $r$-\hspace{0pt}Halbgruppen, so
können wir unser Ergebnis kurz so fassen:\smallskip

Jede lineare $t$-\hspace{0pt}Oberhalbgruppe von~$\mathfrak h$ liefert eine lineare
$r$-\hspace{0pt}Oberhalbgruppe von~$\mathfrak g$.\smallskip

Dieser Satz läßt sich nun umkehren. Ist nämlich eine lineare
$r$-\hspace{0pt}Oberhalbgruppe~$\mathfrak B$ von~$\mathfrak g$ gegeben, so
zeigen wir zunächst, daß $\mathfrak B$~auch $r_a$-Oberhalbgruppe
ist. Dazu sei $d\in(a_1,\dots,a_n)_{r_a}$, \ $a_\nu\in\mathfrak B$
($\nu=1,\dots,n$), d.\ h.\ es gibt ein endliches $r$-Ideal
$\mathfrak c=(c_1,\dots,c_m)_r$, so daß
$d\mathfrak c\subseteq(a_1,\dots,a_n)_r\mathfrak c$.  Da
$\mathfrak B$~linear ist, gibt es unter den Elementen $c_1,\dots,c_m$
eines --~etwa~$c_1$~--, so daß $c_1^{-1}c_\nu\in\mathfrak B$
($\nu=1,\dots,m$). Dann folgt wegen
\[
  d(1,c_1^{-1}c_2,\dots,c_1^{-1}c_m)_r\subseteq(a_1,\dots,a_n)_r(1,c_1^{-1}c_2,\dots,c_1^{-1}c_m)_r\text,
\]
daß $d$~in einem endlichen $r$-Ideal liegt, das aus lauter Elementen
von~$\mathfrak B$ erzeugt wird, also muß~$d$ in~$\mathfrak B$
liegen. \ $\mathfrak B$~ist somit lineare $r_a$-Oberhalbgruppe, woraus
weiter folgt, daß für jedes endliche $r_a$-Ideal
$\mathfrak a=(a_1,\dots,a_n)_{r_a}$ stets $\mathfrak a$ in dem aus
$a_1,\dots,a_n$ erzeugten $s$-Ideal von~$\mathfrak B$ enthalten
ist. \ $\mathfrak a\mathfrak B$ ist also Hauptideal in~$\mathfrak
B$. Wir definieren jetzt eine lineare Oberhalbgruppe~$\mathfrak L$
von~$\mathfrak h$ durch
$\frac{\mathfrak a}{\mathfrak b}\in\mathfrak L$, wenn
\[
  \mathfrak a\mathfrak B\subseteq\mathfrak b\mathfrak B\text.
\]
$\mathfrak L$ ist dann eine lineare $t$-\hspace{0pt}Oberhalbgruppe
von~$\mathfrak h$ und enthält ein Hauptideal~$(a)$ genau dann, wenn
$a$~in~$\mathfrak B$ liegt. Insgesamt erhalten wir unter
Berücksichtigung von §~3:\smallskip

\so{Satz }13. \ \emph{Die linearen\/ $r$-\hspace{0pt}Oberhalbgruppen einer\/
  $r$-\hspace{0pt}abgeschlossenen Halbgruppe entsprechen umkehrbar eindeutig den
  linearen\/ $t$-\hspace{0pt}Oberhalbgruppen der Halbgruppe\/~$\mathfrak h$ der
  ganzen Idealbrüche und damit den\/ $t$-\hspace{0pt}Primidealen
  von\/~$\mathfrak h$.}\smallskip

Als Folgerung dieses Satzes ergibt sich sofort:\smallskip

\so{Satz }14. \ \emph{Jede\/ $r$-\hspace{0pt}abgeschlossene Halbgruppe ist
  Durchschnitt von linearen\/ $r$-\hspace{0pt}Oberhalbgruppen.}\smallskip




\dotfill\medskip

\noindent\hfil\textbf{Literaturverzeichnis.}\nopagebreak\medskip

\dotfill\nopagebreak\smallskip

\small\noindent H. Grell.\nopagebreak

[1] Beziehungen zwischen den Idealen verschiedener Ringe, Math.\ Annalen~\textbf{97}
(1927).

\dotfill\smallskip

\noindent W. Krull.

\dotfill\smallskip


[5] Beiträge zur Arithmetik kommutativer Integritätsbereiche, Beitrag~I, Math.\ Zeitschr.~\textbf{41} (1936).

\dotfill\smallskip

\noindent H. Prüfer.

[1] Untersuchungen über die Teilbarkeitseigenschaften von Körpern, J. reine
angew.\ Math.\ \textbf{168} (1932).

\end{Cadre}

\section{\texorpdfstring{\citealt{Lor1950}}{Lorenzen 1950}: the Fundamentalsatz without valuations.}
\label{sec:1950}

\citealt{Lor1950} takes place in the framework of preordered
noncommutative monoids.  It proposes an analysis of a system of ideals
as an embedding into a semilattice.
\begin{EngCadre}
  In the commutative setting, the systems of ideals are always introduced as
  systems of certain subsets. But if one removes this set-theoretic
  clothing, then the concept of an ideal may be defined quite simply: a
  system of ideals of a preordered set is nothing but an
  embedding into a semilattice.
\end{EngCadre}
This enables Lorenzen to analyse the concept of $r$-\hspace{0pt}supermonoid
of \citealt{Lor1939} as \emph{$r$-allowable monoid} corresponding to
an \emph{$r$-allowable preorder}, i.e.\ a preorder that is the
restriction of a homomorphic image of the system of $r$-ideals. The
reader is invited to return to the part of the letter to Krull dated 6 June 1944
that is translated in the introduction (\cpageref{19440606t}), and
to reflect upon the simplification and clarification w.r.t.\
\cref{w(a+b)} of the definition of a valuation.

In the letter to Krull dated 13 March 1944 (page~\pageref{19440313}), Lorenzen stresses the
importance of these discoveries in his research.
\begin{EngCadre}
  \noindent E.g.\ the insight that a system of ideals is intrinsically
  nothing but a supersemilattice, and a valuation nothing but a linear preorder, strikes me as the most essential result of my
  effort.
\end{EngCadre}

In the noncommutative framework, integral dependence cannot be captured anymore by
passing from the $r$- to the $r_{a}$-system. Therefore Lorenzen
devises a new construction that may be considered as the birth of
dynamical algebra. In the introduction, he describes it for the
system of Dedekind ideals.
\begin{EngCadre}
  Let\label{1950p489t} $I$~be an arbitrary integral domain, $\mathfrak a$~a (finite or
  infinite) ideal of~$I$, $a$~an element of the field of
  fractions~$K$. We ask ourselves under what
  circumstances $a\in\mathfrak aB$ holds for every valuation
  superring~$B$ of~$I$.  For
  arbitrary~$z\in K^*$ and any valuation superring~$B$ of~$I$, always
\[
  z\in B\quad\text{or}\quad z^{-1}\in B\text.
\]
holds. Therefore, whenever $a\in\mathfrak aI[z]$ and $a\in\mathfrak aI[z^{-1}]$ holds, so
does $a\in\mathfrak aB$. Likewise follows: if there are elements
$z_1,\dots,z_n\in K^*$ for which
\[
  a\in\mathfrak aI[z_1^{\pm1},\dots,z_n^{\pm1}]\leqno\qquad3)
\]
holds for each of the $2^n$~combinations of signs, then, for every valuation superring~$B$ of~$I$, holds:
\[
  a\in\mathfrak aB\text.
\]

This condition~3) is also necessary. Namely, if~3) is not fulfilled
for any $z_1,\dots,z_n$, then a simple well-ordering\label{zornintroeng} argument shows
that there is a maximal superring~$\overline{I}$ of~$I$ with the
property that, for all $z_1,\dots,z_n\in K^*$,
\[
a\notin\mathfrak a\overline I[z_1^{\pm1},\dots,z_n^{\pm1}]
\]
holds for at least one combination of signs. Therefore in particular
$a\notin\mathfrak a\overline I$ holds. If this~$\overline I$ were not
a valuation ring, then there would be a $z\in K^*$ with
$z\notin\overline I$ and $z^{-1}\notin\overline I$, therefore there
would be elements $x_1,\dots,x_m$ and $y_1,\dots,y_n$ with
\[
  \begin{aligned}
a&\in\mathfrak a\overline I[z][x_1^{\pm1},\dots,x_m^{\pm1}]\\
a&\in\mathfrak a\overline I[z^{-1}][y_1^{\pm1},\dots,y_n^{\pm1}]
\end{aligned}
\]
for each combination of signs. But from this
\[
a\in\mathfrak a\overline I[z^{\pm1},x_1^{\pm1},\dots,x_m^{\pm1},y_1^{\pm1},\dots,y_n^{\pm1}]
\]
would follow, i.e.\ a contradiction. Therefore $\overline I$~is a valuation ring.

With that, the following new result is gained:

$a\in\mathfrak aB$ holds for every valuation superring~$B$ of~$I$ exactly when there are elements $z_1,\dots,z_n\in K^*$ for which
\[
  a\in\mathfrak aI[z_1^{\pm1},\dots,z_n^{\pm1}]
\]
holds for each combination of signs.

\dotfill\smallskip

\label{end-intro-1950}$a$~is integrally dependent on~$I$ if there are elements $z_1,\dots,z_n$ of the field of fractions that fulfil
$a\in I[z_1^{\pm1},\dots,z_n^{\pm1}]$ for each combination of signs.
\end{EngCadre}

The gist of Lorenzen's new construction is that in a computation in an
integral domain~$I$, one may at any time open two branches concerning
an element~$z$: in one branch, one supposes that~$I$ contains~$z$; in
the other, that~$I$ contains~$z^{-1}$. If a result may be obtained in each branch by a computation,
then there is a computation that yields this result without opening the branches. This construction embodies the
computational content of the proof method described by \citet{krull36}:
see \cpageref{method}.

This leads to a reshaping of Krull's Fundamentalsatz for integral domains into Theorem~26.
We postpone our discussion of this reshaping until the next section, on \citealt{Lor1952}, instead of dwelling on the 1950 version.

But let us stress another discovery reported in this article: the
property of regularity\label{regularity} (see Satz~8 on \cpageref{satz8}). It consists
in requiring of an $\ell$-group that if $a\land xax^{-1}=1$, then
$a=1$. This property may be seen as specifying the amount of commutativity that must be
granted in an $\ell$-group in order to be allowed to open branches in computations as described above. For more on this, see
\citealt{coquandlombardineuwirthkonstanz}.

\subsection{Excerpt.}

Our excerpt is from pages~483--490, 492--494, 496--500, 502--503, 506, 509--512, and~515--518 of ``Über halbgeordnete Gruppen'', \emph{Mathematische Zeitschrift}~52.

\begin{Cadre}\leqnomode\parskip1.5pt plus.5pt minus.5pt\newcommand\dotle{\mathrel{\dot\le}}
  Ein Integritätsbereich~$I$ definiert in der multiplikativen Gruppe~$K^*$
  seines Quotientenkörpers~$K$ eine Teilbarkeitsrelation. [\dots]

  \dotfill\smallskip

  Wie man der Theorie der kommutativen Halbgruppen entnehmen kann,
  sind die Ordnungen, die einen Bewertungsring~$B$ bestimmen, der~$I$
  umfaßt, dadurch gekennzeichnet, daß sie zugleich eine Ordnung der
  Ideale~$\mathfrak a,\mathfrak b,\dots$ von~$I$ ermöglichen, so daß stets gilt:
  \[
    \mathfrak c\le\mathfrak a,\ \mathfrak c\le\mathfrak
    b\implies\mathfrak c\le\mathfrak a+\mathfrak b\text.\leqno\qquad5)
  \]
  Genauer gesagt: die Teilbarkeitsrelationen der Bewertungsoberringe
  von~$I$ sind die Ordnungen von~$K^*$, die durch die Ordnungen des
  Idealsystems von~$I$ induziert werden, die~5) erfüllen. Die hier
  auftretenden Ordnungen des Idealsystems nennen wir die zulässigen
  Ordnungen.

  Wir definieren in §~4 Idealsysteme für beliebige halbgeordnete
  Mengen und insbesondere für halbgeordnete Gruppen. Wir definieren
  auch für diesen allgemeinen Fall die zulässigen Ordnungen eines
  Idealsystems und dann weiter:

  Eine Ordnung einer halbgeordneten Gruppe~$G$, zu der ein
  Idealsystem~$H$ gegeben ist, heißt zulässig (bezgl.~$H$), wenn sie
  durch eine zulässige Ordnung von~$H$ in~$G$ induziert wird. Die
  Teilbarkeitsrelationen der Bewertungsoberringe von~$I$ sind dann also
  genau die (bezgl.\ des \textsc{Dedekind}schen Idealsystems) zulässigen
  Ordnungen von~$K^*$.

  \dotfill\smallskip

  Gegeben sei eine beliebige halbgeordnete Gruppe~$G$ und ein
  beliebiges Idealsystem~$H$ von~$G$. Wann und wie ist die Halbordnung
  von~$G$ als Konjunktion zulässiger Ordnungen von~$G$ darstellbar? [\dots]

  \dotfill\smallskip

\noindent Wie läßt sich eine halbgeordnete Gruppe in
  eine Verbandsgruppe einbetten~-- vorausgesetzt, daß eine solche Einbettung möglich ist?

  \dotfill\smallskip

  Wir erhalten die Antwort auf unsere Frage durch eine genaue
  halbordnungstheoretische Analyse des Idealbegriffs. Im Kommutativen
  werden die Idealsysteme stets als Systeme von gewissen Untermengen
  eingeführt. Beseitigt man aber diese mengentheoretische Einkleidung,
  so läßt sich der Idealbegriff ganz einfach definieren: ein Idealsystem
  einer halbgeordneten Menge ist nichts anderes als eine Einbettung
  in einen Halbverband. (Ich bemerke ausdrücklich, daß man hier statt
  „Halbverband“ nicht „Verband“ setzen darf.)

  \dotfill\smallskip

  Ist die Halbordnung von~$G$ als Konjunktion zulässiger Ordnungen
  von~$G$ darstellbar, so gibt es ein Idealsystem~$H'$ von~$G$, dessen
  $v$-Ideale eine reguläre Verbandsgruppe bilden. [\dots]

  Die hier auftretenden Idealsysteme~$H'$, die in der Arithmetik
  „arithmetisch brauchbar“ heißen, nennen wir auch im allgemeinen Fall
  „brauchbar“.

  \dotfill\smallskip



  Wie in der Arithmetik ergibt sich dann weiter, daß unter den
  brauchbaren Idealsystemen eines ausgezeichnet ist, zu dem alle
  anderen homomorph sind [\dots]


  \dotfill\smallskip

  Zur Beantwortung unserer Ausgangsfragestellung für halbgeordnete
  Gruppen bleibt uns jetzt vor allem folgende Frage:

  Wie\label{wie} läßt sich das ausgezeichnete brauchbare Idealsystem einer
  halbgeordneten Gruppe~$G$ konstruieren --~vorausgesetzt, daß es
  überhaupt existiert.

  Wir setzen also voraus, daß die Halbordnung von~$G$ als Konjunktion
  zulässiger Ordnungen von~$G$ darstellbar ist. Bei der Konstruktion
  des ausgezeichneten brauchbaren Idealsystems wollen wir aber die
  Kenntnis dieser Ordnungen nicht voraussetzen, weil wir gerade mit
  Hilfe des Idealsystems die zulässigen Ordnungen bestimmen wollen.

  Diese Konstruktion wird in der Arithmetik geliefert durch den
  \textsc{Prüfer}schen Begriff der Ganzabhängigkeit von einem Ideal~$\mathfrak a$ des
  gegebenen Integritätsbereiches~$I$.

  Ein Element~$a$ des Quotientenkörpers~$K$ heißt ganz abhängig
  von~$\mathfrak a$, wenn~$a$ einer Gleichung
\[a^n+a_1a^{n-1}+\cdots+a_n=0\leqno\qquad1)\]
mit $a_\nu\in\mathfrak a^\nu$ genügt.

In der Theorie der kommutativen Halbgruppen tritt an die Stelle
von~1) die Bedingung, daß es ein endliches Ideal~$\mathfrak c$ geben soll, so daß
\[
  a\mathfrak c\subseteq\mathfrak a\mathfrak c\leqno\qquad2)
\]
gilt. Auch diese Bedingung 2) läßt sich jedoch nicht auf den
nichtkommutativen Fall ausdehnen.

Um den Grundgedanken unserer neuen Konstruktion, die wir statt~1)
oder~2) verwenden, deutlich zu machen, führe ich sie hier für den
arithmetischen Spezialfall durch. §~6 des Textes bringt die
Durchführung für beliebige halbgeordnete Gruppen.

$I$\label{1950p489}~sei ein beliebiger Integritätsbereich, $\mathfrak a$~ein (endliches
oder unendliches) Ideal von~$I$, $a$~ein Element des
Quotientenkörpers~$K$. Wir fragen uns, wann für jeden
Bewertungsoberring~$B$ von~$I$ gilt: $a\in\mathfrak aB$.  Für beliebiges~$z\in K^*$
und jeden Bewertungsoberring~$B$ von~$I$ gilt stets
\[
  z\in B\quad\text{oder}\quad z^{-1}\in B\text.
\]
Gilt also $a\in\mathfrak aI[z]$ und $a\in\mathfrak aI[z^{-1}]$, so
gilt stets $a\in\mathfrak aB$. Ebenso folgt: Gibt es Elemente
$z_1,\dots,z_n\in K^*$, für die
\[
  a\in\mathfrak aI[z_1^{\pm1},\dots,z_n^{\pm1}]\leqno\qquad3)
\]
für jede der $2^n$~Vorzeichenkombinationen gilt, so gilt für jeden Bewertungsoberring~$B$ von~$I$:
\[
  a\in\mathfrak aB\text.
\]

Diese Bedingung~3) ist auch notwendig. Ist nämlich~3) für kein
$z_1,\dots,z_n$ erfüllt, so zeigt ein einfacher Wohlordnungsschluß\label{zornintrodeu}, daß es
einen maximalen Oberring~$\overline{I}$ von~$I$ gibt, mit der Eigenschaft, daß für
jedes $z_1,\dots,z_n\in K^*$
\[
a\notin\mathfrak a\overline I[z_1^{\pm1},\dots,z_n^{\pm1}]
\]
für mindestens eine Vorzeichenkombination gilt. Insbesondere gilt also
$a\notin\mathfrak a\overline I$. Wäre dieses~$\overline I$ kein Bewertungsring, so gäbe es ein $z\in K^*$ mit $z\notin\overline I$ und $z^{-1}\notin\overline I$, also gäbe es Elemente $x_1,\dots,x_m$ und $y_1,\dots,y_n$ mit
\[
  \begin{aligned}
a&\in\mathfrak a\overline I[z][x_1^{\pm1},\dots,x_m^{\pm1}]\\
a&\in\mathfrak a\overline I[z^{-1}][y_1^{\pm1},\dots,y_n^{\pm1}]
\end{aligned}
\]
für jede Vorzeichenkombination. Hieraus würde aber
\[
a\in\mathfrak a\overline I[z^{\pm1},x_1^{\pm1},\dots,x_m^{\pm1},y_1^{\pm1},\dots,y_n^{\pm1}]
\]
folgen, d.\ h.\ ein Widerspruch. Also ist~$\overline I$ Bewertungsring.

Damit ist das folgende neue Ergebnis gewonnen:

Es gilt genau dann $a\in\mathfrak aB$ für jeden Bewertungsoberring~$B$ von~$I$, wenn es Elemente $z_1,\dots,z_n\in K^*$ gibt, für die
\[
  a\in\mathfrak aI[z_1^{\pm1},\dots,z_n^{\pm1}]
\]
für jede Vorzeichenkombination gilt.

Der\label{einfach} Beweis ist so einfach, daß er sich ohne weiteres in beliebigen
halbgeordneten Gruppen durchführen läßt und hier sowohl eine Konstruktion des ausgezeichneten brauchbaren Idealsystems liefert (vorausgesetzt, daß es existiert), als auch die folgende Frage beantwortet:

Wann ist die Halbordnung einer halbgeordneten Gruppe~$G$ als
Konjunktion zulässiger Ordnungen von~$G$ darstellbar?

Zum Schluß zeigen wir, daß sich diese Ergebnisse mit den bekannten Ergebnissen im Kommutativen in Übereinstimmung befinden. [\dots]
Ferner erhält man für den grundlegenden Begriff der Ganzabhängigkeit
eines Elementes~$a$ von einem Integritätsbereich~$I$ einen neuen
gleichwertigen Ausdruck:

\label{fin-intro-1950}$a$~ist genau dann ganz abhängig von~$I$,
wenn es Elemente $z_1,\dots,z_n$ des Quotientenkörpers gibt, die
$a\in I[z_1^{\pm1},\dots,z_n^{\pm1}]$ für jede Vorzeichenkombination erfüllen.\bigskip

\noindent\hfil§~1.\nopagebreak\smallskip

\noindent\hfil\textbf{Halbordnungen.}\nopagebreak\medskip

\dotfill\smallskip

Für jede Halbordnung benutzen wir außerdem zur Abkürzung:
\[
  \begin{aligned}
    &a\ge b\text{, wenn $b\le a$ gilt,}\\
    &a \equiv b\text{, wenn $a\le b$ und $a \ge b$ gilt,}
  \end{aligned}
\]
\dotfill\smallskip

Die ganzen Elemente einer halbgeordneten Gruppe~$G$ bilden eine
invariante Unterhalbgruppe~$\mathfrak g$ von~$G$. Es gilt nämlich
\[1\ge1\text,\]
für jedes~$a,b\in G$:
\[a\ge1,\ b\ge1\implies ab\ge1\text,\]
für jedes~$a,x\in G$:
\[a\ge1\implies xax^{-1}\ge1\text.\]
Dies bedeutet für die Menge~$\mathfrak g$ der ganzen Elemente
\[1\in\mathfrak g\text,\]
für jedes~$a,b\in G$:
\[a\in\mathfrak g,\ b\in\mathfrak g\implies ab\in\mathfrak g\text,\]
für jedes~$a,x\in G$:
\[a\in\mathfrak g\implies xax^{-1}\in\mathfrak g\text.\]
Wir nennen die invariante Unterhalbgruppe der ganzen Elemente kurz
die (zu~$\le$) zugehörige Halbgruppe.

\so{Satz} 1. \ Die Halbordnungen einer Gruppe~$G$ entsprechen eineindeutig den invarianten Unterhalbgruppen von~$G$.

\dotfill\bigskip

\noindent\hfil§~2.\smallskip

\noindent\hfil\textbf{Verbandsgruppen.}\medskip

\dotfill\smallskip

\noindent Wir nennen zwei Elemente~$a$ und~$b$ einer Verbandsgruppe
teilerfremd, wenn $a \land b \equiv 1$ gilt, und schreiben hierfür
$a\parallel b$. [\dots]

\dotfill\bigskip

\noindent\hfil§~3.\smallskip

\noindent\hfil\textbf{Hauptgruppen.}\medskip

\dotfill\smallskip

Wir nennen diese direkten Produkte von geordneten Gruppen
Vektorgruppen. [\dots]

\label{satz8}\so{Satz} 8. \ In einer Vektorgruppe~$G$ gilt für jedes $a, x \in G$
\[a\parallel xax^{-1}\implies a \equiv 1\text.\]

\so{Beweis.} \ Wir haben nur nachzuweisen, daß diese Bedingung für
geordnete Gruppen gilt. Aus $a > 1 \implies xax^{-1} > 1$ folgt aber sofort
$a > 1 \implies \min (a, x a x^{-1}) > 1$.

Wir nennen jede Verbandsgruppe, die die Bedingung des Satzes~8
erfüllt, regulär. Also ist jede Vektorgruppe eine reguläre Verbandsgruppe.

Dieser Satz läßt sich nicht umkehren, wohl aber gilt, daß jede
reguläre Verbandsgruppe Untergruppe einer Vektorgruppe ist. Dieses
wollen wir im Folgenden beweisen.

Wir definieren zunächst: Jede Untergruppe einer Vektorgruppe
heißt eine Hauptgruppe.

\dotfill\smallskip

Es seien~$H$ und~$H'$ Halbverbandshalbgruppen und $\to$ eine Abbildung von~$H$ in~$H'$, die für jedes~$a, b\in H$, $a',b'\in H'$ mit~$a\to a'$ und $b\to b'$ erfüllt:
\newsavebox{\boximplies}\sbox{\boximplies}{${}\implies{}$}
\begin{align*}
  a\le b&\usebox{\boximplies} a'\le b'\tag*{\qquad1)}\\
  a\land b&\makebox[\wd\boximplies]{\hfil$\to$\hfil}a'\land b'\tag*{\qquad2)}\\
  ab&\makebox[\wd\boximplies]{\hfil$\to$\hfil}a'b'\text.\tag*{\qquad3)}
  \end{align*}
Eine solche Abbildung nennen wir einen Homomorphismus von~$H$.

Jeder Homomorphismus definiert eine Relation~$\dotle$ in~$H$ auf folgende
Weise. Wir setzen fest, daß $a\dotle b$ in~$H$ genau dann gelten soll, wenn
$a'\le b'$ in~$H'$ gilt. $H$~ist bzgl.\ dieser Relation wieder eine Halbverbandshalbgruppe und die Relation~$\dotle$ erfüllt:
\begin{align*}
  \tag*{\qquad \hphantom{I}I)}a\le b&\implies a\dotle b\\
  \tag*{\qquad II)}x\dotle a,\ x\dotle b&\implies x\dotle a \land b\text.
\end{align*}
Es gilt auch die Umkehrung: Ist
$\dotle$~eine Halbordnung von~$H$, die~I)
und~II) erfüllt, so ist~$H$ bzgl.~$\dotle$ eine Halbverbandshalbgruppe (wir
wollen diese mit~$H'$ bezeichnen) und die Abbildung~$a\to a$ von~$H$ auf~$H'$
ist ein Homomorphismus.

\dotfill\smallskip

Aufgrund dieser Gleichwertigkeit betrachten wir statt der Homomorphismen von~$H$ nur die Halbordnungen von~$H$, die I)~und~II) erfüllen. Wir nennen diese Halbordnungen kurz die zulässigen Halbordnungen von~$H$.

Zwei Homomorphismen von~$H$, die dieselbe zulässige Halbordnung
definieren, nennen wir äquivalent. Es entsprechen dann also die Homomorphismen von~$H$ bis auf Äquivalenz eineindeutig den zulässigen
Halbordnungen von~$H$.

Für eine Verbandsgruppe~$G$ entsprechen die zulässigen Halbordnungen von~$G$ aber auch eineindeutig gewissen Unterhalbgruppen von~$G$.
Aus Satz~1 folgt nämlich sofort:

Ist $\mathfrak g$~die zugehörige Halbgruppe einer Verbandsgruppe~$G$, so entsprechen die zulässigen Halbordnungen von~$G$ eineindeutig den invarianten Unterhalbgruppen~$\dot{\mathfrak g}$ von~$G$, die für jedes~$a, b\in G$
\[\text{I)}\quad\mathfrak g\subseteq\dot{\mathfrak g}\qquad\quad
  \text{II)}\quad a\in\dot{\mathfrak g},\ b\in\dot{\mathfrak g}\implies a\land b \in\dot{\mathfrak g}
\]
erfüllen.

Diese Halbgruppen~$\dot{\mathfrak g}$ nennen wir die zulässigen Halbgruppen von~$G$.
Die Homomorphismen einer Verbandsgruppe~$G$ entsprechen also bis
auf Äquivalenz eineindeutig den zulässigen Halbgruppen von~$G$.

Die zulässigen Halbgruppen einer Verbandsgruppe lassen sich nun
durch „Quotientenbildung“ aus der zugehörigen Halbgruppe~$\mathfrak g$ gewinnen.
Dazu haben wir invariante Unterhalbgruppen~$S$ von~$\mathfrak g$ zu betrachten.
Gilt für jedes~$a,b\in\mathfrak g$
\[
  a \le b,\ b \in S \implies a \in S\text,
\]
so nennen wir $S$~vollständig.

\so{Satz} 9. \ Die Homomorphismen einer Verbandsgruppe~$G$ entsprechen
bis auf Äquivalenz eineindeutig den vollständigen, invarianten Unterhalbgruppen der zugehörigen Halbgruppe~$\mathfrak g$.\bigskip

\dotfill\smallskip

\so{Satz} 13. \ Jede reguläre Verbandsgruppe ist Hauptgruppe.

\dotfill\smallskip

Zusammengefaßt erhalten wir also eine Kennzeichnung der Hauptgruppen, die unabhängig ist vom Ordnungsbegriff. Die\label{superlattice} Bedeutung
dieser Kennzeichnung liegt darin, daß sich die Ordnungen einer Gruppe
im allgemeinen nur mit Wohlordnungsschlüssen konstruieren lassen~--
die Oberverbandsgruppen dagegen lassen sich auch ohne Wohlordnung
konstruieren mit Hilfe des Idealbegriffs.

\bigskip

\noindent\hfil§~4.\nopagebreak\smallskip

\noindent\hfil\textbf{Idealsysteme.}\nopagebreak\medskip

Wir definieren den Idealbegriff zunächst für beliebige halbgeordnete
Mengen.

Unter einem Idealsystem~$\mathfrak S$ einer halbgeordneten Menge~$M$ verstehen wir einen minimalen Oberhalbverband von~$M$. Ein Oberhalbverband von~$M$ enthält zu endlich vielen Elementen~$a_{\mu}$ von~$M$ ($\mu = 1 , \dots, m$)
stets den g.\ g.\ T.~$a_1\land\dotsc\land a_m$. Er heißt minimal, wenn er nur aus
Elementen dieser Form besteht. Es gibt dann nämlich keinen echten
Unterhalbverband, der auch noch~$M$ umfaßt.

\dotfill\smallskip

Der Idealbegriff ist für halbgeordnete Gruppen in dieser Form noch
zu unbestimmt. Als Idealsystem einer halbgeordneten Halbgruppe~$H$
wollen wir nur diejenigen Idealsysteme der halbgeordneten Menge~$H$
bezeichnen, in denen durch die Festsetzung
\[
  (a_1\land\dotsc\land a_m)(b_1\land\dotsc\land b_n)\equiv
  a_1b_1\land\dotsc\land a_mb_n
\]
eine Multiplikation definiert werden kann.

\dotfill\smallskip

Die Einführung\label{embedding} der Idealsysteme geschieht, um mit ihrer Hilfe die
Einbettung einer Hauptgruppe~$G$ in eine reguläre Verbandsgruppe
konstruktiv ohne Benutzung der Ordnungen von~$G$ zu ermöglichen.

\dotfill\bigskip

\noindent\hfil§~5.\smallskip

\noindent\hfil\textbf{$r$-Gruppen.}\medskip

\dotfill\smallskip

Damit unsere Untersuchung der halbgeordneten Gruppen auch die
speziellen Fragen der Arithmetik der kommutativen Integritätsbereiche
umfaßt, legen wir daher von jetzt ab eine halbgeordnete Gruppe zugrunde, zu der von vornherein ein festes Idealsystem gegeben ist.
Ist dies das System der $r$-Ideale, so nennen wir die halbgeordnete
Gruppe kurz eine $r$-Gruppe.

In einer $r$-Gruppe erklären wir nur diejenigen Halbordnungen als
$r$-\hspace{0pt}zulässig, die durch eine zulässige Halbordnung des $r$-\hspace{0pt}Idealsystems induziert werden.

Die $r$-\hspace{0pt}zulässigen Halbordnungen lassen sich unabhängig von dem
Begriff der zulässigen Halbordnung eines Halbverbandes kennzeichnen.
Ist nämlich eine beliebige halbgeordnete Gruppe~$G$ gegeben, und ist
$H$~ein Idealsystem von~$G$, so gilt

\so{Satz} 17. \ Eine Halbordnung~$\dotle$ von~$G$ wird genau dann durch
eine zulässige Halbordnung von~$H$ induziert, wenn
\[
  \leqno\qquad1)\enskip a \le b \implies a \dotle b
\]
erfüllt ist und wenn aus $a_1\land\dotsc\land a_m\le a$
in~$H$ für jedes~$x \in G$
folgt
\[
  \leqno\qquad2)\enskip x\dotle a_1,\dots,x\dotle a_m\implies x\dotle a\text.
\]

\dotfill\smallskip

Für die Ordnungen von~$G$ gilt darüber hinaus

\so{Satz} 18. \ Die zulässigen Ordnungen von~$H$ entsprechen eineindeutig den Ordnungen von~$G$, die 1)~und 2) erfüllen.

\dotfill\smallskip

Ist $G$~eine $r$-Gruppe und $\mathfrak g$~die zugehörige Halbgruppe, so entsprechen die $r$-\hspace{0pt}zulässigen Halbordnungen bzw.\ Ordnungen von~$G$ den
invarianten Unterhalbgruppen bzw.\ den linearen, invarianten Unterhalbgruppen~$\dot{\mathfrak g}$ die
\[
  \leqno\qquad1)\enskip \mathfrak g\subseteq\dot{\mathfrak g}
\]
erfüllen, und für die aus $a_1\land\dotsc\land a_m\mathrel{\smash{\underset{r}{\le}}}a$
folgt:
\[
  \leqno\qquad2)\enskip a_1\in\dot{\mathfrak g},\dots,a_m\in\dot{\mathfrak g}\implies a\in\dot{\mathfrak g}\text.
\]
Statt~2) läßt sich auch schreiben:
\[
  \leqno\qquad\phantom{2)}\enskip a_1\in\dot{\mathfrak g},\dots,a_m\in\dot{\mathfrak g}\implies(a_1,\dots,a_m)_r\subseteq\dot{\mathfrak g}\text.
\]

Wir nennen die invarianten Unterhalbgruppen von~$G$, die 1)~und~2)
erfüllen, die $r$-\hspace{0pt}zulässigen Halbgruppen von~$G$.

Eine $r$-Gruppe~$G$ heißt $r$-\hspace{0pt}Hauptgruppe, wenn die Halbordnung von~$G$ nicht nur Konjunktion von Ordnungen von~$G$ ist, sondern sogar
Konjunktion von $r$-\hspace{0pt}zulässigen Ordnungen ist. Das ist also genau dann
der Fall, wenn die zugehörige Halbgruppe Durchschnitt von linearen
$r$-\hspace{0pt}zulässigen Halbgruppen ist.

\dotfill\bigskip

\noindent\hfil§~6.\smallskip

\noindent\hfil\textbf{Die $r$-Abschließung.}\medskip

Wir wollen in diesem~§ nachweisen, daß sich das ausgezeichnete
$r$-\hspace{0pt}brauchbare Idealsystem, das $r_a$-System, einer $r$-Hauptgruppe unabhängig von den Begriffen der vorhergehenden~§§ kennzeichnen läßt.

Es sei $G$~eine beliebige $r$-Gruppe, $\mathfrak g$~die zugehörige
Halbgruppe, $\dot{\mathfrak g}$ eine $r$-\hspace{0pt}zulässige Halbgruppe,
$\mathfrak a = \{a_1,\dots, a_m\}$ eine endliche Untermenge von~$G$
und $a \in G$.

Wir bezeichnen das $r$-Ideal~$(a_1,\dots,a_m)_r$ mit~$\mathfrak a_r$.

Wir definieren das Produkt~$\mathfrak a_r\mathbin.\dot{\mathfrak g}$ als die kleinste Untermenge von~$G$, die sämtliche Elemente~$a\dot g$ mit $a \in \mathfrak a_r$ und $\dot g\in\dot{\mathfrak g}$ enthält und außerdem für jedes $x_1,\dots,x_n \in G$ die Bedingung
\[x_1\in\mathfrak a_r\mathbin.\dot{\mathfrak g},\dots,x_n\in\mathfrak a_r\mathbin.\dot{\mathfrak g}\implies(x_1,\dots,x_n)_r\subseteq\mathfrak a_r\mathbin.\dot{\mathfrak g}\]
erfüllt.

\dotfill\smallskip

[\dots] bezeichnen wir für jedes $x \in G$ mit $\dot{\mathfrak g}(x)_r$ den Durchschnitt aller $r$-\hspace{0pt}zulässigen Halbgruppen, die $\dot{\mathfrak g}$~umfassen und außerdem $x$~enthalten.
$\dot{\mathfrak g}(x)_r$ ist wieder eine $r$-\hspace{0pt}zulässige Halbgruppe, die „$r$-\hspace{0pt}Erweiterung“ von~$\dot{\mathfrak g}$ mit~$x$.

\dotfill\smallskip

\newcommand\ar{\mathrel{\alpha_r}}
Wir definieren dazu eine Relation $\ar$ zwischen den Elementen~$a$
und den Produkten~$\mathfrak a_r\mathbin.\dot{\mathfrak g}$.

Wir setzen $a\ar\mathfrak a_r\mathbin.\dot{\mathfrak g}$, wenn es Elemente
$x_1,\dots, x_n$ von~$G$ gibt mit
$a \in \mathfrak a_r \mathbin. \dot{\mathfrak g} (x_1^{\pm1},\dots,x_n^{\pm1})_r$ für jede der~$2^n$ möglichen Vorzeichenkombinationen.

\label{satz24}\so{Satz} 24. \ $a$ ist genau dann $r$-abhängig von $\mathfrak a_r$ in $\dot{\mathfrak g}$, wenn $a\ar\mathfrak a_r\mathbin.\dot{\mathfrak g}$
gilt.

\dotfill\smallskip

Die\label{radirect} Bedeutung dieses Satzes liegt darin, daß es jetzt möglich ist,
für jede $r$-Hauptgruppe eine direkte Konstruktion des $r_a$-\hspace{0pt}Idealsystems
anzugeben.

\dotfill\smallskip

[\dots] Gilt $a\ar(1)\mathbin.\dot{\mathfrak g}$, so nennen wir $a$~$r$-\hspace{0pt}abhängig von~$\dot{\mathfrak g}$, da $(1)\mathbin. \dot{\mathfrak g} = \dot{\mathfrak g}$ ist.

Liegt jedes von~$\dot{\mathfrak g}$ $r$-\hspace{0pt}abhängige Element in~$\dot{\mathfrak g}$, so nennen wir $\dot{\mathfrak g}$
$r$-\hspace{0pt}abgeschlossen.

\dotfill\smallskip


\label{satz26}\so{Satz} 26. \ Eine $r$-Gruppe ist genau dann $r$-\hspace{0pt}Hauptgruppe, wenn die
zugehörige Halbgruppe $r$-\hspace{0pt}abgeschlossen ist.
\end{Cadre}

\section{\texorpdfstring{\citealt{Lor1952}}{Lorenzen 1952}: the Fundamentalsatz for semilattice domains.}
\label{sec:Lor1952}
\newcommand{\R}{\mathrel{R}}
\newcommand{\rS}{\mathrel{S}}

The first three sections of \citealt{Lor1952} propose a streamlined
version of \citealt{Lor1950} in the more general framework of a
\emph{domain}~$(B,{\preccurlyeq_B},G)$, i.e.\ of a preordered
set~$(B,{\preccurlyeq_B})$ with a 
monoid~$G$ of preorder-preserving operators
on~$B$. $(H,\preccurlyeq_H,G)$ is a \emph{semilattice
  domain} if $(H,\preccurlyeq_H)$~is a semilattice and if the
action of~$G$ preserves meets. A preorder~$\preccurlyeq$ on a
semilattice domain~$H$ is \emph{$\preccurlyeq_H$-allowable} if it is
coarser than~$\preccurlyeq_H$ and if meets (w.r.t.~$\preccurlyeq_H$)
are also meets w.r.t.~$\preccurlyeq$. One can similarly define lattice domains and the corresponding allowability.

A \emph{domain of ideals}~$H_r$ for a domain~$B$ is a minimal
supersemilattice domain. A preorder~$\preccurlyeq$ on~$B$ is
\emph{$r$-allowable} if it is induced by an allowable preorder
of~$H_r$. $B$~is an \emph{$r$-principal domain} if its preorder is a
conjunction of $r$-allowable linear preorders.

Krull's Fundamentalsatz for integral domains becomes a characterisation of an $r$-principal
domain as a domain~$B$ in which behaving as if it were linearly
preordered, i.e.\ in which assuming that certain pairs of
elements~$\alpha_1,\dots,\alpha_e\in B\times B$ are linearly
preordered, does not add new preordered pairs to the preorder. His
formulation of this ``as if'' goes as follows: it is the simultaneous
consideration of $2^e$~different preorders corresponding to
preordering each pair one way or the other (a way described by a combination
of signs $\varepsilon_1,\dots,\varepsilon_e=\pm1$), and each of these
preorders is the
conjunction $R[\alpha_1^{\varepsilon_1},\dots,\alpha_e^{\varepsilon_e}]_r$
of all the $r$-allowable preorders coarser than~$\preccurlyeq_B$ in
which the pairs are preordered in the described way (here below, $R$~is Lorenzen's notation for the preorder~$\preccurlyeq_B$).
\begin{EngCadre}
  \noindent\so{Theorem 1}. \ \emph{An\/ $r$-domain\/~$B$ (w.r.t.\/~$R$ and\/~$G$) is an\/ $r$-principal domain if, and only if, for all pairs\/ $\alpha_1,\dots,\alpha_e$ from\/~$B$, holds:}
  \[\bigcap\limits_{\makebox[0pt][c]{$\scriptstyle\varepsilon=\pm1$}}{}_{\strut\varepsilon_1,\dots,\varepsilon_e}R[\alpha_1^{\varepsilon_1},\dots,\alpha_e^{\varepsilon_e}]_r\subset R\text.\]\vspace{-\belowdisplayskip}
\end{EngCadre}

This final shape of the Fundamentalsatz corresponds to
Satz~26 of~\citealt{Lor1950} (see above) and is
obtained by the arguments developed there. The characterisation given
in Theorem~1 is but a formulation of \emph{$r$-closedness}: compare
the formulation of integral dependence at the end of the introduction
of~\citealt{Lor1950} (\cpageref{fin-intro-1950}, translated on
\cpageref{end-intro-1950}).

Note that the article \citealt{Lor1951} shows that Lorenzen had the
means of a considerably simpler formulation of this ``as if'' by
starting from the logic-free and set-theory-free formal system given
by~$H_r$ and by adding the $e$~axioms corresponding to preordering the
pairs (see, in
this respect,
\citealt{coquandlombardineuwirth19,coquandlombardineuwirthkonstanz}). He may have refrained from doing so in order to stick
to a purely algebraic framework that would suit his potential readers.

In the last section, the property of regularity introduced in \citealt{Lor1950} is given the following symmetrical form for a distributive lattice domain~$(V,\preccurlyeq_V,G)$: that $xa\land yb\preccurlyeq_Vxb\lor ya$ for $a,b\in V$ and $x,y\in G$.
\subsection{Excerpt.}

Our excerpt is from pages~269--274 of ``Teilbarkeitstheorie in Bereichen'',
\emph{Mathematische Zeitschrift}~55.

\begin{Cadre}
  \leqnomode
  \noindent\hfil§~1.\nopagebreak\smallskip

  \noindent\hfil \textbf{Bereiche.}\nopagebreak\medskip

  $B$~sei eine Menge (Elemente $a ,b , \dots$), $R$~eine zweistellige
  Relation in~$B$ und $G$~eine Menge von Operatoren $x,y, \dots$
  von~$B$.

  \so{De{f}inition 1}. \ B heißt ein \emph{Bereich} (bezügl.~$R$ und~$G$),
  wenn für alle~$a$, $b$, $c$ und~$x$ gilt
  \begin{align*}
    &a\R a\tag{1.1}\\
    &a\R b, b\R c\to a\R c\tag{1.2}\\
    &a\R b\to xa\R xb\text.\tag{1.3}
  \end{align*}

Da (1.3) für den Operator 1 (definiert durch $1 \cdot a = a$) und für~$xy$ (definiert durch $(xy) a = x(ya)$) gilt, falls für~$x$ und~$y$, sei im
Folgenden stets angenommen, daß $G$~eine Halbgruppe mit Einselement ist.

$R$~heißt die Halbordnung von~$B$. Statt $a \R b$ werde $a < b$ oder
$b > a$ geschrieben. $a = b$ bedeute $a < b$ und $b < a$ ($=$ braucht nicht
die Identität zu sein).

Gilt $a < b$ oder $a > b$ für alle~$a$ und~$b$, dann heißt~$\R$ eine \emph{Ordnung}.

\so{De{f}inition 2}. \ Ein Bereich~$B$ heißt ein \emph{Halbverbandsbereich},
wenn für alle $a, b$~ein größter gemeinsamer Teiler $a\land b$ existiert, so
daß für alle~$c$ und~$x$ gilt
\begin{gather*}
  c < a \land b \longleftrightarrow c < a , c < b\tag{1.4}\\
  x(a\land b) = x a \land xb\text.\tag{1.5}
\end{gather*}

  \dotfill\smallskip


\so{De{f}inition 3}. \ Eine Halbordnung~$\rS$ eines Halbverbandsbereiches~$B$ (bezügl.~$\R$ und~$G$) heißt \emph{zulässig}, wenn für alle~$a$, $b$, $c$ und~$x$ gilt
\begin{align*}
\tag{1.6}&R\subset S\\
\tag{1.7}&{a \rS b}\to {x a \rS xb}\\
\tag{1.8}&{c\rS a}, {c\rS b} \to {c\rS {a \land b}} \text.
\end{align*}

  \dotfill\smallskip


Die Halbordnung~$\R$ eines Bereichs~$B$ induziert in jeder Untermenge~$B_0$ von~$B$ eine Halbordnung~$R_0$ von~$B_0$. $R$~heißt eine \emph{Fortsetzung} von~$R_0$ auf~$B$. Ist~$B_0$ zulässig bezügl.~$G$ (d.\ h.\ $xB_0\subset B_0$ für
alle~$x$), dann induziert~$G$ eine Operatorenhalbgruppe~$G_0$ von~$B_0$. $B$~heißt dann ein \emph{Oberbereich} von~$B_0$ (bezügl.~$R_0$ und~$G_0$).\bigskip

  \noindent\hfil§~2.\smallskip

  \noindent\hfil \textbf{Ideale.}\medskip

  $B$~sei ein Bereich (bezügl.~$\R$ und~$G$).

  \so{De{f}inition 4}. \ Ein minimaler Halbverbandsoberbereich~$H$
  von~$B$ heißt ein \emph{Idealbereich} von~$B$.

  \dotfill\smallskip

  Die Idealbereiche~$H$ werden durch einen Kennbuchstaben unterschieden. $r$~sei eine Variable für diese Kennbuchstaben.

  Ein Bereich~$B$ mit einem Idealbereich~$H_r$ heiße ein \emph{$r$-Bereich}.

  \so{De{f}inition 5}. \ Eine Halbordnung~$S$ eines $r$-Bereiches~$B$
  heiße \emph{$r$-\hspace{0pt}zulässig}, wenn~$S$ durch eine zulässige Halbordnung
  von~$H_r$ induziert wird.

  \dotfill\smallskip

  \so{De{f}inition 6}. \ Ein $r$-Bereich~$B$ heißt \emph{$r$-\hspace{0pt}Hauptbereich}, wenn
  die Halbordnung~$\R$ von~$B$ Konjunktion $r$-\hspace{0pt}zulässiger Ordnungen
  von~$B$ ist.

  \dotfill\smallskip

\noindent werde für
jede $r$-\hspace{0pt}zulässige Halbordnung~$S$ und jedes Paar $\alpha = a_1, a_2$ aus~$B$ die
$r$-\hspace{0pt}Erweiterung $S[\alpha]_r$ als Konjunktion aller $r$-\hspace{0pt}zulässigen Halbordnungen~$\overline{S}$ von~$B$ mit
\begin{align*}
\tag{2.1}S\subset\overline{S}\\
\tag{2.2}a_1\mathrel{\overline{S}}a_2
\end{align*}
definiert.
Zur bequemeren Formulierung sei ferner gesetzt
\[
  \begin{gathered}
    \hbox to7em{\dotfill}\\
    \alpha^{+1}=a_1,a_2\text,\quad\alpha^{-1}=a_2,a_1\text.
  \end{gathered}
\]
  \dotfill\smallskip

\so{Satz 1}. \ \emph{Ein\/ $r$-Bereich\/~$B$ (bezügl.\/~$R$ und\/~$G$) ist genau dann
    ein\/ $r$-\hspace{0pt}Hauptbereich, wenn für alle Paare\/ $\alpha_1,\dots,\alpha_e$ aus\/~$B$ gilt:}
  \[\bigcap\limits_{\makebox[0pt][c]{$\scriptstyle\varepsilon=\pm1$}}{}_{\strut\varepsilon_1,\dots,\varepsilon_e}R[\alpha_1^{\varepsilon_1},\dots,\alpha_e^{\varepsilon_e}]_r\subset R\text.\]
  \dotfill\bigskip

  \noindent\hfil§~4.\smallskip

  \noindent\hfil \textbf{Reguläre Verbandsbereiche.}\medskip

  \so{De{f}inition} 9. \ Ein Verbandsbereich~$V$ heißt \emph{regulär}, wenn $V$~distributiv ist und für $a, b$ und $x, y$ gilt: $x a \land y b < x b \lor y a$.

  \dotfill\smallskip

  \so{Satz}~3. \ \emph{Die Halbordnung eines Verbandsbereiches\/~$V$ ist genau
dann Konjunktion zulässiger Ordnungen von\/~$V$, wenn\/ $V$~regulär ist.}
\end{Cadre}

\section{\texorpdfstring{\citealt{Lor1953}}{Lorenzen 1953}: the Fundamentalsatz for integral domains as an embedding into a super-\texorpdfstring{$\ell$}{l}-group.}
\label{sec:1953}

The introduction of \citealt{Lor1953} describes his state of the art on Krull's Fundamentalsatz for integral domains:
\begin{itemize}[wide,nosep]
\item\citealt{Lor1950} shows that the problem of representing an
  integral domain~$I$ as intersection of valuation rings may be
  reduced to the problem of constructing a super-$\ell$-group of
  the divisibility group~$G$ of~$I$ such that the meets of elements
  of~$G$ form a homomorphic image of the system of Dedekind ideals.
\item In fact, Satz~3 of \citealt{Lor1952} (see above) establishes
  that the preorder of a regular lattice domain may be represented as
  a conjunction of allowable linear preorders, and this corresponds to an
  intersection of valuation rings in the case of an integral domain.
\end{itemize}

\begin{EngCadre}
  According to the Fundamentalsatz, the integral domains~$I$ that are representable as intersection of valuation rings are characterised by being integrally closed. On the other hand, the representability of~$I$ as an intersection of valuation rings is equivalent to the existence of a lattice-preordered supergroup~$V$ of the multiplicative group~$G$ of the field of fractions~$K$ of~$I$ -- where~$V$ must satisfy the condition that the domain of ideals~$H$ (which consists of all g.c.d.s $a_1\land\dots\land a_n$ with $a_\nu\in G$) contained in~$V$ is a homomorphic image of the domain of Dedekind ideals~$H_d$ of~$I$.

  This last equivalence results from the more general domain-theoretic theorem that the preorder of every regular lattice-preordered domain is a conjunction of allowable total orders -- if one adds that every commutative lattice-preordered group is trivially regular. Therefore two routes to the proof of Krull's Fundamentalsatz are available: one\label{Wege} has to construct, for an integrally closed integral domain, either an intersection representation through valuations or a lattice-preordered supergroup. Both proof possibilities are carried out for commutative groups -- for noncommutative groups, however, and more generally for domains, only the first route has been hitherto practicable.

  The present work describes the second proof course also in the general case. If afterwards one specialises the general method to divisibility in commutative fields~$K$, then the following results. For the construction of a lattice-preordered supergroup it is required only to define when, for arbitrary elements $a_\mu,b_\nu$ out of the multiplicative group~$G$ of~$K$,
\[a_1\land\dots\land a_m\prec b_1\lor\dots\lor b_n\]
is to hold ($\prec$ stands for ``divides'', $\land$~denotes the g.c.d., $\lor$~the l.c.m.).

For every integrally closed integral domain~$I$, one obtains a lattice-preordered supergroup of~$G$ if one defines
\[a_1\land\dots\land a_m\prec b_1\lor\dots\lor b_n\longleftrightarrow1\in{\textstyle\sum\limits_1{}_{\!\raisebox{-2pt}{$\scriptstyle\varkappa$}}}(a_1b_1^{-1},\dots,a_\mu b_\nu^{-1},\dots,a_mb_n^{-1})^\varkappa\text.\leqno(1)\]
For $n=1$, $b_1=b$,  this definition changes into the prüferian definition of integral dependence of~$b$ on $(a_1,\dots,a_m)$
\[a_1\land\dots\land a_m\prec b\longleftrightarrow b^k+c_1b^{k-1}+\dotsb+c_k=0\qquad[c_\varkappa\in(a_1,\dots,a_m)^\varkappa]\text.\leqno(2)\]
Therefore (1)~describes a suitable enlargement of the concept of integral dependence.
\end{EngCadre}

In this last article on the subject, Lorenzen proposes the construction
of an $\ell$-group from a preordered group~$G$ and a semilattice of
ideals~$H_r$; $G$~embeds into this $\ell$-group if and only if it is
$r$-closed. He does so by first constructing the distributive lattice generated by an ``entailment relation'' associated with~$H_r$ (see \citealp{coquandlombardineuwirthkonstanz}) and then proving that it is a lattice-preordered group.

\subsection{Excerpt.}

Our excerpt is from page~15 of ``Die Erweiterung halbgeordneter Gruppen zu
Verbandsgruppen'', \emph{Mathematische Zeitschrift}~58.

\begin{Cadre}
  Nach dem \textsc{Krull}schen Fundamentalsatz [\emph{1}] sind die
  Integritätsbereiche~$I$, die als Durchschnitt von Bewertungsringen
  darstellbar sind, dadurch charakterisiert, daß sie ganz
  abgeschlossen sind. Die Darstellbarkeit von~$I$ als Durchschnitt von
  Bewertungsringen ist andererseits äquivalent [\emph{3}] mit der Existenz
  einer Verbandsobergruppe~$V$ der multiplikativen Gruppe~$G$ des
  Quotientenkörpers~$K$ von~$I$ --~wobei~$V$ der Bedingung genügen muß,
  daß der in~$V$ enthaltene Idealbereich~$H$ (der aus allen
  g.g.T. $a_1\land\dotsc\land a_n$ mit $a_\nu\in G$ besteht) ein
  homomorphes Bild des \textsc{Dedekind}schen Idealbereichs~$H_d$
  von~$I$ ist.

  Diese letztere Äquivalenz ergibt sich aus dem allgemeineren
  bereichstheoretischen [\emph{4}] Satz, daß die Halbordnung jedes regulären
  Verbandsbereiches Konjunktion von zulässigen Ordnungen ist --~wenn
  man hinzufügt, daß jede kommutative Verbandsgruppe trivialerweise
  regulär ist. Daher stehen zum Beweis des \textsc{Krull}schen
  Fundamentalsatzes zwei Wege zur Verfügung: Man hat für einen ganz
  abgeschlossenen Integritätsbereich entweder eine
  Durchschnittsdarstellung durch Bewertungen oder eine
  Verbandsobergruppe zu konstruieren. Beide Beweismöglichkeiten sind
  für kommutative Gruppen durchgeführt --~für nichtkommutative Gruppen
  und allgemeiner für Bereiche ist bisher jedoch nur der erste Weg
  gangbar.

  Die vorliegende Arbeit stellt den zweiten Beweisgang auch im
  allgemeinen Falle dar. Spezialisiert man anschließend die allgemeine Methode auf die Teil%
barkeit in kommutativen Körpern~$K$, so ergibt sich folgendes. Zur Konstruktion
einer Verbandsobergruppe braucht nur definiert zu werden, wann für beliebige
Elemente $a_\mu,b_\nu$ aus der multiplikativen Gruppe~$G$ von~$K$ gelten soll [\emph{5}]
\[a_1\land\dotsc\land a_m\prec b_1\lor\dotsc\lor b_n\]
($\prec$ steht für "`teilt"', $\land$~bezeichnet den g.g.T., $\lor$~das k.g.V.).

Für jeden ganz abgeschlossenen Integritätsbereich~$I$ erhält man eine Ver%
bandsobergruppe von~$G$, wenn man definiert
\[a_1\land\dotsc\land a_m\prec b_1\lor\dotsc\lor b_n\longleftrightarrow1\in{\textstyle\sum\limits_1{}_{\!\raisebox{-2pt}{$\scriptstyle\varkappa$}}}(a_1b_1^{-1},\dots,a_\mu b_\nu^{-1},\dots,a_mb_n^{-1})^\varkappa\text.\leqno(1)\]
Für $n=1$, $b_1=b$ geht diese Definition über in die \textsc{Prüfer}sche Definition der
Ganzabhängigkeit von~$b$ von $(a_1,\dots,a_m)$
\[a_1\land\dotsc\land a_m\prec b\longleftrightarrow b^k+c_1b^{k-1}+\dotsb+c_k=0\qquad[c_\varkappa\in(a_1,\dots,a_m)^\varkappa]\text.\leqno(2)\]
Daher stellt (1) eine zweckmäßige Ausdehnung des Begriffes der Ganzab%
hängigkeit dar.

  \dotfill\medskip

\noindent\hfil\textbf{Literatur.}\smallskip

\small[\emph{1}]~\textsc{Krull}, W.: Allgemeine Bewertungstheorie. J. reine angew.\ Math.~\textbf{167}, 160--196
(1931). -- [\dots]
[\emph{3}]~\textsc{Lorenzen}, P.: Über halbgeordnete
Gruppen. Math.\ Z.~\textbf{52}, 483--526 (1949). -- [\emph{4}]~\textsc{Lorenzen}, P.: Teilbarkeitstheorie in Be%
reichen. Math. Z.~\textbf{55}, 269--275 (1952). -- [\emph{5}]~\textsc{Lorenzen}, P.: Algebraische und logistische
Untersuchungen über freie Verbände. J.\ Symb.\ Log.~\textbf{16}, 81--106 (1951).
\end{Cadre}

\section{A letter from Krull to Scholz from 1953: the
well-ordering theorem.}
\label{sec:letter-from-krull}

\subsection{Extract.}
\label{sec:excerpt}

We translate an extract from
the letter that Krull sent to Scholz on 18 April 1953 (Heinrich Scholz Archive at University and State Library of Münster), edited on \cpageref{19530418}.
\foreignlanguage{german}{Heinrich-Scholz-Archiv} at
\foreignlanguage{german}{Universitäts- und Landesbibliothek Münster}
.

\begin{EngCadre}
In working with the uncountable, in particular with the
well-ordering theorem, I always had the feeling that one uses fictions
there that need to be replaced some day by more reasonable
concepts. But I was not getting upset over it, because I was convinced
that in a careful application of the common ``fictions'' nothing false
comes out, and because I was firmly counting on the man who would some
day put all in order. Lorenzen has now found according to my
conviction the right way [\dots].
\end{EngCadre}

\subsection*{Well-ordering arguments and constructions.}

In his articles, Lorenzen uses well-ordering arguments in several places:
\begin{itemize}[wide,nosep]
\item in \citealt{Lor1939} for the proof of Satz~4 (spelled out in our comment on the proof of Theorem~12 on \cpageref{commenton12});
\item in \citealt{Lor1950} for the proof of Satz~24 (stated in the introduction for the case of integral domains, see \cpageref{zornintrodeu} and the translation on \cpageref{zornintroeng}); for the proof of Satz~13 and of Satz~1 in Teil~II;
\item in \citealt{Lor1952} for the proof of his lemma to Theorem~1, a version of \citealt[Satz~24]{Lor1950}; for the proof of Satz~3, a version of \citealt[Satz~13]{Lor1950}.
\item However, no well-ordering argument appears in \citealt{Lor1953}.
\end{itemize}

Three kinds of uses may be distinguished.
\begin{itemize}[wide,nosep]
\item \citealt[Satz~4]{Lor1939}: Zorn's lemma is used in order to obtain a maximal prime ideal.
\item \citealt[Satz~24]{Lor1950}; \citealt[lemma to Theorem~1]{Lor1952}: a ``more reasonable concept'', which amounts to behaving as if the domain were
linearly preordered, is introduced to replace the ``fiction'' of valuations.
\item \citealt[Satz~13]{Lor1950}; \citealt[Satz~3]{Lor1952}: Zorn's lemma is used to obtain linear preorders. These are ``fictions'' which are replaced by the ``more reasonable concept'' of an $\ell$-group, in which one may compute as if it were linearly preordered.
\end{itemize}

We can retrace Lorenzen's position w.r.t.\ well-ordering arguments through his use of the word \emph{construction} and related words. Let us gather the relevant passages.
\begin{enumerate}[wide,nosep]
\item In Lorenzen's correspondence edited in the appendix, these words appear only in his letter to Krull, dated 25 April 1944 (page~\pageref{19440425}), in which he provides a detailed description of his habilitation, and in his last two letters to Hasse, dated July and September 1963.
\item In the introduction to \citealt{Lor1950}, he writes the following
  (see \cpageref{wie}).
\begin{EngCadre}
  How can a preordered group be embedded into a lattice-preordered
  group -- assuming that such an embedding is possible?

  \dotfill\smallskip

  How can the distinguished brauchbar\footnotemark\ system of ideals of a
  preordered group~$G$ be constructed -- assuming that it actually
  exists?

  We thus assume that the preorder of~$G$ is representable as
  conjunction of allowable linear preorders of~$G$. But we do not
  want to use the knowledge of these linear preorders, because the allowable linear preorders are
  just what we want to determine by the
  aid of the system of ideals.

  In arithmetic, This construction is provided by the prüferian
  concept of integral dependence on an ideal~$\mathfrak a$ of the
  given integral domain~$I$.
\end{EngCadre}\footnotetext{In Lorenzen's approach, the distinguished brauchbar system of ideals is the system of
$r_a$-ideals.}
\item After presenting his ``new construction'', he adds the following
  (see \cpageref{einfach}).
  \begin{EngCadre}
    The proof is so easy that it may be carried out in arbitrary
    preordered groups and here provides a construction of the
    distinguished brauchbar system of ideals (assuming that it exists),
   and also answers the following question:

    When is the preorder of a preordered group~$G$ representable as conjunction of allowable preorders of~$G$?
  \end{EngCadre}
\item After having proved Satz~13, he makes the following comment (see
  \cpageref{superlattice}).
\begin{EngCadre}
  \noindent The significance of this characterisation lies in the fact that, in general, the linear preorders of a group can be constructed only with well-ordering arguments -- whereas the super-lattice-preordered groups can also be constructed without a well-ordering by the aid of the concept of ideal.
\end{EngCadre}
\item He repeats this a few pages later (see \cpageref{embedding}).
\begin{EngCadre}
  Systems of ideals are introduced in order to make possible by their aid the embedding of a principal group~$G$ into a regular lattice-preordered group
  constructively, without use of the linear preorders of~$G$.
\end{EngCadre}
\item Satz~24 is accompanied by the following comment (see
  \cpageref{radirect}).
\begin{EngCadre}
  \noindent The significance of this theorem lies in the fact that it is now
  possible to state, for every $r$-principal group, a direct
  construction of the system of $r_a$-ideals.
\end{EngCadre}
\item\citet[page~273]{Lor1952} uses the following wording, which is the most explicit in the articles studied here.
  \begin{EngCadre}
    For $\varrho$-principal domains~$B$, Theorem~2 provides a
    ``constructive'' definition of the
    domain of $\varrho_a$-superideals~$V_{\varrho_a}$ of~$B$, i.e.\ a definition without use of the
    $\varrho$-\hspace{0pt}allowable linear preorders of~$B$.
\end{EngCadre}
\item Finally, \citet{Lor1953} uses the following formulation (see \cpageref{Wege}).
\begin{EngCadre}
  \noindent [\dots] one has to construct for an integrally closed integral domain either an intersection representation through valuations or a lattice-preordered supergroup.
\end{EngCadre}
\end{enumerate}

As these passages show, Lorenzen holds that a well-ordering
argument provides a construction, but he makes distinctions: such an
argument avers \emph{that} something is possible, but not
\emph{how}; other strategies have to be developed in order to gain
knowledge of the thing over and above its bare existence. He also specifies
for each construction whether it relies on a well-ordering argument or
not, in which latter case it is ``direct'', or ``constructive''. This
adjective and the adverb ``constructively'' have a more restricted use
than the word ``construction'' and mean that the well-ordering
argument or the ``fiction'' of a linear preorder have been avoided.

The distinctions he makes show that, like \citet{krull36}, he thinks
that well-ordering arguments provide only the ``mere existence'' of the
thing constructed. He therefore devises new constructions that bypass
them: systems of ideals, super-lattice-preordered groups, the right to
compute as if everything were linearly preordered. However, we have not
found any
passage in which he explains the problem as clearly and
simply as Krull does, or in which he motivates his work along the
lines of Krull's letter to Scholz.
\[*\ *\ *\]

Lorenzen's reshaping of Krull's Fundamentalsatz for integral domains can be presented as a process in three steps.
\begin{enumerate}[wide,nosep]
\item \citealt{Lor1939} generalises Krull's well-ordering argument from integral domains to preordered cancellative monoids; valuation rings become linear $r$-\hspace{0pt}supermonoids.
\item \citealt{Lor1950} and \citeyear{Lor1952} successively generalise the well-ordering argument
  to preordered noncommutative monoids and to domains, so
  that valuations correspond to $r$-allowable linear preorders; the
  right to compute as if the monoid were linearly preordered is
  embodied in the relation~$\alpha_r$, which is shown to yield finite
  approximations of them on the proviso that they exist; Lorenzen describes a direct construction of an embedding into a lattice-preordered group on
  the proviso that it exists.
\item \citealt{Lor1953} provides the construction of the free lattice-preordered group generated by a regular noncommutative monoid without this proviso by applying Satz 7 of \citealt{Lor1951}, nowadays called the \emph{fundamental theorem of entailment relations} (see \citealp[§~2B]{coquandlombardineuwirth19}).
\end{enumerate}

Let us state a further step forward, which Lorenzen does not make.
\begin{enumerate}[wide,nosep,resume]
\item 
  The condition of integral closedness for an integral domain consists
in the admissibility of the axiom that divisibility is linear for the formal
system generated by the integral domain. The condition of $r$-closedness for a regular noncommutative monoid consists in the admissibility of the axiom that the preorder is linear.
\end{enumerate}

In the articles discussed here, Lorenzen does not question the
existence of a thing provided by a well-ordering argument, but the
nonconstructive nature of its existence fails to provide the
mathematical clarity needed for extending his work to a noncommutative
setting.

The three steps made by Lorenzen are thus genuinely motivated by
his endeavour to clarify mathematics, as he writes in his letter to Krull dated
6 June 1944 (page~\pageref{19440606}, partly translated on
\cpageref{19440606t}), rather than philosophically: viz.,
\begin{enumerate}[wide,nosep]
\item to unveil the order-theoretic nature of the Fundamentalsatz;
\item to devise a counterpart of valuations that works in a noncommutative setting;
\item to bypass valuations by the direct construction of the free lattice-preordered group.
\end{enumerate}

With Step 4 above, dynamical algebra is taking up this process of clarification: see \citet{Lom2016}.

\appendix
\part*{Lorenzen's correspondence with Hasse, Krull, and Aubert, together with some relevant documents.}\addcontentsline{toc}{part}{Lorenzen's correspondence with Hasse, Krull, and Aubert, together with some relevant documents.}

This appendix proposes an edition of the reports on Lorenzen's Ph.D. thesis by Helmut Hasse and Carl Ludwig Siegel, of the known correspondence of Lorenzen with Hasse, Wolfgang Krull, and Karl Egil Aubert, and of a relevant letter from Krull to Heinrich Scholz. It provides evidence for the circumstances in which Lorenzen comes to his insights during the studied period of time.

The reports on Lorenzen's Ph.D. thesis, the letters from Lorenzen to Hasse, and the carbon copies of the letters from Hasse to Lorenzen are in the \foreignlanguage{german}{Universitätsarchiv Göttingen}. The letters from 1938 are in the Wolfgang-Krull-Nachlass at \foreignlanguage{german}{Archiv der Universitäts- und Landesbibliothek} of \foreignlanguage{german}{Universität Bonn}. The letters from Krull and Aubert to Lorenzen, the carbon copies of the letters from Lorenzen to Krull from 1943--1944, to Scholz dated 26 May and 2 June 1944, and to Aubert, as well as the ``\foreignlanguage{german}{Bescheinigung}'' from 1942, are in the \foreignlanguage{german}{Paul-Lorenzen-Nachlass} at \foreignlanguage{german}{Philosophisches Archiv} of \foreignlanguage{german}{Universität Konstanz}. The request of statement to and statement from the Dozentenbundsführer, the Military government of Germany Fragebogen, the report and certificates on Lorenzen's political attitude under national socialism, as well as the notification on Lorenzen's inaugural lecture are in the Universitätsarchiv Bonn. The other letters from Lorenzen to Scholz, the letter from Krull to Scholz, and the carbon copies of the letters from Scholz to Lorenzen are in the \foreignlanguage{german}{Heinrich-Scholz-Archiv} at \foreignlanguage{german}{Universitäts- und Landesbibliothek Münster}.

\section{Synopsis.}

Let us start with a diachronic synopsis of Lorenzen's (L) correspondence with Hasse (H), Krull (K), and Karl Egil Aubert, together with the relevant correspondence between Hasse and Krull edited by \citet{roquette04} (whose dates appear slanted below), of some relevant letters from Lorenzen to Heinrich Scholz, and of some related documents.\smallskip

\begin{liste}{14.02.1938}
\ligne[\slshape 14.02.1938]H suggests to K that L spend some time with K in Erlangen to discuss ideal theory.
\ligne[{\hyperref[19380219]{19.02.1938}}]L pays a visit to K.
\ligne[\textsl{02.03.1938}]K expresses his satisfaction with L.
\ligne[\textsl{07.03.1938}]H joins in this expression of satisfaction.
\ligne[{\hyperref[19380308]{08.03.1938}}]L corrects certain points concerning $v$-ideals and total closedness.
\ligne[{\hyperref[19380313]{13.03.1938}}]L finds the proof method of \citealt[§~3]{Lor1939}.
\ligne[{\hyperref[19380318]{18.03.1938}}]L reports a success that will turn out to be spurious.
\ligne[{\hyperref[19380322]{22.03.1938}}]L provides a glossary for his multiplicative ideal theory.
\ligne[\textsl{03.05.1938}]H sends K the carbon copy of L's Ph.D. thesis and asks for a brief report.
\ligne[\textsl{22.05.1938}]K sends H a very positive report on L's Ph.D. thesis.
\ligne[{\hyperref[19380524]{24.05.1938}}]H writes a very positive report on L's Ph.D. thesis.
\ligne[\textsl{31.05.1938}]H thanks K for his report and joins in his criticism of L's laconic style.
\ligne[{\hyperref[19380602]{02.06.1938}}]Carl Ludwig Siegel writes a vacuous report on L's Ph.D. thesis.
\ligne[{\hyperref[19380609]{09.06.1938}}]L pays a second visit to K.
\ligne[\textsl{21.06.1938}]H emphasises his agreement with K on their appreciation of L.
\ligne[{\hyperref[19380630]{30.06.1938}}]H tells L that he must meet the Dozentenbundsführer [leader of the union of lecturers] and asks a mathematical question.
\ligne[{\hyperref[19380706]{06.07.1938}}]L answers H's question from Erlangen. He reports that he is also working on lattice theory with Gottfried Köthe.
\ligne[\textsl{01.09.1938}]K reports to H that he and L have not been able to repair the defective proofs in L's Ph.D. thesis.
\ligne[\textsl{12.03.1939}]K has L in mind for a position as assistant in Bonn and asks H whether he would agree to let L leave Göttingen.
\ligne[\textsl{17.03.1939}]H agrees, but expresses some apprehension about how to replace L in Göttingen.
\ligne[\textsl{27.03.1939}]K weighs at length the pros and cons of L leaving Göttingen for Bonn.
\ligne[\textsl{30.03.1939}]H agrees with K that L should leave Göttingen for Bonn.
\ligne[\textsl{19.07.1939}]K asks H for L's military address in order to tell him that he will be appointed in Bonn from 1 August on.
\ligne[\textsl{21.07.1939}]H answers that he knows only L's home address in Bad Pyrmont.
\ligne[{\hyperref[19390809]{09.08.1939}}]L thanks H for his years of supervision.
\ligne[{\hyperref[19390811]{11.08.1939}}]H acknowledges L's thanks.
\ligne[{\hyperref[19390906]{06.09.1939}}]War has been declared and L is waiting for his incorporation. He asks H for help in finding a position as a mathematician in the army.
\ligne[\textsl{28.09.1939}]H suggests to K that L help him by composing a manuscript according to his drafts for a projected volume of \emph{Crelles Journal} on groups.
\ligne[{\hyperref[19391011]{11.10.1939}}]L thanks H again on the occasion of the printing of his Ph.D. thesis as an offprint of \emph{Mathematische Zeitschrift}.
\ligne[{\hyperref[19391110]{10.11.1939}}]L submits the manuscript of \citealt{MR0002113} to H as editor of \emph{Crelles Journal}.
\ligne[{\hyperref[19391116]{16.11.1939}}]H acknowledges receipt of L's manuscript and asks a mathematical question.
\ligne[{\hyperref[19391117]{17.11.1939}}]H spells out in detail his question from the day before.
\ligne[{\hyperref[19391213]{13.12.1939}}]L answers H's question.
\ligne[{\hyperref[19400104]{04.01.1940}}]H acknowledges receipt of L's announcement of marriage with Käthe and proposes that L apply for a military position under Alwin Walther by the Baltic sea.
\ligne[{\hyperref[19400109]{09.01.1940}}]L thanks H and tells him that he is waiting for K's permission to apply.
\ligne[{\hyperref[19400110]{10.01.1940}}]L applies for this military position.
\ligne[\textsl{01.02.1940}]K informs H that L has been incorporated.
\ligne[{\hyperref[19400203]{03.02.1940}}]Käthe Lorenzen does so too.
\ligne[{\hyperref[19400204]{04.02.1940}}]L does so himself.
\ligne[\textsl{06.02.1940}]H answers K and tells him that the process of L's candidature sent to Walther might take some time.
\ligne[{\hyperref[19400206]{06.02.1940}}]H answers Käthe Lorenzen along the same lines.
\ligne[{\hyperref[19400209]{09.02.1940}}]L asks H as treasurer of the Deutsche Mathematiker-Vereinigung to waive his member fee.
\ligne[{\hyperref[19400421]{21.04.1940}}]L sends an offprint to H and complains about his soul-killing service.
\ligne[{\hyperref[19400505]{05.05.1940}}]L tells H that his candidature is still on the way and gratefully accepts a proposition by H to join him at the Oberkommando der Kriegsmarine (OKM, Supreme Command of the Navy).
\ligne[{\hyperref[19400509]{09.05.1940}}]H answers that there is a misunderstanding: his proposition concerns a position as a cryptographer under Achim Teubner.
\ligne[{\hyperref[19400521]{21.05.1940}}]L thanks H and tells him that he applies for the position under Teubner.
\ligne[{\hyperref[19400611]{11.06.1940}}]L sends H some news of his application.
\ligne[{\hyperref[19400630]{30.06.1940}}]L reports on his participation in the Battle of France.
\ligne[\textsl{23.10.1940}]K greets L via H.
\ligne[\textsl{14.11.1940}]H tells K that he has not seen L for a while and that he hopes to meet him on 30 November.
\ligne[{\hyperref[19410426]{26.04.1941}}]L reports on his problems with the military hierarchy and on a little progress in group axiomatics made during a three-days arrest.
\ligne[{\hyperref[19410507H]{07.05.1941}}]H tells L what he thinks of his insubordination, stubbornness and unmilitariness, and presents Oswald Teichmüller as the example of a correct attitude. He warns L that he will also be judged in his scientific career according to his military attitude. He will not help L anymore, who should now prove his value as a soldier.
\ligne[{\hyperref[19410507]{07.05.1941}}]Käthe Lorenzen asks H for a good advice regarding her husband's difficulties.
\ligne[{\hyperref[19410517]{17.05.1941}}]L explains to H his attitude towards the military and asks for permission to write to him about his future assignment.
\ligne[{\hyperref[19420102]{02.01.1942}}]K attests that L fulfils the requirements for being admitted to habilitation.
\ligne[{\hyperref[19420318]{18.03.1942}}]L thanks H for his congratulations on the occasion of the birth of L's daughter Jutta and for his indication of Jean Dieudonné's interest in multiplicative ideal theory. He expresses satisfaction with his work as a teacher.
\ligne[{\hyperref[19420402]{02.04.1942}}]H asks L for advice on a manuscript by Hans R. Weber.
\ligne[{\hyperref[19420407]{07.04.1942}}]L answers that Weber's results are already known and can be found in Huntington's work.
\ligne[{\hyperref[19420512]{12.05.1942}}]The director of the Mathematical seminar in Bonn requests a statement of the leader of the union of lecturers on L's ideological and moral prerequisites for being appointed as assistant.
\ligne[{\hyperref[19420601]{01.06.1942}}]Ernst Klapp, the representing leader of the union of lecturers, states that his political attitude is unobjectionable, but that he presents certains deficiencies of character that make his admission to lecturership undesirable: a self-conceit that has also prejudiced his career in the army.
\ligne[{\hyperref[19430507]{07.05.1943}}]On the occasion of a review of \citealt{krull43}, L makes an analysis of the axiomatics of the star-operation.
\ligne[{\hyperref[19430920]{20.09.1943}}]L finds the proof method of \citealt{Lor1950}, the dynamical method in algebra.
\ligne[{\hyperref[19440104]{04.01.1944}}]K discusses the details of the habilitation process with L and projects that it take place at the beginning of the summer term.
\ligne[\textsl{06.02.1944}]K acknowledges receipt of a letter from H dated 23 January 1944, which must contain a strong criticism of the person and the work of~L\@. He defends L's work in lattice theory, but reckons that L deserves a repeated lesson, so that he proposes to suspend the process of L's habilitation.
\ligne[{\hyperref[19440206]{06.02.1944}}]In this letter to L, K expounds several obstacles to his habilitation.
\ligne[\textsl{19.02.1944}]K acknowledges receipt of a postcard by H that enables him to treat L's habilitation in agreement with H.
\ligne[{\hyperref[19440219]{19.02.1944}}]K writes to L that according to H his scientific publications do not suffice for letting him habilitate.
\ligne[{\hyperref[19440313]{13.03.1944}}]L answers that he should be judged on his habilitation manuscript and advocates his research in logic as stemming from the same motivation as in algebra.
\ligne[{\hyperref[19440401]{01.04.1944}}]K has begun to read L's habilitation manuscript and asks a question related to the condition of regularity (see page~\pageref{regularity}).
\ligne[{\hyperref[19440416]{16.04.1944}}]Upon having received a postcard by L, K denounces the lattice-ordered groups as being ``commutatively infected'' and expresses his disappointment that the scope of L's work does therefore not cover the full noncommutative generality. He uses this for reiterating his judgment that L's work is not sufficient for a habilitation.
\ligne[{\hyperref[19440425]{25.04.1944}}]L describes precisely his achievements, emphasises that his proof method is completely different from that of his Ph.D. thesis, and appends another manuscript in the hope of convincing K to change his mind.
\ligne[{\hyperref[19440505]{05.05.1944}}]L writes Scholz an account of his letter to K.
\ligne[{\hyperref[19440526]{26.05.1944}}]L accepts Scholz's proposal to contact Köthe.
\ligne[{\hyperref[19440529]{29.05.1944}}]K praises the clearness of L's letter and expresses the suspicion that the intricate style of his manuscripts might hide their relevance. He does not, however, change his mind.
\ligne[{\hyperref[19440602]{02.06.1944}}]L asks Scholz to send Köthe a copy of his manuscript.
\ligne[{\hyperref[19440606]{06.06.1944}}]L motivates his research as a quest for simplicity and conceptual clarification. He also describes the difficult conditions of his work and career.
\ligne[{\hyperref[194406]{\hfill06.1944}}]L thanks Scholz for his intercession with Köthe.
\ligne[{\hyperref[19440622]{22.06.1944}}]K writes that L's manuscript is a thorough failure and repeats his point with a reference to the requirement of high scientific quality.
\ligne[{\hyperref[19440628]{28.06.1944}}]L writes Scholz an account of K's letter.
\ligne[{\hyperref[19440709]{09.07.1944}}]L thanks Scholz for his encouragement; Köthe has written to him that he is willing to report on his habilitation; however, L sees no way to change K's negative judgment and to proceed with his project of habilitating.
\ligne[{\hyperref[19440716]{16.07.1944}}]K expresses his sympathy with L after the bombing of Wesermünde.
\ligne[{\hyperref[19441001]{01.10.1944}}]Although he is still not satisfied, K admits that L has improved his manuscript. He nevertheless suspects that L has not yet settled the noncommutative case, dealing only with a ``semi-commutative'' one.
\ligne[{\hyperref[19450725]{25.07.1945}}]L asks H whether he has some news on the whereabouts of K, Scholz, Ackermann, and Gentzen.
\ligne[{\hyperref[19450902]{02.09.1945}}]L reports on his political attitude under national socialism.
\ligne[{\hyperref[19450903]{03.09.1945}}]L fills out the Military government of Germany \emph{Fragebogen}, providing a ``chronological record of full-time employment and military service''.
\ligne[{\hyperref[19450906]{06.09.1945}}]Ernst Peschl certifies L's political attitude.
\ligne[{\hyperref[19450914]{14.09.1945}}]A board of examiners (Hellmuth von Weber, Hans Fitting, and Carl Troll) certifies L's political attitude.
\ligne[{\hyperref[19460607]{07.06.1946}}]L writes to Scholz that K agrees that he habilitate at once.
\ligne[{\hyperref[19460608]{08.06.1946}}]L writes to Scholz about his perspectives in Bonn.
\ligne[{\hyperref[19460813]{09.08.1946}}]L is habilitated by giving his inaugural lecture ``On the concept of lattice''.
\ligne[\textsl{15.08.1946}]K writes H on the occasion of L's habilitation and expresses his disagreement with H on the scientific value of L's works in lattice theory and logic. But he shares H's objections to L as person.
\ligne[\textsl{10.09.1946}]K distances himself from ``the conception of mathematics as a pure `theory of structure' in the sense of L (whom [he] tries otherwise to influence vigorously in the opposite direction)''.
\ligne[{\hyperref[19530418]{18.04.1953}}]K praises L in a letter to Scholz.
\ligne[{\hyperref[19530514]{14.05.1953}}]L writes to H on the occasion of his talk on mathematics as science, art, and power \citep{MR0053055}, and asks him to which extent he considers formalisation as a danger. An annotation shows that H discusses this with L on 3 July 1953.
\ligne[{\hyperref[19530605]{05.06.1953}}]H thanks L for his interest and announces that he will shortly pay a visit to Bonn and postpones his answer to that occasion.
\ligne[{\hyperref[19530609]{09.06.1953}}]L invites H to his home in anticipation of the latter’s visit to Bonn.
\ligne[{\hyperref[195907]{\hfill07.1959}}]L gives a sketch of proof of the assertion that the theory of commutative fields is undecidable.
\ligne[{\hyperref[19590801]{01.08.1959}}]H has checked L's sketch but for the logical conclusion.
\ligne[{\hyperref[19600307]{07.03.1960}}]H thanks L for sending a copy of \citealt{MR0111655}
  .
\ligne[{\hyperref[19610627]{27.06.1961}}]H thanks L for sending an offprint of \citealt{MR0147364} and invites him to give a couple of colloquia in Hamburg about his results on the foundations of mathematics.
\ligne[{\hyperref[19610704]{04.07.1961}}]L gratefully accepts and proposes to question how the ``assertions used as `axioms' are to be \emph{proved}''.
\ligne[{\hyperref[19610708]{08.07.1961}}]H organises L's visit to Hamburg.
\ligne[{\hyperref[19610711]{11.07.1961}}]L fixes a last detail of his visit to Hamburg.
\ligne[{\hyperref[19621009]{09.10.1962}}]H congratulates L on his appointment as professor in Erlangen and thanks him for the copy of \citealt{lorenzen62}, expressing the hope of finding the time to read it up to the consistency proof for arithmetic.
\ligne[{\hyperref[196307]{\hfill07.1963}}]L thanks H for sending him a copy of \citealt{MR0153659}. He considers H's construction of the completion of a valuated field as a flaw.
\ligne[{\hyperref[196309]{\hfill09.1963}}]L thanks H for his postcard from Notre Dame. Even though he did not want H to justify himself, he is happy that H agrees that ``when one invokes `all' sequences (e.g.\ of rational numbers), one always means only `sufficiently many'{}''. L illustrates his statement that ``controlling that everywhere really always only `sufficiently many' are used, however, is in [his] opinion not trivial'' with \citealt{MR0057856}.
\ligne[{\hyperref[19780221]{21.02.1978}}]Aubert writes to L with an inquiry about K's attitude toward $t$-ideals (see page~\pageref{t-ideal}), which he considers ``building blocks of general arithmetics. They seem to form the true arithmetical divisors [\dots]''.
\ligne[{\hyperref[19780306]{06.03.1978}}]L answers that he had almost no contact with K in Bonn for ``political reasons (K had prevented [his] habilitation during war)'' and because of his ``foundation-theoretic works that K did not acknowledge as `mathematical' achievements''; he lost track of $t$-ideals when he focussed on divisibility theory in domains. He asks Aubert in return whether this latter theory has been developed by some algebraist.
\ligne[{\hyperref[19780410]{10.04.1978}}]Aubert expands on the importance of $t$-ideals as universal objects and points out the only reference to \citealt{Lor1952} that he is aware of.
\ligne[{\hyperref[19790618]{18.06.1979}}]L thanks Aubert for sending him his article on divisibility theory and points out that it deals only with commutative groups: ``Was the enlargement to noncommutative groups perhaps only an intellectual luxury?''
\end{liste}

\section{The correspondence between Krull and Lorenzen, 1938.}

\selectlanguage{german}

\subsection{19.02.1938. Letter from Lorenzen to Krull.}\input{19380219-1-23-4-Lorenzen-Krull}
\subsection{08.03.1938. Letter from Lorenzen to Krull.}\input{19380308-1-23-6-Lorenzen-Krull}
\subsection{13.03.1938. Letter from Lorenzen to Krull.}\input{19380313-1-23-5-Lorenzen-Krull}
\subsection{18.03.1938. Letter from Lorenzen to Krull.}\input{19380318-1-23-7-Lorenzen-Krull}
\subsection{22.03.1938. Letter from Lorenzen to Krull.}\input{19380322-1-23-8-Lorenzen-Krull}
\subsection{09.06.1938. Letter from Lorenzen to Krull.}\input{19380609-1-23-9-Lorenzen-Krull}

\section{The reports on Lorenzen's thesis.}

\subsection{[24.05.1938. Report by Hasse.]}
\label{19380524}

\begin{tabular}[t]{@{}c@{}}
Mathematisches Institut\\
der Universität\\
Prof.\ Dr.\ Hasse.\\[6pt]
Siegel
\end{tabular}\bigskip\hfill
\begin{tabular}[t]{@{}c@{}}
Göttingen, den 24.\ Mai 1938.\\
Bunsenstraße 3/5
\end{tabular}

Gutachten über die Dissertation Lorenzen.\footnote{\selectlanguage{english}This report is written on the basis of a description given in a letter from Krull to Hasse dated 22 May 1938 as proposed by Hasse in a letter dated 2 March 1938 and gratefully accepted by Hasse in letters dated 7 March and 3 May 1938 \citep[§§ 1.32--1.35]{roquette04}.}\medskip

\noindent Seitdem Dedekind den Idealbegriff in die Arithmetik und Algebra
eingeführt hat, sind vielfach Verallgemeinerungen dieses Be\bs
griffes vorgenommen worden, mit dem Ziel die Struktur beliebi\bs
ger Integritätsbereiche ebenso einfach zu übersehen, wie man
die ganzen Zahlen eines algebraischen Zahlkörpers durch die
Dedekindschen Ideale übersieht. Lorenzen entwickelt eine zu\bs
sammenfassende Theorie aller dieser Idealbegriffe von einem
neuen einheitlichen Gesichtspunkt aus. Es erscheint dem Kenner
fast paradox, dass dabei die Addition, auf der doch der Dede\bs
kindsche Idealbegriff wesentlich beruht, vollständig ausser
Betracht gelassen wird, und nur die Multiplikation gebraucht
wird. Erst durch diese Lorenzensche Auffassung des Idealbegriffs
kommen die Zusammenhänge zwischen der idealtheoretischen und
der bewertungstheoretischen Behandlung des Strukturproblems
allgemeiner Integritätsbereiche klar und abgerundet heraus.
Lorenzen gibt neue verblüffend einfache Beweise für die in dieser Richtung liegenden Sätze von Krull, Prüfer u.\ a.

Es gelingt ihm ferner, die Identität des Prüferschen mit dem
Krullschen Idealsystem zu zeigen und dem Krullschen Hauptsatz
über Durchschnitte von allgemeinen Bewertungsringen (beliebige
nicht-archimedische Wertgruppe) einen entsprechenden Hauptsatz
über spezielle Bewertungsringe (archimedische Wertgruppe) an
die Seite zu stellen. Dieser Satz\footnote{\selectlanguage{english}\label{nonconclusive}The proofs of this and of the next assertion, which one can find in \citealt{Lor1938m}, are not conclusive: see Krull's letter to Hasse dated 1 September 1938 \citep[§~1.38]{roquette04}. As for this assertion, see footnote~\ref{nakayama} on page~\pageref{nakayama}. As for the next assertion, see the counterexample provided by Masayoshi \citet{MR47627,MR72866}.}
  %
  %
  %
  %
  %
  %
gibt also eine idealtheore\bs
tische Charakterisierung derjenigen Integritätsbereiche, die
bewertungstheoretisch durch das Positivsein von Exponenten\bs
bewertungen mit archimedischer Wertgruppe charakterisiert
sind. Ferner gibt Lorenzen auch eine idealtheoretische Cha\bs
rakterisierung der diskreten Bewertungsringe. Diese beiden
Charakterisierungen wurden von den Algebraikern lange gesucht.
Sie stellen wertvolle wissenschaftliche Ergebnisse dar. Auch
sonst finden sich in der Arbeit zahlreiche neue Einblicke in
die Beziehungen zwischen den verschiedenen in der Idealtheo\bs
rie studierten Begriffen.

Die Darstellung ist streng axiomatisch und leider von
äusserster Knappheit. Für den Druck erscheint mir eine etwas
ausführlichere Darstellung erwünscht und ausserdem auch eine
stärkere Anlehnung an die Literatur durch Verweise.

Da die Arbeit wesentlich neue und originelle Gedanken ent\bs
hält und wichtige Ergebnisse bringt, schlage ich die Annahme
mit dem Prädikat
\[\text{\uline{ausgezeichnet}}\]\nopagebreak
vor.\nopagebreak\medskip

\hfill Hasse

\subsection{[02.06.1938. Report by Siegel.]}
\label{19380602}

Eine sorgfältige Nachprüfung der Abhandlung hätte für mich eine Arbeit von
mehreren Monaten bedeutet, da der Text an vielen Stellen kaum zu verstehen ist
und sich wegen der mangelnden Literaturangaben auch nur schwer ergänzen
lässt. Ich muss mich deshalb eines genaueren Urteils über den Wert der Abhand\bs
lung enthalten. Das wenige, was ich mit Mühe habe verstehen können, macht
den Eindruck, als ob der Verfasser jedenfalls über gute mathematische Fähig\bs
keiten verfügt.\medskip

\hfil\hfil\hfil\hfil\hfil\hfil Siegel\hfil 1938 VI 2\qquad

\section{The correspondence between Hasse and Lorenzen, 1938--1942.}

\subsection{30.06.1938. Letter from Hasse to Lorenzen.}\input{19380630-Cod-Ms-H-Hasse-1-1022-Beil-1-Hasse-Lorenzen}
\subsection{06.07.1938. Letter from Lorenzen to Hasse.}\input{19380706-Cod-Ms-H-Hasse-1-1022-1-Lorenzen-Hasse}
\subsection{09.08.1939. Letter from Lorenzen to Hasse.}\input{19390809-Cod-Ms-H-Hasse-1-1022-2-Lorenzen-Hasse}
\subsection{11.08.1939. Letter from Hasse to Lorenzen.}\input{19390811-Cod-Ms-H-Hasse-1-1022-Beil-2-Hasse-Lorenzen}
\subsection{06.09.1939. Letter from Lorenzen to Hasse.}\input{19390906-Cod-Ms-H-Hasse-1-1022-3-Lorenzen-Hasse}
\subsection{11.10.1939. Letter from Lorenzen to Hasse.}\input{19391011-Cod-Ms-H-Hasse-1-1022-4-Lorenzen-Hasse}
\subsection{10.11.1939. Letter from Lorenzen to Hasse.}\input{19391110-Cod-Ms-Hasse-33-3-Lorenzen-Hasse}
\subsection{16.11.1939. Letter from Hasse to Lorenzen.}\input{19391116-Cod-Ms-Hasse-33-3-Hasse-Lorenzen}
\subsection{17.11.1939. Letter from Hasse to Lorenzen.}\input{19391117-Cod-Ms-H-Hasse-1-1022-4-2-Hasse-Lorenzen}
\subsection{13.12.1939. Letter from Lorenzen to Hasse.}\input{19391213-Cod-Ms-H-Hasse-1-1022-4-1-Lorenzen-Hasse}
\subsection{04.01.1940. Letter from Hasse to Lorenzen.}\input{19400104-Cod-Ms-H-Hasse-1-1022-Beil-3-Hasse-Lorenzen}
\subsection{09.01.1940. Letter from Lorenzen to Hasse.}\input{19400109-Cod-Ms-H-Hasse-1-1022-5-Lorenzen-Hasse}
\subsection{10.01.1940. Letter from Lorenzen to Hasse.}\input{19400110-Cod-Ms-H-Hasse-1-1022-6-Lorenzen-Hasse}
\subsection{03.02.1940. Postcard from Käthe Lorenzen to Hasse.}\input{19400203-Cod-Ms-H-Hasse-1-1021-Lorenzen-Hasse}
\subsection{06.02.1940. Letter from Hasse to Käthe Lorenzen.}\input{19400206-Cod-Ms-H-Hasse-1-1021-Beil-Hasse-Lorenzen}
\subsection{04.02.1940. Letter from Lorenzen to Hasse.}\input{19400204-Cod-Ms-H-Hasse-1-1022-7-Lorenzen-Hasse}
\subsection{09.02.1940. Postcard from Lorenzen to Hasse.}\input{19400209-Cod-Ms-H-Hasse-27-1-Lorenzen-Hasse}
\subsection{21.04.1940. Letter from Lorenzen to Hasse.}\input{19400421-Cod-Ms-H-Hasse-1-1022-8-Lorenzen-Hasse}
\subsection{05.05.1940. Letter from Lorenzen to Hasse.}\input{19400505-Cod-Ms-H-Hasse-1-1022-9-Lorenzen-Hasse}
\subsection{09.05.1940. Letter from Hasse to Lorenzen.}\input{19400509-Cod-Ms-H-Hasse-1-1022-Beil-4-Hasse-Lorenzen}
\subsection{21.05.1940. Letter from Lorenzen to Hasse.}\input{19400521-Cod-Ms-H-Hasse-1-1022-10-Lorenzen-Hasse}
\subsection{11.06.1940. Letter from Lorenzen to Hasse.}\input{19400611-Cod-Ms-H-Hasse-1-1022-11-Lorenzen-Hasse}
\subsection{30.06.1940. Letter from Lorenzen to Hasse.}\input{19400630-Cod-Ms-H-Hasse-1-1022-12-Lorenzen-Hasse}
\subsection{26.04.1941. Letter from Lorenzen to Hasse.}\input{19410426-Cod-Ms-H-Hasse-1-1022-13-Lorenzen-Hasse}
\subsection{07.05.1941. Letter from Hasse to Lorenzen.}\input{19410507-Cod-Ms-H-Hasse-1-1022-Beil-5-Hasse-Lorenzen}
\subsection{07.05.1941. Letter from Käthe Lorenzen to Hasse.}\input{19410507-Cod-Ms-H-Hasse-1-1021-Lorenzen-Hasse}
\subsection{17.05.1941. Letter from Lorenzen to Hasse.}\input{19410517-Cod-Ms-H-Hasse-1-1022-15-Lorenzen-Hasse}
\subsection{18.03.1942. Letter from Lorenzen to Hasse.}\input{19420318-Cod-Ms-H-Hasse-1-1022-14-Lorenzen-Hasse}
\subsection{02.04.1942. Letter from Hasse to Lorenzen.}\input{19420402-Cod-Ms-Hasse-33-3-Hasse-Lorenzen}
\subsection{07.04.1942. Letter from Lorenzen to Hasse.}\input{19420407-Cod-Ms-Hasse-33-3-Lorenzen-Hasse}

\section{Documents relating to Lorenzen's career, 1942.}

\subsection{02.01.1942. Attestation by Krull for Lorenzen.}\input{19420102-PL-1-1-136-Krull}
\subsection{12.05.1942. Request of statement from the director of the Mathematical seminar in Bonn.}\input{19420512}
\subsection{01.06.1942. Statement by the Dozentenführer.}\input{19420601}

\section{The correspondence between Krull and Lorenzen, 1943--1944.}

\subsection{07.05.1943. Letter from Lorenzen to Krull.}\input{19430507-PL-1-1-140-Lorenzen-Krull}
\subsection{20.09.1943. Letter from Lorenzen to Krull.}\input{19430920-PL-1-1-139-Lorenzen-Krull}
\subsection{04.01.1944. Postcard from Krull to Lorenzen.}\input{19440104-PL-1-1-149-Krull-Lorenzen}
\subsection{06.02.1944. Postcard from Krull to Lorenzen.}\input{19440206-PL-1-1-145-Krull-Lorenzen}
\subsection{19.02.1944. Letter from Krull to Lorenzen.}\input{19440219-PL-1-1-143-Krull-Lorenzen}
\subsection{13.03.1944. Letter from Lorenzen to Krull.}\input{19440313-PL-1-1-131-Lorenzen-Krull}
\subsection{01.04.1944. Postcard from Krull to Lorenzen.}\input{19440401-PL-1-1-144-Krull-Lorenzen}
\subsection{16.04.1944. Letter from Krull to Lorenzen.}\input{19440416-PL-1-1-142-Krull-Lorenzen}
\subsection{25.04.1944. Letter from Lorenzen to Krull.}\input{19440425-PL-1-1-132-Lorenzen-Krull}
\subsection{29.05.1944. Postcard from Krull to Lorenzen.}\input{19440529-PL-1-1-146-Krull-Lorenzen}
\subsection{06.06.1944. Letter from Lorenzen to Krull.}\input{19440606-PL-1-1-133-Lorenzen-Krull}
\subsection{22.06.1944. Letter from Krull to Lorenzen.}\input{19440622-PL-1-1-141-Krull-Lorenzen}
\subsection{16.07.1944. Postcard from Krull to Lorenzen.}\input{19440716-PL-1-1-148-Krull-Lorenzen}
\subsection{01.10.1944. Postcard from Krull to Lorenzen.}\input{19441001-PL-1-1-147-Krull-Lorenzen}

\section{A postcard from Lorenzen to Hasse, 1945.}

\subsection{25.07.1945. Postcard from Lorenzen to Hasse.}\input{19450725-Cod-Ms-H-Hasse-1-1022-16-Lorenzen-Hasse}

\section{Documents relating to Lorenzen's career, 1945--1946.}

\subsection{02.09.1945. Report by Lorenzen on his political attitude.}\input{19450902-universitaetsarchiv_bonn-lorenzen-bericht_ueber_meine_politische_einstellung}
\subsection{03.09.1945. Military government of Germany Fragebogen: chronological record of full-time employment and military service.}\input{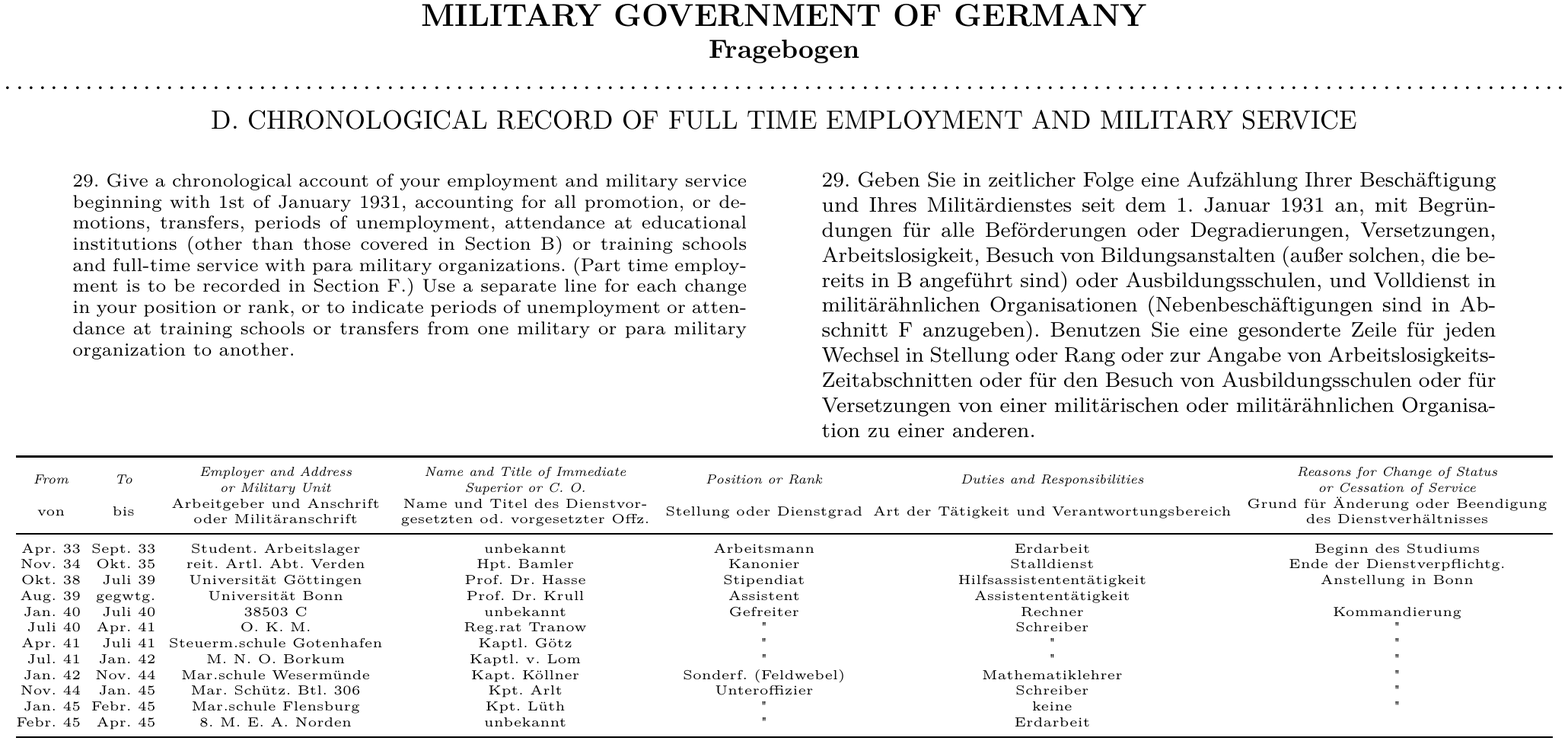}
\subsection{06.09.1945. Certificate by Ernst Peschl on Lorenzen's political attitude.}\input{19450906-universitaetsarchiv_bonn-peschl-lorenzen-gutachten}
\subsection{14.09.1945. Certificate by a board of examiners.}\input{19450914-universitaetsarchiv_bonn-pruefungsauschuss_weber_fitting_troll-lorenzen_gutachten}
\subsection{13.08.1946. Notification on Lorenzen's inaugural lecture of 9 August 1946.}\input{19460813}

\section{A letter from Krull to Scholz, 1953.}

\subsection{18.04.1953. Letter from Krull to Scholz.}\input{19530418-ko-05-0647-Krull-Scholz}

\section{The correspondence between Hasse and Lorenzen, 1953--1963.}

\subsection{14.05.1953. Letter from Lorenzen to Hasse.}\input{19530514-Cod-Ms-H-Hasse-1-1022-17-Lorenzen-Hasse}
\subsection{05.06.1953. Letter from Hasse to Lorenzen.}\input{19530605-Cod-Ms-H-Hasse-1-1022-Beil-6-Hasse-Lorenzen}
\subsection{09.06.1953. Letter from Lorenzen to Hasse.}\input{19530609-Cod-Ms-H-Hasse-1-1022-18-1-Lorenzen-Hasse}
\subsection{July 1959. Letter from Lorenzen to Hasse.}\input{195907-Cod-Ms-H-Hasse-1-1022-18-2-Lorenzen-Hasse}
\subsection{01.08.1959. Letter from Hasse to Lorenzen.}\input{19590801-Cod-Ms-H-Hasse-1-1022-Beil-7-Hasse-Lorenzen}
\subsection{07.03.1960. Letter from Hasse to Lorenzen.}\input{19600307-Cod-Ms-H-Hasse-1-1022-Beil-8-Hasse-Lorenzen}
\subsection{27.06.1961. Letter from Hasse to Lorenzen.}\input{19610627-Cod-Ms-H-Hasse-1-1022-Beil-9-Hasse-Lorenzen}
\subsection{04.07.1961. Letter from Lorenzen to Hasse.}\input{19610704-Cod-Ms-H-Hasse-1-1022-19-Lorenzen-Hasse}
\subsection{08.07.1961. Letter from Hasse to Lorenzen.}\input{19610708-Cod-Ms-H-Hasse-1-1022-Beil-10-Hasse-Lorenzen}
\subsection{11.07.1961. Letter from Lorenzen to Hasse.}\input{19610711-Cod-Ms-H-Hasse-1-1022-20-Lorenzen-Hasse}
\subsection{09.10.1962. Letter from Hasse to Lorenzen.}\input{19621009-Cod-Ms-H-Hasse-1-1022-Beil-11-Hasse-Lorenzen}
\subsection{July 1963. Letter from Lorenzen to Hasse.}\input{196307xx-Cod-Ms-H-Hasse-1-1022-21-Lorenzen-Hasse}
\subsection{September 1963. Letter from Lorenzen to Hasse.}\input{196309xx-Cod-Ms-H-Hasse-1-1022-22-Lorenzen-Hasse}

\section{The correspondence between Aubert and Lorenzen, 1978--1979.}

\selectlanguage{english}
\subsection{21.02.1978. Letter from Aubert to Lorenzen.}\input{19780221-PL-5-Aubert-Lorenzen}
\selectlanguage{german}
\subsection{06.03.1978. Letter from Lorenzen to Aubert.}\input{19780306-PL-5-Lorenzen-Aubert}
\selectlanguage{english}
\subsection{10.04.1978. Letter from Aubert to Lorenzen.}\input{19780410-PL-5-Aubert-Lorenzen}
\selectlanguage{german}
\subsection{18.06.1979. Letter from Lorenzen to Aubert.}\input{19790618-PL-5-Lorenzen-Aubert}

\selectlanguage{english}

\renewcommand\bibname{References.}


\end{document}

%% file: 19380219-1-23-4-Lorenzen-Krull

\label{19380219}%
\hfill Göttingen, 19.2.\nopagebreak

\noindent\hfil Sehr geehrter Herr Professor!\nopagebreak\medskip

Ich danke Ihnen für Ihr freundliches
Interesse, das Sie meiner Arbeit ent\bs
gegenbringen,\footnote{\selectlanguage{english}This interest is subsequent to a letter from Hasse to Krull dated 14 February 1938 \citep[§~1.31]{roquette04}.} und besonders für die
Erlaubnis mit Ihnen persönlich darüber
sprechen zu dürfen.

Da ich aus verschiedenen Gründen erst
am Mittwoch von hier fahren kann,\bruch{}
würde ich mir gestatten, Sie am Donners\bs
tag, dem 24.2., morgens aufzusuchen.
Falls Ihnen diese Zeit ungelegen
sein sollte, können Sie mir ja dann
eine andere Zeit sagen.

Ich freue mich sehr, Sie befragen
zu dürfen, und hoffe mir, daß ich Ihnen
nicht allzu lästig fallen werde.\strut

\hfil%
\begin{tabular}{@{}r@{}}%
  Ihr sehr ergebener\\%
  Lorenzen\rlap.%
\end{tabular}

%% file: 19380308-1-23-6-Lorenzen-Krull

\hfill%
\begin{tabular}{@{}l@{}}%
\label{19380308}%
  Gött., 8.3.\\%
  Feuerschanzengraben 20a%
\end{tabular}\medskip

\hfil Sehr geehrter Herr Professor!\medskip

Ich möchte Ihnen zunächst noch einmal
herzlichst danken für die so sehr freundli\bs
che Aufnahme, die ich bei Ihnen gefunden
habe.\footnote{\selectlanguage{english}This is corroborated by the correspondence between Krull and Hasse dated 2~and~7 March 1938 \citep[\S\S~1.32, 1.33]{roquette04}.} Ich danke auch besonders Ihrer Frau
Gemahlin. Es war für mich eine große
Freude mit Ihnen über meine Arbeit
sprechen zu können.

Was diese anbetrifft, so muß ich lei\bs
der folgendes richtigstellen:\bruch

Sei $\su I$ eine $v$-abg[eschlossene] Halbgruppe, $\su c$ ein
endliches $v$-Ideal, so folgt ($\su I=(1)=1$)%
\[1:\su c^{-1}\su c=(1:\su c^{-1}):\su c=\su c:\su c=1\]%

D.\ h.\ $\su c$ ist $v$-umkehrbar, aber $\su c^{-1}$ braucht nicht
endlich zu sein.\comment{See \citealt[p.~15]{Lor1938m}; \citealt[p.~538]{Lor1939}.} Sind $\su a, \su b$ endliche $v$-Ideale
so folgt $\su a(\su a+\su b)^{-1}+\su b(\su a+\su b)^{-1}=1$ ($v$-Summe).

Damit $\su a\subseteq\su b$ oder $\su b\subseteq\su a$ folgt, muß in~$\su I$ die
$v$-Summe zweier echter $v$-Ideale wieder echt
sein. Ist $\su I$ primär,\comment{\label{p550}See \citealt[p.~32]{Lor1938m}; \citealt[p.~550]{Lor1939}.} so folgt dies nur für
endliche $v$-Ideale. Wenn etwa gelten
würde $(\su a+\su b)_v=(\su a_v\cup\su b_v)_{v_e}$ für beliebiges
$\su a,\su b$ käme man auch durch. Aber so ohne
weiteres nicht.\bruch

Übrigens sind die $v_e$-Ideale schon in
der Arnold\,schen Arbeit\comment{\citealt{arnold29}.} definiert, wie ich
jetzt gesehen habe.

Für weitergehende Sätze besteht die
Schwierigkeit, daß man aus der $v$-Ab\-ge\-schlos\bs
sen\-heit von~$\su I$ nicht auch die von $\su I_S$ schließen
kann, da aus $\alpha\in(\su c\su I_S)_v$ nicht $\alpha s\in(\su c\su I)_v$ \bsp
$s\in S$ folgt (auch nicht bei endlichem~$\su c$!).

Darf ich hier mir vielleicht eine Frage
erlauben: Ist in einer speziellen
Hauptordnung, stets auch $\su I_S$ spezielle H.?

Auf Ihre Anregung hin, habe ich mir
die Arbeit von Akizuki\comment{Most likely \citealt{akizuki35}.}
noch einmal\bruch{}
angesehen, und habe die Kettensätze
für $v$-Ideale untersucht. Es zeigt sich
sofort: die $v$-Maximalbedingung ist
stets mit der schwachen $v$-Min.bed\@.
gleichwertig, denn sind $\su a_n$ \ $n=1,\dots$
irgendwelche Ideale mit $(a)\subseteq\su a_n\subseteq 1$
so folgt zunächst $(a)\subseteq(a)\su a_n^{-1}\subseteq1$ und
für $v$-Ideale folgt aus $\su a_n\subset\su a_{n+1}$ bzw\@.
$\su a_n\supset\su a_{n+1}$ sofort $a\su a_n^{-1}\supset a\su a_{n+1}^{-1}$ bzw\@.
$a\su a_n^{-1}\subset a\su a_{n+1}^{-1}$. \ q.e.d. \ Weiter folgt
aus der $v$-Max.bed., daß kein $v$-Id[eal] eine
unendliche Basis hat und umgekehrt.
(Es kann allgemein aber sein, daß ein Ideal gleichzeitig
eine endliche und eine unendliche Basis hat.)\bruch

Hieraus ergibt sich ein neuer, in
allen Schritten trivialer, Beweis
für den Z.~P.~I.\footnote{``\foreignlanguage{german}{Zerlegungssatz in Primideale}'', see \citealt[12]{krull35}.} aus der schw\@. Min.bed\@.
und der totalen Abgeschlossenheit.

Denn jedes Primideal ist minimal,
also $v_e$-Ideal, also $v$-Ideal wegen der
$v$-Max.bed. Aus $\su a\su a^{-1}\subset1$ würde also 
$(\su a\su a^{-1})_v\subset1$ folgen. \ q.e.d.

Umgekehrt folgt sofort die tot\@. Abg\@. und
die Max.bed.; da jedes Ideal $v$-Ideal ist
also auch die schw\@. Min.bed.\comment{See \citealt[§~5]{Lor1938m}; \citealt[§~2]{Lor1939}.}\bruch

Ich bitte um Entschuldigung, daß dieser
Brief so lang ist, aber ich dachte, es würde
Sie interessieren, denn noch einfacher kann
der Beweis jetzt wohl nicht mehr gemacht
werden.

Ich bitte um eine Empfehlung an
Ihre Frau Gemahlin und um einen freund\bs
lichen Gruß an die beiden Kleinen.\strut\medskip

\hfil%
\begin{tabular}{@{}lll@{}}%
  Ihr\\%
  &sehr ergebener\\%
  &&\llap{Paul }Lorenzen.%
\end{tabular}

%% file: 19380313-1-23-5-Lorenzen-Krull

\hfill%
\begin{tabular}{@{}r@{}}
\label{19380313}%
  Göttingen, 13.3.38\\
  Feuerschanzengraben 20a
\end{tabular}\medskip

\hfil Sehr geehrter Herr Professor!\medskip

Diesmal muß ich wirklich um Entschul\bs
digung bitten, daß ich Sie schon wieder
belästige. Aber es hat sich eine wesent\bs
liche Vereinfachung ergeben, die ich Ihnen
gerne mitteilen möchte.\smallskip

{\newsavebox{\ansbox}%
\newlength{\answd}%
\sbox{\ansbox}{\uline{Definition}:\enskip}%
\setlength{\answd}{\wd\ansbox}%
\parindent-\answd%
\leftskip\answd%
\usebox\ansbox Ein Ideal heißt $\overline{v}$-Ideal, wenn
  es mit zwei $v$-Idealen auch
  deren $v$-Summe enthält.\smallskip

}\noindent\uline{Satz}: Es gibt stets maximale $\overline v$-Ideale.\smallskip\bruch

{\sbox{\ansbox}{\uline{Definition}:\enskip}%
\setlength{\answd}{\wd\ansbox}%
\parindent-\answd%
\leftskip\answd%
\usebox\ansbox Eine Halbgruppe~$\su g$ (Quotienten\bs
  gruppe~$\su G$) heißt primär, wenn
  es keine Quotientenhalbgruppe~%
  $\su g_S$ mit $\su g\subset\su g_S\subset\su G$ gibt.\smallskip

}{\sbox{\ansbox}{\uline{Satz}:\enskip}%
\setlength{\answd}{\wd\ansbox}%
\parindent-\answd%
\leftskip\answd%
\usebox\ansbox Eine primäre total abgeschlossene
Halbgruppe ist linear.\smallskip

}Die Beweise sind ganz einfach, da
mit den $\overline v$-Idealen jetzt wohl die
zweckentsprechenden gefunden sind.\bruch

Besonders bemerkenswert finde ich,
daß alle drei Begriffe des Satzes,
primär, total abgeschlossen, linear
ohne Benutzung eines speziellen Ideal\bs
systemes definiert sind. Es würde
allerdings genügen die $v$-Abgeschlossen\bs
heit zu fordern.\strut\medskip

\hfil%
\begin{tabular}{@{}r@{}}%
  Mit freundlichem Gruß\\%
  Ihr\hphantom{ sehr ergebener}\\%
  sehr ergebener\\%
  Paul Lor\rlap{enzen.}%
\end{tabular}%


%% file: 19380318-1-23-7-Lorenzen-Krull

\hfill%
\begin{tabular}{@{}r@{}}%
\label{19380318}%
  Göttingen, 18.3.\\  
  Feuerschanzengraben 20a%
\end{tabular}\medskip

\hfil Sehr geehrter Herr Professor!\medskip

Endlich kann ich Ihnen mitteilen, daß
jede total abgeschlossene Halbgruppe
spezielle Hauptordnung ist.\comment{\label{nakayama}The argument for this assertion, which one can also find in \citealt{Lor1938m}, is not conclusive. Tadasi \citet{MR14067,MR14068,MR41115} shows that the assertion does not hold in general.}

Sei $\su g$ eine Halbgruppe, $a,b,\cdots$ Elemente
der Quotientengruppe~$\su G$.

Statt $\frac ab\in\su g$ schreibe ich $a\leqslant b$.\smallskip

{\parindent-1em\leftskip1em%
  \uline{Definition}:\enskip Eine Halbgruppe heißt vollstän\bs
  dig, wenn es zu zwei Elementen $a,b$
  einen g.g.T.~$a\lor b$ gibt.\comment{See \citealt[Definition~5]{Lor1939}, reproduced on page~\pageref{1939def5}.}\smallskip

}Die Existenz des k.g.V.\ $a\land b=\frac{ab}{a\lor b}$ folgt.\smallskip\bruch

{\parindent-1em\leftskip1em%
  \uline{Definition}: Eine Untermenge~$\su a$ von~$\su G$
  heißt $g$-Modul, wenn sie zu jedem Ele\bs
  ment~$a$ alle Elemente $b\le a$ und zu
  zwei Elementen~$a,b$ stets $a\lor b$ enthält.
  Ein $g$-Modul~$\su a$ mit $\su a^{-1}\ne0$ heißt $g$-Ideal.\smallskip
  
}Die Ideale $\su a\subseteq\su g$ heißen ganz, die Ideale
$\su a\subset\su g$ echt. Ein ganzes Ideal~$\su p$ heißt Primideal,
wenn die nicht in~$\su p$ liegenden Elemente von\nbsp
$\su g$ eine Halbgruppe bilden.\smallskip

\noindent\uline{Satz}:\enskip Jedes echte Ideal besitzt ein Primoberideal.\smallskip\comment{See \citealt[Satz~4]{Lor1939}.}

{\sbox{\ansbox}{\uline{Satz}:\enskip}%
\setlength{\answd}{\wd\ansbox}%
\parindent-\answd%
\leftskip\answd%
\uline{Satz}:\enskip Die Primidealquotientenhalbgruppen
sind linear.\smallskip\comment{See \citealt[Satz~10]{Lor1939}, reproduced on page~\pageref{1939satz10}.}

\uline{Satz}:\enskip $\su g$ ist Hauptordnung.%
\[\su g=\bigcap_\tau\su L_\tau\quad\su L_\tau\text{ linear}\]\bruch

}Ich bezeichne mit~$a_\tau$ das Element~$a$
als Element von $\su L_\tau$ (den "`Wert"' von~$a$).
Statt~$1_\tau$ einfach~$1$.

Darf~$\su L_\tau$ bei der Darstellung von~$\su g$ nicht
weggelassen werden, so gibt es ein Ele\bs
ment~$a$ mit $a_\tau>1$, $a_\lambda\le1$ für~$\lambda\ne\tau$.
Und zwar gibt es ein solches Element, der\bs
art, daß $a_\tau^{-1}$ nur ein Primoberideal hat,
da man sonst statt~$\su L_\tau$ eine Primideal\bs
quotientenhalbgruppe~$\su L\supset\su L_\tau$ von~$\su L_\tau$ nehmen
könnte.

Besitzt $\su L_\tau$ mehrere Primideale, so gibt
es auch ein Element~$b$ mit $b_\tau\le1$, das meh\bs
rere Primoberideale hat. Dann wird%
\[(1\land b)(1\lor a)^n\text{ ganz für alle }n\]%
$1\lor a$ nicht ganz.\smallskip\bruch

{\parindent-1em\leftskip1em%
  \uline{Satz}:\enskip Jede vollständige total abgeschlossene
  Halbgruppe ist spezielle Hauptordnung.\smallskip

}Ist $\su g$ beliebig total abgeschlossen, so
bilden die ganzen $v$-Ideale eine voll\bs
ständige total abgeschlossene Halbgruppe
womit alles bewiesen ist.

Ich möchte Ihnen, gerade nachdem dieser
Beweis jetzt geglückt ist, nochmals herzlichst
danken für die schönen Tage in Erlangen,
und die ausführlichen, für mich so wert\bs
vollen Besprechungen mit Ihnen.\strut\medskip

\hfil%
\begin{tabular}{@{}l@{}}%
  Mit freundlichem Gruß\\
  \qquad Ihr sehr ergebener\\%
  \qquad\qquad\qquad Paul Lorenzen.
\end{tabular}


%% file: 19380322-1-23-8-Lorenzen-Krull

\hfill%
\begin{tabular}{@{}r@{}}%
\label{19380322}%
  Göttingen, 22.3.\\%
  Feuerschanzengraben 20a%
\end{tabular}\medskip

\hfil Sehr geehrter Herr Professor!\medskip

Vielen Dank für Ihre Karte. Anschließend
das gewünschte Lexikon:

Es sei $\su g$ eine Halbgruppe, $\su G$ die Quotienten\bs
gruppe, $a,b,\dots\in\su G$.\medskip

{\sbox{\ansbox}{\uline{Definitionen}:\ }%
\setlength{\answd}{\wd\ansbox}%
\parindent-\answd%
\leftskip\answd%
\usebox\ansbox $a$ ganz, wenn $a\in\su g$\\[\smallskipamount]%
$a$ Nichteinheit, wenn \begin{tabular}[t]{@{}l@{}}%
                         $a$ ganz\\%
                         $a^{-1}$ nicht ganz%
                       \end{tabular}\\[\smallskipamount]%
  $a\le b$, wenn $\frac ab$ ganz\smallskip
  
}{\sbox{\ansbox}{$\su g$ }%
\setlength{\answd}{\wd\ansbox}%
\parindent-\answd%
\leftskip\answd%
  $\su g$ linear, wenn stets $a\le b$ oder $a\ge b$\smallskip

  $\su g$ einstufig linear, wenn $\su g$ linear und
  es zu jedem ganzen~$a$ und jeder Nichtein\bs
  heit~$b$ einen Exponenten~$n$ mit%
  \[b^n\le a\quad\text{gibt.}\]\bruch

  $\su g$ Hauptordnung, wenn $\su g$ Durchschnitt
  von linearen Oberhalbgruppen ist.\smallskip

  $\su g$ spezielle Hauptordnung, wenn $\su g$ Durchschnitt
  von einstufig linearen Oberhalbgruppen ist.\smallskip

  $\su g$ vollständig, wenn jedes endliche
  $v$-Ideal Hauptideal ist, d.\ h., wenn
  es zu je zwei Elementen~$a,b$ stets ei\bs
  nen g.~g.~T.~$a\lor b$ gibt.\smallskip

  $\su g$ total abgeschlossen, wenn aus $ca^n$ ganz
  für alle~$n$, folgt $a$ ganz.\medskip

}Für Integritätsbereiche entspricht
  \[%
    \begin{tabular}{rcl}%
      linear&---&Bewertungsring\\%
      einstufig linear&---&spezieller Bewertungsring\\%
      total abgeschlossen&---&vollständig ganz abgeschlossen\\\noalign{\bruch}%
    \end{tabular}%
  \]%

  Das Prinzip, das dem Beweis des Satzes:\smallskip

  {\leftskip2\parindent\parindent0pt%
    Jede total abgeschlossene Halbgruppe~$\su g$
    ist spezielle Hauptordnung.\smallskip

  }\noindent zu Grunde liegt, ist der Übergang von
  der gegebenen Halbgruppe~$\su g$ zur Halb\bs
  gruppe~$G$ der ganzen $v$-Ideale von~$\su g$. \bsp
  $G$~muß ebenfalls total abgeschlossen sein,
  und da $G$~vollständig ist, erweist sie sich
  leicht als spezielle Hauptordnung. Aus
  der Durchschnittsdarstellung von~$G$ folgt
  dann die Darstellung v[on]~$\su g$ als spezielle
  Hauptordnung. Ist $\su g$ ein Integritätsbereich
  so sind die linearen Oberhalbgruppen der\bruch{}
  Durchschnittsdarstellung von~$\su g$, sogar Be\bs
  wertungsringe. In Ihrer Terminologie,
  sind es gerade die Bewertungsringe, die
  zum $v$-Idealsyst., als arithmetisch brauch\bs
  barem System gehören.

  Das Unglück der verschiedenen Terminologien
  liegt vor allem in "`vollständig"' "`total ab\bs
  geschlossen"' und "`vollständig ganz abgeschlossen"'.\strut\medskip

\hfil%
\begin{tabular}{@{}l@{}}%
  Mit freundlichem Gruß\\
  \qquad\qquad Ihr sehr ergebener\\%
  \qquad\qquad\qquad Paul Lorenzen.
\end{tabular}


%% file: 19380609-1-23-9-Lorenzen-Krull

\label{19380609}%
\hfill Göttingen, 9.6.\nopagebreak

\noindent\hfil Sehr geehrter Herr Professor!\nopagebreak\medskip

Für Ihren Brief danke ich Ihnen viel\bs
mals. Etwas Besseres konnte ich mir
ja gar nicht wünschen, als daß ich Sie
noch einmal besuchen darf. Daß die
Arbeit so knapp gefaßt ist,\comment{Compare the criticism in Hasse's letter to Krull dated 31 May 1938 \citep[§~1.36]{roquette04}.} lag z.\ T.\ auch
am Zeitmangel, da ich ja dieses Semester
Examen machen wollte. Wenn dieses
nun vorbei ist, werde ich hoffentlich\bruch{}
klarer, ausführlicher und deutlicher be\bs
schreiben können, was ich meine.

Ehe ich komme, gebe ich Ihnen noch
Nachricht.

Mit Ihrem Diskriminantenkriterium\footnote{This probably refers to the ``Allgemeiner Diskriminantensatz'' of \citealt{krull39}.} %
ist es mir in der Tat genau so ergangen,
wie Sie vermuten. Ich wußte nur
noch, daß es den allgemeinst-möglichen
Fall enthielt, wußte aber nicht mehr die
genaue Formulierung. Ob sich der Hasse\bs
sche Beweis verallgemeinern läßt, weiß ich\bruch{}
nicht, ist ja aber auch nicht so wichtig,
da Ihr Resultat vorliegt.

Ich freue mich schon sehr darauf, Sie
besuchen zu dürfen.\strut\medskip

\hfil%
\begin{tabular}{@{}l@{}}%
  Mit freundlichem Gruß\\
  \qquad\qquad Ihr sehr ergebener\\%
  \qquad\qquad\qquad Paul Lorenzen.
\end{tabular}


%% file: 19380630-Cod-Ms-H-Hasse-1-1022-Beil-1-Hasse-Lorenzen
\label{19380630}%
\hfill 30.\ Juni 1938.\nopagebreak

\noindent Prof.\ Dr.\ Hasse.

\hfill\begin{tabular}{@{}l@{}}%
        Herrn\\[6pt]%
        Dr.\ Paul Lorenzen\\[6pt]%
        \uline{z.\ Zt.\ Erlangen.}\\%
        per Adresse Herrn Prof.\ Krull\\%
        Burgbergstr.~53.%
      \end{tabular}\medskip

\noindent Lieber Herr Lorenzen,\nopagebreak\medskip

wie ich eben sehe, ist es erwünscht, wenn Sie sich wegen
Ihres Antrages an das Reichsdozentenwerk hier dem Dozenten\bs
bundsführer vorstellen. Diesen müssen Sie dann auch fragen,
ob Sie sich auch sonst noch bei einem Referenten des Do\bs
zentenbundes vorzustellen haben. Die Sache ist nämlich die,
dass in den Richtlinien von der Vorstellung bei dem Referen\bs
ten die Rede ist. So viel ich weiss, ist aber hier in Göttin\bs
gen gerade ein Wechsel eingetreten. Prof.\ Schriel, der dieses
Amt bisher hatte, gibt es an Prof.\ Mattiat ab. Ich weiss aber
nicht genau, zu welchem Zeitpunkt dieser Übergang erfolgt.
Auch ist mir bekannt, dass Prof.\ Blume, der Dozentenbunds\bs
führer, stets Wert auf persönliche Vorstellung der Leute legt,
über die er zu berichten hat. Es hat natürlich Zeit, wenn
Sie dies nach Ihrer Rückkehr aus Erlangen tun, allerdings
bleibt Ihr Gesuch dann so lange liegen.

\noindent Im Anschluss an meine letzten Bemerkungen im Seminar habe
ich eine Bitte an Sie. Könnten Sie mir wohl mit Hilfe Ihrer
idealtheoretischen Fähigkeiten bestimmen, welches die maxima\bs
len Primideale des folgenden Integritätsbereiches sind: \bruch
$k$~sei ein algebraischer Zahlkörper, $K$~ein Funktionenkörper
einer Unbestimmten über~$k$ als Konstantenkörper, $\mathfrak{O}$~ein Prim\bs
divisor ersten Grades von~$K$ bez.~$k$ und $I$~sei die ganzal\bs
gebraisch abgeschlossene Hülle von~$i[x]$ in~$K$, wo $i$~die
Maximalordnung von~$k$ ist.
Dieser Integritätsbereich~$I$ hängt ersichtlich nur von~$\mathfrak{O}$ ab.
Mit Bestimmung der maximalen Primideale meine ich eine voll\bs
ständige Übersicht über diese, von den Primidealen von~$k$
und den Primdivisoren von~$K$ bez.~$k$ aus.

Für Ihre Arbeit mit Herrn Krull die besten Wünsche.

Ich bitte diesen und sich selbst recht herzlich zu grüssen\strut\medskip

\hfil\hfil Ihr [Hasse]
 

%% file: 19380706-Cod-Ms-H-Hasse-1-1022-1-Lorenzen-Hasse
\label{19380706}%
\hfill Erlangen, 6.7.38.\nopagebreak\smallskip

\noindent\hfil Sehr geehrter Herr Professor!\nopagebreak\medskip

\noindent Für Ihre beiden Briefe herzlichen Dank.
Da ich erst seit gestern hier in Erlangen
bin, konnte ich nicht eher antworten.

Die Referate über die Arbeiten von
Krull will ich gern übernehmen, wenn
dies noch Zeit bis Ende Juli, Anfang
August hat, da ich nicht eher zurück \bruch
sein werde. [Ich möchte nämlich erst
eine Reise durch Bayern anschließen,
das ich gar nicht kenne.]

Falls dies zu spät ist, müßten Sie
mir schon die Formblätter hierher\bs
schicken.

Zu Ihrer mathematischen Anfrage, kann
ich Ihnen -- nach Besprechung mit Herrn
Prof.\ Krull -- nur folgendes mitteilen:
Sei -- nur der Einfachheit halber -- $\Gamma$~der \bruch
ganz rationale Zahlring, $x$~eine Unbestimmte.
Die maximalen Primideale von~$\Gamma[x]$ bestimmen sich
dann eindeutig durch eine Primzahl~$p$
und ein Primelement von~$P_p[x]$. Eine
eindeutige Bestimmung aus den Primzahlen
und den Primelementen von~$P[x]$ ist
dagegen nicht möglich: z.\ B. $(2,x)=(2,x-2)$.

Ich danke Ihnen für die Mitteilung,
dass ich mich zunächst bei Herrn Prof.\ Blume
vorstellen muß. Nach meiner Rückkehr \bruch
werde ich dies sofort tun.

An meiner Arbeit gibt's hier viel
zu tun. Augenblicklich bin ich dabei, für
Köthe\comment{Gottfried Köthe (1905--1989) is professor in Münster from 1937 to 1943, and then until 1946 in Gießen.} einen verbandstheoretischen Exzerpt
zu machen.\medskip

Mit freundlichem Gruß,
auch von Herrn Prof.\ Krull\strut

\hfil\hfil\begin{tabular}{@{}ll@{}}%
  Ihr&\\%
     &sehr ergebener%
\end{tabular}\nopagebreak

\hfill Paul Lorenzen.\par
 

%% file: 19390809-Cod-Ms-H-Hasse-1-1022-2-Lorenzen-Hasse
\label{19390809}%
\hfill\begin{tabular}{@{}l@{}}%
        Bonn, 9.8.39.\\%
        Luisenstr.\ 3%
      \end{tabular}\medskip

\noindent\hfil Sehr geehrter Herr Professor!\nopagebreak\medskip

Es hat mir sehr leid getan, Sie am 6.8.\ %
nicht angetroffen zu haben. Ich hätte Ihnen
so gern nochmal persönlich gedankt für die
Fürsorge und das große Wohlwollen, das \bruch
Sie mir während meiner Lehrjahre
bei Ihnen stets erwiesen haben.

Aber ich hoffe sehr, daß dazu später
einmal Gelegenheit sein wird.\strut\medskip

\hfil\begin{tabular}{@{}c@{}}%
       Mit den ergebensten Grüßen\\%
       Ihr\\%
       \hfil Lorenzen.%
\end{tabular}

%% file: 19390811-Cod-Ms-H-Hasse-1-1022-Beil-2-Hasse-Lorenzen
\label{19390811}%
\hfill 11.8.1939.\nopagebreak

\begin{tabular}{@{}ll@{}}%
  Herrn&\\%
  &\quad Dr.\ Lorenzen%
\end{tabular}\strut\medskip

\hfil\begin{tabular}{@{}l@{}}%
  \uline{\hbox{\so{Bonn}} a.\ Rh.}\\%
  Luisenstrasse 3%
\end{tabular}\medskip

\indent\indent Lieber Herr Lorenzen!\nopagebreak\medskip

Herzlichen Dank für Ihren freundlichen Brief. Es hat mir
sehr leid getan, dass ich Sie vor Ihrer endgültigen Übersiedlung
nach Bonn nicht mehr hier sah. Ich wünsche Ihnen alles Gute für
Ihre neue Stelle. Ich darf ja sicher sein, dass Sie sich dort
sehr wohlfühlen und viele Anregungen für weitere Arbeit bekommen
werden. Die Wahrscheinlichkeit, dass wir uns immer mal wieder\bs
sehen, ist zum Glück gross, denn die mathematische Welt ist doch
recht klein. Sehr dankbar wäre ich Ihnen, wenn Sie mir gelegent\bs
lich aus Ihrem mathematischen Schaffen etwas für Crelles Journal
geben könnten.\medskip

Indem ich Sie bitte, Herrn Krull bestens von mir zu grüssen,

\indent\indent\indent bin ich stets in bester Erinnerung an die Zeiten
gemeinsamer Arbeit,

\hfil Ihr [Hasse]

 

%% file: 19390906-Cod-Ms-H-Hasse-1-1022-3-Lorenzen-Hasse
\label{19390906}%
\hfill\begin{tabular}{@{}l@{}}%
        Bonn, 6.9.39.\\%
        Luisenstr.~3%
      \end{tabular}\medskip

\hfil Sehr geehrter Herr Professor!\medskip

Zunächst möchte ich Ihnen herzlichst dan\bs
ken für Ihren freundlichen Brief. Meine
Tätigkeit in Bonn, die ich mit soviel Hoff\bs
nung aufgenommen habe, erfährt jetzt ein
schnelles Ende. Bisher bin ich zwar noch
nicht einberufen, sondern warte darauf.

Um diese Wartezeit zu verkürzen, \bruch
schreibe ich an Sie. Nun weiß ich
allerdings nicht, wo Sie augenblicklich
sind, aber wäre es Ihnen nicht möglich,
für mich eine meine mathematische Aus\bs
bildung ausnützende Verwendung zu
erreichen? Es ist doch offensichtlich, daß
eine Verwendung der Mathematiker als
Mathematiker nutzbringender für die Ge\bs
samtheit ist als jede andere. Natürlich
liegt diese Verwendung auch sehr in \bruch
meinem eigenen Interesse, hier fallen
eben beide Interessen zusammen.

Welche Institute, Werke oder
militärische Dienststellen in Frage
kommen, weiß ich nicht, aber das aero\bs
dynamische Institut z.\ B.\ wird doch sicherlich
Bedarf an Mathematikern haben.

Es erscheint mir als eine zu billige
Phrase, wollte ich Ihnen meine Dankbarkeit
im Falle Ihrer Hilfe in dieser Angelegenheit \bruch
versichern. Habe ich Ihnen doch schon für
so vieles zu danken. Ich schreibe daher
lieber ganz offen, daß mir sehr viel dran\bs
gelegen ist und daß ich mich an Sie
wende, da ich die Überzeugung habe, daß
Sie mir helfen werden.\medskip

Mit den ergebensten Grüßen
an Ihre Frau Gemahlin und Sie
bin ich stets\strut

\hfil\begin{tabular}{@{}ll@{}}%
       Ihr&\\%
          &Lorenzen.%
     \end{tabular}

%% file: 19391011-Cod-Ms-H-Hasse-1-1022-4-Lorenzen-Hasse
\label{19391011}%
\hfill Bonn, 11.10.39.\nopagebreak

\noindent\hfil Sehr geehrter Herr Professor!\nopagebreak\medskip

Meine Dissertation ist jetzt endlich gedruckt
worden. Ich möchte dies zum Anlaß nehmen,
Ihnen nochmals meinen Dank auszudrücken für
die lange und schöne Zeit, die ich in Göttingen
war.\medskip

Mit den ergebensten Grüßen an Ihre Frau Gemahlin und Sie\strut

\hfil\begin{tabular}{@{}ll@{}}%
       Ihr&\\%
          &Lorenzen.%
     \end{tabular}

%% file: 19391110-Cod-Ms-Hasse-33-3-Lorenzen-Hasse
\label{19391110}%
\hfill Bonn, 10.11.39.\nopagebreak

\hfil Sehr geehrter Herr Professor,\nopagebreak\medskip

Im Zusammenhang mit den gruppenaxiomatischen Unter\bs
suchungen für Herrn Professor Scholz bin ich auf ein Axiomen\bs
system gekommen, das ersichtlich einfacher ist als die bisher
bekannten Systeme. Da dies Ergebnis jedoch völlig aus dem
Rahmen der übrigen Untersuchungen fällt, würde ich es gern für
sich veröffentlichen und möchte Sie daher fragen, ob eine
Veröffentlichung\footnote{\citealt{MR0002113}.} in Crelles Journal möglich ist.\strut\medskip

\hfil\hfil\begin{tabular}{@{}c@{}}%
            Mit den ergebensten Grüßen\\%
            \hfil Ihr\\%
            \hfil\hphantom{Ihr\hspace*{2em}}\rlap{Lorenzen}%
          \end{tabular}

%% file: 19391116-Cod-Ms-Hasse-33-3-Hasse-Lorenzen
\label{19391116}%
\hfill 16.\ Nov.\ 1939.\nopagebreak

\noindent\begin{tabular}{@{}l@{}}%
           Herrn\\[3pt]%
           Dr.\ P.\ \so{Lorenzen}\\[3pt]%
           \uline{\hbox{\so{Bonn}} a.\ Rh.}\\%
           Luisenstr.~3.%
         \end{tabular}\medskip

\noindent Lieber Herr Lorenzen!\nopagebreak\medskip

Herzlichen Dank für Ihre kleine Note über Gruppen\bs
axiome. Ich nehme sie gern in Crelles Journal auf,
sie scheint ja auf demselben Gedanken zu beruhen, der auch bei der
Idealdefinition vorkommt, wo man nur die Subtraktion zu for\bs
dern braucht.

Ich habe noch eine kleine Frage zur algebraischen Funk\bs
tionentheorie. Es handelt sich um die Darstellung der ganzen
Divisoren eines Funktionenkörpers durch ganze homogene Ideale
eines Integritätsbereichs. Es sei $x_1,\dots,x_m$ ein Erzeugen\bs
den-System des Integritätsbereichs der für~$x$ ganzen Elemente
des Körpers, wo $x$~ein beliebiges nicht konstantes Element
ist. Unter welchen Bedingungen entsprechen dann die ganzen
Divisoren des Körpers umkehrbar eindeutig den ganzen Idea\bs
len des durch Homogenisierung der $x_1,\dots,x_m$ entstehenden
Integritätsbereichs? Ich kann zeigen, dass dies richtig
ist, wenn man $x_1,\dots,x_m$ als die von~$1$ verschiedenen Ele\bs
mente der Basis eines Moduls $(\frac{1}{\mathfrak o^m})$ nimmt, wo $m\geqq 3g$ ist. Ich
wüsste aber gern, welcher allgemeine Satz dahintersteckt. Die
Bedingung $m\geqq 3g$ garantiert dafür, dass der Modul eine Nor\bs
malbasis für~$\mathfrak o$ enthält.\medskip

Mit herzlichen Grüssen

\hfil Ihr [Hasse]
 

%% file: 19391117-Cod-Ms-H-Hasse-1-1022-4-2-Hasse-Lorenzen
\noindent\begin{tabular}{@{}c@{}}%
\label{19391117}%
           \bfseries Mathematisches Institut\\%
           \bfseries der Universität%
         \end{tabular}%
\hfill\begin{tabular}{@{}l@{}}%
        \small\textbf{Göttingen, den }17.11.39\\
        \footnotesize\bfseries Bunsenstraße 3/5%
  \end{tabular}\medskip

\indent\indent Lieber Herr Lorenzen!\medskip

Die Frage in meinem gestrigen Brief war etwas knapp
formuliert -- ich diktierte einer nicht-\hspace{0pt}mathematischen Dame in die
Maschine! Ich will daher etwas klarer auseinandersetzen, was ich meine.

Sei $K$~algebr.\ Funktionenkörper einer Unbestimmten über
Konstantenkörper~$\Omega$ vom Geschlecht~$g$ ($\geqq1$), mit einem Primdivisor
ersten Grades~$\mathfrak{o}$. Sei $K^{\mathfrak{o}}$ der Integritätsbereich der für
alle $\mathfrak{x}\neq\mathfrak{o}$ ganzen Elemente aus~$K$. Man folgert dann leicht aus
Riemann-Roch, daß der $\Omega$-Modul $\bigl(\frac{1}{\mathfrak{o}^m}\bigr)$ (der Vielfachen von~$\frac{1}{\mathfrak{o}^m}$ aus~$K$)
für $m\geqq 3g$ den Integritätsbereich~$K^{\mathfrak o}$ ergänzt, d.\ h.\ daß eine
$\Omega$-Basis von~$\bigl(\frac{1}{\mathfrak{o}^m}\bigr)$ ein Erzeugendensystem für~$K^{\mathfrak o}/\Omega$ ist.\smallskip

\uline{Beweis.} Sei $x$ ein nicht-konst.\ Element aus~$K$ mit einer möglichst niedrigen
Potenz $\mathfrak o^n$ im Nenner, so ist $K^{\mathfrak o}$ die ganzalg.-abgeschl.\ Hülle von~$\Omega[x]$ in~$K$.
Eine $\Omega[x]$-Basis von~$K^{\mathfrak o}$ erhält man, indem man zu jeder Restklasse~$i$ mod.~$n$
ein Element~$y_i$ aus~$K$ wählt, dessen Nenner eine möglichst niedrige Potenz~$\mathfrak o^{e_i}$
mit $e_i\geqq0$, $e_i\equiv i\bmod n$ wählt ($y_0=1$). Da $\dim\mathfrak o^e=\dim\mathfrak o^{e-1}+1$ für $e\geqq 2g$
gilt, und $n\leqq g+1$ ist (wegen $\dim\mathfrak o^e\geqq2$ für $e\geqq g+1$), folgt $e_i\leqq3g$. Für $m\geqq3g$
enthält also $\bigl(\frac{1}{\mathfrak{o}^m}\bigr)$ das Erzeugendensystem $x;y_i$ von~$K^{\mathfrak o}/\Omega$.\smallskip

Sei jetzt $1,x_1,\dots,x_r$ ($r=m-g$) eine Basis von~$\bigl(\frac{1}{\mathfrak{o}^m}\bigr)$ ($m\geqq3g$). Ist dann~$\overline{\mathfrak x}\neq0$ ein
algebr.\ Punkt von~$K/\Omega$ (Primdiv.\ der algebr.-abgeschl.\ Konstantenerweiterung~$\overline{K}/\overline{\Omega}$), so
seien%
\[%
  x_1(\overline{\mathfrak x})=\overline{a}_1,\,\dots,\,x_r(\overline{\mathfrak x})=\overline{a}_r%
\]
als die Koordinaten von~$\overline{\mathfrak x}$ bezeichnet. Sie beschreiben~$\overline{\mathfrak x}$ eindeutig, und es gilt%
\[%
  \frac{\overline{\mathfrak x}}{\mathfrak o^m}=(x_1-\overline{a}_1,\dots,x_r-\overline{a}_r)\text,%
\]
wo die Klammer rechts nicht etwa ein Ideal bedeuten soll, sondern nur den
g.\ g.\ T.\ der eingeschlossenen Hauptdivisoren (definiert durch Minima der Exponenten
in der Primdivisorzerlegung in~$\overline{K}/\overline{\Omega}$). Entsprechend hat man auch%
\[%
  \frac{\mathfrak o}{\mathfrak o^m}=(1,x_1,\dots,x_{r-1})\text,%
\]
wenn etwa das Basiselement~$x_r$ (als einziges) den genauen Nenner $\mathfrak o^m$ hat (d.\ h.\bsp
$1,x_1,\allowbreak\dots,x_{r-1}$ Basis von~$\bigl(\frac{1}{\mathfrak{o}^{m-1}}\bigr)$ ist). $\mid^2$\smallskip \bruch

Ferner gilt: \uline{Ist $u$~ein Element aus~$K^{\mathfrak o}$, so stellt sich $u$~als Polynom in $x_1,\dots,x_r$
über~$\Omega$ dar, dessen Grad $\leqq h$ ist, wenn $u$~in~$\bigl(\frac{1}{\mathfrak{o}^{mh}}\bigr)$ liegt.}\smallskip

\uline{Beweis.} Es genügt $\bigl(\frac{1}{\mathfrak{o}^{mh}}\bigr)=\bigl(\frac{1}{\mathfrak{o}^{m}}\bigr)^h$ zu zeigen. Der Beweis durch vollst.\ Induktion
nach~$h$ läuft auf die Feststellung $\bigl(\frac{1}{\mathfrak{o}^{m(h+1)}}\bigr)=\bigl(\frac{1}{\mathfrak{o}^{mh}}\bigr)\bigl(\frac{1}{\mathfrak{o}^{m}}\bigr)$ zurück. Nun zerfallen
die natürlichen Zahlen~$e$ in Nennerzahlen und Lückenzahlen, je nachdem
es in~$K$ ein Element mit dem Nenner $\mathfrak o^e$ gibt oder nicht ($\dim\mathfrak o^e=\dim\mathfrak o^{e-1}+1$
oder $=\dim\mathfrak o^{e-1}$ ist). Es gibt genau $g$~Lückenzahlen, und diese gehören der
Reihe $1,\dots,2g-1$ an. Eine Basis von $\bigl(\frac{1}{\mathfrak{o}^{m(h+1)}}\bigr)$ mod.~$\bigl(\frac{1}{\mathfrak{o}^{mh}}\bigr)$ erhält man,
wenn man zu jeder natürlichen Zahl $e=mh+\mu$ ($\mu=1,\dots,m$)
ein Element aus~$K$ mit dem Nenner~$\mathfrak o^e$ bestimmt. Es genügt dann
zu zeigen, daß jedes solche~$e$ eine Darstellung als Summe einer Nenner\bs
zahl~$\leqq mh$ und einer Nennerzahl~$\leqq m$ hat. Dies leistet aber eine der
$g+1$ Zerlegungen $e=(mh-\gamma)+(\mu+\gamma)$ \ $(\gamma=0,1,\dots,g)$. Ist nämlich $\mu$
Lückenzahl, so ist $\mu+\gamma\leqq3g-1<m$, und unter den $g+1$~Zahlen $\mu+\gamma$
ist mindestens eine Nennerzahl, während die Zahlen $mh-\gamma\geqq m-\gamma\geqq 2g$
sämtlich Nennerzahlen sind. Ist aber $\mu$ Nennerzahl, so genügt $\gamma=0$.\smallskip

Ich setze nun%
\[%
x_1=\frac{X_1}{X_0},\,\dots,\,x_r=\frac{X_r}{X_0}\text.\leqno(1.)%
\]
Dabei sei $X_0$ ein beliebiges Element~$\neq0$ \uuuline{aus~$K$}, wodurch dann $X_1,\dots,X_r$
bestimmt sind. Man kann dann die obigen Darstellungen der alg.\ Punkte\nbsp
$\overline{\mathfrak x}$ von~$K/\Omega$ als g.\ g.\ T.\ auch so schreiben:%
\[%
  \begin{aligned}%
    \overline{\mathfrak x}&=\frac{(0,X_1-\overline{a}_1X_0,\dots,X_r-\overline{a}_rX_0)}{(X_0,X_1,\dots\dots\dots\dots\dots,X_r)}\qquad(\overline{\mathfrak x}\neq0)\\%
    \mathfrak o&=\frac{(X_0,X_1,\dots,X_{r-1},0)}{(X_0,X_1,\dots,X_{r-1},X_r)}\text,%
  \end{aligned}%
\]
unabhängig von der Wahl der Quotientendarstellung~(1.). Nach dem Gaußschen
Satz multiplizieren sich die Klammersymbole $(\mathscr U_0,\mathscr U_1,\dots,\mathscr U_k)$ für den
g.\ g.\ T.\ von Elementen $\mathscr U_0,\mathscr U_1,\dots,\mathscr U_k$ aus~$K$ formal wie die Polynome
$\mathscr U_0t^k+\mathscr U_1t^{k-1}+\dots+\mathscr U_k$ einer Unbestimmten~$t$. Unter Anwendung \uline{dieser}
Multiplikation folgt durch Zusammensetzung der ganzen Divisoren~$\mathfrak a$ von~$K/\Omega$
aus algebr.\ Punkten~$\overline{\mathfrak x}$ von~$K/\Omega$ (vollst.\ Systeme Konjugierter!), daß jeder
ganze Divisor~$\mathfrak a$ von~$K/\Omega$ eine Darstellung der folgenden Form besitzt:%
\[%
\mathfrak a=\frac{(\mathscr A_0,\mathscr A_1,\dots,\mathscr A_{ra})}{(X_0,X_1,\dots,X_r)^a}\text,\mid^3\leqno(2.)%
\]\bruch
wo $a$~der Grad von~$\mathfrak a$ ist, ferner $\mathscr A_0,\dots,\mathscr A_{ra}$ homogene Polynome $a$-ten
Grades \uline{über~$\Omega$} in $X_0,\dots,X_r$ sind.

Von diesen Darstellungen~(2.) aus kann man eindeutig die Distributionen\nbsp
$\mathfrak a(\mathfrak x)$ erklären. Man wähle dazu für jeden Primdivisor~$\mathfrak x$ von~$K/\Omega$
eine besondere Quotientendarstellung~(1.), nämlich so, daß $\mathfrak x$ zu
$(X_0,X_1,\dots,X_r)$ den Beitrag~$\mathfrak x^0$ liefert (\uline{primitiv für~$\mathfrak x$}). Dann liefert
die Festsetzung%
\[%
\mathfrak a(\mathfrak x)=\frac{(\mathscr A_0(\mathfrak x),\dots,\mathscr A_{ra}(\mathfrak x))}{(X_0(\mathfrak x),\dots,X_r(\mathfrak x))^a}\leqno(3.)%
\]
ein bestimmtes Ideal~$\mathfrak a(\mathfrak x)$ aus dem Restklassenkörper~$K(\mathfrak x)$ von~$\mathfrak x$,
und zwar unabh.\ von der Wahl der für~$\mathfrak x$ primitiven Quotientendarstellung~(1.).
Man sieht sofort, daß $\mathfrak a(\mathfrak x)$ distributionsganz ist, daß nämlich der Nenner
von~$\mathfrak a(\mathfrak x)$ im Hauptnenner der Koeffizienten von $\mathscr A_0,\dots,\mathscr A_{ra}$ aufgeht.
Ferner gilt die Regel:%
\[%
\text{aus $\mathfrak a=\mathfrak a_1\mathfrak a_2$\quad folgt\quad$\mathfrak a(\mathfrak x)=\mathfrak a_1\mkern-2mu(\mathfrak x)\mkern2mu\mathfrak a_2\mkern-1mu(\mathfrak x)$,}%
\]
einfach auf Grund der Gültigkeit des Gaußschen Satzes sowohl für die Multipl.\ der
g.\ g.\ T.\ von Hauptdivisoren als auch für die Multipl.\ der g.\ g.\ Teiler von Hauptidealen.

Um den Weilschen Hauptsatz zu beweisen, muß man dann noch wissen,
daß die Definition~(3.) im Sinne der Distributionsgleichheit von der Wahl der
besonderen zuvor konstruierten Darstellung~(2.) unabhängig ist,
daß nämlich folgendes gilt: Ist%
\[\mathfrak a=\frac{(\mathscr U_0,\mathscr U_1,\dots,\mathscr U_s)}{(X_0,X_1,\dots,X_r)^k}\leqno(2'.)\]
irgendeine Darstellung von~$\mathfrak a$, wo $\mathscr U_0,\mathscr U_1,\dots,\mathscr U_s$ homog.\ Polynome in $X_0,X_1,\allowbreak\dots,X_r$
vom Grade~$k$ sind, so ist auch
\[\mathfrak a(\mathfrak x)\doteq\frac{(\mathscr U_0(\mathfrak x),\dots,\mathscr U_s(\mathfrak x))}{(X_0(\mathfrak x),\dots,X_r(\mathfrak x))^k}\qquad\text{(distributionsgleich!)}\leqno(3'.)\]
Ich kann das zwar beweisen, indem ich den alten Schluß mit der algebraischen
Gleichung (Charakt.\ der Nenner durch ganzalgebraische Abhängigkeiten) anwende; dazu
brauche ich den Satz auf S.~2 oben.
Besser würde mir jedoch der Weilsche Schluß aus \emph{Actualités}\comment{\citealt{weil35}.} gefallen, der
darauf beruht, daß die Klammersymbole in diesen Darstellungen auch als
Ideale (lineare Komposita!) aufgefaßt werden dürfen, daß also die Gesamtheit
der~$\mathscr U$ mit $(\mathscr U_0,\mathscr U_1,\dots,\mathscr U_s)\divides\mathscr U$ durch eine homogene lineare Darstellung $\mathscr U=\mathscr F_0\mathscr U_0+\dots+\mathscr F_s\mathscr U_s$
gekennzeichnet ist. Wie ist das genau zu formulieren? Und wie ist es zu
beweisen? Das sind meine Fragen. $\mid^4$

Natürlich wüßte ich gerne eine möglichst allgemeine Regel über die Zuordnung
der ganzen Divisoren~$\mathfrak a$ von~$K/\Omega$ zu homogenen Idealen, möglichst
unabhängig von der besonderen Wahl meiner Erzeugenden $x_1,\dots,x_r$.

Können Sie mir da helfen?\strut\medskip

\hfil\hfil\begin{tabular}{@{}r@{}}%
       Herzlichst Ihr\\%
       Hasse%
     \end{tabular}


%% file: 19391213-Cod-Ms-H-Hasse-1-1022-4-1-Lorenzen-Hasse
\label{19391213}%
\hfill Bonn, 13.12.39.\nopagebreak

\noindent\hfil Sehr geehrter Herr Professor,\nopagebreak\medskip

Ich muß Sie sehr um Entschuldigung bitten, daß ich Ihnen
jetzt erst antworte. Da ich nämlich gerade stark in anderen
Fragen drin war, wollte mir durchaus nichts einfallen zu Ihrer
Frage. So hab ich schließlich aus der Not eine Tugend gemacht
und glaube jetzt Ihre Frage dahin beantworten zu müssen, daß
eine Zurückführung der Divisorenteilbarkeit in~$K$\kern1pt\textsuperscript{1)} auf die Ideal\bs
teilbarkeit in~$\Omega[X_0,\dots,X_r]$ nicht möglich ist mit einem~$X_0\in K$.
Denn die ganzen homogenen Polynome in $X_0,\dots,X_r$ sind dann keine
ganzen Divisoren, sondern gebrochene Hauptdivisoren.

Wenn man dagegen $X_0$ als neue Unbestimmte einführt, so gilt zu\bs
nächst für rationale Funktionenkörper einer Unbestimmten
\begin{enumerate}[1)]
\item Die homogenen Elemente von~$\Omega[X_0,X_1]$ entsprechen
  umkehrbar eindeutig den ganzen Divisoren von~$K$,
\item Die Idealsummenbildung entspricht der Bildung des
g.\ g.\ T.\ von Divisoren.
\end{enumerate}
Für algebraische Funktionenkörper sind in~1) die Elemente durch
homogene Ideale (im Sinne der Weilschen Äquivalenz) zu er\bs\bruch
setzen, während~2) ungeändert bleibt. Für Funktionenkörper
mehrerer Unbestimmten sind in~1) und~2) die Ideale durch $v$-Ideale
zu ersetzen (und außerdem eine Transzendenzbasis von~$K$ auszu\bs
zeichnen!)

Gleichzeitig übersende ich Ihnen die Korrektur des ,,Axiomen\bs
systems für Gruppen``,\comment{See Lorenzen's letter to Hasse dated 10 November 1939 on page~\pageref{19391110}.} für dessen Aufnahme in Ihr Journal ich
Ihnen nochmals herzlichst danken möchte.\strut\medskip

\hfil\begin{tabular}{@{}lll@{}}%
       Mit den ergebensten Grüßen&&\\[6pt]%
                                 &Ihr&\\%
                                 &&Lorenzen%
     \end{tabular}\nopagebreak

\vspace*{\skip\footins}\nopagebreak\footnoterule\nopagebreak{\footnotesize\rule{0pt}{\footnotesep}\textsuperscript{1)} Ich nehme die Bezeichnungen aus Ihrem Brief, den ich beifüge.}


%% file: 19400104-Cod-Ms-H-Hasse-1-1022-Beil-3-Hasse-Lorenzen
\label{19400104}%
\hfill 4.1.1940.\nopagebreak

\noindent\begin{tabular}{@{}l@{}}%
           Herrn\\%
           Dr.\ Paul \so{Lorenzen}\\[6pt]%
           \uline{\hbox{\so{Bonn}} a.\ Rh.}\\%
           Luisenstr.~3.%
         \end{tabular}\medskip

\noindent Lieber Herr Lorenzen!\medskip

Ganz überraschend kam mir heute Ihre Vermählungs\bs
anzeige nicht. Man spricht in Göttingen schon davon.
Nehmen Sie meinen herzlichsten Glückwunsch zu diesem,
für Ihr ferneres Leben so bedeutsamen Schritt.

Auch unabhängig von diesem freudigen Anlass wollte
ich Ihnen dieser Tage schreiben, und zwar in Verfolg
Ihrer damaligen Anfrage nach einer kriegerischen Betäti\bs
gung. Ich erhielt nämlich kürzlich eine Anfrage von Pro\bs
fessor Walther (Darmstadt).\footnote{\selectlanguage{english}Alwin Walther (1898--1967) is a mathematician and engineer at Universität Darmstadt who participates in the V2 program at Peenemünde.} Dieser hat den Auftrag, Mathe\bs
matiker für eine militärische Dienststelle an der Ostsee
namhaft zu machen, für die er selbst tätig ist. Die An\bs
stellung dort ist mit der Reklamation vom Heeresdienst
verbunden. Die Bezahlung ist gut, Sie würden in Ihrem
Falle, unter Berücksichtigung Ihres jetzigen verheirate\bs
ten Standes -- nach Abzug der Kürzungen -- RM~418.74 im
Monat erhalten. Falls Sie sich dafür zur Verfügung stel\bs
len wollen, bitte ich Sie, mir den beiliegenden Bogen
ausgefüllt zurückzusenden, anderenfalls leer.\strut\medskip

\begin{tabular}{@{}lll@{}}%
  Mit besten Grüssen und Wünschen&&\\%
  &&Ihr [Hasse]%
\end{tabular}

Anlage.
 

%% file: 19400109-Cod-Ms-H-Hasse-1-1022-5-Lorenzen-Hasse
\noindent
\label{19400109}%
Dr.\ Paul Lorenzen\hfill\begin{tabular}[t]{@{}l@{}}%
                                   Bonn, den 9.\ Januar 1940\\%
                                   Luisenstr.~3%
                                 \end{tabular}\medskip

\noindent Sehr geehrter Herr Professor,\medskip

\noindent diesesmal habe ich Ihnen wiederum für vieles zu danken.
Und zwar zunächst für Ihren schönen Brief vom 18.12.
Ich wünschte, ich fände die Zeit, um alles völlig zu ver\bs
stehen. (Dazu brauche ich ja oft leider sehr viel Zeit!)

\noindent Dann danke ich Ihnen herzlichst auch im Namen meiner Frau
für Ihre Glückwünsche zu unserer Hochzeit. Augenblicklich
wohnen wir in einer möblierten Wohnung, wobei zur Vervoll\bs
kommnung unseres Glückes eigentlich nur fehlt, daß der Krieg
aufhört.

\noindent Damit komme ich zum dritten, für das ich Ihnen danken
möchte, nämlich für die Übersendung des Fragebogens von Herrn
Professor Walther. Da Herr Professor Krull zurzeit in Bonn
ist, habe ich ihn zunächst gefragt, ob es mir -- als
Assistent -- überhaupt möglich sei, zuzusagen. Herr Professor
Krull hat sich die Antwort für einige Tage vorbehalten, da
er der Meinung ist, es sei evtl.\ zunächst der Kurator zu
fragen.

Ich rechne aber damit, sehr bald zusagen zu können und ver\bs
bleibe bis dahin mit den ergebensten Grüßen\strut\medskip

\hfil\hfil\hfil\begin{tabular}{@{}ll@{}}%
                 Ihr&\\%
                    &Paul Lorenzen\\%
               \end{tabular}
 

%% file: 19400110-Cod-Ms-H-Hasse-1-1022-6-Lorenzen-Hasse
\noindent
\label{19400110}%
Dr.\ Paul Lorenzen\hfill\begin{tabular}[t]{@{}l@{}}%
                                   Bonn, den 10.\ Januar 1940\\%
                                   Luisenstr.~3%
                                 \end{tabular}\medskip

\noindent Sehr geehrter Herr Professor,\medskip

\noindent nach erneuter Rücksprache mit Herrn Prof.\ \so{Krull} habe
ich die Erlaubnis bekommen, Ihnen den Fragebogen zuzusen\bs
den.

\noindent Allerdings würde mich das Math.\ Seminar hier nur sehr un\bs
gern gehen lassen, da -- zumindestens für das laufende Tri\bs
mester -- ziemlich dringend zwei Assistenten benötigt wer\bs
den. Falls Herr Prof.\ \so{Walther} mich also für seine
Dienststelle gebrauchen könnte, möchte ich ihn bitten,
sich dazu an Herrn Prof.\ Krull als den Direktor des Math.\bsp
Seminars zu wenden. Da mit meiner Einberufung zum Wehr\bs
dienst doch wohl bald zu rechnen ist, würde Herr Prof.\ Krull
in diesem Fall wahrscheinlich zustimmen, falls nicht meine
Reklamation für die hiesige Universität sich inzwischen er\bs
möglichen läßt.

\noindent Indem ich Ihnen für die viele Mühe, die Sie sich schon
wieder für mich gemacht haben (ich hoffe nur, daß es nicht
noch mehr wird) nochmals herzlichst danke, verbleibe ich
mit den ergebensten Grüßen\strut\medskip

\hfil\hfil\hfil\begin{tabular}{@{}ll@{}}%
                 Ihr&\\%
                    &Paul Lorenzen\\%
               \end{tabular}

%% file: 19400203-Cod-Ms-H-Hasse-1-1021-Lorenzen-Hasse
\begin{tabular}[t]{@{}l@{}}%
\label{19400203}%
  Lorenzen\\%
  Bonn, Luisenstr.~3%
\end{tabular}%
\hfill\begin{tabular}[t]{@{}l@{}}%
        Herrn\\%
        Professor Hasse,\\%
        Göttingen\\%
        Bunsenstr.~2\\%
        Math.\ Institut%
      \end{tabular}\medskip

\hfill Bonn, 3.II.40\medskip

\noindent Sehr geehrter Herr Professor,\medskip

\noindent da ich nicht weiß, ob mein Mann Ihnen
schon geschrieben hat, daß er seit dem
19.I.\ als Soldat eingezogen ist u.\ gestern
die Beitragsrechnung für die Deutsche Mathe\bs
matiker Vereinigung kam, möchte ich Ihnen
doch noch einmal hiervon Mitteilung
machen. Wahrscheinlich wird Ihnen mein
Mann aber selbst noch davon schreiben, \bruch
da der Transport jetzt
beendet ist.\medskip

\noindent Mit ergebenen Grüßen\nopagebreak

\hfil\hfil\hfil Frau K.\ Lorenzen.
 

%% file: 19400206-Cod-Ms-H-Hasse-1-1021-Beil-Hasse-Lorenzen
\label{19400206}%
\hfill 6.2.40.

\noindent
\begin{tabular}{@{}l@{}}
  Frau\\%
Dr.\ \so{Lorenzen}\\%
\uline{\hbox{\so{Bonn} }a.\ Rh.}\\%
Luisenstr.~3.%
\end{tabular}\medskip

\noindent Sehr verehrte Frau Lorenzen!\medskip

Von der Einberufung Ihres Mannes erfuhr ich auch
durch einen Brief von Professor Krull.\footnote{\selectlanguage{english}Letter dated 1 February 1940 \citep[§~1.72]{roquette04}.} Es bleibt nun
abzuwarten, ob die Heeresdienststelle an der Ostsee,
für die sich Ihr Mann gemeldet hatte, seine Anstellung
weiter betreibt und eine Freistellung erreicht. Jeden\bs
falls scheint diese Stelle es nicht sehr eilig zu haben.
Herr Eichler,\footnote{\selectlanguage{english}The mathematician Martin Eichler (1912--1992) will join Walther's team in Peenemünde.} ein Assistent von hier, der sich kurz vor
Ihrem Mann dort meldete, und der auch schon persönlich
in Berlin mit dem Leiter der Stelle verhandelte, sollte
an sich zum 1.~Februar eingestellt werden, hat aber bis
heute keine Nachricht erhalten. Von Ihrem Mann direkt
habe ich noch nichts gehört.

Mit freundlichen Grüssen\nopagebreak

\hfill Ihr sehr ergebener [Hasse]
 

%% file: 19400204-Cod-Ms-H-Hasse-1-1022-7-Lorenzen-Hasse
\label{19400204}%
\hfill Alsheim, 4.2.\medskip

\hfil Sehr geehrter Herr Professor,\medskip

So leid es mir tut, muß ich Sie -- da Sie
damals so liebenswürdig waren, bei der Anfrage
an Herrn Prof.\ Walther an mich zu denken --
nochmals in dieser Sache belästigen. Ich bin
nämlich mittlerweile eingezogen worden.
(Wahrscheinlich war das Wehrmeldeamt der
Meinung, daß es mir nach meiner Heirat
entschieden zu gut ging -- -- ich muß gestehen,
daß diese Meinung nahe lag.)

Alsheim ist ein kleines Nest in der Nähe von
Worms. Wir sind recht gut in Privatquartie\bs
ren untergebracht und beschäftigen uns \bruch
meist mit Schneeschaufeln. Ich finde das
zwar nicht schön, aber es wird nun einmal
gerade dieser "`Einsatz"' von mir gefordert.

Ich möchte Sie nun bitten, falls es nötig
sein sollte, Herrn Prof.\ Walther mitzu\bs
teilen, daß ich am 19.1.40 einberufen bin
und jetzt unter der Feldpostnummer~38503~c
zu erreichen bin.\medskip

Mit den ergebensten Grüßen an Ihre Frau Gemahlin und Sie verbleibe ich\strut

\hfil\hfil\hfil\begin{tabular}{@{}ll@{}}%
                 Ihr&\\%
                    &Paul Lorenzen%
               \end{tabular}
 

%% file: 19400209-Cod-Ms-H-Hasse-27-1-Lorenzen-Hasse
\begin{tabular}[t]{@{}l@{}}
\label{19400209}%
  Gefr.\ Lorenzen\\%
  38503~c%
\end{tabular}%
\hfill\begin{tabular}[t]{@{}l@{}}
        Herrn\\%
        Prof.\ Dr.\ H.\ Hasse\\%
        Göttingen\\%
        Bunsenstr.~3~-~5%
      \end{tabular}\medskip

\hfill 9.2.\medskip

\hfil Sehr geehrter Herr Professor,\medskip

Diesmal wende ich mich an Sie als Schatz\bs
meister der D. M. V. und bitte Sie
mir den Beitrag für~1940 zu erlassen.\medskip

\hfil Heil Hitler!\strut\medskip

\hfil\hfil\begin{tabular}{@{}lll@{}}%
            Ihr&&\\%
               &sehr ergebener&\\%
               &&Paul Lorenzen.%
          \end{tabular}
 

%% file: 19400421-Cod-Ms-H-Hasse-1-1022-8-Lorenzen-Hasse
\label{19400421}%
\hfill21.4.\medskip

\hfil Sehr geehrter Herr Professor,\medskip

ich freue mich sehr Ihnen aus
dem kriegerischen Westen mal
etwas so Unpolitisches wie dieses
Separatulum\comment{The only candidate for this offprint seems to be \citealt{lorenzen39b}.} zusenden zu können.
Für Ihre Vermittlung an Herrn
Prof.\ Walther bin ich Ihnen noch
täglich von Herzen dankbar. Es hat
sich zwar noch nichts entschieden, ist
aber meine einzige Hoffnung in  \bruch
meinem augenblicklichen geist\bs
tötenden Zustand, dessen evtl.\bsp
lange Dauer ja das Schlimmste ist.\medskip

Mit den ergebensten Grüßen an Ihre Frau Gemahlin und Sie verbleibe ich\strut

\hfil\hfil\hfil\begin{tabular}{@{}ll@{}}%
                 Ihr&\\%
                    &Paul Lorenzen.%
               \end{tabular}
 

%% file: 19400505-Cod-Ms-H-Hasse-1-1022-9-Lorenzen-Hasse
\label{19400505}%
\hfill Hamburg, 5.5.40.\smallskip

\hfil Sehr geehrter Herr Professor,\medskip

vorgestern bekam ich Ihren freundlichen Brief, über
den ich mich sehr gefreut habe und für den ich
mich vielmals bei Ihnen bedanken möchte.

Muß ich Ihnen doch so dankbar dafür sein, daß
Sie auch in Ihrer neuen Tätigkeit gleich wieder an
mich gedacht haben, der ich mich hier mit Pferden,
Stiefeln, Karabinern u.\ ä.\ herumzuquälen habe.

Ich hoffe, daß Sie sich allmählich in Ihre
jetzige Beschäftigung hineinfinden und sich die
Reminiszenzen an Ihre Staatsexamensängste
dabei verlieren werden. Wenn Ihnen genügend
freie Zeit bliebe (was aber wohl kaum zu \bruch
hoffen ist?)\ könnten Sie ja mit Herrn
Rohrbach und Herrn Kochendörffer beinah eine
"`Filiale"' des math.\ Seminars Göttingen
aufrechterhalten.

Bezgl.\ meiner Bewerbung für die
"`Dienststelle an der Ostsee"', die sich in Peene\bs
münde befindet, hat sich noch nichts weiter
ereignet. Auf eine direkte Anfrage von
mir habe ich noch keine Antwort.

\looseness=-1 Daher gehe ich mit Freuden auf
Ihre Anregung ein. Sie fragen, ob mir die
Beschäftigung mit reichlich kniffligen Fragen
der praktischen Physik liegen würde. Nun wage
ich es natürlich nicht, zu entscheiden, in wie\bs
weit ich in Ihrer Forschungsgruppe zu \bruch
gebrauchen wäre und muß es Ihnen daher
überlassen, darüber zu urteilen.

Aber ich darf vielleicht soweit gehn
und behaupten, daß mir jede theoretische Tä\bs
tigkeit ausgesprochen mehr liegt als das,
was ich augenblicklich zu tun habe.

Sie werden daher ermessen können, wie
sehr froh ich wäre, wenn Sie eine Anforderung
des A. A. erreichen könnten.\medskip

Zum Schluß erlaube ich mir noch, Sie um
Grüße an die Herrn Rohrbach, Köthe und Kochendörffer
zu bitten und verbleibe\strut\medskip

\hfil\hfil\begin{tabular}{@{}c@{}}%
            mit den ergebensten Grüßen an Sie\\%
            Ihr%
          \end{tabular}\nopagebreak

\hfill Paul Lorenzen.\par

%% file: 19400509-Cod-Ms-H-Hasse-1-1022-Beil-4-Hasse-Lorenzen
\begin{tabular}[t]{@{}l@{}}%
\label{19400509}%
Oberkommando der Kriegsmarine\\%
M Wa Stb F\\%
Berlin W 35, Tirpitzufer 60-62
\end{tabular}\hfill Berlin, 9.5.1940\medskip

\noindent Lieber Herr Lorenzen,\medskip

\noindent Besten Dank für Ihren freundlichen Brief vom 5.5.40.
Ich verstehe daraus, dass ich mich wohl nicht genügend
klar ausgedrückt habe. In meiner Forschungsgruppe habe
ich leider für Sie keinen Platz. Das Referat "`Mathe\bs
matik"' ist bereits besetzt, übrigens auch durch einen
früheren Göttinger (Schüler von Münzner). Dagegen hatte
ich an eine Verwendung entsprechender Art innerhalb der
Marine gedacht, wie sie Rohrbach und Kochendörffer im
Auswärtigen Amt haben. Es handelt sich dabei um eine
besondere Tätigkeit, die mathematische Fähigkeiten voraus\bs
setzt, und die etwa dem Lösen von Kreuzworträtseln ver\bs
gleichbar ist. Sie können sich danach wohl ungefähr
denken, was verlangt wird. Wie ich nun heute erfahren
habe, hat die Marine in der Tat Verwendung für einen
Mathematiker in dieser Stelle. Wenn Sie sich dafür inter\bs
essieren, so bewerben Sie sich doch bitte unter Berufung
auf mich (hinter meinen Namen OKM, M Wa Stb F setzen!) bei
Herrn Korvettenkapitän Teubner,\footnote{\selectlanguage{english}Achim Teubner (1905--1945) is officer at the Marinenachrichtendienst (Naval intelligence service).} Oberkommando der Kriegs\bs
marine, Berlin W~35, Skl(B). Bei dieser Bewerbung müssen
Sie sich als Mathematiker für den Dienst in der Dienst\bs
stelle dieses Kapitäns anbieten, ohne näher auf die Art
des Dienstes selbst einzugehen. Sie müssen ferner genau
mitteilen, welchem Truppenteil Sie jetzt angehören.

Ob es dann gelingt, Sie dort freizubekommen, ist aller\bs
dings eine besondere Frage. Ich habe jedenfalls Kapt.\ Teubner
ein ausführliches Gutachten über Sie gegeben und Ihre
besondere Geeignetheit für diese Verwendung hervorgehoben.
In der Tat glaube ich, dass von allen Gebieten der Mathe\bs
matik die Algebra am besten für die fragliche Arbeit zu
gebrauchen ist. Bitte lassen Sie mich laufend wissen, was
in dieser Sache geschieht, damit ich gegebenenfalls von
meiner zentralen Stelle aus helfend eingreifen kann.\medskip

Mit bestem Gruss und in der Hoffnung, Sie demnächst
hier zu sehen, Ihr [Hasse]\par
 

%% file: 19400521-Cod-Ms-H-Hasse-1-1022-10-Lorenzen-Hasse
\begin{tabular}[t]{@{}l@{}}%
\label{19400521}%
  Gefr.\ Lorenzen\\%
  38503 C%
\end{tabular}\hfill 21.5.40\medskip

\hfil Sehr geehrter Herr Professor,\medskip

Ihr liebenswürdiger und für mich so erfreuliche
Brief erreichte mich mitten in meiner
ersten Schlacht, bei der unsre Batterie aber
keine Verluste erlitten hat.

Augenblicklich sind wir -- für wenige Tage --
in Ruhestellung, in der aber leider uns der
Dienst nicht das schöne Wetter genießen läßt.

Herrn Korvettenkapitän Teubner habe
ich geschrieben und mich für seine Dienststelle
als Mathematiker beworben.

Wie sehr ich mich freuen würde, wieder
mich mathematisch betätigen zu können, ist \bruch
wohl überflüssig zu betonen.

Aber ich kann ja leider nichts dazu tun
als zu warten.

Von meiner Truppe habe ich nur die Feld\bs
postnummer angegeben, da alles weitere ver\bs
boten ist. Das O~K~M wird daraus ja sicherlich
auch alles andere ermitteln können.\medskip

Indem ich Ihnen von ganzem Herzen
danke für Ihre Bemühungen verbleibe ich
mit den ergebensten Grüßen

\hfil\qquad Ihr\nopagebreak

\hfill Paul Lorenzen.\par
 

%% file: 19400611-Cod-Ms-H-Hasse-1-1022-11-Lorenzen-Hasse
\begin{tabular}[t]{@{}l@{}}%
\label{19400611}%
  Gefr.\ Lorenzen\\%
  38503 C%
\end{tabular}\hfill%
\begin{tabular}[t]{@{}l@{}}%
  Leutnant z.\ See\\%
  Prof.\ Dr.\ Hasse\\%
  O~K~M \ M Wa Stb F\\%
  Berlin W 35\\%
  Tirpitzufer 60-62%
\end{tabular}\medskip

\hfill11.6.\hspace*{3em}\medskip

\hfil Sehr geehrter Herr Professor,\medskip

\noindent vor wenigen Tagen habe ich vom O~K~M
(3.~Abt.\ Skl.) schon die Aufforderung erhal\bs
ten meinen ausführlichen Lebenslauf ein\bs
zuweisen. Mein jetziger Batteriechef wäre
mit meiner Abkommandierung einver\bs
standen. Diese beiden Tatsachen zusammen
erleichtern mir augenblicklich meinen ziemlich
ungemütlichen Aufenthalt in einem selbst ge\bs
grabenen Erdloch von $\np{0,8}\cdot\np[m]{1,5}$ (hoch \np[m]{1,5})
Aber man darf sich dadurch (d.\ h.\ durch dieses
dauernde Granatenkrachen) nicht erschüttern
lassen -- hoffentlich kann ich mich ja bald
sinnvoller beschäftigen.\medskip

Mit den ergebensten Grüßen\nopagebreak

\hfil Ihr

\hfil\hfil Paul Lorenzen
 

%% file: 19400630-Cod-Ms-H-Hasse-1-1022-12-Lorenzen-Hasse
\label{19400630}%
\hfill Fxx Vogesen, 30.6.\medskip

\noindent\hfil Sehr geehrter Herr Professor,\medskip

\noindent darf ich Ihnen, nachdem wir "`unsern"' Feldzug 
siegreich beendet haben, von hier, aus ruhiger
Stellung, die ergebensten Grüße übermitteln.

Falls mit dieser Besatzungszeit, die hoffent\bs
lich nicht lange andauert, der Krieg für uns
tatsächlich beendet sein sollte, so wird man
nicht sagen können, daß es irgendwie schlimm
gewesen sei. Der Durchbruch durch die Maginot\bs
linie hat sich innerhalb zweier Tage vollzogen
und danach hatten wir dauernd zu marschieren
um hinter den Franzosen herzukommen, bis sie
sich dann kampflos schließlich in den Vogesen
ergeben haben.

Ich vermute aber stark, daß sich dieser
geringe Widerstand hauptsächlich aus den vorange\bs
gangenen deutschen Erfolgen erklärt.

Von Herrn Korvettenkapitän Teubner habe
ich, nachdem ich meinen Lebenslauf eingereicht
habe, nichts mehr gehört. Ich würde mich über
meine Abkommandierung natürlich -- so lange
dieser Krieg dauert (und der kann ja immer
noch sehr lange sein) -- immer noch sehr freuen,
da ich mich sehr nach einer "`geistigeren"' Be\bs
schäftigung sehne, als es hier die Besetzung eines
winzigen französischen Dorfes ist.\strut\medskip

\begin{tabular}{@{}ll@{}}%
  Mit den ergebensten Grüßen verbleibe ich stets&\\%
                                                &Ihr%
\end{tabular}

\hfill Paul Lorenzen.\par
 

%% file: 19410426-Cod-Ms-H-Hasse-1-1022-13-Lorenzen-Hasse
\label{19410426}%
\hfill Gotenhafen, 26/4\nopagebreak\medskip

\noindent\hfil Sehr geehrter Herr Professor,\medskip

Wie sich jetzt herausgestellt hat, war meine
unglückliche Stimmung, in der ich Sie in Berlin
verlassen habe, vollauf berechtigt -- obwohl doch
alles so günstig aussah.

Ich empfinde es eigentlich als aufdringlich, Sie
dauernd mit klagenden Berichten über mein Er\bs
gehen zu behelligen -- ich kann  mich nur damit ent\bs
schuldigen, daß Sie sich bisher so wohlwollend meiner
angenommen haben, und daß es eine menschliche Schwäche ist,
ein solches Wohlwollen zu benutzen. 

Ehe ich meinen Bericht fortsetze, bitte ich zu\bs
nächst diese "`offene Bemerkung"' zu entschuldigen.

Der Kommandeur der Steuermannschule hat ent\bs
schieden, daß ich aus militärischen Gründen als Lehrer
an seiner Schule ungeeignet bin. Dazu hat ihn im
wesentlichen die Beurteilung veranlaßt, die meine
letzte Dienststelle über mich geschrieben hat, in der
nämlich steht, daß ich -- kurz und gut -- unmilitärisch
sei. Hinzu kam noch persönliches Pech, das mich veran\bs
laßt hat, mich nichtsahnend hier am Montag morgen \bruch
zum Dienst zu melden, was in Berlin völlig
korrekt gewesen wäre, hier aber dazu führte, daß
ich 3~Tage Arrest bekam, da ich schon seit Sonn\bs
abend im Standort war.

Ich darf an dieser Stelle in meinem Bericht
etwas Mathematisches einschalten, worüber nach\bs
zudenken ich erfreulicherweise während der 3~Tage
endlich Gelegenheit hatte:

Nach meiner Notiz über Gruppenaxiomatik,
die Sie in Crelle aufgenommen haben, ist
eine Menge~$\mathfrak G$, in der eine stets ausführbare Ver\bs
knüpfung $a:b$ definiert ist, eine Gruppe, wenn%
\begin{enumerate}[1)]%
\item $(a:c):(b:c)=a:b$%
\item Es gibt ein~$a$, zu dem für alle~$c$ ein~$b$ existiert mit $a:b=c$%
\end{enumerate}
erfüllt sind. Es war zu vermuten, daß sich 1),~2)
zusammenfassen lassen zu einer -- naturgemäß kom\bs
plizierten -- Formel. Dies hat sich bestätigt:

Nennt man die Elemente~$x$ und $(x:x):x$ invers
zueinander und gilt für inverse Elemente~$a,a'$ und~$b,b'$ stets
\[a:(b:(a':((b':d):(c:d))))=c\]
so ist~$\mathfrak G$ eine Gruppe.\footnote{Compare \citealt{MR0011686}.}\bruch

Vermutlich erinnert Sie diese Zusammensetzung
sehr an die Allüren der «~\emph{Fundamenta mathematica}~»,
was sich vielleicht als Nachwirkung der polnischen
Vorinsassen der Arrestzelle erklären läßt.

Nachdem diese 3~Tage herum waren, hat man
mir erklärt, daß ich abkommandiert würde;
man weiß aber nicht wohin, da dieses der 2.~Admi\bs
ral der Nordseestation bewerkstelligt. Dennoch bin
ich jetzt also genau so weit, wie zu Anfang des
Krieges und möchte Sie daher auch genau wie
damals bitten, mir zu einer Beschäftigung zu
verhelfen. Irgendwelche Ansprüche zu stellen, steht mir
natürlich nicht zu, und ich möchte es auch nicht: mir ist
jede Tätigkeit irgendwo in einem Amtszimmer recht.
Nur als braver Ehemann, füge ich den persönlichen Wunsch
meiner Frau hinzu -- der allerdings auch der meine ist --
daß sie sich sehr nach einer ruhigeren Zeit sehnt, wozu
z.\ B. Kiel und Wilhelmshaven gar nicht geeignet
wären und wegen der Wohnungsverhältnisse auch
Gotenhafen und sogar Berlin nicht besonders. \bruch

Ihre Behauptung, daß Gotenhafen "`widerlich"' sei,
ist wirklich nicht zu kraß formuliert. Wir sind
hier nur behelfsweise untergekommen, und es
scheint auch unmöglich zu sein, eine Wohnung zu
finden, was sich ja aber auch erübrigt hat durch
den Lauf der Dinge, den ich versucht habe, Ihnen
darzustellen.\medskip

Mit den ergebensten Grüßen auch im Namen meiner Frau\strut

\hfil\hfil\begin{tabular}{@{}ll@{}}%
            \rlap{verbleibe ich}&\\%
            Ihr&\\%
                                &Paul Lorenzen.%
          \end{tabular}
 

%% file: 19410507-Cod-Ms-H-Hasse-1-1022-Beil-5-Hasse-Lorenzen
\begin{tabular}[t]{@{}l@{}}%
\label{19410507H}%
Korv.Kapt.\ Prof.\ Dr.\ Hasse\\%
OKM, M Wa Stb F\\%
Berlin W 35, Bissingzeile 13%
\end{tabular}\hfill Berlin, den 7.5.1941\medskip

\indent\indent Lieber Herr Lorenzen,\medskip

Ihr Brief hat mich sehr nachdenklich gemacht, und ich
halte es für richtig, Ihnen ganz offen zu schreiben, was
ich denke.

Ich bin zunächst sehr enttäuscht darüber, dass Sie sich
auf dem Kommando, das ich Ihnen vermittelt hatte, eine
derartig ungünstige militärische Beurteilung zugezogen haben.
Wenn ich mich damals dazu entschlossen hatte, mich für Ihre
Verwendung zu geistiger Arbeit einzusetzen, so geschah das
selbstverständlich in der Erwartung, dass Sie dieser meiner
Empfehlung Ehre machen würden, nicht nur durch Ihre dienst\bs
lichen Leistungen -- ganz gleich welcher Art diese sein
würden~\hbox{--,} sondern auch durch Ihr ganzes soldatisches Ver\bs
halten. Dass Sie sich selbst als unmilitärisch empfinden
und nicht die geringste Begeisterung für das Soldatsein auf\bs
bringen, ist in meinen Augen keine Entschuldigung, sondern
im Gegenteil genau das, was mir nicht gefällt. Es gibt viele
andere junge Wissenschaftler, die da eine durchaus andere und
gesundere Auffassung haben, wie etwa Teichmüller, der seine
eben erfolgte Abkommandierung zu der gleichen Tätigkeit wie
Dr.\ Franz\footnote{The mathematician Wolfgang Franz (1905--1996) works in the Cipher Department of
the Supreme Command of the Wehrmacht.} als eine Herabsetzung empfindet und viel lieber wie\bs
der zu seiner Truppe in Norwegen zurückkehren würde. Jeder
von uns muss heute seine persönliche Bequemlichkeit und seine
eigenen Wünsche zurückstellen und sich in das grosse Ganze
willig einfügen. Menschen, die sich dem entziehen wollen,
können wir heute nicht gebrauchen, und sie gelten heute mit
Recht nichts. Sie sind jung und körperlich kräftig. Sie
haben für wissenschaftliche Arbeit noch ein langes Leben vor
sich. Ich an Ihrer Stelle wäre glücklich, wenn ich in vorder\bs
ster Front mit dabei sein könnte, wo immer um unsere Zukunft
gekämpft wird. Und wenn Sie schon von sich aus keine Begeiste\bs
rung für das Soldatsein aufbringen können, dann vergessen Sie
doch bitte nicht, dass die Haltung, die Sie heute einnehmen,
für alle Zukunft bei Ihrer Beurteilung ganz entscheidend mit\bs
spricht. Gerade von der Intelligenz unseres Volkes muss man \bruch
mit vollstem Recht erwarten, dass sie auf Grund der vertieften
Einsicht in die harte Notwendigkeit und tiefste Berechtigung
dieses Kampfes mit allem ihren Denken und Fühlen bei der
kämpfenden Truppe ist und in ihrer Einsatzbereitschaft allen
anderen als leuchtendes Vorbild vorangeht. Es ist grundfalsch
zu sagen, dass dies mit Ihrem späteren Beruf als Wissenschaftler
nichts zu tun hat. Im späteren Leben werden eben die Menschen
nicht nur nach ihren beruflichen Leistungen gemessen, sondern
es wird der ganze Mann gewogen. Und es ist voll gerecht\bs
fertigt, wenn ein Mann, der in dieser höchsten Bewährungszeit
zu leicht befunden wurde, auch für später diesen Stempel auf
sich trägt.

Noch ist es für Sie Zeit, den Eindruck, den Sie bei mir
und wohl auch schon bei anderen durch Ihre bisherige Haltung
erzeugt haben, zu entkräften. Gerade weil ich glaube, dass
Sie später einmal in der Wissenschaft vorwärts kommen können,
mit Ihren schönen Gaben und Ihrem scharfen Blick für das
Wesentliche in komplizierten begrifflichen Zusammenhängen,
möchte ich heute Ihren Wunsch nach einer stillen ruhigen
Amtszimmertätigkeit, weit entfernt vom Getöse des Krieges,
in keiner Weise begünstigen, ganz abgesehen davon, dass ich
es aus den genannten Gründen auch vor mir selbst gar nicht
verantworten kann, noch einmal bei einer militärischen
Dienststelle für Sie einzutreten. Tun Sie zunächst einmal
als Soldat Ihre Pflicht -- und mehr als das~--, wohin auch
immer man Sie jetzt stellt. Dann werden Sie später auch für
sich selbst das schöne Gefühl haben, sich die Segnungen des
kommenden Friedens wirklich verdient zu haben und nicht bloss
durch andere mit dem Einsatz ihres Lebens haben verdienen zu
lassen.\strut\medskip

\begin{tabular}{@{}rr@{}}%
  Mit besten Grüssen, auch an Ihre Frau,&\\%
  Ihr&[Hasse]%
\end{tabular}
 

%% file: 19410507-Cod-Ms-H-Hasse-1-1021-Lorenzen-Hasse
\hfill\begin{tabular}{@{}l@{}}%
\label{19410507}%
  Gotenhafen, 7.5.41\\%
  Leuthenstr.~21 b.\ Gaspersen%
\end{tabular}\medskip

\noindent Sehr geehrter Herr Professor,\medskip

\noindent Sie haben durch den letzten Brief
meines Mannes ja schon von unserem
"`Unglück"' erfahren. Gewiß hat es
meinen Mann nicht unverdient ge\bs
troffen, da er wohl wirklich gerade\bs
zu typisch unmilitärisch ist, aber
es ist ja nicht damit zu ändern,
daß man ihn deswegen bestraft u.\bsp
in häßlicher Art beschimpft. Aber
hier legt man leider mehr Wert auf
militärische Erziehung, während wir \bruch
uns in Berlin schon darauf gefreut
hatten, daß das hier mehr in den
Hintergrund treten würde. Ich habe
schon oft versucht, meinen Mann
auf sein unmilitärisches Verhalten
aufmerksam zu machen, obwohl ich
nicht unbedingt dafür bin, daß
er sich grundlegend ändert, aber für
Kriegsdauer wäre ich schon damit
einverstanden. Das Schlimmste ist,
daß sich dies alles, wie man an
anderen Mathematikern sieht, noch
nicht einmal mit seiner Wissen\bs
schaft entschuldigen läßt; und doch
ist es ja die einzige Entschuldi\bs\bruch
gung, denn es gibt ja wohl nur sehr
wenige Stunden an einem Tage, in denen
mein Mann nichts Mathematisches im
Kopf hat. Nun erscheint es mir recht
bedauerlich -- um nicht ungerecht
zu sagen -- daß er unter Unter\bs
offiziersaufsicht Flure fegt, Striche
zieht oder ähnliche Tätigkeiten ausübt.
Da ich nun täglich sehe, wie gerade\bs
zu unglücklich mein Mann sich hier\bs
bei fühlt, wage ich es auch persön\bs
lich -- ohne Wissen meines Mannes --
mich an Sie zu wenden. Ich weiß
ja gar nicht, wie weit es Ihnen möglich
sein wird, an diesem Zustand etwas
zu ändern, aber ich vermute doch, daß Sie \bruch
unter Berücksichtigung der geschilderten
Verhältnisse vielleicht schon mit einem
guten, "`richtigen"' Rat helfen können.
Ist es wohl nicht möglich, daß mein
Mann (den die Marine ja nun auf
Kriegsdauer beschäftigen muß u.\ den
sie wegen Seekrankheit nicht auf ein
Schiff stecken kann) irgendwo arbeiten
könnte, wo er einen vernünftigen Vor\bs
gesetzten hat, der nicht \uline{nur} das Mili\bs
tärische an einem Menschen gelten läßt.

\noindent Es ist natürlich eine sehr große Hoff\bs
nung mit wenig Aussicht darauf, daß
diese Möglichkeit wirklich gefunden
werden kann, aber Sie werden das
sicherlich am besten beurteilen können,
ob es überhaupt noch lohnt, \bruch
sich an diese Hoffnung zu klammern.
Ich denke mir auch, daß Sie am ehesten
Gelegenheit haben werden, eine vernünf\bs
tige, nutzbringende Beschäftigung
herauszufinden, denn Sie reisen von
dem zentralen Berlin aus ja in ganz
Europa herum.

\noindent Vor allem aber vertraue ich darauf,
daß Sie  bereit sind, meinem Mann
zu helfen, da Sie bisher schon sozu\bs
sagen sein "`guter Engel"' gewesen
sind -- soweit ich aus Erzählungen
weiß schon vor dem Kriege, u.\ wie ich
aus Erfahrung weiß auch im Kriege.

\noindent Ich darf Ihnen versichern, wie dankbar
ich Ihnen hierfür bin, ohne fürchten
zu müssen, daß Sie diesen Dank für \bruch
eine Phrase halten. Und ich
freue mich auch, Ihnen meinen
Dank nun auch ganz persönlich
aussprechen zu können.\medskip

\noindent Mit freundlichen Grüßen bin ich\strut

\hfil\begin{tabular}{@{}ll@{}}%
  Ihre&\\%
  &sehr ergebene\\[3pt]%
  \strut\rlap{Käthe Lorenzen.}&%
\end{tabular}
 

%% file: 19410517-Cod-Ms-H-Hasse-1-1022-15-Lorenzen-Hasse
\begin{tabular}[t]{@{}l@{}}
\label{19410517}%
  F. B. Gfr.\ Lorenzen\\%
  2.~Komg.\ Strm.schule%
\end{tabular}\hfill Gotenhafen, 17/5\medskip

\hfil Sehr geehrter Herr Professor,\medskip

Es fällt mir jetzt wesentlich schwerer, Ihnen
zu schreiben, als das letzte Mal, wo es sich um
rein militärische Dinge handelte, während ich dieses
Mal gezwungen bin, mich persönlich vor Ihnen zu
rechtfertigen. Ich kann es nämlich unmöglich er\bs
tragen, daß Sie aus meinem letzten Brief, in
dem wohl steht, daß ich gern auf Kriegsdauer
einem Amtszimmer zugewiesen würde, heraus\bs
lesen, ich wolle mich drücken.

Ich darf Ihnen ganz kurz meine militärische
Vergangenheit schildern: Nach dem freiwilligen
Arbeitsdiensthalbjahr habe ich 1933/34 freiwillig ak\bs
tiv  gedient und bin dabei ein Jahr lang im Pferdestall
so behandelt worden, daß ich für militärischen Arbeits\bs\bruch
dienst, Exerzieren, Ehrenbezeigung erweisen u.\ ä.\bsp
keine Begeisterung mehr übrig habe.

Ich sollte das nicht so laut sagen, aber da
ich eben so sicher fühle, daß dies nichts mit
einem Sich-Drücken-Wollen zu tun hat, unterlasse
ich's leider doch selten.

Beim Durchbruch durch die Maginot-Linie habe
ich es bestätigt bekommen, daß ich keine Angst
vor dem Einsatz des Lebens habe. Dort bin ich
nicht als "`zu leicht"' befunden. Trotzdem bin ich
dort gern weggegangen, denn meine Hauptbeschäf\bs
tigung außer Warten war wieder Pferdepflegen
und in Ruhetagen Exerzieren. Ich wollte mich
wirklich gerne nützlicher betätigen.

Wie es mir im O. K. M. ergangen ist, wissen
Sie selbst. Ich bin viermal bei meinem Vorgesetzten
gewesen, um zu sagen, daß ich nun das Abschreiben
begriffen habe -- wovon man aber keine Notiz \bruch
genommen hat. Daraufhin habe ich -- törichterweise
wie immer -- Kameraden gegenüber nicht verhehlt,
daß ich auch für das Abschreiben keine Begeisterung
mehr übrig hätte. Und daher stammt die schlechte
Beurteilung.

Ich werde es wohl kaum noch lernen, mich zu
verstellen; aber ich sehe ein, daß ich dann
kein Recht habe, anders behandelt werden zu
wollen, als ein Drückeberger. Da ich jetzt in
der Funk-Beobachter-Laufbahn bin, ist es so gut
wie sicher, daß ich erst zu meiner Stamm-Abteilung
kommandiert werde, wo ich zu warten habe, (d.\ h.\bsp
Kasernenreinigen und Exerzieren) bis bei einem
Marine-Nachrichten-Offizier oder in Berlin ein Funk-\relax
Beobachter-Gefreiter gebraucht wird.

Obwohl -- wie Sie schreiben -- nach dem
Kriege noch viel Zeit sein wird, Mathematik \bruch
zu treiben, ist es doch für mich ein bitteres
Gefühl, gerade in diesen Kriegsjahren völlig
nutzlos sein zu müssen.

Sehr geehrter Herr Professor, ich hoffe mit Zu\bs
versicht, Sie werden es mir nicht verargen,
daß ich Ihnen so ausführlich und offen meine
Einstellung geschildert habe. Der Notwendigkeit,
für Deutschland das zu tun, was man von mir
verlangt, gleichgültig, ob ich dafür geeignet 
bin oder nicht, werde ich mich selbstverständlich
fügen. Es wäre mir eine große Freude,
wenn Sie mir gestatten würden, Ihnen demnächst
zu schreiben, wohin ich abkommandiert werde.

Mit den ergebensten Grüßen auch von meiner Frau\medskip

\hfil Ihr

\hfil\hfil\hfil Paul Lorenzen.
 

%% file: 19420318-Cod-Ms-H-Hasse-1-1022-14-Lorenzen-Hasse
\hfill\begin{tabular}{@{}l@{}}%
\label{19420318}%
        Wesermünde, 18/3\\%
        4/I Marineschule%
\end{tabular}\medskip

\hfil Sehr geehrter Herr Professor,\medskip

für Ihren freundlichen Glückwunsch möchten
meine Frau und ich Ihnen unsern herzlichsten
Dank aussprechen.

Daß unser Töchterlein ebenfalls Jutta
heißt, bitte ich nicht als Plagiat auf\bs
fassen zu wollen --~-- soweit ich weiß, hat
meine Frau davon unabhängig diesen Namen gewählt.

Ich bin Ihnen sehr dankbar für Ihren
Hinweis, daß \emph{J. Dieudonné} sich für die
multiplikative Idealtheorie interessiert.\footnote{See \citealt{Die1941}.} Es
wäre sehr schön, wenn eine "`Kollaboration"' \bruch
sich anbahnen ließe, --~-- jedenfalls werde
ich Herrn \emph{Dieudonné} bald einmal schreiben.

Mit meiner Tätigkeit hier an der
Schule bin ich sehr zufrieden, da ich hier,
zum ersten Mal im Kriege, das Gefühl
habe, nicht überflüssig zu sein.\medskip

Mit den ergebensten Grüßen verbleibe ich

\hfil Ihr\nopagebreak

\hfill Paul Lorenzen.\par
 

%% file: 19420402-Cod-Ms-Hasse-33-3-Hasse-Lorenzen
\begin{tabular}[t]{@{}l@{}}%
\label{19420402}%
  Korv.Kapt.\ [Prof.\ Dr.\ Hasse]\\%
  Berlin-Wannsee\\%
  Am Sandwerder 5%
\end{tabular}\hfill\begin{tabular}[b]{@{}l@{}}%
                     \llap[Betr.\ Ms.\ u.\ Korr.\\%
                     Hans R. Weber, München]%
                   \end{tabular}

\hfill Berlin, den 2.4.1942

\hfil\begin{tabular}{@{}llll@{}}%
       Herrn&&&\\%
            &Dr.\ P. Lorenzen&&\\%
            &&4/I M S&\\%
            &&&Wesermünde%
\end{tabular}\medskip

\hfil Lieber Herr Lorenzen,\hfil\medskip

\noindent Besten Dank für Ihren freundlichen Brief. Ich sende Ihnen bei\bs
liegend eine Korrektur, über die ich gerne Ihren Rat erbitten
möchte. Ich hatte damals das Ms.\ nur mit grossem Zögern Herrn
Perron zuliebe angenommen, der es empfahl. Nachdem nun aber
der Satz 6, der mich hauptsächlich zur Annahme bestochen hatte,
weggefallen ist, glaube ich nicht, dass die Arbeit noch irgend\bs
etwas Neues oder Bedeutendes bringt. Wie denken Sie darüber?\strut\medskip

\begin{tabular}{@{}ll@{}}%
  Mit besten Grüssen an Sie und Ihre Frau&\\%
                                         &Ihr [Hasse]%
\end{tabular}
 

%% file: 19420407-Cod-Ms-Hasse-33-3-Lorenzen-Hasse
\begin{tabular}{@{}l@{}}%
\label{19420407}%
  Wesermünde, 7/4\\%
  4/I Marineschule%
\end{tabular}\medskip

\hfil Sehr geehrter Herr Professor,\medskip

Das beiliegende Manuskript, das ich heute erhielt,
kann ich leider wirklich nicht zur Veröffentlichung
empfehlen, denn das Axiomensystem~I-III ist
schon etwa 1905 von Huntington\footnote{\selectlanguage{english}This probably refers to \citealt{huntington05}.} aufgestellt. Garver\footnote{\selectlanguage{english}This probably refers to \citealt{garver34}.}
hat außerdem inzwischen bewiesen, daß~I überflüssig
ist. Schließlich läßt sich~III noch wesentlich
abschwächen, worüber allerdings noch nichts veröf\bs
fentlicht ist.\footnote{\selectlanguage{english}Soon after, Lorenzen will submit \citealt{MR0011686} to remedy this.}

Falls Herr Weber sich für diese Dinge interessiert,
stehe ich ihm jederzeit nach Möglichkeit zur Verfü\bs
gung. Vielleicht sind Sie so liebenswürdig Herrn We\bs
ber dies mitzuteilen, wenn Sie das Ms.\ zurück\bs
senden sollten.\strut\medskip

\begin{tabular}{@{}ll@{}}%
  Mit den ergebensten Grüßen verbleibe ich Ihr&\\%
                                              &Paul Lorenzen%
\end{tabular}
 

%% file: 19420102-PL-1-1-136-Krull

\begin{tabular}{@{}c@{}}%
\label{19420102}%
  Mathematisch-naturwissenschaftliche\\%
  Fakultät\\%
  der Rheinischen Friedrich-Wilhelms-Universität\\[\medskipamount]%
  J.-Nr. 397%
\end{tabular}\medskip%

\hfill
Bonn, den 2.\ Januar 1942\medskip

\hfil Bescheinigung\medskip

Dem M.Gefr.\ Dr.\ Paul \so{Lorenzen} wird bescheinigt,
dass er die Bedingungen erfüllt, die von der Reichshabilitations\Bs
Ordnung für die Zulassung zur Habilitation gefordert werden. Falls
Dr.\ Lorenzen imstande ist eine Arbeit rechtzeitig vorzulegen, die
von der Fakultät als Habilitationsschrift angenommen wird, kann
das Habilitationsverfahren noch in diesem Semester abgeschlossen
werden.\strut\medskip

\hfill
\begin{tabular}{@{}c@{}}
  Der \so{Dekan}\\
  Krull
\end{tabular}

%% file: 19420512
\noindent BUK -- Az Lorenzen\hfill Bonn, den 12.\ Mai 1942\hspace*{1cm}\medskip
\label{19420512}%

\leavevmode\llap{1) }\begin{tabular}[t]{@{}l@{}}%
            An den\\%
            Herrn Rektor der Universität\\%
            \hfil\uline{\hbox{\so{Bonn}}}\\[2pt]%
            -- mit 2 Durchschlägen --
          \end{tabular}\medskip

Der Direktor des Mathematischen Seminars
beantragt mit den \uline{gegen gefl.\ Rück\-gabe} bei\bs
liegenden Unterlagen Dr.\ Paul \so{Lorenzen}
vom 1.\ Mai 1942 ab zum wissenschaftlichen Assistenten zu ernennen. Ich bitte,
die Stellungnahme\comment{See \citet[pages 174--181]{MR1991149} for the signification of the \emph{Dozentenschaft} reports: ``Frequently, the \emph{Dozentenführer} had no idea about an 
individual and had to ask a politically trusted member of his discipline for a report,
which would then be passed on verbatim''. It is plausible that the views expressed in the subsequent statement are Krull's.} des Dozentenschaftsleiters
und Dozentenbundsführers zu den weltanschau\bs
lichen und charakterlichen Voraussetzungen des
Vorgeschlagenen herbeizuführen und selbst zu
dem Anstellungsantrage zusammenfassend Stel\bs
lung zu nehmen.\medskip

\hfil I.\ A.\nopagebreak

\leavevmode\llap{2) }\uline{Wv.\ 25.~Mai 1942}.


%% file: 19420601
\noindent\begin{tabular}[b]{@{}c@{}}%
\label{19420601}%
           \large Dozentenschaft der\\%
           \large Universität Bonn%
         \end{tabular}\hfill{\footnotesize\begin{tabular}[t]{@{}l@{}}%
                              \normalsize\textbf{Bonn, }den 1.\ Juni 1942\\%
                              Poppelsdorfer-Schloß\\%
                              Fernruf 6294$$
                            \end{tabular}}\medskip

\noindent\begin{tabular}{@{}l@{}}%
           Jahrbuch-Nr.~112\\%
           \small bei allen Antworten angeben
         \end{tabular}\medskip

\noindent
\begin{tabular}{@{}l@{}}%
  An\\%
  Se.\ \begin{tabular}[t]{@{}l@{}}%
            Magnifizenz den Herrn Rektor der\\%
            Rhein.\ Friedr.-Wilhelms-Universität,%
          \end{tabular}\\[2pt]%
  \uline{\hbox{\so{Bonn}}}%
\end{tabular}\medskip

\noindent \uline{Betr.\ Dr.\ P. \hbox{\so{Lorenzen}}, Ihr Schr.\ v.\ 16.5.42 Nr.~983}\medskip

\noindent Dr.\ Lorenzen ist ein vielversprechender Mathematiker mit Begei\bs
sterung für seine Wissenschaft und großem Arbeitseifer.

\noindent Soweit bekannt, ist die politische Einstellung einwandfrei.
Es soll jedoch nicht verschwiegen werden, daß L. wohl einige Cha\bs
raktermängel aufzuweisen scheint, die es z.\ Zt.\ jedenfalls nicht
erwünscht erscheinen lassen, ihn etwa zur Dozentur zuzulassen.
L. neigt zu starker Selbstüberschätzung, was ihm offenbar auch
in seiner Laufbahn bei der Wehrmacht geschadet hat.

\noindent Angesichts seiner großen Befähigung ist es daher nur erwünscht,
wenn er jetzt Gelegenheit erhält als Assistent zu zeigen, wie
weit er sich einzufügen vermag, und wie weit er sich für die
Hochschullaufbahn eignet.

\noindent Der Antrag des Herrn Direktors des Mathematischen Seminars wird
daher be\-für\-wortet.\strut\medskip

\noindent\hfil\begin{tabular}{@{}c@{}}
                Heil Hitler!\\%
                Klapp.\footnotemark\\%
                Dozentenführer i.\ V.%
              \end{tabular}\footnotetext{\selectlanguage{english}Ernst Klapp (1894--1975) is professor of agricultural sciences at the University of Bonn.}
 

%% file: 19430507-PL-1-1-140-Lorenzen-Krull

\label{19430507}%
\hfill
\begin{tabular}{@{}l@{}}%
Wesermünde-Lehe, den 7.~Mai 1943\\%
Hafenstr.~92%
\end{tabular}\medskip%

\indent\indent Sehr geehrter Herr Professor,\medskip

\noindent von den "`Fortschritten"' erhielt ich Ihren Beitrag~VIII\comment{\citealt{krull43}.} über die
$\varLambda$-\hspace{0pt}Operationen zum Referat.\comment{The \emph{\foreignlanguage{german}{Jahrbuch über die Fortschritte der Mathematik}} ceases to appear with the volume on the year 1942.} Leider habe ich nun Ihren Beitrag\nbsp
I nicht zur Hand, dessen Kenntnis dazu nötig wäre. Insbesondere
weiß ich nicht mehr, ob und in welchem Zusammenhang der folgende
Satz darin enthalten ist:\smallskip

\noindent Ist $w$~eine arithmetisch brauchbare Operation eines Integritäts\bs
bereiches~$R$ (Quotientenkörper~$K$), so entsprechen die arithmetisch
brauchbaren $'$-\hspace{0pt}Operationen, für die stets $\su a_w\subseteq\su a'$ gilt, einein\bs
deutig denjenigen Quotientenringen~$M'$ des Funktionalringes~$M$
mit $M'\cap K=R$.\smallskip

Dieser Satz läßt sich direkt mit der Halbgruppe der ganzen Ideal\bs
brüche beweisen. Denn der Übergang von~$\su a_w$ zu~$\su a'$ ist nichts
anderes als eine homomorphe Abbildung der $w$-\hspace{0pt}Ideale auf die $'$-\hspace{0pt}Ideale.

Da in Beitrag~VIII Ihre Bemerkungen zu den Definitionsfor\bs
meln der $'$-\hspace{0pt}Operationen (insbesondere Anmerkung~21) leider durch
Druckfehler entstellt sind, möchte ich die Abhängigkeiten zwischen
den Formeln:%
\[%
  \begin{aligned}[t]%
    &1)\enskip\su a\subseteq\su a'\\%
    &2')\enskip(\su a'+\su b')'=(\su a+\su b)'\\%
    &3)\enskip(a\su a)'=a\su a'\\%
    &4)\enskip\su R=\su R'%
  \end{aligned}\quad%
  \begin{aligned}[t]%
    &2)\enskip\su a\subseteq\su b'\rightarrow\su a'\subseteq\su b'\\%
    2''&\textrm a)\enskip\su a\subseteq\su b\rightarrow\su a'\subseteq\su b'\quad2''\textrm b)\enskip\su a''=\su a'\\%
    &3')\enskip(\su a'\cdot\su b')'=(\su a\cdot\su b)'\\%
    &4')\enskip(a)=(a)'%
  \end{aligned}%
\]%
hier darlegen.\comment{Compare \citealt[§~4]{Lor1950}.}\smallskip

\noindent I. Unter Voraussetzung von $1)$ ist $2)$ gleichwertig mit $2''a)$ und
$2''b)$.\smallskip\bruch

Beweis: \begin{tabular}[t]{@{}rll@{}}
  Es gelte $1)$ und $2)$. Dann gilt&$\su a\subseteq\su b\rightarrow\su a\subseteq\su b'$&nach $1)$\\%
  und&$\su a\subseteq\su b\rlap{$'$}\rightarrow\su a'\subseteq\su b'$&nach $2)$\\%
  also&$\su a\subseteq\su b\rightarrow\su a'\subseteq\su b'$,&\\%
  ferner gilt&$\su a'\subseteq\su a''$&nach $1)$\\%
  und&$\su a'\subseteq\su a'\rightarrow\su a''\subseteq\su a'$&nach $2)$\\%
  also&$\su a'=\su a''$&\\%
\multicolumn3l{Es gelte $2''a)$ und $2''b)$. Dann gilt}\\%
&$\su a\subseteq\su b'\rightarrow\su a'\subseteq\su b''$&nach $2''a)$\\%
  also&$\su a\subseteq\su b'\rightarrow\su a'\subseteq\su b'$&nach $2''b)$%
\end{tabular}\smallskip

\noindent II. Unter Voraussetzung von $1)$ ist $2)$ gleichwertig mit $2')$.\smallskip

\noindent Beweis: \begin{tabular}[t]{@{}rll@{}}%
         Es gelte $1)$ und $2)$. Dann gilt&$\su a'\subseteq(\su a+\su b)'$&\\%
         und&$\su b'\subseteq(\su a+\su b)'$&nach $2''a)$\\%
         also&$\su a'+\su b'\subseteq(\su a+\su b)'$&\\%
         und&$(\su a'+\su b')'\subseteq(\su a+\su b)'$&nach $2)$\\%
         \multicolumn3l{Es gelte $1)$ und $2')$. Dann gilt $\su a''=(\su a'+\su a')'=(\su a+\su a)'=\su a'$}\\%
         \multicolumn3l{also $2''b)$. Ferner gilt}\\%
         \multicolumn2l{\hspace{9em}$\su a\subseteq\su b\rightarrow\su a+\su b=\su b$}&\\%
         \multicolumn2l{\hspace{9em}$\hphantom{\su a\subseteq\su b}\rightarrow(\su a+\su b)'=\su b'$}&\\%
         \multicolumn2l{\hspace{9em}$\hphantom{\su a\subseteq\su b}\rightarrow(\su a'+\su b')'\subseteq\su b'$}&nach $2')$\\%
         \multicolumn2l{\hspace{9em}$\hphantom{\su a\subseteq\su b}\mathrel{\dot\rightarrow}\su a'+\su b'\subseteq\su b'$}&nach $1)$\\%
         \multicolumn2l{\hspace{9em}$\su a\subseteq\su b\rightarrow\su a'\subseteq\su b'$}&%
       \end{tabular}\smallskip
       
\noindent III. Unter Voraussetzung von $1)$ und $2)$ ist $3)$ gleichwertig mit $3')$.\smallskip

\noindent Beweis:
\begin{tabular}[t]{@{}rlr@{}}%
  \multicolumn{3}{l}{Es gelte $1)$, $2)$ und $3)$. Dann gilt für $b\in\su b$}\\%
  &$\su a'\cdot b\subseteq(\su a\cdot b)'\subseteq(\su a\cdot\su b)'$&\\%
  &$\su a'\cdot\su b\subseteq(\su a\cdot\su b)'$&\llap{nach $3)$ und $2''a)$}\\%
  \multicolumn3l{Genau so beweist man $\su a'\cdot\su b'\subseteq(\su a'\cdot\su b)'$}\\%
  \multicolumn1l{also gilt}&$\su a'\cdot\su b'\subseteq(\su a\cdot\su b)''=(\su a\cdot\su b)'$&nach $2''b)$\\%
  &$(\su a'\cdot\su b')'=(\su a\cdot\su b)'$&\llap{nach $2)$ und $2''a)$}\\%
  \multicolumn3l{Es gelte $1)$ und $3')$. Dann gilt}\\%
  \multicolumn2c{$a\cdot\su a'\subseteq(a)'\cdot\su a'\subseteq((a)'\su a')'=(a\su a)'$}&\llap{nach $1)$ und $3')$}\\\noalign{\bruch}%
  Genau so beweist man&$%
                        \begin{aligned}[t]%
                          a^{-1}(a\su a)'&\subseteq(a^{-1}a\su a)'\\%
                          a^{-1}(a\su a)'&\subseteq\su a'\\%
                          (a\su a)'&\subseteq a\su a'%
                        \end{aligned}$&%
\end{tabular}\smallskip
       
\noindent IV. Unter Voraussetzung von $3)$ ist $4)$ gleichwertig mit $4')$.\smallskip

\noindent Beweis:
\begin{tabular}[t]{@{}rl@{}}%
\multicolumn{2}{l}{Es gelte $3)$ und $4)$. Dann gilt}\\%
  \multicolumn{2}{c}{$(a)'=(a\cdot\su R)'=a\su R'=a\cdot\su R=(a)$}\\%
  Es gelte $4')$. Dann gilt&$\su R=(1)=(1)'=\su R'$%
\end{tabular}\smallskip

\noindent Eine $'$-\hspace{0pt}Operation läßt sich also durch~$1)$, $2')$, $3')$,~$4)$ definieren,
wobei~$2')$ und~$3')$ nichts anderes als die Homomorphie aussagen.
Sind daher die $w$-\hspace{0pt} und die $'$-\hspace{0pt}Operationen arithmetisch brauchbar,
so hat man eine homomorphe Abbildung der Halbgruppe der $w$-\hspace{0pt}Ideale
auf die Halbgruppe der $'$-\hspace{0pt}Ideale, und damit auch eine homomorphe
Abbildung der Halbgruppe~$\su g_w$ der ganzen $w$-\hspace{0pt}Idealbrüche auf die
Halbgruppe~$\su g'$ der ganzen $'$-\hspace{0pt}Idealbrüche.	\ $\su g_w$~ist vollständig
(nach der Terminologie in Prüfer und meiner Dissertation, ich würde
jetzt lieber sagen: "`$\su g_w$	ist Verbandshalbgruppe"'). Die Homomorphis\bs
men einer Verbandshalbgruppe entsprechen aber eineindeutig \barre[deren]{ihren}
Quotientenhalbgruppen. (Das ist Ihr Satz aus Beitrag~I, daß jeder
Oberring eines Hauptidealringes stets Quotientenring ist, denn die
Homomorphismen einer Verbandshalbgruppe entsprechen eineindeutig den
Oberhalbgruppen, die mit~$a$ und~$b$ auch~$a\lor b$ enthalten --~vgl.\ in
meiner Dissertation S.\barre[~545]{55}, Absatz~1~--).\comment{\citealt[first paragraph on page~545]{Lor1939}, reproduced on page~\pageref{1939p545} and spelled out on page~\pageref{1939p545t}.}

Ich wäre Ihnen dankbar, wenn Sie mir mitteilen würden, ob
Ihnen eine arithmetisch brauchbare $'$-\hspace{0pt}Operation bekannt ist, die
keine $\varLambda$-\hspace{0pt}Operation ist. Es müßte dann $\su g'$ aus $\su g_w$ durch ein multi\bs
plikativ abgeschlossenes System~$S$ entstehen, das nicht nur Ideale,
sondern auch Idealbrüche enthält.\strut\medskip

\begin{tabular}{@{}lll@{}}%
  \strut Mit den ergebensten Grüßen verbleibe ich\\%
  &Ihr\\%
  &&Lo[renzen]%
\end{tabular}

%% file: 19430920-PL-1-1-139-Lorenzen-Krull

\label{19430920}%
\hfill
\begin{tabular}{@{}l@{}}%
Wesermünde-Lehe, den 20.~September 1943\\%
Hafenstr.~92%
\end{tabular}\medskip%

\noindent Sehr geehrter Herr Professor,\ajout{ Krull}\medskip

\noindent ich möchte Ihnen mitteilen, daß sich die Identität des
$r_a$-\hspace{0pt}Idealsystems mit dem $r_b$-\hspace{0pt}Idealsystem auch für beliebige
halbgeordnete Gruppen mit beliebigem $r$-\hspace{0pt}Idealsystem bewei\bs
sen läßt.

\noindent Die Beweismethode liefert für die kommutativen Integritäts\bs
bereiche folgende Vereinfachung:

\noindent Es sei $I$~ein Integritätsbereich (Quotientenkörper~$K$),
$\su a$~ein Dedekindsches $I$-Ideal, $x\in K$.\smallskip

\noindent Hilfssatz: Ist $a\in K$ ganz abhängig von~$\su aI[x]$ und~$\su aI[x^{-1}]$,
so ist~$a$ ganz abhängig von~$\su a$.\smallskip

\noindent Beweis: Es gibt endliche Ideale~$\su e_1$ und~$\su e_2$ mit
\[%
  \begin{aligned}%
    \su e_1a&\subseteq\su e_1\su aI[x]\\%
    \su e_2a&\subseteq\su e_2\su aI[x^{-1}]%
  \end{aligned}%
\]
Also gilt für~$\su e=\su e_1\su e_2$ und geeignete $n_1$, $n_2$%
\[%
  \begin{aligned}%
    \su ea&\subseteq\su e\su a(1,x,\dots,x^{n_1})\\%
    \su ea&\subseteq\su e\su a(1,x^{-1},\dots,x^{-n_2})%
  \end{aligned}%
\]
und daher für $n=n_1+ n_2$%
\[\su ea(1,x,\dots,x^n)\subseteq\su e\su a(1,x,\dots,x^n)\]

\noindent \uline{Satz:} Ist $a$ nicht ganz abhängig von~$\su a$, so gibt es einen
Bewertungsoberring~$B$ von~$I$ mit $a\notin\su aB$.\smallskip

\noindent Beweis: Es gibt einen maximalen Oberring~$B$ von~$I$, für den
$a$ nicht ganz abhängig von $\su aB$ ist (Wohlordnungsschluß). \bsp
Aus $B\subset B[x]$ und $B\subset B[x^{-1}]$ würde folgen, daß $a$ ganz abhängig\bruch{}
von $\su aB[x]$ und $\su aB[x^{-1}]$ ist, also von~$\su aB$ (Hilfssatz). \bsp
Also gilt $x\in B$ oder $x^{-1}\in B$.\smallskip

\noindent In Schiefkörpern gilt der Hilfssatz nicht für die übliche
Ganzabhängigkeit. Nennt man aber ein Element~$a$\smallskip%

{\leftskip2.5\parindent\parindent0pt%
  \leavevmode\llap{1. }$d_0$-abhängig von~$\su a$, wenn $a\in\su a$%
  
  \leavevmode\llap{2. }$d_{n+1}$-abhängig von~$\su a$, wenn $a$ $d_n$-abhängig von~$\su aI[x]$ und~$\su aI[x^{-1}]$%

  \leavevmode\llap{3. }$d$-abhängig von~$\su a$, wenn $a$ $d_n$-abhängig von~$\su a$ für
mindestens ein~$n$,\smallskip%

}\noindent so gilt der Hilfssatz für die $d$-Abhängigkeit. Also gilt auch
der Satz für die $d$-Abhängigkeit. Die Umkehrung des Satzes
ist trivial.

Indem ich hoffe, daß Sie diese Mitteilung interessiert
hat, verbleibe ich mit den ergebensten Grüßen\medskip

\hfil[Ihr Lorenzen]\comment{The following notes were written on this carbon copy by hand. \smallskip

  Für kom[mutative] Gr[uppen]:%
\[%
  a\mathrel{\alpha_r}\su a\mid\su g\prec\bigcurlyvee_{\su e}\su ea\subseteq\su e\su a\su g%
\]

Beweis durch Ind[uktion]: $a\mathrel{\alpha_r^0}\su a\mid\su g\prec a\in\su a\su g$%
\[%
  \begin{aligned}%
    a\mathrel{\alpha_r^{n+1}}\su a\mid\su g&\prec\bigcurlyvee_xa\mathrel{\alpha_r^n}\su a\mid\su g(x)_r\land a\mathrel{\alpha_r^n}\su a\mid\su g(x^-)_r\\%
    \text{Ind[uktions]vor[aussetzung]}\quad&\prec\bigcurlyvee_{\makebox[0pt][c]{$\scriptstyle\su e_1,\su e_2$}}\su e_1a\subseteq\su e_1\su a\su g(x)_r\land\su e_2a\subseteq\su e_2\su a\su g(x^-)_r\\%
    &\prec\bigcurlyvee_{\makebox[0pt][c]{$\scriptstyle\su e=\su e_1\su e_2$}}\su ea\subseteq\su e\su a(1,\dots,x^{n_1})_r\land\su ea\subseteq\su e\su a(1,\dots,x^{-n_2})_r\\%
    &\prec\bigcurlyvee_{\makebox[0pt][c]{$\scriptstyle n=n_1+n_2$}}\su ea(1,\dots,x^n)_r\subseteq\su e\su a(1,\dots,x^n)_r%
  \end{aligned}%
\]}\par
 

%% file: 19440104-PL-1-1-149-Krull-Lorenzen

\noindent\begin{tabular}[t]{@{}l@{}}%
\label{19440104}%
           Krull\\%
           Bonn\\%
           Kaiser-Friedrichstr.~18%
         \end{tabular}\hfill%
         \begin{tabular}[t]{@{}l@{}}%
           Sonderführer (M)\\%
           Dr.\ Paul Lorenzen\\%
           Wesermünde\\%
           Hafenstr.~92%
         \end{tabular}\medskip%

\hfill Bonn, 4.1.44\medskip

\indent\indent Lieber Herr Lorenzen!\medskip

\noindent Durch den stellvertretenden Dekan, Herrn v.\ Antropoff erfuhr
ich von Ihrem Habilitationsgesuch. Es freut mich sehr, dass Sie
so bald nach Ihrem Unteroffizierskurs Ihre Arbeit abschliessen
konnten. Dumm ist es nur, dass ich nun das Exemplar, das Sie
mir zuschickten, nicht bekommen habe. Es liegt sicher in Greifswald,
während ich hier zunächst festgehalten bin, da ich jetzt nach Ab\bs
lauf meines "`Genesungsurlaubs"' nochmals --~oder richtiger zum
ersten Mal, denn ich war bisher dauernd in ambulanter Be\bs
handlung~-- ins Lazarett muss. Vielleicht bitte ich mir da
einfach vom Dekan ein Exemplar aus. Allerdings glaube ich
nicht, dass Sie schon damit rechnen können, Ende dieses Se\bs
mesters sich zu habilitieren. Ich muss Ihre Arbeit immer\bs
hin ganz sorgfältig durchsehen, um für die Richtigkeit des In\bs
halts garantieren zu können, Sie wissen ja selber, wie leicht
einem bei einer Arbeit gelegentlich ein Fehler unterläuft.\comment{This is an allusion to Lorenzen's nonconclusive proofs, see footnote~\ref{nonconclusive} on page~\pageref{nonconclusive}.}
Ausserdem weiss ich ja nicht, wie stark mich eventuell die
Lazarettbehandlung in anspruch nehmen wird. -- Also rich\bs
ten Sie sich am besten gleich für eine Habilitation Anfang des\bruch{}
SS ein. Die Themen fürs Kolloquium, die Sie
an Herrn Besselhagen\footnote{Erich Bessel-Hagen (1898--1946) is professor of mathematics in Bonn.} geschickt haben, waren
ganz passend, zum mindesten angesichts der
Tatsache, dass Sie zur Zeit alle Ihre Wissen\bs
schaft nebenher arbeiten müssen. Ich würde am
liebsten über die Gentzenschen Sachen Sie sprechen
hören. --~Übrigens, wüssten Sie einen Herrn ausser\bs
halb Bonn, der Ihre wissenschaftlichen Leistungen begut\bs\bruch
achten könnte? Wir holen hier grundsätzlich bei jeder Habilitation solche Gutachter
von auswärts ein. Mit den besten
nachträglichen Neujahrswünschen Ihr\medskip

\hfill Wolfg.\ Krull.


%% file: 19440206-PL-1-1-145-Krull-Lorenzen

\noindent\begin{tabular}[t]{@{}l@{}}%
\label{19440206}%
           Krull\\%
           Bonn\\%
           Kaiser-Friedrichstr.~18%
         \end{tabular}\hfill%
         \begin{tabular}[t]{@{}l@{}}%
           Sonderführer Dr.\\%
           Paul Lorenzen\\%
           23 Wesermünde\\%
           Hafenstr.~92%
         \end{tabular}\medskip%

\hfill Bonn, 6.2.44.\nopagebreak\medskip

\indent\indent Lieber Herr Lorenzen!\nopagebreak\medskip

\noindent Leider habe ich die Reste meiner früheren Sonderabzüge
irgendwie verkramt. Sobald ich sie aber wiederfinde, be\bs
kommen Sie das gewünschte Separat. -- Nun aber zu Ihrer
Habilitationsschrift.\comment{On the same day, Krull answers a letter from Hasse dated 23~January 1944 \citep[§~1.89]{roquette04} that must contain a harsh criticism of Lorenzen and a refusal to write a report on Lorenzen's habilitation. Krull takes Lorenzen's defense, ``at least scientifically'', but concludes that ``he is still so self-confident that a repeated lesson would do no harm at all''.} Von Greifswald habe ich bisher nur
die Anmerkungen, nicht aber das Manuskript selber be\bs
kommen. Hier auf dem Dekanat bekam ich auf mehr\bs
fache Anfrage immer nur die Antwort, es sei bisher von
Ihnen nichts eingegangen. Das beunruhigt mich doch leb\bs
haft. -- Im übrigen, können Sie mir ausser Hasse und
Scholz keinen anderen Referenten angeben. Hasse liegt Ihrer
eigentlichen Richtung ("`Verbände"') doch ziemlich fern,
und Scholz ist halbwegs Philosoph. Dabei sind
die auswärtigen Gutachter für mich, vor allem nach
dem Standpunkt, den ich als Dekan immer eingenommen\bruch{}
habe, \uline{sehr} wichtig. Ferner bitte ich Sie
um eine Zusammenstellung Ihrer bisherigen Veröffentlichungen.\strut\medskip

\hfil\begin{tabular}{@{}ll@{}}%
  Mit den besten Grüßen&\\%
&Ihr Wolfgang Krull.%
\end{tabular}


%% file: 19440219-PL-1-1-143-Krull-Lorenzen
\label{19440219}%

\hfill Bonn, 19.2.44\nopagebreak\medskip

\noindent\hfil Lieber Herr Lorenzen!\nopagebreak\medskip

\noindent Leider\comment{On the same day, Krull writes an answer to a postcard by Hasse \citep[§~1.90]{roquette04}, with whom he agrees on how to answer to Lorenzen.} muss ich mit meinem Brief Ihnen vermutlich eine ernsthafte
Enttäuschung bereiten. Erstens: Ihre Habilitationsschrift nebst Bewerbung
ist tatsächlich hier nicht eingetroffen. Wenn ich Ihnen früher in anderm
Sinne schrieb, so nur deshalb, weil ich mich ohne Rückfrage bei dem De\bs
kanat auf einen Brief verliess, den mir Besselhagen zur Verfügung ge\bs
stellt hatte, und in dem Sie von Ihrem Schritt beim Dekanat als von
einer Selbstverständlichkeit sprachen. Zweitens: Es ist in gewissem Sinne
ein Glück, dass Ihr Gesuch verlorengegangen ist, wenigstens wenn Sie, wie
ich hoffe, noch ein Exemplar Ihrer Arbeit in Händen haben. Prof.\bsp
Hasse, auf dessen Gutachten Sie ja selber besonderen Wert legten, hat
sich keineswegs so geäussert, wie Sie es wohl erwarteten. Er ist der Ansicht,
dass Sie bisher noch sehr wenig mathematische Leistungen von wirklichem Belang
aufzuweisen haben, und mit dieser Stellungnahme in der Hand ist es mir
unmöglich, Ihre Habilitation im Augenblick zu befürworten, zumal da ich
Hasses Bedenken beim Betrachten der Liste Ihrer bisherigen Veröffentlichungen
sehr gut verstehe. Auch für mich zählt von diesen außer der Dissertation ei\bs
gentlich nur die im Druck befindliche Arbeit über den Verfeinerungssatz, und\bruch{}
diese in gewissem Sinne nur halb, da mir die logistische Seite zu fern liegt.
Natürlich bestände noch die Möglichkeit, dass Sie mir als Habilitationsschrift
nicht nur eine solide, sondern eine ganz aussergewöhnliche Leistung vorlegen, der
gegenüber alle Bedenken verstummen. Aber das halte ich nach dem, was
ich bisher von Ihrer Arbeit weiss nicht für wahrscheinlich, und so kann
ich Ihnen im Augenblick nur den einen Rat geben: Betrachten Sie Ihr Ha\bs
bilitationsgesuch als nicht erfolgt, legen Sie zunächst mir Ihre Arbeit
vor, -- (hoffentlich ist das Manuskript, das Sie nach Greifswald sandten,
dort angekommen), -- und wenn mir die Arbeit doch allein für die Habili\bs
tation noch nicht auszureichen scheint, so arbeiten Sie eben weiter. Suchen
Sie vor allem Ihre Basis zu erweitern und auch unter den reinen Mathe\bs
matikern bekannt zu werden. Prof.\ Scholz kommt als Philosoph mit
mathematischen Interessen für mich nur als Mitbegutachter, nicht als
Hauptreferent inbetracht, das müssen Sie immer berücksichtigen.

Es tut mir leid, wenn ich Ihnen mit diesem Brief, wie schon einmal
früher, eine Enttäuschung bereiten muss, aber glauben Sie mir, es ist
in Ihrem eigenen Interesse.\strut\medskip

\hfil\begin{tabular}{@{}ll@{}}%
  Mit den besten Grüssen&\\%
&Ihr Wolfgang Krull.%
\end{tabular}


%% file: 19440313-PL-1-1-131-Lorenzen-Krull

\noindent Dr.\ Paul Lorenzen\hfill%
\label{19440313}%
\begin{tabular}[t]{@{}l@{}}%
Wesermünde-Lehe, den 13.3.44\\%
Hafenstr.~92%
\end{tabular}\medskip

\noindent Sehr geehrter Herr Professor\ajout{ Krull},\medskip

ich bitte zu entschuldigen, daß ich auf Ihren Brief erst
heute antworten kann. Durch die gegenwärtige starke dienst\bs
liche Inanspruchnahme sind die neuen Exemplare meiner Ar\bs
beit erst jetzt fertig geworden. Hoffentlich gelangt die
Arbeit diesmal aber nun auch wirklich und endlich in Ihre
Hände -- darf ich Sie wohl bitten, mir den Eingang be\bs
stätigen zu wollen?

Allerdings werde ich befürchten müssen, daß meine Arbeit
jetzt eine ungünstigere Aufnahme finden wird, als sie vor
einem Vierteljahr gefunden hätte, da Sie mir schreiben, daß
meine Habilitation unter normalen Umständen kaum möglich sein
wird.

Ihrer Meinung, die auch die Meinung von Herrn Prof.\ Hasse ist,
daß meine bisherigen Veröffentlichungen ``keine Leistungen von
wesentlichem mathematischen Belang'' sind, möchte ich durchaus
nicht widersprechen. Hierin bin ich völlig Ihrer Meinung, ein
anderes Urteil wäre wohl auch kaum möglich.

Mit meinem Anliegen, mich zu habilitieren, stütze ich
mich ja aber nicht auf diese bisherigen Veröffentlichungen,
sondern auf die Habilitationsarbeit.

Soweit ich weiß und aufgrund der "`Bescheinigung"', die
Sie mir vor zwei Jahren ausstellten, annehmen mußte, bestehen
keine Bestimmungen, die die Habilitation von vorausgegangenen
Veröffentlichungen abhängig machen. Es ist der Nachweis wissen\bs
schaftlicher Tätigkeit zu erbringen -- und da ist es nun ja
nicht meine Schuld, daß ich die letzten vier Jahre hierfür --\nbsp
ganz allein auf mich gestellt~-- nur die dienstfreien Abende
zur Verfügung hatte. (Daß ich diese voll ausgenutzt habe,
werde ich aber sagen dürfen.) Sämtliche Examensbestimmungen,
die für die Kriegszeit erlassen sind, gehen darauf hinaus,
daß die notwendig entstehenden Härten für den Examinanden nach
Möglichkeit auszugleichen sind. Mir scheint daher mein Wunsch,
mich jetzt -- 6~Jahre nach meiner Promotion -- zu habilitie\bs\bruch
ren, durchaus gerechtfertigt, da ich mich normalerweise doch
schon vor mindestens 3~Jahren hätte habilitieren können.

Da ich nicht weiß, wie lange der Krieg noch dauert, und
wie die Arbeitsbedingungen für mich werden, müßte ich, wenn
Sie meine Habilitation jetzt ablehnen, die beiliegende Arbeit
zurück erbitten --~und weiterarbeiten, wie Sie mir schrei\bs
ben. Allerdings muß ich dann noch damit rechnen, daß meine
Arbeit völlig unnütz ist, weil Sie sie nicht gelten lassen
werden. Im Anschluß an eine algebraische Untersuchung über 
orthokomplementäre Halbverbände versuche ich jetzt, den Zu\bs
sammenhang dieser Fragen mit der Widerspruchsfreiheit der
klassischen Logik herauszubekommen. Diesen Versuch möchte ich
machen selbst auf die Gefahr hin, daß das evtl.\ Ergebnis von
Ihnen gar nicht gezählt wird --~weil ich nicht umhin kann,
solche Fragen als Fragen von wesentlichem mathematischem Be\bs
lang zu empfinden (und in dieser Auffassung der Logistik 
darf ich mich sogar auf Hilbert berufen).

Ich möchte Sie daher bitten, mir sagen zu wollen, ob es
unmöglich ist, sich als Mathematiker zu habilitieren, wenn man
etwa entschlossen ist, die eigene Forschungsarbeit an mathe\bs
matisch-logistische Dinge zu wenden.

Diese Voraussetzung trifft allerdings auf mich noch nicht 
einmal zu, da ich selber eigentlich viel mehr an der alge\bs
braischen Seite der Beweistheorie interessiert bin als an der
rein logischen. Genau dasselbe Interesse würde ich auch daran
haben, die Fragestellungen der algebraischen Geometrie be\bs
grifflich zu durchdringen, in dem Sinne, wie es etwa die bei\bs
liegende Arbeit für die multiplikative Idealtheorie versucht.
Aber auch auf diesem Gebiete fürchte ich, daß meine Auffassung
von der Ihrigen abweicht. Z.~B.~erscheint mir die Erkenntnis,
daß ein Idealsystem eigentlich nichts anderes als ein Ober\bs
halbverband und eine Bewertung nichts anderes als eine Ordnung
ist, als das wesentlichste Ergebnis meiner Bemühung. In die\bs
sem Sinne läßt sich sogar der Inhalt von §~6 als auf rein halb\bs
ordnungstheoretischen Tatsachen beruhend erkennen, worauf ich\bruch{}
aber in der Arbeit nicht näher eingegangen bin.

Es bleibt mir auch an dieser Stelle nur übrig, zu fra\bs
gen, ob ich mit einer solchen Auffassung habilitationsfähig
bin oder nicht.

Wenn das nicht der Fall sein sollte, so sehe ich eigent\bs
lich nicht, was ich anderes machen sollte, als mich nach einem
anderen Beruf umzusehen, denn es wird schwer sein, meine
Auffassung über den Sinn der mathematischen Forschung zu än\bs
dern.

Es wird jedoch ebenfalls sehr schwer für mich sein, mein
Habilitationsgesuch als nicht geschehen zu betrachten, nach\bs
dem ich meine Absicht, mich zu habilitieren, ja nicht nur
Ihnen gegenüber geäußert habe.

Ich möchte nichts unternehmen, was Ihrem Ratschlag wider\bs
spricht, muß aber doch um Verständnis bitten, daß ich mich
nicht einfach damit abfinden kann, die Habilitation auf spä\bs
ter zu verschieben, da ich entweder in den nächsten Jahren
gar nicht werde arbeiten können, oder aber damit rechnen muß,
daß meine Arbeit nicht gezählt wird.\medskip

\hfil Mit den ergebensten Grüßen verbleibe ich

\hfil Ihr\hfil

\hfil\hfil\hfil\hfil\hfil Paul Lorenzen\medskip

Für die Übersendung der beiden Separata bin ich Ihnen sehr zu Dank verbunden.


%% file: 19440401-PL-1-1-144-Krull-Lorenzen

\noindent\begin{tabular}[t]{@{}l@{}}%
\label{19440401}%
           Reg.Rat Prof.\ Krull\\%
           Greifswald\\%
           Mar[ine]obs.%
         \end{tabular}\hfill%
         \begin{tabular}[t]{@{}l@{}}%
           Herrn Dr.\\%
           Paul Lorenzen\\%
           23 Wesermünde-Lehe\\%
           Hafenstr.~92\\%
         \end{tabular}\medskip

\hfill Greifswald, 1.4.44\nopagebreak\medskip

\noindent\hfil Lieber Herr Lorenzen!\hfil\medskip

Auf Ihren ausführlichen Brief vom 13.3.\ möchte ich Ihnen
heute noch nicht im Einzelnen antworten. Zunächst nur das eine,
dass ich Ihre Arbeit ohne jedes Vorurteil, zum mindesten ohne
jedes ungünstige, zu lesen angefangen habe. Indessen bin ich
schon auf S.~8 auf eine Schwierigkeit gestossen, die ich Sie bitten muss,
mir aufzuklären, da ich durch andere Untersuchungen (Korrelations\bs
theorie) zu stark inanspruch genommen bin. Sie schreiben bei
Satz~4:\comment{See \citealt[Satz~4]{Lor1950}.} "`Aus $c\equiv a\cdot b^{-1}$ und $a\parallel b$ folgt $a\equiv b\cdot c$"'. D.\ h.\bsp
aber anders ausgedrückt: "`Aus $a\land b\equiv1$ folgt $a\equiv bab^{-1}$"',
und diese Tatsache, die wenn richtig, unbedingt eine Formulie\bs
rung als Satz verdient hätte, (weil hier ein sehr starkes kom\bs
mutatives Element ins Nichtkommutative hineinkommt), konnte
ich jedenfalls aus dem Handgelenk heraus nicht beweisen. An\bs
dernfalls scheint der weitere Wortlaut Ihres Textes zu zeigen, dass
keineswegs nur ein Schreibfehler vorliegt, dass Sie nicht etwa
$c\cdot b$ meinten und $b\cdot c$ schrieben. --~Also bitte, klären Sie
mir diesen Punkt auf!\medskip

\hfil Mit den besten Grüßen\hfil Ihr Wolfgang Krull.

%% file: 19440416-PL-1-1-142-Krull-Lorenzen

\label{19440416}%
\hfill Greifswald, 16.4.44\nopagebreak\medskip

\noindent\hfil Lieber Herr Lorenzen!\hfil\medskip

\noindent Vielen Dank für Ihre Karte! Dieser Punkt wäre also geklärt. Also die Tatsache, dass aus $a\land b\equiv%
1$ stets $a\cdot b\equiv b\cdot a$ folgt, sollte unbedingt als Satz formuliert werden. Nicht nur, weil erst dadurch
der Beweis von Satz~4 ganz in Ordnung kommt. Es handelt sich auch darum, zu zeigen, wie stark die
Verbandsgruppen "`kommutativ infiziert sind"'. -- Und hier wären wir nun an einem
Punkt, über den ich etwas ausführlicher werden muss. Ihre Arbeit hat mich insofern etwas ent\bs
täuscht, als sie doch anscheinend fürs Nichtkommutative wesentlich weniger bringt, als ich mir
gedacht hatte. Im Verbandsbegriff steht eben ein so starkes kommutatives Moment, dass bei der
Theorie Ihrer Verbandsgruppen die Beweise für kommutativ und nichtkommutativ völlig gleich
laufen. Andererseits bedeutet die Beschränkung auf Verbandsgruppen vom Nichtkommutativen aus
gesehen offenbar eine sehr starke Beeinträchtigung der Allgemeinheit. -- Unter diesen Umstän\bs
den hätte ich aber Ihre Arbeit auch ohne die Stellungnahme von Hasse nicht als Habilitations\bs
schrift empfehlen können. Dazu bewegt sie sich doch gar zu sehr im Gedankenkreise Ihrer
Dissertation. Aber natürlich bedeutet  dieses Urteil keineswegs, dass ich Ihre Arbeit "`nicht als
mathematische Leistung zähle"'; und wenn Sie mich mit weiteren Veröffentlichungen belehren
würden, dass Ihre Methoden doch auch fürs Nichtkommutative grössere Bedeutung haben,
als es mir bis jetzt scheint, so sollte mich das nur freuen. Das wäre eine Arbeits\bs
richtung, in der ich Sie lieber sehen würde als in der nach der Logistik hin. Ich halte
die Spezialisierung nach der logistischen Seite hin nun einmal für nicht unbedenklich, so
fern es mir liegt, bestreiten zu wollen, dass auch auf diesem Gebiete Bedeutendes gelei\bs
stet werden kann (man denke nur etwa an Gödel oder Gentzen). -- Also wenn Sie glau\bs\bruch
ben, wirklich wesentlich neue Gedanken zur Beweistheorie zu haben, so lassen Sie sich ja nicht
abhalten, sie auszuarbeiten. -- Was die Tatsache angeht, dass Sie schon an verschiedenen Stel\bs
len über Ihre baldige Habilitation geredet habe, so war das eine Unvorsichtigkeit
von Ihnen, die mir eigentlich unverständlich ist. Über solche Dinge redet man doch erst
dann, wenn man seiner Sache absolut sicher ist. Und wo es möglich ist, wie in Ihrem Falle,
legt man doch zu allererst die zukünftige Arbeit ganz persönlich dem voraussichtlichen
Referenten vor, und macht die weiteren offiziellen Schritte erst im Einverständnis mit
ihm! Wenn Sie so vorgegangen wären, wären Ihnen alle Peinlichkeiten erspart
geblieben. Aber natürlich helfe ich Ihnen auch jetzt gerne so viel ich kann. Ihrer Dienst\bs
stelle gegenüber wird es Ihnen ja nicht schwer fallen, irgend einen formalen Grund für
die Hinausschiebung Ihrer Habilitation zu erfinden. Unserm Dekan gegenüber brau\bs
chen Sie einfach von der Sache nicht mehr zu reden oder zu schreiben; dann wird er
nicht mehr daran denken. Und Besselhagen und Peschl kann ich ja sagen, ich hätte ge\bs
wünscht, dass Sie Ihre Untersuchungen über das bisher Vorliegende, zunächst zu Ver\bs
öffentlichende hinaus und weiter aufs Nichtkommutative ausdehnten. Sie sehen, ich
helfe Ihnen gerne, wo ich kann. -- Was übrigens die Form der Veröffentlichung der
bisherigen Ergebnisse angeht, so müssen wir uns darüber ein anderes Mal unterhalten.
Ich glaube, Sie müssen die Sache im Sinne der Entlastung von allen irgendwie ent\bs
behrlichen Hilfsbegriffen sehr stark kürzen, wenn Sie wollen, dass ein gewisser Kreis
von Fachgenossen sich die Arbeit wirklich genauer ansieht. Aber wie gesagt, davon das
nächste Mal!\medskip

\hfil Mit den besten Grüßen\hfil Heil Hitler!\hfil

\hfill Ihr Wolfgang Krull.


%% file: 19440425-PL-1-1-132-Lorenzen-Krull

\label{19440425}\hfill%
\begin{tabular}{@{}l@{}}%
Wesermünde, den 25.4.44\\%
Hafenstr.~92%
\end{tabular}\hspace{3em}\medskip

\noindent Sehr geehrter Herr Professor,\medskip

\noindent Ihren Brief habe ich vor einigen Tagen erhalten. Ich danke
Ihnen für Ihre Freundlichkeit, mir in den Schwierigkeiten,
die durch Ihr ablehnendes Urteil entstehen, helfen zu
wollen. Allerdings bitte ich zuvor, mir erlauben zu wollen,
mich gegen das Urteil, daß sich meine Behandlung der nicht\bs
kommutativen Gruppen ``zu sehr in dem Gedankenkreis meiner
Dissertation'' befinde, zu verteidigen.

Die Sätze der \S\S~1-3\comment{See \citealt[§§~1--3]{Lor1950}.} sind allerdings tatsächlich im
Anschluß an meine Dissertation entstanden. Ich habe damals
die nichtkommutative Theorie jedoch liegen lassen, weil mir
diese Übertragung vom Kommutativen aufs Nichtkommutative
nicht interessant schien (mit Ausnahme der Regularitätsbe\bs
dingung) -- und vor allem, weil es aussichtslos war,
die wesentlichen Bestandteile der kommutativen Theorie:
die Konstruktionen des $a$-\hspace{0pt}Idealsystems und der Gruppe der
Idealbrüche auch im Nichtkommutativen durchzuführen.

Die Prüfersche Definition des $a$-Systems liefert im Nicht\bs
kommutativen nämlich kein Idealsystem. Andererseits kann man
im Nichtkommutativen aus einer Halbgruppe mit $ac=bc\Rightarrow a=b$
keine Quotientengruppe konstruieren.

Erst nachdem mir -- vor etwa 4~Jahren -- endlich klar wur\bs
de, daß ein Idealsystem nichts anderes als ein Halbverband
ist, ergab sich plötzlich die Möglichkeit, die Konstruk\bs
tion der Gruppe der Idealbrüche durch etwas ganz Neues zu
ersetzen: nämlich durch eine iterierte Anwendung des Ideal\bs
begriffs. Jede minimale Oberverbandsgruppe einer halbgeord\bs
neten Gruppe~$G$ ist das $v$-\hspace{0pt}Idealsystem eines Idealsystems von~$G$.
Das ist der Inhalt von~\S~4.\comment{See \citealt[§~4]{Lor1950}.}

In der Hoffnung, auf dieser Grundlage die Theorie, sowie
einige Ansätze über geordnete Gruppen, während eines Arbeits\bs
urlaubs ausbauen zu können, habe ich Sie damals um Zulassung
zur Habilitation gebeten.\bruch

Ich war mir bewußt, daß es erwünscht gewesen wäre, wenn ich
ein Thema beantwortet hätte, dessen Fragestellung nicht in
die Richtung der Dissertation fiel -- hätte ich auch nur
einmal wenigstens einige Literatur zur Einarbeitung in ein
neues Problem zur Verfügung gehabt, so hätte ich das bestimmt
vorgezogen.

Unter den gegebenen Umständen schien mir jedoch meine Bitte um
Zulassung zur Habilitation nicht unangemessen zu sein.

Nachdem damals dieser Arbeitsurlaub sich wegen meiner Versetzung nach
hier nicht verwirklichen ließ, hat sich nun inzwischen herausgestellt,
daß der kommutative Aufbau:\smallskip%

{\leftskip2.25\parindent\parindent0pt%
\leavevmode\llap{1) }Definition des $a$-\hspace{0pt}Idealsystems nach Prüfer%

\leavevmode\llap{2) }Konstruktion der Idealbrüche (bzw.\ der Funktionale)%

\leavevmode\llap{3) }Identitätsbeweis des $a$- und $b$-Systems mit Hilfe der
  Idealbrüche bzw.\ Funktionale\smallskip%

  }\noindent im Nichtkommutativen wieder durch eine ganz andere Methode er\bs
setzt werden muß. Zur Erläuterung dieser neuen Methode bitte
ich Sie, die beiliegende Bemerkung durchblättern zu wollen.
Diese Bemerkung ist vorläufig nicht zur Veröffentlichung be\bs
stimmt, sie soll nur versuchen, Ihnen darzulegen, daß diese
Methode, deren Grundgedanke in §~6\comment{See \citealt[§~6]{Lor1950}.} meiner Arbeit enthalten ist,
durchaus verschieden von den Methoden meiner Dissertation\insere{ ist}.
Denn ich werde sagen dürfen, daß die Erkenntnis, daß solche
Sätze wie die bewertungstheoretischen Fundamentalsätze mit
rein halbordnungstheoretischen Mitteln zu erhalten sind, durch\bs
aus nicht im Gedankenkreis meiner Dissertation zu finden ist.

Als dritte "`neue Methode"' hat meine Arbeit den Satz~3 der
beiliegenden Blätter aufzuweisen, dessen Anwendung auf halb\bs
geordnete Gruppen der Teil~2 meiner Arbeit durchführt.

Aufgrund dieser Ergebnisse --~und gerade weil die ent\bs
scheidenden Methoden von meiner Dissertation völlig abweichen
trotz der gleichen Fragestellung~-- habe ich vor gut einem
halben Jahre Ihnen mitgeteilt, daß ich glaubte, die Arbeit nun
abschließen zu können. Wenn ich die Angelegenheit dann etwas
überstürzte, so bitte ich das entschuldigen zu wollen, weil
mir ausgerechnet mein militärischer Kursus dazwischen kam.

Mit den ergebensten Grüßen und nochmaligem Dank für die
Mühe, die Sie sich in meiner Sache geben, verbleibe ich

\hfill[Ihr Paul Lorenzen]\comment{\label{19440505}\selectlanguage{english}{In a letter to Scholz dated 5 May 1944
  , Lorenzen gives an account of this writing, }``\selectlanguage{german}die darlegt, daß gewisse Fundamentalsätze der Krullschen Bewertungstheorie nicht nur für halbgeordnete Gruppen gelten, sondern sogar für beliebige halbgeordnete Mengen.
  
    ``Allerdings zweifle ich daran, ob gerade dieses Ergeb\bs
    nis, das der neuen Methode meiner Arbeit zu danken ist,
    sich dazu eignet, die Meinung von Herrn Prof.\ Krull zu meinen
    Gunsten zu wenden.

    ``Ich versuche, diese Fatalität in stoischem Sinne zu
    ertragen, es ist mir dabei aber ein schöner Trost, Ihnen
    darüber schreiben zu dürfen.''

    \selectlanguage{english}\label{19440526}On 26 May 1944, he addresses the following letter to Scholz.
    \selectlanguage{german}\input 19440526-PL-1-1-134-Lorenzen-Scholz}\par%
 

%% file: 19440526-PL-1-1-134-Lorenzen-Scholz

\hfill Wesermünde, 26.5.44\hspace*{2em}

Sehr geehrter Herr Professor,\smallskip

\noindent Dieses darf ich Ihnen zunächst sagen, daß ich immer davon
überzeugt bin, daß "`man in Münster"' wirklich an mich denkt
-- und nicht nur denkt.

\noindent Trotzdem hätte ich es gern vermieden, Ihnen mit meiner
Habilitationsangelegenheit explizit lästig zu fallen, weil
ich annehmen muß, daß Sie genug anderes zu tun haben. Auf
Ihren Brief vom 22/5, kann ich jetzt aber Ihre Hilfe nicht
mehr ausschlagen -- und darf Ihnen also versichern, wie sehr
mir eine solche Hilfe gelegen kommt.

\noindent Herr Prof.\ Krull hat mir noch nicht geantwortet, vielleicht
ist mein letzter Brief gar nicht angekommen.

Herrn Prof.\ Köthe habe ich schon in Würzburg versucht, meine
Theorie vorzutragen -- es war aber kaum ausreichend Zeit.

\noindent Ich lege noch eine "`Bemerkung"' bei, die eine Fortführung
des \S~6 meiner Arbeit ist, ohne jedoch diese vorauszusetzen.
Da diese Bemerkung rein halbordnungstheoretisch ist, wäre
ich Ihnen auch für Ihr Urteil sehr dankbar.

Würden Sie Herrn Prof.\ Köthe wohl mitteilen, daß mein Vortrag
in Würzburg den Inhalt von~\S~3 zum Gegenstand hatte? Die
Neuerungen gegenüber der kommutativen Theorie liegen in~\S~4
und~\S~6.

Ich werde für \uline{jedes} Urteil dankbar sein, denn wie sollte
ich sonst lernen, was man in dieser Welt als "`gut"' bezeichnet.\smallskip

\hfil Mit den ergebensten Grüßen verbleibe ich

\hfil Ihr [Lorenzen]\smallskip

Die fehlenden Anmerkungen werde ich Ihnen in den nächsten
Tagen nachsenden. %

%% file: 19440529-PL-1-1-146-Krull-Lorenzen

\label{19440529}%
\hfill\begin{tabular}[t]{@{}l@{}}%
        Sonderführer Dr.\\%
        Paul Lorenzen\\%
        23 Wesermünde\\%
        Hafenstr.~92%
      \end{tabular}\medskip%

\noindent Lieber Herr Lorenzen!\hfill Greifswald, 29.5.44\medskip

\noindent Endlich komme ich dazu, Ihren Brief vom 25.4.\ zu beantworten. Ich
wollte mir zunächst Ihr beigelegtes Manuskript ansehen, und dazu
fand ich erst jetzt während der Pfingsttage Zeit. Am besten hat mir
Ihr Brief selber gefallen, da sagen Sie am klarsten, worauf es ankommt.
An Ihrem Manuskript ist wieder das störende, dass nicht zu sehen ist,
was für einen Vorteil man aus der Verallgemeinerung des Fundamental\bs
satzes der Bewertungstheorie auf bel.\ Halbordnungen und aus Ihrer Umformung
des Kriteriums für "`ganz abgeschlossen"' gewinnt. So bleibt immer
das Gefühl, ob nicht der umständliche Weg in keinem rechten Ver\bs
hältnis zum Endergebnis steht, und das ist das gleiche bei diesem
nicht für die Veröffentlichung bestimmten Manuskript ebenso wie bei
Ihrer geplanten Habilitationsschrift. Ich halte es nun durchaus für
möglich, dass das mehr ein Mangel der Darstellung ist, Sie wissen
irgendwie aus Ihren Ideen nicht das zu machen, was man aus
ihnen herausholen könnte. Aber das ist leider ein Punkt, den man
schriftlich kaum richtig klären kann, wir müssten uns ein\bs
mal sehr gründlich über alle Einzelheiten aus\uline{sprechen}, und da\bs
zu ist gerade augenblicklich leider keine Möglichkeit gegeben. Ich
bedauere das vor allem deshalb, weil Sie nach meiner Ansicht den
Inhalt Ihrer ursprünglich als Habilitationsschrift gedachten\bruch{}
Untersuchungen veröffentlichen sollten, weil ich aber
befürchte, dass in der vorliegenden Form die
Arbeit nicht geeignet ist, Sie, wie es doch sein
sollte, einem gewissen Kreis von Fachgenossen
wirklich bekannt zu machen. Es wäre da
zunächst eine sehr gründliche Umarbeitung
nötig (starke Straffung, ev.\ Zweiteilung, an
anderer Stelle wieder Erweiterung), und dabei
würde ich Sie gerne beraten. Aber wir werden
uns eben damit abfinden müssen, dass im
Augenblick eine derartige Zusammenarbeit, die
mündliche Besprechungen erfordert, nicht
möglich ist. Seien Sie aber überzeugt, dass ich
Ihnen wirklich gerne helfen möchte.\medskip

\hfil Mit den besten Grüßen\hfil Heil Hitler!

\hfill Ihr Wolfgang Krull.\comment{\foreignlanguage{english}{In a letter dated 2 June 1944\label{19440602}, Lorenzen gives an account of this letter to Scholz and writes:}
  ``Da er von meiner Habilitation aber nichts mehr schreibt,
  will ich es noch ein letztes Mal versuchen, um seine Zu\bs
  stimmung zu bitten. Aber wird es was nützen? Darf ich Sie
  trotzdem bitten, Herrn Prof.\ Köthe das `Hasse-Exemplar' ein\bs
  schließlich der Bemerkung mit den beiliegenden Anmerkungen
  zu schicken. Das neue Exemplar stelle ich ganz zu Ihrer Ver\bs
  fügung.''

  \foreignlanguage{english}{Then, in an undated letter\label{194406}, he writes:}
  ``Für Ihre so schnelle Vermittlung zu Herrn Prof.\
  Köthe bin ich Ihnen sehr dankbar, --~-- ebenso auch für
  Ihr Gedenken bei unserm Angriff, bei dem unsere
  Wohnung immerhin so beschädigt wurde (wenn sie
  auch noch integrierbar sein wird), daß zunächst
  meine Frau mit Jutta in ein Dorf der Sächsischen
  Schweiz abgereist ist, und ich ein `möbliertes
  Zimmer' erworben habe.

  ``Ich werde Herrn Prof.\ Köthe Genaueres schreiben,
  werde ihm aber darin recht geben müssen, daß es
  leicht möglich ist, daß Herr Prof.\ Krull nicht erbaut
  sein wird --~-- sondern eine Einmischung jedem Beteiligten
  sehr übelnehmen wird. Daher werde ich Herrn Prof.\
  Köthe dankbar sein, wenn er Ihnen oder mir sein
  Urteil mitteilen wird.''}\par
 

%% file: 19440606-PL-1-1-133-Lorenzen-Krull

Dr.\ Paul Lorenzen\hfil\hfil(23)\quad%
\label{19440606}%
\begin{tabular}[t]{@{}l@{}}%
Wesermünde, den 6.6.44\\%
Hafenstr.~92%
\end{tabular}\medskip

\noindent Sehr geehrter Herr Professor,\medskip

\noindent für Ihre Karte von 29.5.\ danke ich Ihnen sehr, obwohl sie
mir eine sehr ernste Enttäuschung bereitet, indem als we\bs
sentlicher Mangel meiner Arbeit jetzt die Umständlichkeit
des Verfahrens bezeichnet wird.

Dabei ist gerade die Vereinfachung und Klärung der Be\bs
weismethoden das eigentliche Hauptziel meiner Arbeit. Ich
habe nicht versucht, die Sätze der multiplikativen Ideal\bs
theorie um jeden Preis zu verallgemeinern, auch um den Preis
einer Komplizierung -- sondern mir liegt im Gegenteil nur
daran, die Grundgedanken der Beweismethoden in ihrer letzten
Einfachheit zu erkennen. Wenn ich z.\ B.\ den Bewertungsbegriff
ersetze durch einen ``Homomorphismus eines Halbverbandes in
eine geordnete Menge'', so sehe ich darin nämlich eine be\bs
griffliche Vereinfachung, und nicht etwa eine Komplizierung.
Denn die Einführung des Bewertungsbegriffs (z.\ B.\ die zunächst
willkürliche Dreiecksungleichung) rechtfertigt sich nur durch
den späteren Erfolg, der Homomorphiebegriff trägt dagegen
seine Berechtigung in sich selbst. Ich würde sagen, daß der
Homomorphismus in eine Ordnung der "`reine Begriff"' ist, der
dem Bewertungsbegriff zugrunde liegt. Und der zugrundeliegende,
reine Begriff scheint mir unbestreitbar der einfachere zu
sein.

Wenn ich mich in diesem Punkte irren sollte, so bitte ich
Sie aufs dringendste darum, es mir sagen zu wollen, weil es
nämlich bei meiner ganzen mathematischen Arbeit bisher immer
mein Bestreben war, diese zugrundeliegenden, reinen Begriffe
selbst in ihrer einfachen und durchsichtigen Klarheit ans
Licht zu bringen.

In dieser Tendenz der begrifflichen Klärung unterscheidet
sich die Arbeit ebenfalls grundsätzlich von meiner Dissertation.
Daß in dieser Klärung, in dem Verständnis der inneren Bedeu\bs
tung, wie Sie es einmal nennen, die vordringlichste Aufgabe\bruch{}
liegt, das ist mir nämlich erst in den letzten Jahren wirk\bs
lich bewußt geworden.

Der Wert einer solchen begrifflichen Erkenntnis liegt
m.\ E.\ vor allem in sich selbst. Nur in zweiter Linie kommt die
Vereinfachung in Frage, die sich ergibt, wenn man die reinen
Begriffe auf den ursprünglichen Spezialfall anwendet. (Z.\ B.\bsp
wird der Beweis des bewertungstheoretischen Fundamentalsatzes
fast trivial: Gilt für kein $z_1,\dots,z_n$ \ $a\in\su a\su I[z_1^{\pm1},\dots,z_n^{\pm1}]$,
so gibt es einen maximalen Oberring~$\su B$ mit dieser Eigenschaft.
Aus $\su B\subset\su B[z^{\pm1}]$ folgt%
\[%
a\in\su a\su B[z^{+1},x_1^{\pm1},\dots,x_r^{\pm1}]\text{ und }%
a\in\su a\su B[z^{-1},y_1^{\pm1},\dots,y_s^{\pm1}]%
\]
also $a\in\su a\su B[z^{\pm1},x_1^{\pm1},\dots,y_s^{\pm1}]$. Widerspruch!)\comment{See \citealt[p.~489]{Lor1950}, reproduced on page~\pageref{1950p489} and translated on page~\pageref{1950p489t}.}

Ebenso ist es eigentlich nur ein Nebenergebnis, daß die
Sätze jetzt auch im Nichtkommutativen gelten. Im Vordergrund
steht stets die Erkenntnis der reinen Begriffe.

Ich bin natürlich weit davon entfernt, zu behaupten, daß
diese Auffassung der Mathematik die richtige sei, aber daß
sie eine berechtigte Auffassung ist, werde ich behaupten dürfen.

Wenn allerdings an dieser Auffassung meine Habilitation
scheitern sollte, so will ich mich bemühen, mir eine gegen\bs
teilige Auffassung zu eigen zu machen, soweit das möglich ist.

Ich darf mich mit diesem Anliegen meiner Habilitation noch
einmal an Sie wenden, da Ihnen ja nichts daran gelegen sein
kann, mich einfach vor ein unabänderliches und auswegloses
Nein zu stellen.

Die Gründe, die Sie meiner Habilitation entgegenhalten,
daß ich nichts Wesentliches bisher veröffentlicht habe, daß
meine Arbeit in Richtung der Dissertation liegt und daß sie
ungerechtfertigt umständlich ist, kann ich --~mit Ausnahme
des letzten~-- nicht leugnen. Aber ich darf Sie bitten, auch
die Bedingungen zu berücksichtigen, unter denen ich stehe:
daß es mir nicht möglich war, ein neues Arbeitsgebiet ohne Lite\bs
ratur und Anregung (und ohne ausreichend Zeit) wirklich zu
erschließen. Wenn ich die Gewißheit hätte, den Krieg abwarten
zu können, so brauchte mir jetzt nicht so viel an meiner Habi\bs
litation gelegen zu sein. Ich würde die zuversichtliche Hoff\bs
nung haben, es nach dem Kriege zu schaffen --~aber wer weiß,
ob und wann wir normale Nachkriegszeiten erleben werden. Die
Berufsausbildung, das Erreichen eines Abschlußes wird ja darum\bruch{}
bei allen übrigen, wenn irgend möglich so gefördert. Es ist
mir ja auch nicht nur praktisch unmöglich, einen anderen
Beruf zu ergreifen, ich sehe vor allem aufgrund meiner Ver\bs
anlagung keine Möglichkeit dazu, da eigentlich alle meine
Gedanken und Bestrebungen sich ausschließlich auf die mathe\bs
matische Erkenntnis richten.

Als eine besondere Härte muß ich die Nichtabge\bs
schlossenheit meiner äußeren Berufsausbildung deshalb empfin\bs
den, weil mir dadurch auch in meinem gegenwärtigen Dienst
mehrere Möglichkeiten abgeschnitten sind, die an die Bedingung
der Habilitation geknüpft sind.

Was den zweiten Hinderungsgrund anbetrifft, daß meine
Arbeit der Dissertation gegenüber nicht neu ist, so darf ich
hier noch einmal wiederholen, daß die entscheidenden Methoden
der Dissertation ($a$-Ideale für Halbgruppen, Idealbrüche,
$t$-Ideale) weder explizit noch implizit in meiner Arbeit eine
Rolle spielen.

Ich bitte Sie darum, dieses Ihren Gründen gegenüber halten
zu wollen, und verbleibe mit den ergebensten Grüßen\medskip

\hfill Ihr [Paul Lorenzen]


%% file: 19440622-PL-1-1-141-Krull-Lorenzen

\label{19440622}%
\hfill Greifswald, 22.6.44.\qquad

\noindent\hfil Lieber Herr Lorenzen!\hfil\medskip

\noindent Es fällt mir nicht ganz leicht, Ihren letzten Brief zu beantworten, denn ich muss Ihnen
voraussichtlich noch einmal eine Enttäuschung bereiten. Zuerst: Seien Sie überzeugt, dass
ich recht habe, wenn ich Ihre grosse, 53~Seiten lange Arbeit in ihrer derzeitigen Form für
gründlich verfehlt erkläre. Ich glaube Ihnen gerne, dass wirklich neue Gesichtspunkte
gegenüber Ihrer Dissertation drin stecken und ich erkenne auch Ihr Streben nach
äusserster begrifflicher Klarheit grundsätzlich durchaus an. Aber gerade eine
solche Klarheit, die sich auch in einer durchsichtigen Form der Darstellung äussern müsste,
vermisse ich in Ihrer Arbeit durchaus. Und es langt nicht, dass Sie sich selbst über
etwas klar sind, Sie wollen und müssen es auch den andern ebenso klar machen,
\uline{darauf} kommt es an. Schon das letzte Mal schrieb ich Ihnen, im Einzelnen müsse
man sich über diese Dinge mündlich aussprechen. Wenn ich aber Ihnen schriftlich einen
Rat geben soll, so kann es nur etwa der sein: Versuchen Sie doch einmal den Inhalt
Ihrer Arbeit auf 20-25~Seiten zusammen zudrängen, aber so, dass Sie an den
wirklich wesentlichen Stellen sich nicht scheuen deutlich zu sagen, was Sie wollen und
nicht etwa hinter ein paar formalen Rechnungen die Gedanken verstecken. Ich denke,
das müsste gehen, wenn es Ihnen auch zunächst unmöglich erscheint, und so
kämen Sie dann vielleicht zu einer Arbeit, die wirklich geeignet wäre, Ihren
Namen bekannter zu machen. Das wäre die nächste Aufgabe, die ich Ihnen
stellen möchte, denn ich will natürlich nicht, dass Ihnen das, was Sie sich da
alles überlegt haben, verloren geht. Und dann später ein neues, wenn mög\bs\bruch
lich auf speziellere Anwendungen (etwa das Nichtkommutative) zugeschnittenes Problem
angepackt! Denn mit dieser Arbeit allein lasse ich Sie noch nicht zur Habilitation
zu, dabei bleibe ich, -- wie ich überzeugt bin, in Ihrem eigenen Interesse. Die Habili\bs
tation ist kein Abschluss einer Laufbahn, eine "`abgeschlossene Hochschulausbildung"',
wie sie z.\ B.\ auch für die höhere Wehrmachtsbeamtenlaufbahn gefordert wird, haben
Sie als Promovierter schon längst. Und --~entweder geht der Krieg einigermassen an\bs
ständig aus, wie wir alle hoffen, dann werden Sie auch Ihre Arbeit für die
Habilitation als Assistent in Ruhe in Angriff nehmen können,~-- oder, ja ich glaube,
Sie sind sich nicht klar, dass es dann mit einer Hochschullaufbahn für Sie ohnehin alle
wäre. Ein à la Baisse-Spekulieren kann und darf es heute bei niemandem
geben. --~Schliesslich, um es nochmal zu sagen: Die Habilitation ist nicht ein Ab\bs
schluss, sondern ein neuer Anfang, gewissermassen als Gesellenstück für eine Laufbahn,
die zu betrachten letzten Endes nur dann Sinn tut, wenn man einigermassen Gewiss\bs
heit hat, dass man es irgendwann auch zum Meister bringt. Und so weit sind
Sie eben noch nicht, wobei Ihnen gerne zugegeben sei, dass Sie unter den Kriegsbedin\bs
gungen so gehemmt sind, dass niemand Ihnen deswegen einen Vorwurf machen dürfte.
Aber es wäre unverantwortlich von mir, Ihnen jetzt schon das Tor zu öffnen, wo ich
mir noch längst kein genügend sicheres Bild über Ihre zukünftige Weiterentwicklung
machen kann.\medskip

\hfil Mit den besten Grüssen und Heil Hitler!

\hfill Ihr Wolfgang Krull.\comment{\foreignlanguage{english}{In a letter to Scholz dated 28 June 1944\label{19440628}
  , Lorenzen writes:}
    ``Heute habe ich die endgültige Absage von Herrn
    Prof.\ Krull erhalten, da `meine Arbeit die begriffliche
    Klarheit durchaus vermissen läßt', und es `unver\bs
    antwortlich wäre, mir schon jetzt das Tor (zum Dozenten)
    zu öffnen'.

    ``Herrn Prof.\ Köthe habe ich geschrieben, er möchte
    zunächst mir sein Urteil mitteilen, da Herr Prof.\ 
    Krull in der Tat vermutlich nicht erbaut sein würde
    über eine direkte Einmischung.''

    \foreignlanguage{english}{On 9 July\label{19440709}
    , he writes:}
    ``ich darf Ihnen versichern, wie dankbar
    ich Ihnen dafür bin, daß Sie sich so um mich kümmern.
    Für Ihre ermutigenden Worte danke ich Ihnen beson\bs
    ders. Wenn ich auch nicht im geringsten das Gefühl
    habe, durch das Urteil von Herrn Prof.\ Krull `umge\bs
    worfen' zu sein, so sind mir Ihre Zeilen doch sehr
    wohltuend gewesen. Und gerade kam auch ein ebenso
    wohltuender Brief von Herrn Prof.\ Köthe, der sich
    damit einverstanden erklärt, Herrn Prof.\ Krull gegenüber
    als Referent zu fungieren, und mir sogar ein positives
    Urteil zusagt. Auch für die Worte von Herrn Prof.\
    Peschl bin ich sehr dankbar.

    ``Trotzallem muß ich die Absicht, mich zu habili\bs
    tieren, zunächst fallen lassen, da ich es für völlig
    ausgeschlossen halte, daß sich jetzt noch eine Meinungs\bs
    änderung bei Herrn Prof.\ Krull vollziehen wird (sodaß
    ich Herrn Prof.\ Köthe gar nicht als Referenten nennen
    können werde) und an eine neue Arbeit ist nicht
    zu denken, weil Herr Prof.\ Krull darauf besteht,
    die vorliegende Arbeit erst gründlich umzuarbeiten
    -- obwohl mir diese jetzt gründlich verleidet ist,
    und ich eigentlich nur den Wunsch habe, endlich einmal
    etwas anderes zu machen, als immer irgendwelche
    Halbordnungen. Wie ich aus diesem Dilemma
    herauskommen werde, weiß ich noch nicht.

    \dotfill\smallskip

    ``Zum Schluß bitte ich noch, Herrn Prof.\ Köthe die Anmerkungen
    zu meiner Arbeit schicken zu wollen --~-- ich habe keine mehr,
    und meine Frau ist ja leider auch nicht da. Prof.\ Köthe
    hat um die Anmerkungen gebeten.''
  }\hspace{2em}


%% file: 19440716-PL-1-1-148-Krull-Lorenzen

\noindent\begin{tabular}[t]{@{}l@{}}%
\label{19440716}%
           Reg.Rat Prof.\ Krull\\%
           Greifswald\\%
           Mar[ine]obs.%
         \end{tabular}\hfill%
         \begin{tabular}[t]{@{}l@{}}%
           Herrn Dr.\\%
           Paul Lorenzen\\%
           23 Wesermünde\\%
           Elbestr.~42\kern1pt\raise2pt\hbox{I}.\\%
         \end{tabular}\medskip

Lieber Herr Lorenzen!\hfill Greifswald, 16.7.44\medskip

\noindent Vor allem diesmal die Versicherung, dass es mir aufrichtig
leid getan hat, dass Sie jetzt auch vom Bombenspuk be\bs
troffen wurden. Ein Glück wenigstens, dass offenbar Ihre
Angehörigen keinen Schaden genommen haben. Hoffentlich
ist Ihnen nicht allzuviel von eigenen Sachen zerstört worden.
Oder hatten Sie überhaupt wirklich gewohnt? Aus der
Bemerkung, Sie hätten sich jetzt ein möbliertes Zimmer
genommen, glaube ich leider auf das Gegenteil schliessen
zu müssen. -- Zu Ihrer Arbeit möchte ich Sie immer
wieder daran erinnern, dass Sie stets bedenken müssen,
dass Sie jetzt im Kriege nicht unter normalen Bedin\bs
gungen schaffen, und dass Sie sich also nicht deprimie\bs
ren lassen dürfen, wenn Sie noch nicht so weit gekom\bs
men sind, als Sie sich vorgestellt hatten. Persönlich
bedaure ich es vor allem, dass ich mit Ihnen nicht mündlich
öfter Ihre Untersuchungen durchsprechen kann. Nur auf\bruch{}
diese Weise könnte ich Ihnen wirklich
helfen, wie ich es gerne täte. Aber
augenblicklich ist da eben leider
keine Gelegenheit dazu.\strut\medskip

\hfil\begin{tabular}{@{}l@{}}%
       Mit den besten Grüssen\\%
       Ihr Wolfgang Krull.%
\end{tabular}


%% file: 19441001-PL-1-1-147-Krull-Lorenzen

\noindent\begin{tabular}[t]{@{}l@{}}%
\label{19441001}%
           Reg.Rat Prof.\ Krull\\%
           4 Greifswald\\%
           Mar[ine]obs.%
         \end{tabular}\hfill%
         \begin{tabular}[t]{@{}l@{}}%
           Sonderführer Dr.\\%
           Paul Lorenzen\\%
           23 Wesermünde\\%
           Marinefachschule\\%
         \end{tabular}\medskip

\hfill Greifswald, 1.10.44\nopagebreak

Lieber Herr Lorenzen!\medskip

\noindent Heute bin ich endlich dazu gekommen, Ihre ausführliche Inhalts\bs
skizze durchzulesen. Entschuldigen Sie also bitte, dass ich Ihnen
den Empfang erst heute bestätige. Ihre eigentliche Schrift werde
ich Ihnen in den nächsten Tagen zusenden. In den letzten Wochen
hatte ich es unter den verschiedensten Abhaltungen einfach vergessen.
Was Sie mir diesmal zuschickten, war entschieden ein Fort\bs
schritt. Ich habe nur immer noch das Gefühl, dass man die
Sache noch wesentlich kürzer und klarer sagen kann. Aber
das sind eben Dinge über die man sich mündlich aussprechen
müsste. Schriftlich gibt es zu leicht Missverständnisse.
Ich würde mich aber auch an Ihrer Stelle im Augenblick gar\bs
nicht an diese Aufgabe der Herstellung eines "`wirksamen"'
Manuskripts halten. Arbeiten Sie doch lieber etwas über das
bisherige hinaus! Was mir immer noch zweifelhaft erscheint,
ist die Frage, wie weit Sie wirklich das \uline{Nicht}kommutative
erledigt haben. Mein Eindruck war, dass bei Ihnen auch im\bruch{}
Nichtkommutativen sehr viel an "`Vertausch\bs
barkeitsforderungen"' drin steckt, so dass man
hierbei von einem "`halbkommutativen"'
Fall reden könnte. Wenn es Ihnen gelänge
diese Bedenken zu zerstreuen, schiene
mir das ein grosser Fortschritt.\medskip

Mit
den besten
Grüssen

Heil Hitler!\qquad Ihr Wolfgang Krull.


%% file: 19450725-Cod-Ms-H-Hasse-1-1022-16-Lorenzen-Hasse
\begin{tabular}{@{}l@{}}%
\label{19450725}%
  Dr.\ Paul Lorenzen\\%
  Bad Pyrmont\\%
  Bahnhofstr.~8\\%
  German%
\end{tabular}\hfill\begin{tabular}{@{}l@{}}%
                     Herrn\\%
                     Prof.\ Dr.\ H. Hasse\\%
                     Göttingen\\%
                     Math.\ Sem.\ d.\ Univers.%
\end{tabular}\medskip

[Antw.\ 1.8.45]\hfill 25/7\medskip

Sehr geehrter Herr Professor,\medskip

\noindent In der Hoffnung, daß Sie das Kriegsende
glücklich überstanden haben, erlaube ich
mir, Ihnen zu schreiben. Ich bin zur Zeit
hier in Bad Pyrmont gut aufgehoben --~-- von
Bonn habe ich noch keine Nachricht. Dürfte
ich Sie um eine Mitteilung bitten, wenn Ihnen
über den Verbleib von Herrn Prof.\ Krull etwas
bekannt ist. Auch von Herrn Prof.\ Scholz, Herrn
Ackermann und Herrn Gentzen, denen ich noch im \bruch
Kriege geschrieben habe wegen
einer neuen Methode für Widerspruchsfreiheitsbe\bs
weise  habe ich keinerlei
Nachricht.\medskip

Mit den ergebensten
Grüßen verbleibe ich\strut

\hfill\begin{tabular}{@{}ll@{}}%
        Ihr&\\%
           &Paul Lorenzen\\%
\end{tabular}
 

%% file: 19450902-universitaetsarchiv_bonn-lorenzen-bericht_ueber_meine_politische_einstellung
\noindent\begin{tabular}{@{}l@{}}%
\label{19450902}%
           Dr.\ Paul \so{Lorenzen}\\[2pt]%
           Wissenschaftl.\ Assistent am\\%
           Mathemat.\ Seminar der Universität \so{Bonn}.
         \end{tabular}\medskip

\noindent\hfil\uuline{Bericht über meine politische Einstellung}\medskip

Ich bin 1915 geboren und ging noch zur Schule, als Hitler
an die Macht kam. In meinem Elternhaus bin ich politisch liberal
erzogen worden und in den letzten Jahren antinationalsozialistisch,
da mein Vater Freimaurer war. Ich selbst missbilligte am National\bs
sozialismus am schärfsten die chauvinistischen und antisemitischen
Tendenzen.

Nach dem Abitur habe ich zunächst ein halbes Jahr Arbeitsdienst
abgeleistet. Im Herbst 1933 begann ich mein Studium und musste
mich -- ehe die Immatrikulation möglich war -- zum Eintritt in
die SA und in den NSDStB melden. Wegen meines passiven Widerstandes
gegen den SA-Dienst bin ich dort nie Scharführer oder Ähnliches
geworden.

Vom Herbst 1934 bis 1935 leistete ich ein Jahr aktiven
Militärdienst ab, da die zweijährige Dienstpflicht drohte. Seit
dieser Zeit bin ich entschiedener Antimilitarist und bin auf Grund
dieser Einstellung während des Krieges erst im letzten halben Jahr
Unteroffizier geworden.

Nach der Dienstzeit setzte ich mein Studium fort. 1936 wur\bs
de ich volljährig und trat aus der evangelischen Kirche aus.
Dieser Austritt geschah nicht aus politischen Gründen, sondern
auf Grund der Überzeugung, dass ich, wenn ich wahrhaftig sein
wollte, mich nicht als gläubigen Christen bezeichnen konnte --
und daher der Kirche nicht angehören dürfte.

1937 musste ich als Mitglied der SA in die Partei ein\bs
treten. Auch dort habe ich niemals ein Amt innegehabt. 1939
wurde ich Assistent am Mathematischen Seminar in Bonn und An\bs
fang 1940 zum Heer eingezogen. Später aber kam ich zur Marine,
wo ich seit 1942 als Mathematiklehrer in Wesermünde an der
Marineschule tätig war, bis ich im Januar 1945 an die Marine\bs
schule Flensburg versetzt wurde.

Meine Frau, die ich 1939 heiratete, stammt aus einer streng
kirchlichen Familie, war aktiv tätig für die Bekennende Kirche --
und daher selbstverständlich Gegnerin des Nationalsozialismus.
Sie war kein Parteimitglied.

Als Zeugen für meine Gegnerschaft führe ich Herrn Professor
Scholz in Münster an, ferner Frau Dr.\ K. Wolff, Bonn, Luisenstr.~3,
bei der ich seit 1939 wohnte. Frau Wolff ist die Witwe eines
jüdischen Arztes.\medskip

\noindent Bonn, den 2.\ September 1945

\hfil\hfil\hfil\hfil Lorenzen\par
 

%% file: 19450903-universitaetsarchiv_bonn-military_government_of_germany_fragebogen-D


  
\newsavebox{\employment}
\sbox{\employment}{\fontsize{5.3}{6.36}\selectfont\sffamily\setlength\tabcolsep{2pt}%
  \begin{tabular}{@{}rrccccc@{}}\toprule
    \multicolumn1c{\emph{From}}&\multicolumn1c{\emph{To}}&\begin{tabular}{@{}c@{}}\emph{Employer and Address}\\\emph{or Military Unit}\end{tabular}&\begin{tabular}{@{}c@{}}\emph{Name and Title of Immediate}\\\emph{Superior or C. O.}\end{tabular}&\emph{Position or Rank}&\begin{tabular}{@{}c@{}}\emph{Duties and}\\\emph{Responsibilities}\end{tabular}&\begin{tabular}{@{}c@{}}\emph{Reasons for Change of Status}\\\emph{or Cessation of Service}\end{tabular}\\
    \multicolumn1c{von}&\multicolumn1c{bis}&\begin{tabular}{@{}c@{}}Arbeitgeber und Anschrift\\oder Militäranschrift\end{tabular}&\begin{tabular}{@{}c@{}}Name und Titel des Dienstvor-\\gesetzten od.\ vorgesetzter Offz.\end{tabular}&\begin{tabular}{@{}c@{}}Stellung oder\\Dienstgrad\end{tabular}&\begin{tabular}{@{}c@{}}Art der Tätigkeit und\\Verantwortungsbereich\end{tabular}&\begin{tabular}{@{}c@{}}Grund für Änderung oder Been-\\digung des Dienstverhältnisses\end{tabular}\\\midrule
    Apr.~33&Sept.~33&Student.\ Arbeitslager&unbekannt&Arbeitsmann&Erdarbeit&Beginn des Studiums\\
    Nov.~34&Okt.~35&reit.\ Artl.\ Abt.\ Verden&Hpt.\ Bamler&Kanonier&Stalldienst&Ende der Dienstverpflichtg.\\
    Okt.~38&Juli 39&Universität Göttingen&Prof.\ Dr.\ Hasse&Stipendiat&Hilfsassistententätigkeit&Anstellung in Bonn\\
    Aug.~39&gegwtg.&Universität Bonn&Prof.\ Dr.\ Krull&Assistent&Assistententätigkeit&\\
    Jan.~40&Juli 40&38503 C&unbekannt&Gefreiter&Rechner&Kommandierung\\
    Juli 40&Apr.~41&O. K. M.&Reg.rat Tranow&"&Schreiber&"\\
    Apr.~41&Juli 41&Steuerm.schule Gotenhafen&Kaptl.\ Götz&"&"&"\\
    Jul.~41&Jan.~42&M. N. O. Borkum&Kaptl.\ v.\ Lom&"&"&"\\
    Jan.~42&Nov.~44&Mar.schule Wesermünde&Kapt.\ Köllner&Sonderf.\ (Feldwebel)&Mathematiklehrer&"\\
    Nov.~44&Jan.~45&Mar. Schütz.\ Btl.~306&Kpt.\ Arlt&Unteroffizier&Schreiber&"\\
    Jan.~45&Febr.~45&Mar.schule Flensburg&Kpt.\ Lüth&"&keine&"\\
    Febr.~45&Apr.~45&8.~M. E. A. Norden&unbekannt&"&Erdarbeit&\\\bottomrule
  \end{tabular}}
\newlength{\demipage}
\setlength\demipage{.5\wd\employment}
\addtolength\demipage{-3.5mm}
\[\makebox[0pt]{%
      \parbox[t]{\demipage}{\selectlanguage{english}\scriptsize\noindent\label{19450903}29. Give a chronological account of your employment and military service beginning with 1st of January 1931, accounting for all promotions, or demotions, transfers, periods of unemployment, attendance at educational institutions (other than those covered in Section B) or training schools and full-time service with para military organizations. (Part time employment is to be recorded in Section F.) Use a separate line for each change in your position or rank, or to indicate periods of unemployment or attendance at training schools or transfers from one military or para military organization to another.}\hspace{7mm}%
      \parbox[t]{\demipage}{\selectlanguage{german}\footnotesize\noindent29. Geben Sie in zeitlicher Folge eine Aufzählung Ihrer Beschäftigung und Ihres Militär\bs
dienstes seit dem 1.\ Januar 1931 an, mit Begründungen für alle Beförderungen oder
Degradierungen, Versetzungen, Arbeitslosigkeit, Besuch von Bildungsanstalten (außer
solchen, die bereits in~B angeführt sind) oder Ausbildungsschulen, und Volldienst in
militärähnlichen Organisationen (Nebenbeschäftigungen sind in Abschnitt~F anzugeben).
Benutzen Sie eine gesonderte Zeile für jeden Wechsel in Stellung oder Rang oder zur
Angabe von Arbeitslosigkeits-\hspace{0pt}Zeitabschnitten oder für den Besuch von Ausbildungs\bs
schulen oder für Versetzungen von einer militärischen oder militärähnlichen Organisa\bs
tion zu einer anderen.}}\]
\[\makebox[0pt]{\usebox\employment}\]\par 

%% file: 19450906-universitaetsarchiv_bonn-peschl-lorenzen-gutachten
\noindent\begin{tabular}[t]{@{}l@{}}%
\label{19450906}%
           Prof.\ Dr.\ Ernst Peschl\\%
           Bonn, Arndtstr.~2.
         \end{tabular}\hfill Bonn, den 6.9.45.

\begin{flushright}%
  An den Prüfungsausschuß der Universität
  \mbox{in Hdn von Herrn Prof.\ H. von Weber}\strut\medskip

  \begin{tabular}{l}%
    \uline{\hbox{\so{Bonn}}}\\%
    Sternwarte.%
  \end{tabular}%
\end{flushright}

\noindent Betrifft: \begin{tabular}[t]{@{}l@{}}%
            Politische Einstellung des Herrn\\%
            \uline{Dr.\ Paul Lorenzen}, geb. 24.3.15,\\%
            pl.\ Assistent am Math.\ Seminar d.\ Univ.%
          \end{tabular}\medskip

Ich kenne Herrn Dr.\ Lorenzen vor allem aus der
Zeit vom August~39 bis zu seiner Einziehung zur
Wehrmacht Anfang~40. In dieser Zeit hatte ich
ausgiebig Gelegenheit mich mit ihm eingehend zu
unterhalten, da ich ihn fast täglich sah. Aber auch vor
dieser Zeit (36/37) wie auch nachher traf ich ihn aus
Anlaß der Jahresversammlungen der Deutschen Mathe\bs
matiker-Vereinigung und hatte dabei ebenfalls längere
Unterhaltungen mit ihm. Ich glaube ihn daher
eingehendst zu kennen, zumal er seiner ganzen Ver\bs
anlagung nach ein offener absolut aufrichtiger Charak\bs
ter ist.

Er ist ein sehr kritischer Mensch und hatte vom ersten
Augenblick seiner Fühlungnahme mit mir an stärkste
Ausdrücke (meist sehr sarkastischer Art) der völligen
Ablehnung des Nationalsozialismus und der Person
Hitlers und aller seiner Trabanten geäußert. Dies
entsprach auch seiner ganzen inneren Haltung. Des\bs
gleichen lehnte er jede Form von Militarismus
scharf ab.

Auch sein mir bekannter Kirchenaustritt im Jahre~37
hat gar nichts mit politischen Motiven zu tun. Es war
lediglich die aufrichtige Schlußfolgerung aus der Haltung
eines philosophischen Agnostizismus, die sich nach
langem ehrlichen Ringen um erkenntnistheoretische
Fragen bei ihm herausgebildet hatte. Jede Feindselig\bs
keit gegen religiöse Überzeugungen ist ihm völlig
fremd. Im Gegenteil hat er eine Frau aus streng
religiöser Familie der evangelischen Bekenntniskirche
geheiratet.\medskip

\hfill Prof.\ Dr.\ Ernst Peschl.\comment{Ernst Peschl (1906--1986) is professor of mathematics in Bonn. See \citet[pages 461--462]{MR1991149} on his attitude towards national socialism.}\par
 

%% file: 19450914-universitaetsarchiv_bonn-pruefungsauschuss_weber_fitting_troll-lorenzen_gutachten
\label{19450914}\noindent\hfil In der Prüfungssache\nopagebreak

\noindent des Assistenten am Mathematisch-Naturwissenschaftlichen Seminar
der Universität Bonn Dr.\ Paul \so{Lorenzen}

\noindent erstattet der Prüfungsausschuss bestehend aus den Professoren
von Weber, Fitting und Troll in seiner Sitzung vom 14.9.1945
folgendes\medskip

\noindent\hfil\so{Gutachten.}\medskip

Lorenzen, geboren 1915, trat mit Beginn seines Studiums der
SA und dem NSDStB bei. 1937 wurde er in die Partei übernommen.
Er ist ein Schüler von Professor Scholz in Münster. Seine Haltung
war eindeutig antinationalsozialistisch, wie auch das beiliegende
Zeugnis von der Witwe des jüdischen Arztes Wolff bezeugt. Für die
Partei hat er sich niemals betätigt.

Lorenzen ist nur formales Parteimitglied. Der Ausschuss befür\bs
wortet seine Belassung in seiner bisherigen Stellung.\nopagebreak\medskip

\hfil von Weber\nopagebreak\medskip

\uline{1 Anlage}
 

%% file: 19460813
\noindent\label{19460813}
\hspace*{-9pt}\begin{tabular}[t]{@{}c@{}}
  \bfseries\large Mathematisch-naturwissenschaftliche\\%
  \bfseries\large Fakultät\\%
  \bfseries der Rheinischen Friedrich-Wilhelms-Universität\\[7pt]%
  J.-Nr.~864
\end{tabular}\hfill\llap{\textbf{Bonn, den }13. August 1946}

\hfill
\begin{tabular}{@{}l@{}}
  An\\%
  den Herrn Oberpräsidenten der Nordrheinprov.\\[2pt]%
  \uline{\hbox{\so{Düsseldorf}}}\\[1pt]
  durch Se.\ Magnifizenz, den Herrn Rektor\\%
  der Rheinischen Friedr.-Wilh.-Universität,%
\end{tabular}\medskip

Namens der Mathematisch-naturwissenschaftlichen Fakultät der Rheinischen Fried\-rich-\hspace{0pt}Wilhelms-\hspace{0pt}Universität in Bonn teile ich ergebenst mit, dass die Mathe\-ma\-tisch-\hspace{0pt}naturwissenschaftliche Fakultät Herrn%
\begin{center}
  Dr.\ Paul \so{Lorenzen}\comment{\foreignlanguage{english}{In a letter to Scholz dated 7 June 1946\label{19460607}
    , Lorenzen writes:}
    ``Herr Prof.\ Krull
    ist seit kurzem glücklicherweise zugelassen und hat sich
    mit meiner sofortigen Habilitation vollkommen einver\bs
    standen erklärt.''

    \foreignlanguage{english}{The day after,\label{19460608} he writes to him
    :}
    ``Welche Aussichten dann hier
    bestehen, darüber weiß ich nichts, aber
    ich vermute, daß ich hier irgendwie
    einen Lehrauftrag bekommen werde. Für dieses
    Semester habe ich -- ausnahmsweise -- einen
    Auftrag für eine Algebravorlesung bekommen.''}
\end{center}
nach am 9.8.1946 gehaltener Antrittsvorlesung über:%
\[\text{"`Über den Verbandsbegriff"'}\]
als Privatdozent für Mathematik zugelassen und die venia legendi verlie\bs
hen hat.\medskip

\hfill
\begin{tabular}{@{}c@{}}
  gez.\ Reichensperger\\%
  Dekan%
\end{tabular}\hspace*{2cm}
 

%% file: 19530418-ko-05-0647-Krull-Scholz
\label{19530418}%
\hfil Lieber Herr Scholz!\hfil\hfil\hfil\hfil\hfil Bonn, 18.4.53.\medskip

Was Sie mir über Lorenzen schrieben, hat mich sehr gefreut. Ich schätze
seine Arbeiten zur Begründung der Analysis ausserordentlich hoch. Ich hatte
bei dem Arbeiten mit dem Überabzählbaren, insbesondere mit dem Wohlord\bs
nungssatz, immer das Gefühl man benutzt da Fiktionen, die eines Tages
durch vernünftigere Begriffsbildungen ersetzt werden müssen. Aber ich regte mich
darüber nicht auf, weil ich überzeugt war, dass bei vorsichtiger Anwendung
der geläufigen "`Fiktionen"' nichts Falsches herauskommt, und weil ich mit Sicher\bs
heit auf den Mann rechnete, der eines Tages alles in Ordnung brächte. Lorenzen
hat nun nach meiner Überzeugung den richtigen Weg gefunden, und das ist
allein schon eine Leistung, für die er ein Ordinariat verdiente. -- Es ist mir\bruch{}
nun eine grosse Befriedigung, dass auch Sie, --~wenn auch unter einem
etwas anderen Blickwinkel als dem meinen,~-- von Lorenzen so viel halten.\comment{Compare Krull's letter to Hasse dated 15~August 1946 \citep[§~1.95]{roquette04}, in which he already expresses his appreciation of Lorenzen.}
Dass es für die mathematische Logik und Grundlagenforschung an der nö\bs
tigsten Stellenzahl fehlt, ist auch meine Überzeugung. Ich bin also sehr
gerne bereit, Sie bei einer Aktion, die die Vermehrung der einschlägigen Stellen
anstrebte, aufs wärmste zu unterstützen. In Bonn selbst wird sich aller\bs
dings im Augenblick in dieser Hinsicht nichts tun lassen, da wir gerade um
unser Extraordinariat bzw.\ Ordinariat für Angewandte Mathematik kämp\bs
fen, und man nur schrittweise vorgehen kann. Aber schreiben Sie mir doch
bitte, wie Sie sich das Aktionsprogramm vorstellen. -- Es tut mir sehr leid,
dass Sie wieder so lange und so schwer unter Magenkrämpfen litten. Hoffent\bs
lich bleiben Sie jetzt recht lange davon verschont!\medskip

Mit herzlichen Grüssen\quad Ihr Wolfgang Krull

%% file: 19530514-Cod-Ms-H-Hasse-1-1022-17-Lorenzen-Hasse

\begin{tabular}[b]{@{}c@{}}%
\label{19530514}%
  Prof.\ Dr.\ P. Lorenzen\\%
  Bonn, Luisenstr. 3%
\end{tabular}\hfill\begin{tabular}[t]{l}%
                     \llap[Mein R.\\%
                     bespr.\ 3.7.53]%
                   \end{tabular}\hfill Bonn, den 14.5.53\hspace*{3em}\medskip

\noindent Sehr verehrter, lieber Herr Hasse,\nopagebreak\medskip

\noindent aus Freude und Dankbarkeit darüber, dass Sie in Ihrem Vor\bs
trag über Mathematik als Wissenschaft, Kunst und Macht\footnote{\citealt{MR0053055}.}
auch einmal die menschliche Seite der Mathematik --\nbsp
die Empfindungen und Beweggründe, die uns gerade an diese
Wissenschaft vor allen anderen binden (obwohl Sie betonen,
dass Sie nur Ihre subjektive Einstellung darstellen, haben
Sie doch zugleich für viele andere gesprochen) -- ins
Bewusstsein heben, erlaube ich mir, Ihnen zu schreiben.

\noindent Das von Ihnen zuerst behandelte Motiv: Mathematik als
reinste Wissenschaft, als Prototyp unvergänglicher Erkennt\bs
nis, das ist wohl auch im Verlauf der bisherigen Geschichte
-- seit Thales und den Pythagoräern -- das wirksamste Motiv
gewesen. Gegenwärtig scheint sich dieses Motiv allerdings
in einer besonderen Gefahr zu befinden. Die Geometrie,
die doch bisher -- wie ja offiziell auch heute noch --
immer zur Mathematik gerechnet wurde, hat seit dem vorigen
Jahrhundert den Charakter der absoluten Sicherheit verloren.
Die Mathematik hat sich daher bezüglich der Geometrie auf
Implikationen der Form: "`wenn die und die Axiome gelten,
dann gilt auch \dots{}"' zurückgezogen, also auf Aussagen, die
kaum noch geometrisch zu heissen verdienen. Wenn ich Sie
nicht missverstanden habe, verstehen Sie unter Mathematik
den engeren Begriff, der nur Logik, Arithmetik und Analysis
umfasst. Neben diesen konkreten Teilen der Mathematik würden
dann die abstrakten Teile -- also die axiomatischen Theorien
der Algebra und Topologie -- eine sekundäre Rolle spielen. 
Will man die Geometrie, eventuell zusammen mit der Mechanik
nicht als ein eigenes Fach konstituieren, so wird man sie
mit zur theoretischen Physik rechnen müssen. Gegenwärtig
wollen nun manche sogar noch der Mathematik im engeren Sinne
die unumstössliche Gültigkeit absprechen, also z.\ B.\ auch die
Arithmetik auf Implikationen: "`wenn die Peano-Axiome gelten,
dann gilt auch \dots{}"' reduzieren.

\noindent Dann gäbe es keine Mathematik als Wissenschaft mehr, es gäbe
nur noch formales Ableiten aus Axiomen, die aufgrund pragma\bs
tisch-empirischer Kriterien gewählt werden. Zieht man schliess\bs
lich die Logik mit in diesen Relativierungsprozess hinein,
so würden auch noch die Regeln des formalen Ableitens aus
Axiomen nur durch aussermathematische Kriterien zu
gewinnen sein.

\noindent Ich wage, zu vermuten, dass es vor allem diese Gefahr ist,
die Sie als Formalisierung eine inhärente Tragik der Mathe\bs
matik nennen. Hier würde ich mir nun gern die Frage erlauben:
ist diese Formalisierung zu unterscheiden von der "`Formali\bs
sierung"', die in Ihrer Betrachtung über die Bedeutung einer
klaren und prägnanten Bezeichnungsweise als eine Gefahr dar\bs
gestellt wird? Die Existenzaussage auf p.~22 würde nach Er\bs
setzung der logischen Partikeln durch Symbole, nämlich von%
\[\begin{tabular}{@{}llc@{}}%
\rlap{wenn, so}\hphantom{es gibt genau ein}&durch&$\rightarrow$\\%
und&durch&$\land$\\%
es gibt&durch&$\biglor$\\%
ist&durch&$\varepsilon$%
\end{tabular}\]%
lauten:%
\[%
  p\mathrel\varepsilon\mathrm{prim}\land p=4n+1\ \rightarrow\ \biglor\nolimits_{\!\!x,y}\ p=x^2+y^2\text.
\]

Zur "`Formalisierung"' der Eindeutigkeitsaussage könnte%
\[\begin{tabular}{@{}llc@{}}%
es gibt genau ein&durch&$\biglor\limits^1$%
\end{tabular}\]%
ersetzt werden, und wir erhielten dann%
\[%
  p\mathrel\varepsilon\mathrm{prim}\land p=4n+1\ \rightarrow\ \mathop{\biglor\limits^1}\nolimits_{\!\!x,y}x<y\land p=x^2+y^2\text.%
\]

\noindent Wenn man die logischen Symbole nur als prägnante Abkürzung
für die umgangssprachlichen Partikel benutzt -- ihnen also
alle Inhaltlichkeit lässt, und die Aussagen jetzt nicht
plötzlich als bedeutungslose Zeichenreihen behandelt --
dann scheinen mir genau dieselben Gründe für diese Symbolik zu
sprechen, die uns veranlasst haben, von:%
\[\begin{tabular}{@{}ll@{}}%
  &\llap{"`}$p$ ist die Summe der Quadrate von~$x$ und~$y$"'\\%
  zu&$p=x^2+y^2$%
\end{tabular}\]%
überzugehen. Könnte man der letzten Formel vorwerfen, dass
sie dem Fluss der suggestiven Sprache Gewalt antut?

\noindent Dadurch, dass die logische Symbolik von Hilbert zu dem Zweck
gebraucht ist, eine inhaltliche Theorie in einen axiomatischen
Formalismus zu verwandeln, scheint mir die Möglichkeit
eines inhaltlichen Gebrauches der logischen Symbole zu kurz
gekommen zu sein.

\noindent Ich bitte Sie herzlich, diese Bemerkungen nicht als Pro\bs
paganda für die Logistik auffassen zu wollen --~ganz im
Gegenteil, es ist ja gerade die Logistik, die das inhaltliche
Denken verkennt -- ich würde mich aber freuen, wenn Sie
diesen Brief als einen Ausdruck meiner Dankbarkeit annehmen
wollen, dafür, dass Sie die geistigen Schmerzen eines
"`deutschen Aufsatzes"' nicht gescheut haben, um einmal zur
Besinnung auf unser Tun aufzurufen.\strut\medskip

\begin{tabular}{@{}rl@{}}%
  Mit den ergebensten Grüssen bin ich&\\%
  Ihr\qquad&\\%
  &Paul Lorenzen%
\end{tabular}


%% file: 19530605-Cod-Ms-H-Hasse-1-1022-Beil-6-Hasse-Lorenzen
\hfill 5.\ Juni 1953\medskip
\label{19530605}%
 
\noindent
\begin{tabular}{@{}l@{}}%
Prof.\ Dr.\ H. Hasse\\%
Ahrensburg i.\ H.\\%
Hamburgerstr.~43%
\end{tabular}\medskip

Lieber Herr Lorenzen,\medskip

\indent\indent es ist sehr freundlich von Ihnen, dass Sie solches Interesse an
meinem kleinen Büchlein genommen haben. Es fehlt mir jetzt mitten im
Semesterbetrieb leider die Zeit, mich zu den von Ihnen gemachten Bemer\bs
kungen zu äussern. Da ich aber voraussichtlich noch im Laufe dieses
Semesters nach Bonn kommen werde, können wir dann vielleicht einmal
mündlich darüber sprechen.\strut\medskip

\hfil\begin{tabular}{@{}c@{}}%
       Mit herzlichen Grüssen\\%
       Ihr [Hasse]%
     \end{tabular}
 

%% file: 19530609-Cod-Ms-H-Hasse-1-1022-18-1-Lorenzen-Hasse
\label{19530609}%
\indent\indent Sehr verehrter, lieber Herr Hasse,\medskip

\noindent für die freundliche Aufnahme meines Briefes zu Ihrem
"`Büchlein"' danke ich Ihnen sehr. Wenn Ihr Plan,
nach Bonn zu kommen, sich verwirklichte, wäre das
sehr schön. Wenn Ihre Zeit nicht allzu knapp sein
wird, darf ich Sie vielleicht für eine kurze Zeit
bitten, bei uns sein zu wollen.\strut\medskip

\hfil\begin{tabular}{@{}rl@{}}%
       Mit den ergebensten Grüßen&\\%
       bin ich stets Ihr&\\%
                        &\llap{Paul} Lorenzen.%
\end{tabular}\strut\medskip

\noindent
\begin{tabular}{@{}l@{}}
  Bonn, 9.6.53\\%
  Luisenstr.~3%
\end{tabular}

%% file: 195907-Cod-Ms-H-Hasse-1-1022-18-2-Lorenzen-Hasse
\hfill [P. Lorenzen, Juli 1959]\medskip
\label{195907}%

\noindent\uline{Sätze}\quad [Hasse J\,f\,M~152 (1923)]

\[\begin{aligned}%
 p\equiv1\bmod 4,\ q\equiv1\bmod\xcancel{2}\smash{\llap{\raise.7\baselineskip\hbox{4}}},\footnotemark\ \Bigl(\frac qp\Bigr)=-1&\enskip\rightarrow\enskip\bigland_{\text{rat}}\!_{\raisebox{-3pt}{$\scriptstyle t$}}\mathrel.\Phi(p,q,t)\leftrightarrow t\text{ $p$-$q$-ganz.}\\%
 p\equiv3\bmod 4,\ \rlap{$q=2$}\hphantom{q\equiv1\bmod\xcancel{2}\smash{\llap{\raise.7\baselineskip\hbox{4}}},\footnotemark\ \Bigl(\frac qp\Bigr)=-1}&\enskip\rightarrow\enskip\bigland_{\text{rat}}\!_{\raisebox{-3pt}{$\scriptstyle t$}}\mathrel.\Phi(1,p,t)\leftrightarrow t\text{ $p$-$q$-ganz.}%
\end{aligned}\]\addtocounter{footnote}{-4}\footnotetext{This correction is by Hasse, see the following letter.}
\[\Phi(r,s,t)\leftrightharpoons\biglor_{\text{rat}}{}_{\raisebox{-3pt}{$\!\!\scriptstyle x,y,z$}}x^2+ry^2-sz^2=2+rst^2\text.\]

\noindent\uline{Folgerungen}\nopagebreak

\[\bigland\nolimits_p\biglor\nolimits_{\!\!q,m,n}\bigland\nolimits_t\mathrel.\Phi(m,n,t)\enskip\leftrightarrow\enskip t\text{ $p$-$q$-ganz.}\]

Für eine Klasse~$K$ von Paaren~$m,n$ gilt also
\begin{enumerate}[(1)]
\item $\bigland_K\!_{\raisebox{-3pt}{$\scriptstyle m,n$}}\Phi(m,n,0)$
\item $\bigland_K\!_{\raisebox{-3pt}{$\scriptstyle m,n$}}\bigland_{\text{rat}}\!_{\raisebox{-3pt}{$\scriptstyle t$}}\mathrel.\Phi(m,n,t)\enskip\rightarrow\enskip\Phi(m,n,t+1)$
\item $\bigland_{\text{rat}}\!_{\raisebox{-3pt}{$\scriptstyle u$}}\mathrel.\bigland_K\!_{\raisebox{-3pt}{$\scriptstyle m,n$}}\Phi(m,n,u)\enskip\rightarrow\enskip u\text{ ganz}$.
\end{enumerate}

\noindent\uline{Satz}\quad $\bigland_{\text{rat}}\!_{\raisebox{-3pt}{$\scriptstyle u$}}\mathbin.u\text{ ganz}\enskip\leftrightarrow\enskip\bigland_{\text{rat}}\!_{\raisebox{-3pt}{$\scriptstyle r,s$}}\mathbin.\Phi(r,s,0)\land\bigland\nolimits_t{\mathbin.\Phi(r,s,t)\rightarrow\Phi(r,s,t+1)\mathbin.}\rightarrow\Phi(r,s,u).$
[J. Robinson J\,S\,L~14 (1949)]\smallskip

\noindent\uline{Satz}\quad Die Theorie der kommutativen Körper ist unentscheidbar.
 

%% file: 19590801-Cod-Ms-H-Hasse-1-1022-Beil-7-Hasse-Lorenzen
\hfill1.\ August 1959\medskip
\label{19590801}%

Lieber Herr Lorenzen,\medskip

Heute endlich bin ich dazu gekommen, mich mit den Kriterien über
quadratische Darstellungen zu beschäftigen.

Es ist alles in Ordnung, nur dass es im Falle $p = 1$ mod.~$4$ auch
$q = 1$ mod.~$4$ (statt nur mod.~$2$) heissen muss,
was wohl auf einem Schreibfehler Ihrerseits beruht.

Die Kriterien lassen sich allerdings nicht unmittelbar meiner
Dissertation entnehmen, in der ja von Ganzzahligkeit nicht die Rede
ist. Dass aus der Darstellbarkeit die Ganzheit von~$t$ für~$p,q$ folgt,
ergibt sich ohne weiteres durch Übergang zu der zugeordneten
ganzzahligen primitiven homogenen Gleichung und Kongruenzbetrachtung
mod.~$p,q$ (wobei für $q = 2$ sogar mod.~$4$ oder gar~$8$ zu rechnen ist).
Dass umgekehrt für $p,q$-ganze~$t$ wirklich Darstellbarkeit besteht,
beruht darauf, dass in meiner Dissertation die ternäre Darstellbarkeit
auf binäre lokale \uline{Nicht}darstellbarkeiten zurückgeführt ist. Schliesst
man also indirekt, so führt die Annahme der ternären Nichtdarstellbar\bs
keit auf binäre lokale Darstellbarkeiten, und daraus kann dann wie
vorher durch Übergang zu den ganzzahligen primitiven homogenen Glei\bs
chungen hier auf die Nichtganzheit von~$t$ für $p,q$ geschlossen werden.

Was die von Robinson gezogenen logistischen Folgerungen betrifft,
so verstehe ich sie leider nicht ganz und habe hier auch keine Möglich\bs
keit, das JSL einzusehen.

Mit besten Grüssen und Ferienwünschen, auch für Ihre Frau und
Tochter,\medskip

\hfil Ihr [Hasse]\hfil
 

%% file: 19600307-Cod-Ms-H-Hasse-1-1022-Beil-8-Hasse-Lorenzen
\label{19600307}%
\hfill 7.\ März 1960\nopagebreak\medskip

\noindent
\begin{tabular}{@{}lc@{}}%
  Herrn Prof.\ Dr.\ P. Lorenzen&\\%
                               &\uline{Kiel}\\%
                               &Philosophisches Seminar d.\ Universität%
\end{tabular}\medskip

\noindent Lieber Herr Lorenzen,\medskip

recht herzlich möchte ich mich für die Zusendung Ihres kleinen
Büchleins über die Entstehung der exakten Wissenschaften\footnote{\citealt{MR0111655}.} bedanken.
Ich habe es in den letzten Tagen mit grosser Freude gelesen.

\hfil Mit besten Grüssen,
auch an die verehrte Gattin,\nopagebreak\medskip

\hfil\hfil Ihr [Hasse]
 

%% file: 19610627-Cod-Ms-H-Hasse-1-1022-Beil-9-Hasse-Lorenzen
\label{19610627}%
\hfill 27.\ Juni 1961\nopagebreak\medskip

\noindent Lieber Herr Lorenzen,\medskip

haben Sie zunächst recht herzlichen Dank für die Zusendung
Ihres Sonderabdrucks "`Ein didaktisches Konstruktivitätskriterium"'.\footnote{In fact ``Ein dialogisches Konstruktivitätskriterium'', \citealt{MR0147364}.}
Wir haben schon seit längerer Zeit die Absicht, Sie um einen
ausführlichen Vortrag über Ihre Ergebnisse zur Grundlegung der
Mathematik zu bitten. In diesem Sommersemester liess sich das
deshalb schlecht einrichten, weil wir laufend auswärtige, ins\bs
besondere überseeische Gäste hatten und dazu auch noch einige
Mathematiker aus Gründen der Besichtigung für eine evtl.\ Berufung
einladen mussten. Ich möchte Sie heute fragen, ob es Ihnen mög\bs
lich wäre, uns im Wintersemester die Freude Ihres Besuchs und
eines Vortrags über Ihre Ergebnisse zu machen. Es würde uns daran
liegen, dass dieser Vortrag uns Ihre Gedankengänge möglichst in
einer uns allen wirklich verständlichen Form nahebringt. Wir
würden Ihnen dafür auch gern zwei aufeinanderfolgende Kolloquiums\bs
tage (immer Dienstag 16~Uhr) zur Verfügung stellen, wo Sie dann
jeweils maximal 60~Minuten Vortragszeit hätten. Selbstverständ\bs
lich würden wir Ihnen Ihre Reise- und Aufenthaltskosten erstatten
mit DM~150,-- je Vortragstag. Ich schreibe schon heute, weil jetzt
noch fast alle Dienstage des Wintersemesters zur Verfügung stehen,
ausgenommen lediglich der erste Dienstag, nämlich der 7.11., an
dem unser P. Jordan über ein Problem aus der Theorie der nicht\bs
kommutativen Verbände etwas erzählen wollte.

Sie würden uns wirklich eine sehr grosse Freude machen,
wenn Sie auf diese unsere Bitte eingingen.

\hfil Mit herzlichen Grüssen

\hfil\hfil Ihr [Hasse]
 

%% file: 19610704-Cod-Ms-H-Hasse-1-1022-19-Lorenzen-Hasse
\label{19610704}%
{\footnotesize\begin{tabular}[t]{@{}c@{}}%
  \normalsize Prof.\ P. Lorenzen\\%
  \small Kiel\\%
  Clausewitzstraße~14%
\end{tabular}}\hfill \raise1pt\hbox{4}/\lower1pt\hbox{7}~61\medskip

Sehr verehrter, lieber Herr Hasse,\medskip

\noindent für Ihre liebenswürdige Anfrage und Einladung
danke ich Ihnen herzlichst.

Ich gehe in der Tat gern auf dieses Anerbieten
ein, weil ich meine gegenwärtige "`Philosophie"' ja nur
umständehalber als Philosophie vertrete, während es
doch nichts als eine mathematische Bemühung um
die Frage ist, ob denn tatsächlich (wie man heute meint)
die Arithmetik und Mengenlehre eine axiomatische
Theorie -- wie die Geometrie -- sei. Ist diese Frage
zu verneinen, so bleibt zu fragen, wie denn die Sätze,
die als "`Axiome"' benutzt werden, zu \uline{beweisen} sind.

Leider sind diese Fragen z.\ Zt.\ mit vielerlei
unnötigem Ballast beladen, sodaß ich zur Verständlich\bs
keit von dem Angebot zweier Kolloquien gerne Gebrauch \bruch
mache. Am ersten Tage würde ich vorschlagen ,,Logik
und Arithmetik`` zu behandeln, am zweiten Tage
"`Geometrie und Mengenlehre"'.

Da ich mir die Tage noch aussuchen darf,
schlage ich Dienstag, den 5.~Dezember und Dienstag,
den 12.~Dezember vor.

\hfil Mit herzlichen Grüßen stets\nopagebreak\medskip

\hfill Ihr P. Lorenzen\par
 

%% file: 19610708-Cod-Ms-H-Hasse-1-1022-Beil-10-Hasse-Lorenzen
\label{19610708}%
\hfill 8.7.1961\nopagebreak\medskip

\indent\indent Lieber Herr Lorenzen,\medskip

Haben Sie recht herzlichen Dank für Ihre Zusage, über die
wir uns alle sehr freuen. Wir sind mit Ihren Vorschlägen hin\bs
sichtlich Daten und Themen gern einverstanden, haben davon Vormer\bs
kung genommen und hoffen nur, daß nichts dazwischenkommt.

Für den Fall, daß wir für Sie Quartier besorgen sollen, geben
Sie uns bitte rechtzeitig Nachricht. Oder wollen Sie jedesmal
am Abend wieder nach Kiel zurückfahren?

Hit herzlichen Grüßen, auch von meiner Frau und an die Ihre,\medskip

\hfil Ihr\rlap{ [Hasse]}
 

%% file: 19610711-Cod-Ms-H-Hasse-1-1022-20-Lorenzen-Hasse
\label{19610711}%
{\footnotesize\begin{tabular}[t]{@{}c@{}}%
  \normalsize Prof.\ P. Lorenzen\\%
  \small Kiel\\%
  Clausewitzstraße~14%
\end{tabular}}\medskip

Sehr verehrter, lieber Herr Hasse,\medskip

\noindent herzlichen Dank für die Bestätigung.

Da ich von Hamburg ja noch gegen 23\textsuperscript{\,\underline{00}} nach
Kiel kommen kann, werde ich sicherlich abends
zurückfahren und werde also kein Quartier brauchen.

\hfil Mit herzlichen Grüßen stets\nopagebreak

\hfil\hfil Ihr P. Lorenzen\nopagebreak\medskip

Kiel \raise1pt\hbox{11}/\lower1pt\hbox{7} 61


%% file: 19621009-Cod-Ms-H-Hasse-1-1022-Beil-11-Hasse-Lorenzen
\label{19621009}%
\hfill Prof.\ Dr.\ H. Hasse\nopagebreak\medskip

\begin{tabular}{@{}l}%
  Herrn Prof.\ Dr.\ P. Lorenzen\\%
  Erlangen\\%
  Philosophisches Institut d.\ Universität%
\end{tabular}

\hfill9.10.1962\medskip

\noindent Lieber Herr Lorenzen,\nopagebreak\medskip

Mit meinen herzlichen Glückwünschen zur Ernennung in Erlangen
verbinde ich meinen besten Dank für die Zusendung Ihres neuesten
Opus "`Meta-Mathematik"'.\footnote{\citealt{lorenzen62}.} Nach genauer Lektüre der Einleitung und
Gewinnung eines Überblicks über den Gang der Handlung habe ich
dies sicher von hoher Warte mit größter Klarheit und bemerkens\bs
wertem didaktischen Geschick geschriebene Werk neben der Scholz\bs
schen "`Metaphysik"' in meinem Bücherschrank abgestellt, mit dem
Vorsatz, wenigstens einmal bis zum Widerspruchsfreiheitsbeweis
der Arithmetik im Einzelstudium vorzudringen. Vielleicht finde
ich dazu nach meiner Emeritierung die Zeit.

Für heute recht herzliche Grüße
und gute Wünsche von Haus zu Haus\nopagebreak\medskip

\hfil Ihr\rlap{ [Hasse]}
 

%% file: 196307xx-Cod-Ms-H-Hasse-1-1022-21-Lorenzen-Hasse
\label{196307}%
{\footnotesize\begin{tabular}[t]{@{}c@{}}%
  \normalsize Prof.\ P. Lorenzen\\%
  \small 852 Erlangen\\%
  Saalestraße~1%
\end{tabular}}\hfill\begin{tabular}[t]{@{}l@{}}%
                     \llap[Beantw.~21.7.1963\\%
                     Kritik anerkannt.\\%
                     Aber in solchem Werk kein\\%
                     Platz für Grundlagenfragen]%
                   \end{tabular}\medskip

Sehr verehrter, lieber Herr Hasse,\medskip

\noindent Ihr großes Buch\comment{\citealt{MR0153659}.} kam hier vor einigen Tagen an: welch
ein vorteilhafter "`Austausch"' ist das für mich, da Sie ja
nur das Metamathematikbändchen erhalten haben!

So versuche ich noch, meinen herzlichen Dank dazuzu\bs
geben. Daß die Zahlentheorie die Königin der mathe\bs
matischen Disziplinen ist, ist mir durch Ihr Buch
wieder einmal deutlich geworden. Welche Wohltat,
daß am Anfang kein willkürliches Axiomensystem
steht: die Gegenstände, die man zu erkennen trachtet,
sind vielmehr von vornherein eindeutig durch Konstruktion
gegeben. [Die axiomatischen Begriffe wie Körper, Gruppe, \dots\bsp
dienen stets nur als ein Mittel, um die Beziehungen zwischen
den konstruierten Gegenständen -- vor allem also Beweiszusammen\bs
hänge -- klar und deutlich zu erfassen] \bruch

Bei meiner notorischen Kritik am Cantorschen Mengen\bs
begriff werden Sie es mir verzeihen, daß ich die ,,Kon\bs
struktion`` der vollständigen Hülle eines bewerteten
Körpers als einen Schönheitsfehler betrachte: es werden
ja "`alle"' $\varphi$-konvergenten Folgen~$(a_n)$ betrachtet. Wie
läßt sich diese Unmenge aber "`konstruieren"'?

M. E. kommt man zu den gewünschten Resultaten
über die Fortsetzungen von Bewertungen, wenn man statt
"`aller"' konvergenten Folgen immer nur "`genügend viele"'
konvergente Folgen zur Erweiterung des bewerteten Körpers
benutzt.

Diese "`puristische"' Kritik wäre wohl der Forderung
zu vergleichen, etwa beim Beweis des Dirichletschen Einheiten\bs
satzes \uline{ohne} den \uline{vollen} Körperraum (und Logarithmenraum)
auszukommen: hier braucht man ersichtlich nicht \uline{alle}
reellen Zahlen, sondern könnte sich mit Abschätzungen
behelfen.

Daher möchte ich -- statt solcher Kritik -- lieber meinen herzlichsten Dank
noch einmal wiederholen.\strut\medskip

\hfill\begin{tabular}{@{}l@{}}%
        Mit den ergebensten Grüßen\\%
        stets\\%
        \hphantom{stets\qquad}Ihr P. Lorenzen%
\end{tabular}
 

%% file: 196309xx-Cod-Ms-H-Hasse-1-1022-22-Lorenzen-Hasse
\label{196309}%
{\footnotesize\begin{tabular}[t]{@{}c@{}}%
  \normalsize Prof.\ P. Lorenzen\\%
  \small 852 Erlangen\\%
  Saalestraße~1%
\end{tabular}}\hfill\begin{tabular}[t]{@{}l@{}}%
                     \llap[Gedankt\\%
                     durch Postkarte\\%
                     26.9.63]%
                   \end{tabular}\hfill September 1963\medskip

Sehr verehrter, lieber Herr Hasse,\medskip

\noindent in der Hoffnung, daß Sie dieser Brief nach guter
Rückkehr von der Amerikareise antrifft, möchte ich
Ihnen vielmals für die Karte aus Morris Inn, Notre Dame
(ich habe dort 1957 einmal übernachtet) danken.

Meine "`Kritik"' war durchaus nicht so gemeint,
daß es einer "`Rechtfertigung"' Ihrerseits bedurft hätte:
das Ziel des Buches, wirklich zu den Höhen der Zahlen\bs
theorie zu führen rechtfertigt -- in unserer gegenwärtigen
Situation -- vollkommen die Benutzung der traditionellen
Fundamente.

Trotzdem freue ich mich natürlich sehr über Ihre
Zustimmung, daß dann, wenn z.\ B. von "`allen"' Folgen (etwa
rationaler Zahlen) geredet wird, immer nur "`genügend viele"' \bruch
gemeint ist. Die Kontrolle, daß überall in der Tat
immer nur "`genügend viele"' gebraucht werden, ist
allerdings m.\ E. nicht trivial (In der Math.\ Zeitschr.~59,
1953, habe ich einen Beweis dafür dargestellt, daß in jedem
abzählbaren bewerteten Körper~$K$ ein \uline{abzählbarer} Oberkörper~$\overline{K}$
existiert, derart daß sich die Bewertung von~$K$ auf~$\overline{K}$ fortsetzen
läßt und daß für jeden algebraischen Oberkörper~$M$ von~$\overline{K}$
\uline{höchstens} \uline{eine} Bewertung von~$M$ existiert, die die Bewertung von~$\overline{K}$ 
fortsetzt).

Vor einiger Zeit erhielt ich die Neuauflage Ihrer Höheren
Algebra~I.\comment{\citealt{MR0158893}.} Daß die determinantenfreie Algebra jetzt zu
konstruktiven Entscheidbarkeits- und Berechenbarkeitsverfahren
führt, ist eine wesentliche Verbesserung.

Mir fiel auf (aber auch das soll keine Kritik sein, weil
es dazu zu unwichtig ist), daß für Körper und Gruppen immer
die Eindeutigkeit der Division (bzw.\ Subtraktion) gefordert
wird. In Satz~14 beweisen Sie aber die Eindeutigkeit des Eins\bs
elementes, nämlich: \(\begin{aligned}[t]%
                      \mathscr CE_A&=\mathscr C\text{ für alle }\mathscr C\\%
                      F_B\mathscr C&=\mathscr C\text{ für alle }\mathscr C%
\end{aligned}\)\\[2pt]
also $F_B=F_BE_A=E_A$ (jedes~$F_B$ ist jedem~$E_A$ gleich, d.\ h.\ alle~$E_A$ und~$F_B$ sind dasselbe Element).

Ebenso beweisen Sie Satz~15 die Eindeutigkeit des Reziproken.\strut\medskip

\hfil\begin{tabular}{@{}ll@{}}%
  Mit ergebensten Grüßen&\\%
                        &\rlap{stets Ihr P. Lorenzen}%
\end{tabular}
 

%% file: 19780221-PL-5-Aubert-Lorenzen
\label{19780221}%
\noindent\hfil \scshape university of oslo\nopagebreak

\noindent\begin{tabular}[t]{@{}l@{}}%
institute of mathematics\\%
\scriptsize p.\ o.\ box 1053 -- blindern, oslo 3%
\end{tabular}%
\hfill%
blindern, \upshape February 21, 1978\medskip

\noindent Dear Professor Lorenzen,\medskip

I am aware of the fact that you are now mainly working
as a philosopher and not as a mathematician. But since the question
I am going to ask you is really of a historical nature (although
connected with technical mathematics) I am nevertheless taking the
liberty to write to you. I believe that you are now (after Krull's
death) the only person who could possibly enlighten me on the questions
which I touch upon below. Hence I would be very thankful for any
comment that you would be willing to give me concerning this letter.

For some time I have been puzzled by the fact that so
few people know about $t$-ideals (``divisorial ideals of finite
character''). Specialized books on the topic of general arithmetics or
divisibility theory (multiplicative ideal theory, divisors, etc.)\
do not even mention this notion (like the books of Bourbaki, Gilmer,
Larsen--McCarthy and Fossum).\footnote{\selectlanguage{english}\citealt{bourbaki72,gilmer72,larsenmccarthy71,fossum73}.} This in spite of the fact that there
is in my opinion an overwhelming evidence pointing in favour of
$t$-ideals as the building blocks of general arithmetics. They seem
to form the true arithmetical divisors with nice properties which
are shared neither by the $d$-ideals nor by the $v$-ideals.

It is a puzzle to me how the basic arithmetical virtues of $t$-ideals
have escaped the notice of the specialists --~let alone the general
mathematical public. The only book which seems to give a fair
treatment of them is Jaffard's monograph
``\foreignlanguage{french}{Les systèmes
  d'idéaux}''.\footnote{\citealt{Jaf1960}. See
  \citet[Chapter~13]{bosbach15} for another fair treatment.}  And
Jaffard's main source when it comes to $t$-ideals is certainly your
own 1939 paper.

I believe it is high time to make some propaganda for
the $t$-ideals exhibiting them to a more general public in a somewhat
more pedagogical and attractive way than has been done hitherto.
I have hence decided to write an article on this subject with the 
tentative title ``Divisorial ideals of finite character''.\footnote{\citealt{aubert79}.} It is
going to be partly expository and partly historical, but will also
contain some more technical results of my own, going slightly
further than the results which can be found in your 1939 paper and
in Jaffard's monograph.

In connection with the history of the subject I am
especially puzzled by finding that even Krull did not seem to have
grasped the relevance of $t$-ideals --~at least if we are to judge
from his published works. I have not been able to spot a single
reference to $t$-ideals in his papers or in his `Idealbericht'.\footnote{\citealt{MR0229623}.}
I may of course have overlooked some part of his work, but I
cannot really see what that should be. I have checked rather
carefully his `\foreignlanguage{german}{Idealbericht}', his long series of papers on the
arithmetics of integral domains\footnote{\citealt{krull36,krull36a,krull37,krull38,krull37b,krull38b,krull39,krull39b,krull43b,krull43}.} as well as his papers on Krull
domains (`\foreignlanguage{german}{Endliche diskrete Hauptordnungen}').\footnote{\citealt{krull51,krull52,krull53}.} In doing so I have
made a couple of strange observations. In a short historical and
bibliographical note on $v$-ideals on page~121 of his `\foreignlanguage{german}{Idealbericht}'
he refers among other things to Arnold's paper ``\foreignlanguage{german}{Ideale in kommutativen
Halbgruppen}'' from 1929\footnote{\citealt{arnold29}.} as if this paper were concerned with $v$-ideals.
To describe the content of that paper Krull writes in a parenthesis:
($v$-Ideale in ``\foreignlanguage{german}{Halbgruppen}''). But Arnold's paper is concerned with
$t$-ideals and not with $v$-ideals!\footnote{Compare Lorenzen's letter to Krull dated 8 March 1938 (page~\pageref{19380308}).} My conjecture is simply that Krull
was mainly interested in the arithmetics of integral domains and
much less concerned with monoids. So he may not have studied Arnold's
paper very attentively -- in fact so superficially that the real
content of that paper may have escaped him.

On the other hand Krull makes extensive reference to
your 1939 paper (and in particular to the `finite character property'
of ideal systems) in the last paper (No~VIII from~1943) of his long
series of papers on the arithmetics of integral domains. And even on
a first and rather superficial reading of your paper one could hardly
avoid to see the crucial role played by the $t$-ideals. Only a glimpse
should be enough for a man like Krull in order to recognize this~--
the $t$-ideals lying at the cross-road of two of his pet topics
``\foreignlanguage{german}{Gruppensätze}'' and ``\foreignlanguage{german}{Endliche diskrete Hauptordnungen}'': An integral
domain is a Krull ring if and only if the fractional $t$-ideals form
a group (under $t$-multiplication). (Satz~7 in your 1939 paper is
equivalent to this result.) How could possibly this and many other
basic arithmetical facts concerning $t$-ideals have escaped Krull's 
keen eye for the essentials?

The $t$-ideals put the Krull rings and the Dedekind rings
on an equal footing describing the Krull rings simply as those
rings which are ``\foreignlanguage{german}{Dedekind}'' relative to the $t$-ideals. In the same way
the $t$-ideals also put the divisor class group on an equal footing
with the ideal class group. They are just particular instances
($r = t$, $r = d$) of a general notion of an ``$r$-class group''. Not to
mention that the unique factorization domains (monoids) are just
those domains (monoids) where every $t$-ideal is principal. (See the
footnote on page~543 in your 1939~paper.) In addition to this comes
the fact that (in contradistinction to the case of $v$-ideals) the
$t$-ideals are of finite character and hence allow the use of Zorn's
lemma, allow a nice theory of localization, etc. And this is very far
from exhausting the nice properties and the possible applications
of $t$-ideals.

My question to you concerns the possible additional
information you might have as to Krull's attitude to -- or knowledge
of -- $t$-ideals. Since you and Krull shared some of the same arithmetical
interests surrounding the notion of a $t$-ideal and at the same time
being colleagues for a number of years in Bonn -- I fancy that you
may have had discussions with him on this particular topic and that
you may possibly still have some recollections from such discussions.
If you know of some written material which points in the direction
of an awareness of the relevance of $t$-ideals on Krull's part, I would
of course appreciate this still more.

I am asking you these questions because I want my
historical remarks accompanying the above-mentioned projected paper
to be as accurate as possible. This in itself is of course a detail,
but I do not think it is just a detail or a matter of taste when
Bourbaki as well as most other writers are founding their exposition of
divisors, divisor class groups, Krull domains, etc.\ on the notion
of a $v$-ideal and not on the notion of a $t$-ideal. But it is not only
the $v$-ideals which on many occasions should be replaced by $t$-ideals --
it is even more so in the case of $d$-ideals. And this touches upon
an essential point of more far-reaching consequences: It seems
inappropriate to force the heavily additive Dedekind notion of an
ideal onto the solution of purely multiplicative problems of
divisibility --~just because there already happens to exist a well\Bs
developed commutative algebra surrounding this additive ideal concept.
I would advocate another and more natural procedure, namely to use
the $t$-ideals which are generally better suited for these purely
multiplicative problems and develop the corresponding necessary
piece of commutative algebra as adapted to this ideal concept
(like for instance the localization technique, the seeds of which
may already be found in your 1939 paper).

You may perhaps have lost all your former interest in
$t$-ideals, but in that case I still hope that the slight leanings
towards the history and philosophy of science of the present letter
may catch your interest -- or at least a sufficient interest to
make a reaction possible!\strut\medskip

\hfill%
\begin{tabular}{@{}l@{}}
Sincerely yours,\\[\medskipamount]%
Karl Egil Aubert%
\end{tabular}\medskip%

\noindent P.S. We once met in Bonn in 1954. I also knew Krull from that time
as well as from more recent encounters. But at that time I had not
really grasped the importance of $t$-ideals, missing my opportunity
to ask him these questions directly.


%% file: 19780306-PL-5-Lorenzen-Aubert
\label{19780306}%
\hfill  6.3.1978\nopagebreak

\hfill 2323\qquad\qquad\qquad

\noindent Prof.\ Dr.\ Paul Lorenzen\strut\medskip

\noindent\begin{tabular}{@{}l@{}}%
Herrn\\%
Prof.\ Karl Egil Aubert\\%
Institute of Mathematics\\%
University of Oslo\\%
P.\ O.\ Box 1053\\[2pt]%
Blindern, Oslo 3%
\end{tabular}\medskip%

\noindent Lieber Herr Aubert,\medskip

\noindent auf Ihren Brief hin habe ich erst einmal in meiner Dissertation
(das sind ja 40~Jahre her) nachsehen müssen, was $t$-Ideale sind.
In meiner Habilitationsarbeit (Über halbgeordnete Gruppen, Math.\ 
Z. 52, 1949) treten sie nur noch einmal auf (p.~486): für nicht\bs
kommutative Verbandsgruppen treten ``prime vollständige invariante
Unterhalbgruppen'' an die Stelle von $t$-Primidealen.

\noindent Durch die anschließende Verallgemeinerung der Teilbarkeitstheorie
von Gruppen auf Bereiche (Math.\ Z. 55, 1952) habe ich die $t$-Ideale
wohl ganz aus den Augen verloren -- und bin dann ja auch immer mehr
zu Grundlagenfragen der Logik und Mathematik gekommen, bis ich
jetzt schließlich ein sog.\ "`Wissenschaftstheoretiker"' geworden bin.

\noindent Meine Zusammenarbeit mit Krull war in der ganzen Zeit unseres Zu\bs
sammenseins schlecht. Ich wurde August~39, einen Monat vor Kriegs\bs
beginn, sein Assistent. Erst seit 1946 waren wir wieder in Bonn
zusammen, aber da war ich Dozent und wir waren fast ohne Kontakt.
Das hatte politische Gründe (Krull hatte meine Habilitation während
des Krieges verhindert)\comment{Lorenzen qualifies the reasons why he maintained almost no contact with Krull when they were both back in Bonn as ``political'' and not as personal. I see this as an indication that he was aware that the prevention of his habilitation by Krull (in agreement with Hasse) was connected with his lack of submissiveness under national socialism.} und lag außerdem an meinen grundlagentheo\bs
retischen Arbeiten, die Krull nicht als "`mathematische"' Leistungen
anerkannte.

\noindent Es kann also durchaus sein, daß er die Rolle der $t$-Ideale, im Unter\bs
schied zu den $v$-Idealen, nicht deutlich genug gesehen hat\comment{See however Krull's letter to Hasse dated 22~May 1938 \cite[§~1.35]{roquette04}.} -- wir
werden kaum jemals darüber gesprochen haben. Sofern ich damals
überhaupt noch Algebra gemacht habe, war das ja diese ``Teilbar\bs
keitstheorie in Bereichen'' (vgl.\ auch Math.\ Z.~58, 1953) oder
Theorie der Korrespondenzen (Math.\ Z.~60, 1954).\footnote{\citealt{Lor1952,Lor1953,lorenzen54}, respectively.}

\noindent Obwohl das kaum eine befriedigende Antwort auf Ihre Fragen ist
-- es ist aber so ge\-we\-sen -- möchte ich eine Gegenfrage stellen:
Kennen Sie Algebraiker, die die Teilbarkeitstheorie in Bereichen
weiter entwickelt haben?\strut\medskip

\hfill\hfill\hfill\hfill%
\begin{tabular}{@{}l@{}}%
  Mit besten Grüßen\\%
  \phantom{Mit }Ihr\\%
  \phantom{Mit besten Grüßen}\llap{Lo}[renzen]\\%
\end{tabular}\par 

%% file: 19780410-PL-5-Aubert-Lorenzen
\noindent\hfil\scshape universitetet i oslo
\label{19780410}%

\noindent\begin{tabular}{@{}l@{}}%
matematisk institutt\\%
  \scriptsize\begin{tabular}{@{}l@{}}%
               avd.~a: matematikk\\%
               avd.~b: mekanikk\\%
               avd.~c: statistikk og forsikringsmatematikk%
             \end{tabular}%
\end{tabular}\hfill%
\begin{tabular}{@{}l@{}}%
blindern, oslo 3, \upshape April 10, 1978\\%
\scriptsize telefon * 46 68 00%
\end{tabular}\upshape\medskip

\noindent Dear Professor Lorenzen,\nopagebreak\medskip

Many thanks for your letter of March~6. As you say
yourself, your reply does not give a really positive answer
to my specific question about Krull's eventual awareness of
the relevance of $t$-ideals. In spite of this, your ``negative''
information is also a piece of information for me.

I have little time now to work on the projected paper
on $t$-ideals which I mentioned in my first letter to you. But
I hope to be able to finish it during the summer and I will
then send you a copy of it. I am more and more convinced that
there is a need for such a paper, although it will hardly be
more than a reworking and a slight continuation of some of
your own ideas. One essential idea which is intimately linked
with $t$-ideals --~and which I did not mention in my last letter
to you~-- is the notion of a `Lorenzen group' (in Jaffard's
terminology). This is the purely multiplicative generalization
of a `Kronecker function ring' and is related to $t$-ideals and
$r$-valuations in the same way as the Kronecker function ring
is related to $d$-ideals and Krull valuations. What the concept
of a Lorenzen group shows, is that the $t$-system in a way has
a `universal property' reducing (in certain situations) the
study of $r$-systems to the study of the more concrete and well\Bs
behaved $t$-system in the corresponding Lorenzen group (which is
lattice-ordered and hence `$t$-Bézout').

As to your ``\foreignlanguage{german}{Gegenfrage}'' I will, generally speaking, say
that there have not been any remarkable innovations in divisibi\bs
lity theory since your own contributions --~which so far repre\bs
sent the end-point of a long-lasting German arithmetical traditi\bs
on. Your 1939~paper is widely cited, but I have somehow the
impression that only very few mathematicians have fully grasped 
the content of it. One sign of this is precisely the almost
total disregard of $t$-ideals in recent research literature as
well as in the books which lately have been published on multi\bs
plicative ideal theory.

As to your more specific question concerning further
developments of ``\foreignlanguage{german}{Teilbarkeitstheorie in Bereichen}'' I must confess
that I virtually know of no such contribution --~except for \uline{one}
paper ``\foreignlanguage{german}{Über verbandsgeordnete Vektorgruppen mit Operatoren}''
by I. Fleischer, published in Mathematische Nachrichten B~72
(1976).\footnote{\citealt{fleischer76}.} The opening sentences of this paper seem to confirm
the scarcity of contributions to this subject: ``\foreignlanguage{german}{Es sind schon
mehr als 20 Jahre vergangen, seitdem Lorenzen seine grundlegende
Arbeit ``Über halbgeordnete Gruppen'' veröffentlicht hat. Dagegen
scheint die Fortsetzung, Lorenzen [6] (``Teilbarkeitstheorie in
Bereichen'') völlig unbeachtet geblieben zu sein, -- leider, denn
viele der nachfolgenden Abhandlungen ließen sich ja in ihren
Rahmen zwangslos eingliedern.}''

It is hard to say why your papers on divisibility theory
have not quite exercised the influence which I think they deserve.
One reason may be that they are written in a style which many
mathematicians find hard to penetrate. This latter remark also
applies to Jaffard's monograph ``\foreignlanguage{french}{Les systèmes d'idéaux}'' which
is really a very fine book but which is written in a style and
uses a terminology which may have prevented many from reading
it who otherwise could have been attracted by its rich content.

I have for some time thought about writing a comprehen\bs
sive treatise on ideal systems which should not only be limited
to the arithmetical traditions of the theory --~but rather be a
kind of `generalized commutative algebra’ where the notion of
an ideal system is also seen in connection with ideal concepts
occurring outside of general arithmetics (like for instance in
general ring theory, differential algebra, lattice theory, etc.).
I have a lot of unpublished material on this topic, but I seem
to refrain from the rather laborious enterprise of putting it
all together.\strut\medskip

\hfill%
\begin{tabular}{@{}c}
  With my best regards\\[\medskipamount]%
  Karl Egil Aubert%
\end{tabular}

%% file: 19790618-PL-5-Lorenzen-Aubert
\label{19790618}%
\hfill 18.6.1979\nopagebreak

\hfill 2323\qquad\qquad\nopagebreak

\noindent Prof.\ Dr.\ Paul Lorenzen\medskip

\noindent \begin{tabular}{@{}l}%
Herrn\\%
Prof.\ Dr.\ K. E. Aubert\\%
Matematisk Institutt\\%
Universitetet i Oslo\\%
Postboks 1053\\[2mm]%
Blindern -- Oslo 3\\[2mm]%
Norwegen%
\end{tabular}\medskip%

\noindent Lieber Herr Aubert,\nopagebreak\medskip

\noindent herzlichen Dank für die Zusendung Ihrer Arbeit über die Teil\bs
barkeitstheorie.\footnote{\selectlanguage{english}\citealt{aubert79}, published subsequently as \citealt{aubert83}.} Ich habe ja diese Dinge fast ganz aus den
Augen verloren, es freut mich aber selbstverständlich sehr,
daß die multiplikative Idealtheorie bei Ihnen solche sorg\bs
fältige Neubearbeitung erfährt.

\noindent Es fällt mir auf, daß Sie sich auf kommutative Gruppen be\bs
schränken. War die Ausdehnung auf nichtkommutative Gruppen
vielleicht nur ein intellektueller Luxus?\strut\medskip

\hfill\hfill\hfill\hfill%
\begin{tabular}{@{}l@{}}%
  Mit herzlichen Grüßen\\%
  \phantom{Mit }stets\\%
  \phantom{Mit stets}Ihr\\%
  \phantom{Mit stets Ihr }Lo[renzen]\\%
\end{tabular}